\documentclass[UTF-8,reqno]{amsart}
\usepackage{enumerate}
\usepackage{mhequ}
\usepackage{amsthm,amsmath,amssymb,url,color, booktabs,nccmath}
\usepackage[left=3.2cm,right=3.2cm,top=4cm,bottom=4cm]{geometry}
\usepackage{mathrsfs}
\usepackage{enumitem,dsfont}
\usepackage{subfigure}
\usepackage[graphicx]{realboxes}
\usepackage{bbm}
\usepackage{subcaption}
\usepackage{times}\usepackage{mathrsfs}
\usepackage{amsfonts,amssymb,esint}
\usepackage{graphics}
\usepackage{enumerate}
\usepackage{mathtools}
\usepackage{xfrac}
\usepackage{tikz}
\usepackage{color}
\usepackage[colorlinks=true]{hyperref}
\hypersetup{
    linkcolor=blue,          
    citecolor=red,        
    filecolor=blue,      
    urlcolor=cyan
}

\definecolor{darkergreen}{rgb}{0.0, 0.5, 0.0}

\numberwithin{equation}{section}
\def\theequation{\arabic{section}.\arabic{equation}}
\newcommand{\be}{\begin{eqnarray}}
\newcommand{\ee}{\end{eqnarray}}
\newcommand{\ce}{\begin{eqnarray*}}
\newcommand{\de}{\end{eqnarray*}}
\newtheorem{theorem}{Theorem}[section]
\newtheorem{lemma}[theorem]{Lemma}
\newtheorem{proposition}[theorem]{Proposition}
\newtheorem{Examples}[theorem]{Example}
\newtheorem{corollary}[theorem]{Corollary}

\newtheorem{definition}[theorem]{Definition}
\newtheorem{conjecture}[theorem]{Conjecture}

\theoremstyle{definition}
\newtheorem{remark}[theorem]{Remark}
\newcommand{\RR}{\mathring R}
\newcommand*{\curl}{\ensuremath{\mathrm{curl\,}}}

\newcommand{\norm}[1]{\left\|#1\right\|}

\newcommand{\Rbar}{\mathring {\overline R}}
\newcommand{\RSZ}{\mathcal R}

\newcommand{\Proj}{\ensuremath{\mathbb{P}}}
\newcommand{\les}{\lesssim}

\DeclareMathOperator{\supp}{supp}

\def\T{\mathbb{T}^{3}}

\def\[{{\Big[}}
\def\]{{\Big]}}
\def\<{{\langle}}
\def\>{{\rangle}}
\def\({{\Big(}}
\def\){{\Big)}}

\def\bx{{\mathbf{x}}}
\def\tr{\mathrm {tr}}

\def\dif{{\mathord{{\rm d}}}}

\def\={&\!\!=\!\!&}

\def\cF{{\mathcal F}}

\def\cL{{\mathcal L}}

\def\mN{{\mathbb N}}

\def\mP{{\mathbb P}}
\def\mQ{{\mathbb Q}}
\def\mR{{\mathbb R}}

\def\mT{{\mathbb T}}

\def\mZ{{\mathbb Z}}

\def\bP{{\mathbf P}}
\def\bQ{{\mathbf Q}}
\def\bE{{\mathbf E}}
\def\1{{\mathbf{1}}}

\def\geq{\geqslant}
\def\leq{\leqslant}

\def\div{\mathord{{\rm div}}}

\def\[{{\Big[}}
\def\]{{\Big]}}
\def\<{{\langle}}
\def\>{{\rangle}}
\def\({{\Big(}}
\def\){{\Big)}}

\def\bx{{\mathbf{x}}}
\def\tr{\mathrm {tr}}

\def\dif{{\mathord{{\rm d}}}}

\def\={&\!\!=\!\!&}
\def\bt{\begin{theorem}}
\def\et{\end{theorem}}
\def\bl{\begin{lemma}}
\def\el{\end{lemma}}
\def\br{\begin{remark}}
\def\er{\end{remark}}
\def\bx{\begin{Examples}}
\def\ex{\end{Examples}}
\def\bd{\begin{definition}}
\def\ed{\end{definition}}
\def\bp{\begin{proposition}}
\def\ep{\end{proposition}}
\def\bc{\begin{corollary}}
\def\ec{\end{corollary}}
\def\bcon{\begin{conjecture}}
\def\econ{\end{conjecture}}

\def\geq{\geqslant}
\def\leq{\leqslant}

\def\div{\mathord{{\rm div}}}

\def\Id{\textrm{Id}}

\def\bP{{\mathbf P}}

 \def\R{\mathbb R}
 \def\R{\mathbb R}    
\def\N{\mathbb N}  
   
\def\<{\langle} \def\>{\rangle}

\allowdisplaybreaks

\begin{document}

\title[Non-uniqueness of Lagrangian Trajectories for  Euler]{Non-uniqueness of (Stochastic) Lagrangian Trajectories for Euler Equations}

\author{Huaxiang L\"u}
\address[H. L\"u]{Academy of Mathematics and Systems Science,
Chinese Academy of Sciences, Beijing 100190, China}
\email{lvhuaxiang22@mails.ucas.ac.cn }

\author{Michael R\"ockner}
\address[M. R\"ockner]{ Fakult\"at f\"ur Mathematik, Bielefeld Universit\"at, D 33615 Bielefeld, Germany, and  Academy of Mathematics and Systems Science,
Chinese Academy of Sciences, Beijing 100190, China}
\email{roeckner@math.uni-bielefeld.de}
 
\author{Xiangchan Zhu}
\address[X. Zhu]{State Key Laboratory of Mathematical Sciences, Academy of Mathematics and Systems Science,
Chinese Academy of Sciences, Beijing 100190, China}
\email{zhuxiangchan@126.com}

\thanks{Research  supported   by National Key R\&D Program of China (No. 2022YFA1006300, 2020YFA0712700) and the NSFC (No. 12426205, 12090014, 12288201). The financial support by the DFG through the CRC 1283 ``Taming uncertainty and profiting
 from randomness and low regularity in analysis, stochastics and their applications'' is greatly acknowledged.  Funded by the Deutsche Forschungsgemeinschaft (DFG,German Research Foundation)-Project-ID 317210226-SFB 1283"
}

\begin{abstract}
 We are concerned with the (stochastic) Lagrangian trajectories associated with Euler or Navier-Stokes equations. First, in the vanishing viscosity limit, we establish sharp non-uniqueness results for positive solutions to transport equations advected by weak solutions of the 3D Euler equations that exhibit kinetic energy dissipation with $C_{t,x}^{1/3-}$ regularity. As a corollary, in conjunction with  the superposition principle, this yields the non-uniqueness of associated (deterministic) Lagrangian trajectories.
 Second,   in dimension $d\geq2$,  for any $\frac{1}{p}+\frac{1}{r}>1$ or $p\in(1,2),r=\infty$, we construct solutions to the  Euler or Navier-Stokes equations in the space $L_t^rL^p\cap L_t^1W^{1,1}$, demonstrating that the associated  (stochastic)  Lagrangian trajectories are not unique.  Our result is sharp in 2D in the sense that: (1) in the stochastic case, for any vector field $v\in C_tL^p$ with $p>2$, the associated  stochastic Lagrangian trajectory associated with $v$ is unique (see \cite{KR05}); (2) in the deterministic case, the LPS  condition guarantees that for any weak solution $v\in C_tL^p$ with $p>2$ to the Navier-Stokes equations, the associated  (deterministic) Lagrangian trajectory is unique. Our result is also sharp in dimension $d\geq2$ in the sense that for any divergence-free vector field $v\in L_t^1W^{1,s}$ with $s>d$, the associated  (deterministic) Lagrangian trajectory is unique (see \cite{CC21}).
\end{abstract}


\date{\today}

\maketitle

\tableofcontents

\section{Introduction}
In this paper, we study the  Lagrangian trajectories associated with weak solutions to the incompressible Euler equations on the torus $\mathbb{T}^d = \mathbb{R}^d / \mathbb{Z}^d$ for $d \geq 2$: 
\begin{align}
  \frac{\dif}{\dif t} X_t &= v(t, X_t),\ t\in[0,T], \label{eq:ode} \\
  X_0 &= x, \notag
\end{align}
where  $T>0,v: [0, T] \times \mathbb{T}^d \to \mathbb{R}^d$ is a weak solution to the incompressible deterministic Euler equations on $[0,T]\times \mT^d$:
\begin{align}
     \partial_t v + \div(v \otimes v) + \nabla \pi &= 0,\label{eq:euler} \\
     \div v &= 0, \notag
\end{align}
where $\pi$ denotes the pressure field associated with the fluid.  
 In this context, $X: [0, T] \to \mathbb{T}^d$ represents the trajectory of a particle in an incompressible, non-viscous fluid.  In this paper, we focus on the  Lagrangian trajectory  in the following sense.
 
\bd\label{def:intcurve}
 Let $v : [0, T ] \times \mathbb{T}^d \to \R^d$ be a Borel map  and $x\in\mT^d$. We say that  an absolutely continuous function $X_{\cdot}^x\in AC([0, T ]; \mathbb{T}^d )$  is a  (deterministic) Lagrangian trajectory of $v$ starting from $x$ if   for  all $t \in [0, T ]$,
\begin{align*}
    X_t^x=x+\int_0^tv(s,X_s^x)\dif s.
\end{align*}
\ed

It is now widely believed that both the 3D Euler equations and their corresponding Lagrangian trajectories may develop singularities. 
The famous Onsager conjecture \cite{Ons49}  says that the threshold regularity for energy conservation  of weak solutions to the Euler equations \eqref{eq:euler} is $\frac{1}{3}$, which comes from Onsager's  attempt to explain the primary mechanism of energy dissipation in turbulence. Consequently, he implicitly suggested that H\"older continuous weak solutions of the Euler equations could provide an appropriate mathematical description of turbulent flows in the inviscid limit. In recent years, this conjecture has been fully proven using the technique of convex integration (see \cite{Ise18,BDLSV19}). For  a further discussion, we refer to Section \ref{sec:tur}  below.  We also emphasize that $\frac{1}{3}$ regularity is also related to the dissipative length scales in Kolmogorov’s 1941
 (K41) phenomenological theory of turbulence \cite{Kol41} (see \cite{HPZZ23,Nov23}).

On the other hand, a fundamental insight into the Lagrangian origin of turbulent scalar dissipation is the concept of spontaneous stochasticity, which was predicted by Lorenz \cite{Lor63,Lor69}. This notion suggests that multiscale fluid flows can inherently lose their deterministic nature and become intrinsically random.  The pioneering mathematical work of Bernard, Gawedzki, and Kupiainen \cite{BGK97} examined Kraichnan's turbulence model \cite{Kra68}, where the advected velocity is represented as a Gaussian random field with white-noise correlation in time. They demonstrated that, due to the spatial roughness of the advected field, Lagrangian trajectories become non-unique and stochastic in the limit of infinite Reynolds number, even for a fixed initial particle position.
We refer to \cite{DE17,DE17b,ED18,TBM20,MR23,JS24} for recent studies on spontaneous stochasticity.

From the above, studying non-uniqueness of Lagrangian trajectories associated with Euler equations is of primary interest for turbulence.
 In this work, we first adopt an "Eulerian“ perspective, which provide a more robust framework for understanding non-uniqueness phenomena. Specifically, we examine the evolution of the superposition of a family of  Lagrangian  trajectories starting from $x \in \mT^3$. This approach naturally leads us to study the following transport equation:
\begin{align}\label{eq:tpe}
\partial_t \rho+\div(v\rho)=0,\\
\rho(0)=\rho_0.\notag
\end{align}
Here, $\rho:[0,T]\times \mT^d\to \mR^d$ represents the probability density of the particles.

Our first main result of the paper is the following sharp non-uniqueness for the transport equations  advected by weak solutions to Euler equations with H\"older regularity up to the Onsager exponent $\frac13$.

\bt\label{thm:con:1/3}
 Let $\beta ,\tilde{\beta}\in(0, 1)$ and $T>0$ be fixed.\\
 (1). Let $\beta+2\tilde{\beta}>1$.
 For any divergence-free vector field $v\in C^\beta([0,T]\times \mathbb{T}^3)$, the transport equation \eqref{eq:tpe}  has a unique solution  $\rho\in C^{\tilde{\beta}}([0,T]\times \mathbb{T}^3)$, which  conserves the kinetic energy.\\
 (2). Let  $0<\beta+2\tilde{\beta}<1,0<\beta<\frac13$
 and $e:[0,T]\to\mathbb{R}$ be a strictly positive smooth function. Then  there exists a solution $v\in C^\beta([0,T]\times \mathbb{T}^3)$ to the 3D Euler equation \eqref{eq:euler} satisfying
 \begin{align*}
     \int_{\mT^3}|v(t,x)|^2\dif x=e(t),
 \end{align*}
 such that there is a  non-constant  probability density $\rho  \in C^{\tilde{\beta}}([0,T]\times \mathbb{T}^3)$  solving the transport equation \eqref{eq:tpe}  with initial data $\rho_0 = 1$ the  unit volume. Then \eqref{eq:tpe}  admits at least two probability density solutions in the space $C^{\tilde{\beta}}([0,T]\times \mathbb{T}^3)$ by noticing that $\overline \rho(t)\equiv 1$ is always a solution. 
\et

 We note that non-uniqueness of solutions to \eqref{eq:tpe} is not entirely new, as it also follows from combining two existing results:  the work  \cite[Theorem 1.1]{BSW23} on the anomalous dissipation of passive scalar advected by  solutions of Euler equations with H\"older regularity $\frac13-$, and the work  \cite[Theorem 1.3]{Row24} on the relationship between anomalous dissipation and non-uniqueness of Lagrangian trajectories. In contrast to the
 prior results,  we establish the
  H\"older regularity  for solutions to the transport equations, and even obtain a optimal regime for non-uniqueness except for the endpoint $\beta+2\tilde{\beta}=1$. While the aforementioned approaches relied  on quantitative homogenization  and convex integration applied solely to the Euler equation, our work introduces a   convex integration iteration that simultaneously addresses both the transport equations and the Euler equations.

The first part of Theorem \ref{thm:con:1/3} is established using commutator estimates, following the lines of \cite{CET94}. The second part of Theorem \ref{thm:con:1/3} is our main result and  will be established using convex integration method in Section \ref{sec:Euler:C1/3} below. The detailed construction using the convex integration method is shown in Section \ref{sec:proof13}. 

The connection between the Eulerian description (governed by the transport equation) and the Lagrangian description (involving particle trajectories) is established through the superposition principle (see \cite[Theorem 3.2]{Amb08} in the case of $\mR^d$. For the case of manifolds as our $\mT^d$, see \cite[Section 7.2]{Tre14}). Specifically, the superposition principle allows us to express the solution $\rho$ as a superposition of  time marginal laws of probability  measures $\bQ^x$ supported on Lagrangian trajectories  started at $x\in\mT^3$.  By applying the superposition principle with $\beta<\frac13$ close to $\frac13$, along with a decreasing function $e(t)$, we establish the following result: there exist dissipative solutions to the 3D Euler equations below the Onsager  exponent $\frac13$ for which deterministic Lagrangian trajectories are non-unique. 

\begin{corollary}\label{thm:nonode1/3}
  Let $T>0$ and $0<\beta<\frac13$ be fixed. There exists a  solution  $v\in C^\beta([0,T]\times \mT^3)$   to the Euler equation \eqref{eq:euler} which dissipates the kinetic energy, such that non-uniqueness of (deterministic) Lagrangian trajectories holds in the sense that:\\
There is a measurable  set $A(v) \subset \mathbb{T}^3$ with positive Lebesgue measure such that for every $x \in A(v)$ there are at least two  trajectories of $v$ starting at $x$. 
\end{corollary}



 Moreover,  it is also mathematically interesting to consider more general cases, where the fluid dynamics incorporates viscous effects, and the Lagrangian trajectories are subject to external stochastic perturbations. In the following, we shall show some sharp non-uniqueness  results of stochastic Lagrangian trajectories associated  with the Euler equation, and even with the Navier-Stokes equation in general dimension $d\geq2$, i.e. solutions to
 \begin{align}
  \dif X_t &= v(t, X_t) \dif t + \sqrt{2\kappa} \dif W_t, \ t\in[0,T],\label{eq:sde} \\
  X_0 &= x, \notag
\end{align} 
where $\kappa \in (0, 1]$ be fixed, $X: [0, T] \to \mathbb{T}^d$ is the stochastic process  representing the stochastic particle trajectory, and $W_t$ is a  standard
$\mathbb{T}^d$-valued Brownian motion  defined on a probability space $(\Omega, \mathcal{F},(\cF_t)_{t\geq0}, \mathbf{P})$. See below for a precise definition.  
The drift term $v: [0, T] \times \mathbb{T}^d \to \mathbb{R}^d$ is a weak solution to the incompressible deterministic  Navier-Stokes or Euler equations on $[0,T]\times \mT^d$:
 \begin{align}
     \dif v+\div(v\otimes v)\dif t-\nu\Delta v\dif t&+\nabla{\pi}\dif t=0,\label{eq:ns}\\
\div v&=0, \notag
\end{align}
where $\nu\in[0,1]$, $\pi$ denotes the pressure field associated with the fluid. When  $\nu=0$, \eqref{eq:ns} is the Euler equation \eqref{eq:euler}.  Now we focus on the stochastic Lagrangian trajectories  in the following sense.

\bd\label{def:intcurvesde}
 Let $v : [0, T ] \times \mathbb{T}^d \to \R^d$ be a Borel map and  $W$ be a given $\mT^d$-valued Brownian motion on the probability space $(\Omega,\cF,(\cF_t)_{t\geq0},\bP)$. We say that    an $(\cF_t)_{t\geq0}$-adapted map $X_\cdot^x\in C([0, T ]; \mathbb{T}^d ),\bP$-a.s.  is a stochastic Lagrangian trajectory of $v$ starting from $x$ if 
\begin{align*}
    X_t^x=x+\int_0^tv(s,X_s^x)\dif s+\sqrt{2\kappa}W_t,\ \ t \in [0, T ], \ \ \mathbf{P} {\rm -a.s.}
\end{align*}
\ed

\bd We say that uniqueness in law holds for the stochastic Lagrangian   trajectories if for any two sets of stochastic Lagrangian trajectories $\{X^x\}_{x\in\mT^d}$ and $\{\overline X^x\}_{x\in\mT^d}$ (which may be defined on different probability spaces), we have $\bP\circ (X^x)^{-1}=\overline\bP\circ (\overline X^x)^{-1}$ for a.e. $x\in\mT^d$.
\ed

\br
The SDE in Definition \ref{def:intcurvesde} is meant in the sense of solving the corresponding martingale problem. This is sufficient, since we shall  only consider the law of the stochastic Lagrangian trajectory.
\er

 In the study of SDEs, there is evidence that  a suitable stochastic noise may provide a regularizing effect on deterministic ill-posed problems. 
 On the whole space $\R^d$, for vector fields that are only Lebesgue-integrable, specifically $v\in L^r_tL^p :=L^r([0,T];L^p)$ satisfying  $\frac dp+\frac2r<1$, Krylov and the second named author \cite{KR05}   established  strong existence and uniqueness of stochastic Lagrangian   trajectories  to \eqref{eq:sde}, i.e. adapted to the given Brownian motion.  By an analogous argument, this result is also valid on the torus.  However, beyond the condition $\frac{d}p+\frac 2r<1$,  the question whether  stochastic Lagrangian trajectories are pathwise unique-or, more weakly, unique in law-is not completely understood.   We refer to Section \ref{sec:pre_sde} for recent results on  this problem.
 
 On the other hand, in the study of hydrodynamics, there is a deep connection between  SDEs \eqref{eq:sde} and  the Navier-Stokes
equations \eqref{eq:ns} with $\nu=\kappa>0$. When $v$ is  a smooth solution to the Navier-Stokes equations, Constantin and Iyer \cite{CI08} provided the following stochastic representation:
\begin{align} 
v(t, x) = \mP_{\rm H} \bE[\nabla^T(X^x_t)^{-1}v_0((X^x_t)^{-1})],\label{intro:prob_pres}
\end{align}
where $\mP_{\rm H}$ is the Leray projection and $v_0$ is the initial condition. Conversely, if
$v$ is smooth and $(v,X^x)$ solves \eqref{eq:sde}
and \eqref{intro:prob_pres}, then $v$ is also   a solution to the Navier-Stokes equations with initial condition $v_0$. The  condition  $\frac dp+\frac2r< 1$ is   special case of the famous Ladyzhenskaya-Prodi-Serrin (LPS) condition  $\frac dp+\frac2r\leq  1$,  which provides a sufficient condition for the regularity and uniqueness  of weak solutions to the Navier-Stokes equations. Within the LPS condition, both the NS equation and the corresponding   Lagrangian trajectories are  unique. Beyond this condition, using  the convex integration method, sharp non-uniqueness of weak solutions to the Navier-Stokes equations has been shown in \cite{BVb,BCV18,CL22,CL23}.
 
 In summary, whether from the perspective of SDEs or from the perspective of  fluid dynamics, studying the well/ill-posedness of stochastic Lagrangian trajectories with  hydrodynamic drifts beyond the  LPS condition is of  theoretical interest. 
 Our second main result is to show the non-uniqueness in law of the stochastic Lagrangian trajectories of weak solutions to the Navier-Stokes or Euler equations  \eqref{eq:ns} in the supercritical regime, i.e. $\frac dp+\frac2r>1$.  To state the main result,  for $d\geq2$ we define $\mathcal{A}:=\mathcal{A}_1\cup \mathcal{A}_2$, where
\begin{align*}
    \mathcal{A}_1:&=\left\{(p,r,s)\in[1,\infty]^3:\frac1p+\frac1r>1,1<s<d\right\},\\
    \mathcal{A}_2:&=\left\{(p,r,s)\in[1,\infty]^3:\frac1p+\frac1r\leq 1,1<p<2,\frac1d+\frac12\frac{1-\frac1r-\frac1p}{\frac12-\frac1r}<\frac1s<1\right\}.
\end{align*}

\bt\label{thm:nonuniode+sde}
 Let $d\geq2,\kappa\in[0,1],\nu\in[0,1]$. For any triple $(p,r,s)\in \mathcal{A}$, there exists a divergence-free
vector field  $v\in L^r([0,T];L^p)\cap L^2([0,T]\times\mathbb{T}^d)\cap L^1([0, T]; W^{1,s})\cap C([0,T];L^1)$  which is a  solution to the Navier-Stokes or 
Euler equations  \eqref{eq:ns}, such that  the law of stochastic Lagrangian trajectories of $v$ is not unique in the sense that:\\
There is a measurable $A(v) \subset \mathbb{T}^d$ with positive Lebesgue measure such that for every $x \in A(v)$ there are at least two stochastic  trajectories of $v$ starting at $x$, admitting distinct laws satisfying 
$\mathbf{E}[\int_0^T|v(s,X^x_s)|\dif s] <\infty$.

  Moreover, when $r=\infty$, the solution $v$ is continuous in $L^p$-norm, i.e. $v\in C([0,T];L^p) =:C_tL^p$.
\et

As a corollary, in the stochastic case $\kappa>0$, by taking $s>1$ close to 1 in Theorem \ref{thm:nonuniode+sde}, we show that sharp non-uniqueness in law holds for stochastic Lagrangian trajectories of weak solutions to the Navier-Stokes or Euler equations \eqref{eq:ns}, within a range of supercritical regimes:
 
\bc\label{thm:nonlaw}
  Let $d\geq2,\kappa\in(0,1],\nu\in[0,1],p,r\in[1,\infty]$ be fixed.  There exists  a solution
  \begin{align*}
v\in
\begin{cases}
L^r([0,T];L^p)\cap L^2([0,T]\times \mT^d)\cap L^1([0, T]; W^{1,1})\cap C([0,T];L^1),&\ {\rm if}\ \ \frac1p+\frac1r>1,\\
 C([0,T];L^p)\cap L^2([0,T]\times \mT^d) \cap L^1([0, T]; W^{1,1}),&\ {\rm if}\ \ 1<p<2,
\end{cases}
\end{align*}
  to the Navier-Stokes or Euler equations   \eqref{eq:ns}, such that the law of stochastic Lagrangian trajectories of $v$ is not unique in the sense that:\\
 There is a measurable $A(v) \subset \mathbb{T}^d$ with positive Lebesgue measure such that for every $x \in A(v)$ there are at least two  trajectories of $v$ starting at $x$, admitting distinct laws satisfying $\mathbf{E}[\int_0^T|v(s,X^x_s)|\dif s] <\infty$. 
\ec

\begin{figure}
    \centering
   
\begin{tikzpicture}[scale=5]
    \draw[->] (0,0) -- (1.1,0) node[right] {$\frac1p$};
    \draw[->] (0,0) -- (0,1.1) node[above] {$\frac1r$};

    \draw[red,thick] (0,1/2) -- (1/2,0);

    \fill[red,opacity=0.3] (0,1/2) -- (1/2,0) -- (0,0) -- cycle;
    \fill[green,opacity=0.3] (0,1) -- (1/2,1/2) -- (0.5,0) --(1,0)-- (1,1) -- cycle;

    \node at  (0.7,0.6) {¬!};
    \node at (0.1,0.2) {!};
      \node at  (0.27,0.4) { $\frac 1p+\frac1r= \frac12$};

    \draw[dotted] (0,1) --(1/2,1/2)-- (1/2,0);

     \draw (1/2,1/2) node[right] {$(\frac12,\frac12)$} ;
    \draw (1,0) node[below] {$1$} -- +(0,0.02);
    
     \draw (0,1/2) node[left] {$\frac12$} -- +(0.02,0);
    \draw (0,1) node[left] {$1$} -- +(0.02,0);
    \draw (0.5,0) node[below] {$\frac12$} -- +(0,0.02);
    
    \filldraw[black] (1/2,1/2) circle (0.03mm);
\end{tikzpicture}

    \caption{State of (non-)uniqueness to SDE \eqref{eq:sde} for  vector fields $v\in
L^r_tL^p$ in  the 2D case.\\
Red area: in the subcritical case, the SDE admits a unique strong solution.\\
Red line: in the critical case
, the well-posedness of  SDE remains open in $d=2$.\\
Green area: in the supercritical case, our main result shows the non-uniqueness in law holds in this range.}
\label{tu:sde}
\end{figure}
In Figure \ref{tu:sde}, we summarize the (non-)uniqueness in law
results on SDE \eqref{eq:sde} for  vector fields $v\in
L^r_tL^p$ in the case $d=2$. Our result implies that for any $p<2$, there exists a solution $v\in C_tL^p$ to  the Navier-Stokes or Euler equations such that non-uniqueness in law of stochastic Lagrangian trajectories of $v$ holds, which  is sharp  in the sense  that for  any vector field  $v\in C_tL^p$, $p>2$, \eqref{eq:sde} admits a unique strong solution  as recalled above. Our result covers the optimal range, except for the endpoint $p=2$.
 We also notice that since $\div v=0$, \eqref{eq:sde} admits  an invariant measure $ {\rho}\equiv1$. We obtain non-uniqueness in law even when starting from this invariant measure.

In contrast to this, in the deterministic case $\kappa=0$, despite the absence of regularization-by-noise phenomena, the LPS criteria guarantee that any weak solution $v\in C_tL^p$ with $p>d$ to the Navier-Stokes equations is  Leray and regular, which still ensures  the uniqueness of the associated  (deterministic) Lagrangian trajectory.  Taking $r=\infty$ in Theorem \ref{thm:nonuniode+sde}, the  following corollary  presents a sharp counterexample  $v\in C_tL^{2-}\cap L^1_tW^{1,1+}$ for which the Lagrangian trajectories are non-unique, which covers the optimal range, except for the endpoint $p=2$.

\bc\label{thm:nonuniode}
  Let $d\geq2,\nu\in[0,1]$. For any $1< p<2,\frac1p+\frac1s>1+\frac1d$, there exists a divergence-free
vector field  $v\in  C([0,T];L^p)\cap L^2([0,T]\times\mathbb{T}^d)\cap L^1([0, T]; W^{1,s})$  which is a  solution to the Navier-Stokes or Euler equations  \eqref{eq:ns}, such that   the non-uniqueness of (deterministic) Lagrangian trajectories holds in the sense that:\\
There is a measurable $A(v) \subset \mathbb{T}^d$ with positive Lebesgue measure such that for every $x \in A(v)$ there are at least two  trajectories of $v$ starting at $x$.
\ec

Furthermore, for any $s<d$, there is a solution to the Navier-Stokes or Euler equations  \eqref{eq:ns} in the  space 
 $C_tL^1\cap L^2_{t,x}\cap L^1_tW^{1,s}$
for which the Lagrangian trajectories are non-unique. Our result is sharp  in comparison to  the result of Caravenna and Crippa \cite[Corollary 5.2]{CC21}, which states that  for any divergence-free vector field  $v\in L^1_tW^{ 1,s}$ with  $s > d$, then for a.e. $x \in\mathbb{T}^d$, there is a unique trajectory  of $v$ starting at $x$. Our result covers the optimal range except for the endpoint $p=d$.  We also remark that Bru\'e,  Colombo and De Lellis \cite{BCDL21}  provided a nice divergence-free counterexample  $   v\in C_t(W^{1,d-}\cap L^{\infty-})$ using convex integration.

As before, we adopt an "Eulerian" perspective and establish the non-uniqueness of solutions to the corresponding Fokker-Planck equations on  $[0,T] \times \mathbb{T}^d$:
\begin{align}\label{eq:fpe}
\partial_t \rho - \kappa \Delta \rho + \div(v \rho) &= 0, \\
\rho(0) &= \rho_0. \notag
\end{align}

\bt\label{thm:non_pde_23}
Let $\kappa,\nu\in[0,1]$. Then for any triple $(p,r,s)\in\mathcal{A}$,  there exists a solution $v\in L^r([0,T];L^p)\cap L^2([0,T]\times \mT^d)\cap L^1([0, T]; W^{1,s})\cap C([0,T];L^1)$ to  \eqref{eq:ns}, and a non-constant density $\rho  \in  L^r([0,T];L^p)\cap L^2([0,T]\times\mT^d)\cap C([0,T];L^1)$ satisfying \eqref{eq:fpe} with initial data $\rho_0 = 1$. 

Moreover, when $r=\infty$, the solutions $v,\rho$ are also continuous in  $L^p$-norm.
\et

This result is proved in Section \ref{cogpss2}, while the detailed construction via the convex integration method is given in Section \ref{proof:prop2}. Then Theorem \ref{thm:nonuniode+sde} follows    by the superposition principle proved in  \cite[Section 7.2]{Tre14}.

\subsection{Convex integration and Onsager conjecture}\label{sec:tur}
 The rigid part of Onsager’s conjecture was established using commutator estimates, as demonstrated in \cite{CET94, CCFS08}. The flexible part was proven using the convex integration method by Isett in \cite{Ise18} for the 3D case, with further constructions of strictly dissipative solutions discussed in \cite{BDLSV19}. 
 This convex integration technique was first introduced to fluid dynamics by De Lellis and Sz\'ekelyhidi Jr. \cite{DLS09, DLS10, DLS13}, leading to numerous groundbreaking results. 
 For the $L^2$-based Sobolev scale,  Buckmaster,  Masmoudi,  Novack, and  Vicol \cite{BMNV23} constructed non-conservative weak solutions of the 3D  Euler equations  in $C_tH^{1/2-}$, and  then refined by   Novack and Vicol \cite{NV23}
 to the class $C^0_ t(H^{1/2-}\cap L^{\infty})$. By interpolation, such solutions belong  to $C_tB^{1/3-}_{3,\infty}$, which can be seen as a proof of  the $L^3$-based intermittent Onsager theorem. 
 Notably, the 2D   Onsager conjecture was addressed by Giri and Radu \cite{GR24} through convex integration with Newton iteration.    We refer to \cite{Cho13,TCP,Buc15,BDLS16,DSJ17,DLK22,Ise22,BHP23,BC23,GKN23,GKN24,BM24a,BM24b}  for more  results on the Euler equations.
We mention that the convex integration method also  led to a breakthrough to the  non-uniqueness of weak solutions to  the Navier-Stokes equations, see for example \cite{BVb,BCV18,CL22,CL23,MNY24a,MNY24b,CZZ25}.    We refer interested readers to the comprehensive reviews \cite{BV, BV21,DLS22} for more details and references.

Very recently, Bru\'e, Colombo and Kumar \cite{BCK24b} introduced a new “asynchronization” idea for the building blocks in the convex integration  to  derive non-uniqueness of weak solutions in $L^\infty_tL^2$ with vorticity in $L^\infty_tL^{1+}$, which solves a longstanding
open problem of non-uniqueness of the 2D Euler equations with integrable vorticity.

We also note that the convex integration method has been successfully applied to the stochastic fluid dynamics (see \cite{BFH20,HZZ20,Yam22a,Yam22b, HZZ,HZZ21b,LZ23,HLP24,HZZ19,LRS24,Pap24,HZZ22,LZ23b,LZ24} and 
references therein).

\subsection{Previous results on the ODE level}
If the advected vector field $v$  is Lipschitz continuous, classical theorems ensure the uniqueness of Lagrangian trajectories starting from any $x\in\mT^d$. For less regular vector fields  on the whole space $\mR^d$,  DiPerna and Lions \cite{DL89} proved the uniqueness of trajectories in the class of regular Lagrangian flows under suitable Sobolev regularity conditions. 
  Moreover, for a divergence-free vector field  $v\in L^1_tW^{1,s}$ and  every $\rho_0 \in  L^p$ with $\frac1p+\frac1s\leq1$,
there exists a unique weak solution $\rho \in L_t^\infty L^p$ to the transport equation \eqref{eq:tpe}. An analogous version of this result holds on the torus  $\mT^d$.  Later, the DiPerna-Lions theory  was  extended by Ambrosio \cite{Amb08} to the bounded variation case $v\in L^1_t(BV)$. 
However, it is not clear whether  for almost every $x$, there is a unique trajectory of $v$ starting at $x$.  
 We refer to  \cite{ DL08,Amb17,CC21,BCDL21} for more related uniqueness results beyond the DiPerna-Lions condition.

 Concerning fluid equations on the whole space $\mR^3$, if $v$ is a Leray solution to the 3D Navier-Stokes equations as defined in \eqref{eq:ns} and  if it is obtained through approximation with $v(0)\in H^{1/2}(\mR^3)$, Robinson and Sadowski \cite{RS09a,RS09b} proved the existence and uniqueness of Lagrangian trajectories for a.e.
 initial point $x\in\R^3$. Then Galeati \cite{Gal24} refined it to the case $v(0)\in L^2(\mR^3)$.

  Regarding non-uniqueness,  to the best of our knowledge, the literature currently features several distinct approaches. The first approach is Lagrangian, which involves using the degeneration of the flow map to show non-uniqueness at the ODE level, and we refer to \cite{DL89,Dep03,YZ17,ACM19,DEIJ22,Kum24,Pap23,BCK24,MS24b}. These constructions are usually  quite specific  and play an important role in understanding the fluid dynamics. 
The second approach is Eulerian, by  demonstrating non-uniqueness directly at the PDE level. Crippa, Gusev, Spirito, and Wiedemann \cite{CGSW15} were the first to utilize convex integration to derive non-uniqueness results for transport equations \eqref{eq:tpe}  on torus. Subsequently, numerous significant breakthroughs concerning Sobolev vector fields were achieved through convex integration in \cite{MS18,MS19,MS20,BCDL21,CL21,CL22b,PS23}.
A more recent approach, developed by Armstrong and Vicol \cite{AV25}, employs quantitative homogenization techniques. Subsequently, Burczak, Sz\'ekelyhidi Jr, and Wu \cite{BSW23} successfully combined this methodology with convex integration to demonstrate anomalous dissipation of passive scalar advected by weak solutions of the Euler equations.

\subsection{Previous results on the SDE level}\label{sec:pre_sde}
 As mentioned earlier, the regularization-by-noise phenomenon plays a significant role in the well-posedness of SDEs, which  comes from the effect of  the Laplacian in the Fokker-Planck equations.
 
   On the whole space, when $v$ is a bounded measurable function, Veretennikov \cite{Ver80} proved the uniqueness of probabilistically strong solutions. When $v\in L^r_tL^p$, Krylov and the second named author  in \cite{KR05} established the existence and uniqueness of strong solutions to \eqref{eq:sde} in the class $\int_0^T|v(s,X_s)|^2\dif s<\infty,\ \mathbf{P}$-a.s., under the  condition $\frac dp+\frac2r<1.$  Moreover, the square integrability condition can be removed, see for example \cite[Lemma 3.4]{Hao23}. For more well-posedness  results in the subcritical case, we refer to \cite{Zha05,Zha11,Zha16,XXZZ20,RZ21a}.

In the critical case $\frac{d}{p} + \frac{2}{r} = 1$,  on the whole space, Krylov \cite{Kry20} proved the strong
well-posedness  of SDEs in the case  $ v\in L^d (\mR^d )$, which is a significant
progress on this topic.  The second named author and Zhao \cite{RZ21} showed that for any vector field $v\in C_tL^d$ or $v\in L^r_tL^p,r,p\in(2,\infty)$, \eqref{eq:sde} admits a unique strong solution within a class satisfying a Krylov-type estimate. They \cite{RZ23} also proved weak uniqueness 
with  divergence-free $v \in L^\infty_t L^d$ within a class satisfying a Krylov-type estimate.  We refer to \cite{Nam20,Kry20b,Kry20c} for further results.

Beyond the LPS condition $\frac{d}{p}+\frac2r\leq1$, the known results in the supercritical regime are  very limited. 
 When the vector field is not divergence-free, in \cite{BFGM19}, there is a counterexample showing that \eqref{eq:sde} may not have weak solutions  if  $v$ is in the Lorentz space $L^{d,\infty}(\R^d)$.  Then Zhao \cite[Theorem 5.1]{Zha19}  constructed a  divergence-free vector field $v \in L^p(\R^d) + L^\infty(\R^d),p\in(\frac d2,d),d\geq3$, such that weak uniqueness  fails.  Given an additional divergence-free property on the drift $v$, Zhang and Zhao \cite{ZZ21} established weak existence and uniqueness in the sense of approximation for  $\frac{d}{p} + \frac{2}{r} \leq 2$.
 Recently, Galeati \cite{Gal24}  demonstrated strong existence and pathwise uniqueness for every initial $x\in\R^3$  for the SDEs \eqref{eq:sde}, under the assumption that $v$ is a Leray solution to 3D Navier-Stokes equations obtained through approximation with divergence-free $v(0)\in H^{1/2}(\mR^3)$.   It is also worth mentioning that very recently, for $d\geq1$ and any $\frac{d}{p}+\frac2r>1,p>d$, Galeati and Gerencs\'er in \cite{GG25} constructed an example $v\in L_t^rL^p(\mR^d)$  such that non-uniqueness in law holds for \eqref{eq:sde} when starting from  $x = 0$.  We also refer to \cite{Gal23, BG23,HZ23,GP24,HRZ24} for more results in the supercritical regime.

\subsection{Ideas of the proof}\label{sec:idea}

In this paper, to demonstrate the non-uniqueness of stochastic Lagrangian trajectories, in view of the superposition principle, we work at the PDE level to construct a fluid field such that the transport equation admits two solutions. To achieve this, we concurrently apply the convex integration method to the fluid equation and the transport equation simultaneously.  However, when considering two different scales $C_{t,x}^0$ or $L^1_tW^{1,s}$, we face several challenges, which require to address them using distinct strategies tailored to each scale's characteristics.

\subsubsection{Ideas in the $C_{t,x}^0$-scales}
We apply the convex integration method to transport equations and the Euler equations simultaneously. At each step $q\in\mN_0$, we construct a pair $(v_q, \rho_q, \mathring{R}_q, M_q)$  satisfying the following system:
\begin{align}
\partial_t \rho_q+\div(v_q \rho_q)&=-\div M_q,\notag\\
\partial_tv_q+\div(v_q\otimes v_q)+\nabla\pi_{q}
&=\div\mathring{R}_q,\ \ \div v_q=0,\label{eq:idea:rhov}
\end{align}
where $\mathring{R}_q$  is  a trace-free symmetric matrix, and $M_q$ is  a vector field. Here, $\mathring{R}_q$ and $M_q$ converge  to 0 and $(v_q,\rho_q)$  converge  to  a weak solution to the transport equation and the Euler equations respectively.

At each iterative step, we need to construct new perturbations to simultaneously cancel two stress terms and obtain smaller residual stress terms. 
However, a new  velocity perturbation constructed to cancel one of the stress terms will also  influence  the other stress term.   To overcome  this,  inspired by the work of Isett \cite{Ise22},  we define the perturbations $(w_{q+1}+\overline w_{q+1},\theta_{q+1})$ as a sum of highly oscillatory  Mikado flows such that  the support of $ (w_{q+1},\theta_{q+1})$ and that of $\overline w_{q+1}$ are disjoint  by choosing disjoint building blocks.
 Here  the perturbation $(w_{q+1},\theta_{q+1})$ is used to  cancel the stress term $ M_q$. As we can see, it produces a new  term $w_{q+1}\otimes w_{q+1}$ in the Euler equation.
Then we construct the perturbation $\overline w_{q+1}$ to cancel the 
stress term $\mathring{ R}_q$ and the low frequency  part, which comes from $w_{q+1}\otimes w_{q+1}$.
More precisely, 
\begin{align}
  w_{q+1}\theta_{q+1}&\sim  M_{q}+({\rm high\ frequency\ error}),\notag\\   \overline w_{q+1}\otimes  \overline w_{q+1}&\sim -\mathring{{R}}_{q}-\int_{\mathbb{T}^3}w_{q+1}\otimes w_{q+1} \dif x+({\rm high\ frequency\ error}). \label{eq:intro:ww}
\end{align}

Once we have the above relations,  since the support of $ (w_{q+1},\theta_{q+1})$ and that of $\overline w_{q+1}$ are disjoint, it follows that
\begin{align*}
 ( w_{q+1}+\overline w_{q+1})\theta_{q+1}&\sim  M_{q}+({\rm high\ frequency\ error}),\\   ( w_{q+1}+\overline w_{q+1})\otimes ( w_{q+1}+\overline w_{q+1})&\sim -\mathring{R}_{q}+ \mathbb{P}_{\neq0}(w_{q+1}\otimes w_{q+1} )+({\rm high\ frequency\ error}).\ \ 
\end{align*}
Here we denote $\mathbb{P}_{\neq0} f := f -\fint f\dif x$.
So, the principle part of  the oscillation errors has  been canceled, while  the high frequency errors and the high frequency part of the product $w_{q+1}\otimes w_{q+1} $ are small in $C^0$-norms. 

To achieve the  regime $\beta+2\tilde{\beta}<1$ and the Onsager regime up to  $1/3$ at the same time, we need to choose the parameters  carefully. Let us heuristically show the typical  estimates of the Nash errors $\div^{-1}((w_{q+1}+\overline w_{q+1})\cdot  \nabla    v_{q})$ and  $\div^{-1}((w_{q+1}+\overline w_{q+1})\cdot  \nabla    \rho_{q})$ from  the Euler equations and  the transport equation respectively. 
We assume that the frequencies grow hypergeometrically $\lambda_q=\lambda_{q-1}^b$ for some  $b>1$, and the perturbations obey the following  H\"older scaling:
$$\|w_{q+1}\|_{C^0}+\|\overline w_{q+1}\|_{C^0}\leq \lambda_{q+1}^{-\beta},\ \ \|\theta_{q+1}\|_{C^0}\leq \lambda_{q+1}^{-\tilde{\beta}},$$
for $\beta,\tilde{\beta}>0$.
In view of \eqref{eq:intro:ww}, we need to ensure
\begin{align*}
   \| \mathring{R}_{q}\|_{C^0}\leq  \lambda_{q+1}^{-2\beta},\ \  \| M_{q}\|_{C^0}\leq  \lambda_{q+1}^{-\beta-\tilde{\beta}}.
\end{align*}
Then we need to verify   the above bounds at the level $q+1$. As usual   in the convex integration method, the  iterate $(v_q,\rho_q)$ is a  sum of building blocks with frequency   not bigger  then $\lambda_q$. Then by choosing $\lambda_{q+1}\gg\lambda_q$, 
 we notice that the inverse divergence $\div^{-1}$ will give us a factor $\lambda_{q+1}^{-1}$, while the gradient will give us a factor  $\lambda_{q}$.   Then two Nash errors   can be estimated by
\begin{align*}
\| \div^{-1}( (w_{q+1}+\overline w_{q+1})\cdot  \nabla    v_{q})\|_{C^0}&\les
  \frac{\|w_{q+1}+\overline w_{q+1}\|_{C^0}\| v_{q}\|_{C^1}}{\lambda_{q+1}}\les\frac{\lambda_{q+1}^{-\beta}\lambda_{q}^{-\beta}\lambda_q}{\lambda_{q+1}}  \leq \lambda_{q+2}^{-2\beta},\\
 \|  \div^{-1}((w_{q+1}+\overline w_{q+1})\cdot  \nabla   \rho_{q})\|_{C^0}&\les
  \frac{\|w_{q+1}+\overline w_{q+1}\|_{C^0}\| \rho_{q}\|_{C^1}}{\lambda_{q+1}}\les\frac{\lambda_{q+1}^{-\beta}\lambda_{q}^{-\tilde{\beta}}\lambda_q}{\lambda_{q+1}}  \leq \lambda_{q+2}^{-\beta-\tilde{\beta} }.
\end{align*}
To ensure the  validity of the iteration, we need 
\begin{align*}
    1-b-\beta-\beta b<-2\beta b^2,\ 1-b-\tilde{\beta}-\beta b<-\beta b^2-\tilde{\beta} b^2,
\end{align*}
i.e.
\begin{align*}
(b-1)(\beta(2b+1)-1)<0,\ (b-1)(\tilde\beta(b+1)+\beta b-1)<0,    
\end{align*}
which is satisfied by choosing $b>1$ close to 1, 
$\beta<\frac13$ and $\beta+2\tilde{\beta}<1$.

Moreover, before the perturbation step,  we  apply the gluing step introduced by Isett in \cite{Ise18} to achieve the above regularities. Specifically, we combine exact solutions to the Euler equations, following a method similar to that in \cite{BDLSV19}. Correspondingly, we also glue together exact solutions to the transport equation, ensuring that the associated  errors $\overline M_q$ are supported in disjoint temporal intervals. Here we need to estimate the difference between the glued solution  $\rho_i$ and the original solution to the transport equation  $ \rho_q$ in certain negative order spaces. Since $\rho_i$ and $ \rho_q$  are both  scalars, we can not apply the Biot-Savart law as done in \cite{BDLSV19}. Instead, we utilize the inverse divergence operator to define $y_i=\div^{-1}\rho_i,y_q=\div^{-1}y_q$, which satisfies an equation of the form:
\begin{align*}
    \div[(\partial_t+v_l\cdot\nabla)(y_i-y_q)]=\div[\cdot\cdot\cdot].
\end{align*}
Since there is no suitable left-inverse of the $\div$ operator, we use the identity $\nabla\div=\Delta+\curl\curl$ and the fact that $\curl(y_i-y_q)=0$ to derive that 
\begin{align*}
    (\partial_t+v_l\cdot\nabla)(y_i-y_q)=\Delta^{-1}\nabla\div[\cdot\cdot\cdot]-\Delta^{-1}\curl\curl [v_l\cdot\nabla(y_i-y_q)].
\end{align*}
By combining various analytic identities (see Lemma \ref{lem;identity} below),  we obtain that $y_i- y_q$ is bounded by  estimates on the transport equation together with Gronwall's inequality.  We refer to Section \ref{sec:gluing:rho} for more details.

\subsubsection{Ideas in the $L^1_tW^{1,s}$-scales}\label{sec:idea2}
When we consider the $L^1_tW^{1,s}$-scales or $L^r_tL^p$-scales for SDEs, to prove Theorem \ref{thm:non_pde_23},  we once again apply convex integration. At  each step $q\in\mN_0$, we need to deal with a similar system as \eqref{eq:idea:rhov} with an extra $\Delta$ in both equations. 
As in the previous analysis, to eliminate the stress terms, we construct new perturbations of the form $(w_{q+1}+\overline w_{q+1},\theta_{q+1})$, where $(w_{q+1},\theta_{q+1})$ is to cancel the  stress term $ M_q$, while $\overline w_{q+1}$ is to cancel the stress term $\mathring{ R}_q$. 
To address the dissipative Laplacian term, we enhance more intermittency in the building blocks  by introducing generalized intermittent space-time jets, which are inspired by \cite[Section 3]{LZ23}, \cite[Section 4]{BCDL21} for the spatial direction and by \cite[Section 4.2]{CL21} for the temporal direction.

However, because of the existence of extra intermittency, the previous
method breaks down. Roughly speaking, we aim to construct a perturbation $w_{q+1}=\sum_\xi a_{(\xi)}W_{(\xi)},\theta_{q+1}=\sum_\xi a'_{(\xi)}\Theta_{(\xi)}$ with a large oscillation parameter $\mu$. On the one hand,   similar to  the procedure in  \cite[Section 4.3]{BVb}, to cancel the additional oscillation errors in  the Navier-Stokes equations arising from
 \begin{align*}
 \div(\mathbb{P}_{\neq0}(w_{q+1}\otimes w_{q+1} ))
 \sim \frac1{\mu}\partial_t(\sum_{\xi}a_{(\xi)}^2|W_{(\xi)}|^2\xi) +({\rm high\ frequency\ error}),
\end{align*}
 we need to introduce a temporal corrector and choose  $\mu> {r}_\parallel^{-\frac12} {r}_\perp^{-\frac{d-1}{2}}$
to control this corrector. However, when dealing with the transport  equation, the intermittent building blocks should satisfy the equation
\begin{align*}
     \partial_t\Theta_{(\xi)}+\mu {r}_\parallel^{\frac12} {r}_\perp^{\frac{d-1}{2}}\div (W_{(\xi)}\Theta_{(\xi)})=0
 \end{align*}
as seen in \eqref{eq:ptthe+}. This condition requires $\mu={r}_\parallel^{-\frac12} {r}_\perp^{-\frac{d-1}{2}}$, which leads to a contradiction.

To address this problem, at the beginning of the iteration, we iterate two equations at different scales. Heuristically speaking, we assume that the frequencies grow hypergeometrically $\lambda_q=\lambda_{q-1}^b$ for some  $b$ large enough, and  that the perturbations obey the following $L^2$ scaling:
$$\|\overline w_{q+1}\|_{L^2}\leq \lambda_{q+1}^{-\beta},\ \ \|w_{q+1}\|_{L^2}+\|\theta_{q+1}\|_{L^2}\leq \lambda_{q+2}^{-{\beta}},$$
and correspondingly
\begin{align*}
   \| \mathring{R}_{q}\|_{L^1}\leq  \lambda_{q+1}^{-2\beta},\ \  \| M_{q}\|_{L^1}\leq  \lambda_{q+2}^{-2\beta}.
\end{align*}
Then the undesired product $w_{q+1}\otimes w_{q+1}$  automatically satisfies the desired estimates for $\mathring{R}_{q+1}$ since $$\|w_{q+1}\otimes w_{q+1}\|_{L^1}\lesssim \|w_{q+1}\|_{L^2}^2\les \lambda_{q+2}^{-2\beta}.$$

Finally, we emphasize that although we obtain a spatial range similar to \cite{BCDL21}, the building blocks and the corresponding estimates in our work are different. Specifically, in \cite[Section 4]{BCDL21}, the authors constructed building blocks satisfying a certain $L^p$-normalization property with  some $p>0$, and derived solutions in the space $C_t(W^{1,d-}\cap L^{\infty-})$. In contrast  to that, in the present paper, our construction of building blocks is limited to the $L^2$-normalization property, which is crucial for  solutions to the  Navier-Stokes equations. By introducing temporally intermittent jets and  a careful choice of parameters, we ultimately achieve a spatial range similar to \cite[Section 4]{BCDL21}, at the cost  that the time integrability of  the solutions is only $L^1$.

\subsection{ Organization of the paper.}
This paper is organized as follows.  First, in Section \ref{sec:Euler:C1/3}, we establish our first main PDE result, Theorem \ref{thm:con:1/3}, via the convex integration method. Subsequently, by applying the superposition principle, we derive Corollary \ref{thm:nonode1/3}. 
The   implementation of the main convex integration procedure is presented in Section \ref{sec:proof13}.  Then, Section \ref{cogpss2} is devoted to the main results  in  $L^1_tW^{1,s}$-based  scales: Theorem \ref{thm:nonuniode+sde} and  the two Corollaries \ref{thm:nonlaw} and \ref{thm:nonuniode}. We also work at the PDE level by demonstrating non-uniqueness  for Fokker-Planck equations advected by solutions to the Navier-Stokes or Euler equations, as stated in Theorem \ref{thm:non_pde_23}. The  implementation of the main convex integration procedure is presented in Section \ref{proof:prop2}. 
In Appendix \ref{app:b}, we collect some technical tools used in the construction.
  In Appendix \ref{s:appA.3} and Appendix \ref{s:appA.5}  we  provide the building blocks and some auxiliary estimates used in the constructions  in Section \ref{sec:proof13} and Section \ref{proof:prop2} respectively.

\noindent{\bf Notations.} Let $T>0$, $\mN_{0}:=\mN\cup \{0\}$. Throughout the manuscript, we write  $\mT^d = \mR^d/\mZ^d$ for the $d$-dimensional flat torus. 
We define the natural projection ${\rm Pr}:\mR^d\to\mT^d$ by ${\rm Pr}(x)=x-[x]$, where $[x]$ is the integer part of $x$,  which is continuous. 
A $\mathbb{T}^d$-valued Brownian motion is seen as 
the natural projection of  $\mathbb{R}^d$-valued Brownian motion onto  $\mathbb{T}^d$.
We refer to \cite{Hsu02} for a comprehensive treatment of the Brownian motion  on general manifolds.
We use the following notations.
\begin{itemize}
    \item 
We employ the notation $a\lesssim b$ if there exists a constant $c>0$ such that $a\leq cb$. 
\item   Given a Banach space $E$ with a norm $\|\cdot\|_E$, we write $C_tE=C([0,T];E)$ for the space of continuous functions from $[0,T]$ to $E$, equipped with the supremum norm. For $p\in [1,\infty]$ we write $L^p_tE=L^p([0,T];E)$ for the space of $L^p$-integrable functions from $[0,T]$ to $E$, equipped with the usual $L^p$-norm.

\item  
For $\alpha\in(0,1)$ we  define $C^\alpha_tE$ as the space of $\alpha$-H\"{o}lder continuous functions from $[0,T]$ to $E$, endowed with the norm $\|f\|_{C^\alpha_tE}=\sup_{s,t\in[0,T],s\neq t}\frac{\|f(s)-f(t)\|_E}{|t-s|^\alpha}+\sup_{t\in[0,T]}\|f(t)\|_{E},$   and write  $C_t^\alpha$ in the case when $E=\mathbb{R}$. 
 
\item  We use $L^p$ to denote the set of  standard $L^p$-integrable functions on $\mathbb{T}^d$.
For $s>0$, $p>1$ we set $W^{s,p}:=\{f\in L^p; \|(I-\Delta)^{\frac{s}{2}}f\|_{L^p}<\infty\}$ with the norm  $\|f\|_{W^{s,p}}=\|(I-\Delta)^{\frac{s}{2}}f\|_{L^p}$.

\item For $N\in \N_0 $, $C^N$ denotes the space of $N$-times differentiable functions equipped with the norm
	$$
	\|f\|_{C^N}:=\sum_{\substack{|\alpha|\leq N, \alpha\in\N^{d}_{0} }}\| D^\alpha f\|_{ L^\infty_x}.$$ Similarly, if the norm is taken in space-time, we use $C^{N}_{t,x}$.
	For $N\in \mathbb{N}_0$ and $\kappa \in (0,1)$, $C^{N+\kappa}$ denotes the subspace of $C^N$ whose $N$-th derivatives are $\kappa$-H\"{o}lder continuous, with the norm
	\begin{align*}
		\|f\|_{C^{N+\kappa}}:=	\|f\|_{C^N}+\sum_{\substack{|\alpha|= N,  \alpha\in\N^{d}_{0} }}[D^\alpha f]_{C^\kappa},
	\end{align*}
	where  $	[f]_{C^\kappa_{x}} :=\sup_{\substack{x\neq y, x,y\in \mathbb{T}^d }} \frac{ \left|  f(x)- f(y)\right| }{|x-y|^\kappa}$ is the H\"{o}lder seminorm.

\item 
We define  the projections $\mathbb{P}_{=0} f := \fint_{\mathbb{T}^d} f\dif x$, and  $\mathbb{P}_{\neq0} f := f -\fint_{\mathbb{T}^d} f\dif x$.
\item  For a matrix $R$, we denote its traceless part by $\mathring{R} := R- \frac{1}{d} \tr (R )\Id$.
\item 
We denote  the Lebesgue measures on $\mT^d$  by $\cL^d$.
\end{itemize}

\section{Construction of non-unique solutions in $C_{t,x}^0$ scales}
\label{sec:Euler:C1/3}
In this section, we    work on the PDE level to  demonstrate the non-uniqueness of solutions to the transport equations advected by Euler flows, as stated in Theorem \ref{thm:con:1/3} (2). Specifically, for any fixed $\beta,\tilde{\beta}>0$  satisfying $0<\beta+2\tilde{\beta}<1,0<\beta<\frac13$, we apply the convex integration method simultaneously to the Euler equations \eqref{eq:euler} and the transport equations \eqref{eq:tpe}. We construct a solution to the 3D Euler equations in the space $C^\beta([0,T]\times\mT^d)$, which satisfies $e(t)=\|v(t)\|_{L^2}^2$ for a prescribed energy profile, such that  the associated transport equation admits a non-constant, positive solution in $C^{\tilde{\beta}}([0,T]\times\mT^d)$  and  with
initial data $\rho_0 = 1$.
By selecting  $\beta$  sufficiently close to $1/3$ and choosing a decreasing energy profile $e(t)$, we  establish  Corollary \ref{thm:nonode1/3}, with the help of the superposition principle.

Without loss of generality, we assume $T=1$. When considering the $C_{t,x}^0$ scales, we will primarily focus on spaces equipped with the supremum norm. In this section and  the subsequent Section \ref{sec:proof13}, we use the notation $\|\cdot\|_{B}=\|\cdot\|_{C_tB}$ for any Banach space $B$.

The convex integration iteration is indexed
by a parameter $q \in \mathbb{N}_0$. We define the frequency parameter $\{\lambda_q\}_{q\in\mathbb{N}_0}\subset\mathbb{N}$ which diverges to $\infty$, and the amplitude parameters $\{\delta_q,\tilde{\delta}_q\}_{q\in\mathbb{N}_0}\subset (0,1]$
 which are decreasing to 0 by
\begin{align*}
\lambda_q&=   a^{(b^q)}, \  \ \delta_q=\lambda_q^{-2\beta},\ \ 
\tilde{\delta}_q=\lambda_q^{-2\tilde{\beta}},\ q\geq0,
\end{align*}
where $a>1$ is a  large parameter and $b>1$ is close to $1$.  Here we recall that $\beta,\tilde{\beta}>0$ are given in Theorem \ref{thm:con:1/3} satisfying $0<\beta<1/3$ and $0<\beta+2\tilde{\beta}<1$. In the following, without loss of generality, we additionally assume  $\beta\leq \tilde{\beta}$.
 In addition, we have
\begin{align}\label{ieq:ab2}
\sum_{q\geq1} \tilde{\delta}_q^{1/2}\lesssim  \frac1{a^{b\tilde{\beta}}-1}< \frac1{3}
\end{align}
by choosing $a$ large enough in terms of $b$ and $\tilde{\beta}$.

At each step $q\in\mathbb{N}_0$, a pair $(v_q,\rho_q, \mathring{R}_q,M_q)$ is constructed solving the following systems on $\mT^d$:
\begin{align}
\partial_t \rho_q+\div(v_q \rho_q)&=-\div M_q,\label{e:euler_trans}\\
\partial_tv_q+\div(v_q\otimes v_q)+\nabla\pi_{q}
&=\div\mathring{R}_q,\ \
\div v_q=0,\label{e:euler_reynolds}
\end{align}
where $\mathring{R}_q$ is assumed as a trace-free symmetric matrix, and $M_q$ is a vector field.

To take into account the initial condition, we require that $\rho_q=1$ near the origin, more precisely,  that $\rho_q=1$ on $[0,T_q]$, where 
\begin{align}
    T_q := \frac13 - \sum_{ 1\leq r\leq q} \tilde{\delta}_r^{1/2}. \label{def:tq}
\end{align} 
Applying  \eqref{ieq:ab2} we obtain $0< T_q\leq\frac13$. Here we define $\sum_{1\leq r\leq 0}:=0.$ 

Our main iteration on the approximate solution $(v_q,\rho_q, \mathring{R}_q,M_q)$  reads as follows:
\begin{proposition}\label{prop:1/3}
 Let
 $\beta\in(0,\frac13),\tilde{\beta}\in[\beta,\frac{1-\beta}{2})$ and  $1<b<\frac{1-\tilde{\beta}}{\beta+\tilde{\beta}}.$
  Let $e:[0,1]\to \mR$ be a strictly positive function satisfying $|e'(t)|\leq1$. Then there exists a choice of parameters $\alpha\in(0,1)$ and $a>1$ such
that the following holds true: 
 let $(v_q,\rho_q, \mathring{R}_q, M_q)$ be a solution to \eqref{e:euler_trans}-\eqref{e:euler_reynolds} satisfying   $\int\rho_q\dif x=1$,
\begin{align}
\norm{v_q}_{C^0} & \leq M\sum_{i=0}^q\delta_i^{1/2},\ \ \norm{\rho_q}_{C^0} \leq 2+\sum_{i=0}^{q}{\tilde{\delta}}_{i}^{1/2},\label{e:v_q_0}\\
\norm{v_q}_{C^1}&\leq M \delta_q^{1/2}\lambda_q,\ \ \norm{\rho_q}_{C^1}\leq  {\tilde{\delta}}_{q}^{1/2}\lambda_q,\label{e:v_q_inductive_est}\\
\|\mathring R_q\|_{C^0}&\leq  \delta_{q+1}\lambda_q^{-3\alpha},\ \ \norm{M_q}_{C^0}+\frac{1}{\lambda_{q}}\norm{M_q}_{C^1}\leq  {\delta}_{q+1}^{1/2}\tilde{\delta}_{q+1}^{1/2}\lambda_q^{-3\alpha}\label{e:R_q_inductive_est},\\
&\rho_q-1=M_q=0\ {\rm on}\ [0,T_q],\label{bd:mq=rhoq=0}
\end{align}
where $M>1$ is a universal geometric constant.  Moreover, for any $t\in[0,1]$
		\begin{align}\label{estimate:energy}
			\delta_{q+1} \lambda_{q}^{-\alpha/3}\leq e(t) - \| v_q(t)\|_{L^2}^2  \leq \delta_{q+1}.
		\end{align}
Then there exists 
$(v_{q+1},\rho_{q+1}, \mathring{R}_{q+1}, M_{q+1})$ which solves \eqref{e:euler_trans}-\eqref{e:euler_reynolds}, satisfies \eqref{e:v_q_0}-\eqref{estimate:energy} at the level $q+1$ and 
 \begin{align}
\norm{v_{q+1}-v_q}_{C^0}
\leq M\delta_{q+1}^{1/2}\label{e:v_diff_prop_est},\ \
\norm{\rho_{q+1}-\rho_q}_{C^0}
\leq \tilde{\delta}_{q+1}^{1/2}.
\end{align}
\end{proposition}
As noted in \cite{BDLSV19}, the wiggle room provided by the factor $\lambda_q^{-3\alpha}$ is beneficial during the gluing step. Additionally, the extra bound on $M_q$ in the $C^1$ norm is utilized to establish the time regularity of $\rho$, see the proof of Theorem \ref{thm:con:1/3} below.

\begin{proof}[Proof of Theorem \ref{thm:con:1/3}]

We start with the proof of the first part  of  the theorem by in  a similar way as \cite{CET94}. Specifically, we   prove that  for $\beta+2\tilde{\beta}>1$ and vector field $v\in C_{t,x}^\beta$, any weak solution $\rho\in C_{t,x}^{\tilde{\beta}}$ to the transport equation \eqref{eq:tpe} conserves energy.  Then the uniqueness follows from the linearity.

  Let $\rho_l:=\rho*\varphi_l$ be  the spatial mollification of $\rho$ with a length scale $l$. Similarly, we define $\rho_l$ and $(v\rho)_l$.
 Then $\rho_l$ satisfies for any $t\in[0,1]$
 \begin{align*}
     \int_{\mathbb{T}^3}|\rho_l(x,t)|^2\dif x- \int_{\mathbb{T}^3}|\rho_l(x,0)|^2\dif x=2\int_0^t \langle (v\rho)_l,\nabla\rho_l\rangle_{L^2}\dif s.
 \end{align*}
   Since \begin{align*}
      \langle v_l\rho_l,\nabla\rho_l\rangle _{L^2}  \equiv0,
   \end{align*}
we have
\begin{align*}
     \int_{\mathbb{T}^3}|\rho_l(x,t)|^2 \dif x- \int_{\mathbb{T}^3}|\rho_l(x,0)|^2 \dif x=2\int_0^t\langle (v\rho)_l-v_l\rho_l,\nabla\rho_l\rangle _{L^2}\dif s.
\end{align*}
Using the commutator estimate Lemma \ref{p:CET} we have
\begin{align*}
   \left| \|\rho_l(t)\|_{L^2}^2- \|\rho_l(0)\|_{L^2}^2\right|  \lesssim  \|(v\rho)_l-v_l\rho_l\|_{C_{t,x}^0}\|\rho_l\|_{C_{t,x}^1}\lesssim l^{\beta+2\tilde{\beta}-1}\|v\|_{C_{t,x}^\beta}\|\rho\|_{C_{t,x}^{\tilde{\beta}}}^2.
\end{align*}
   Thus, as $\beta+2\tilde{\beta}> 1$, the right hand side converges to zero as $l\to0$. We conclude the proof with the observation that $\rho_l$ converges to $\rho$ in $C_tL^2$.

For the second part, since $0<\beta<\frac13$, it suffices  to prove the statement in the case where $\beta\leq\tilde{\beta}<1-2\beta$ and $T=1$.  Without loss of generality, by  the same argument as in \cite[Proof of Theorem 1.1]{BDLSV19}, we may assume that 
\begin{align*}
    \inf_{t\in[0,1]}e(t)\geq \delta_1\lambda_0^{-\alpha/3},\ \ \sup_{t\in[0,1]}e(t)\leq \delta_1,\ \ \sup_{t\in[0,1]}e'(t)\leq1.
\end{align*}

By choosing $\alpha>0$ small enough, there exists $\gamma\in\mathbb{Q}$ such that $3\alpha+(\tilde{\beta}+\beta)b<\gamma<1-\tilde{\beta}$. Then we define $\lambda=\lambda_0^{\gamma}\in2\mN$ by choosing a suitable $a$ large enough, and  start the iteration from 
$$\rho_0(t,x)=1+\frac{\sin\lambda \pi x_1}2\chi_0(t),\  v_0=0,\ \mathring{R}_0=0,\ M_0(t,x)=\partial_t\chi_0(t)\frac{\cos \lambda \pi x_1}{2\lambda\pi}(1,0,0)^T,$$
where $x=(x_1,x_2,x_3)$ and $\chi_0$ is a smooth function with $\chi_0(t)=0$ on $[0,\frac13]$, $\chi_0(t)=1$ on $[\frac23,1]$.

 Then  $\int\rho_0\dif x=1$, and by choosing $a$ large enough to absorb the constant, we have
\begin{align*}
   \| \rho_0\|_{C^0}\leq 2,\ \  \| \rho_0\|_{C^1}\lesssim\lambda\leq \tilde{\delta}_0^{1/2}\lambda_0,\ \ \|M_0\|_{C^0}+\frac{1}{\lambda}\norm{M_0}_{C^1}\lesssim \frac{1}{\lambda}\leq{\delta}_{1}^{1/2}\tilde{\delta}_{1}^{1/2}\lambda_0^{-3\alpha},
\end{align*}
which implies that \eqref{e:v_q_0}-\eqref{bd:mq=rhoq=0}  hold at the level  $q=0$. \eqref{estimate:energy} is automatically satisfied by the assumption on $e(t)$ and the fact that $v_0=0$.

Applying Proposition \ref{prop:1/3} iteratively, we obtain a sequence of fields $(v_q,\rho_q,\mathring{R}_q,M_q)$  satisfying \eqref{e:v_q_0}-\eqref{estimate:energy} and converging  to a weak solution $(v,\rho)$ to \eqref{eq:euler}, \eqref{eq:tpe}  by \eqref{e:R_q_inductive_est}. Using the estimate \eqref{e:v_diff_prop_est} yields for any $\beta'<\beta$:
\begin{align*}
\sum_{q=0}^{\infty} \norm{v_{q+1}-v_q}_{C^{\beta'}} &
\lesssim \; \sum_{q=0}^{\infty} \norm{v_{q+1}-v_q}_{C^0}^{1-\beta'}\norm{v_{q+1}-v_q}_{C^1}^{\beta'}
\lesssim \; \sum_{q=0}^{\infty} \delta_{q+1}^{\frac{1-\beta'}2}\left(\delta_{q+1}^{1/2}\lambda_{q+1}\right)^{\beta'}
\lesssim  \sum_{q=0}^{\infty} \lambda_{q+1}^{\beta'-\beta}.
\end{align*}
Hence we obtain that $v\in C^0_tC^{\beta'}_x$. By the same argument and \eqref{e:v_diff_prop_est}  we obtain that $\rho\in C^0_tC^{{\beta}''}_x$ for any ${\beta}''<\tilde{\beta}$. 
The time regularity of $v$ is obtained by the same argument as in \cite{Ise18} or \cite{BDLSV19}. We  thus obtain $v\in C_xC_t^{\beta'}$ and then $v\in C^{\beta'}([0,1]\times \mathbb{T}^3)$ for arbitrary $\beta'<\beta$.

For the  time regularity of $\rho$, using \eqref{e:v_q_inductive_est}, \eqref{e:R_q_inductive_est} above, we immediately have  for $q\in\mN_0$
\begin{align*}
  \|  \partial_t \rho_q\|_{C^0}\les\|v_q \|_{C^0}\|\nabla\rho_q\|_{C^0}+\|\div M_q\|_{C^0}\lesssim {\tilde{\delta}}_{q}^{1/2}\lambda_q+ \lambda_q{\delta}_{q+1}^{1/2}\tilde{\delta}_{q+1}^{1/2}\lesssim {\tilde{\delta}}_{q}^{1/2}\lambda_q,
\end{align*}
which together with interpolation and \eqref{e:v_diff_prop_est} implies that
\begin{align*}
 \sum_{q=0}^{\infty}  \|\rho_{q+1}-\rho_q\|_{C_x^0C_t^{\beta''}}&\les \sum_{q=0}^{\infty} \| \rho_q-\rho_{q+1}\|_{C^0}^{1-\beta''} \|\partial_t{\rho}_{q}-\partial_t{\rho}_{q+1}\|_{C^0}^{\beta''}\\
  &\les\sum_{q=0}^{\infty} \tilde{\delta}_{q+1}^{\frac{1-{\beta}''}2}\left(\tilde{\delta}_{q+1}^{1/2}\lambda_{q+1}\right)^{{\beta}''} \lesssim  \sum_{q=0}^{\infty} \lambda_{q+1}^{{\beta}''-\tilde{\beta}}.
\end{align*}
 As a consequence,  $\rho_q\to \rho\in  C_x^0C_t^{\beta''}$. Thus we obtain $\rho\in  C^{\beta''}([0,1]\times \mathbb{T}^3)$ for any $\beta''<\tilde{\beta}$.

Moreover,  \eqref{bd:mq=rhoq=0} ensures that $\rho(t) \equiv 1$ for every $t$ sufficiently close to 0,  and
\eqref{estimate:energy} ensures that $e(t)=\|v(t)\|_{L^2}^2$ for any $t\in[0,1]$.

Clearly by \eqref{ieq:ab2} and  \eqref{e:v_diff_prop_est} we have
\begin{align*}
    \|\rho-\rho_{0}\|_{C^0}\leq \sum_{q=0}^\infty \|\rho_{q+1}- \rho_q\|_{C^0}\leq \sum_{q=0}^\infty \tilde{\delta}_{q+1}^{1/2}\leq\frac1{3}, 
\end{align*}
  which implies that 
$\rho$ is nonnegative on $\mathbb{T}^3$:
\begin{align*}
\inf_{t\in [0,1]}\rho \geq \inf_{t\in [0,1]} \rho_0 -
 \|\rho-\rho_{0}\|_{C^0} \geq\frac12-\frac1{3}>0,
\end{align*}
and  $\rho$ does not coincide with the solution which is constantly 1, as
\begin{align*}
  \|\rho - 1\|_{C^0}\geq \|1 - \rho_0\|_{ C^0} -
  \|\rho-\rho_{0}\|_{C^0}\geq\frac12-\frac1{3}>0.
\end{align*}
\end{proof}

\begin{remark}
   If one applies the same argument as in \cite[(2.16)]{BDLSV19} to derive the time regularity for $\rho$, one can only achieve the time regularity up to $\tilde{\beta}/(1+\tilde{\beta}-\beta) < \tilde{\beta}$, due to the low regularity of $v$.
\end{remark}

 In Theorem \ref{thm:con:1/3} we choose $\beta<\frac13$ close enough to $\frac13$  and choose $e(t)$ to be decreasing. Then the proof of Corollary \ref{thm:nonode1/3} can be established using an argument analogous to  \cite[Theorem 1.3]{BCDL21} (see also the proof of Theorem \ref{thm:nonuniode+sde} detailed in Section \ref{cogpss2}).

\section{Proof of Proposition ~\ref{prop:1/3}}\label{sec:proof13}
The proof follows a series of main steps. First, we fix some necessary parameters and proceed with a mollification step in Section~\ref{sec:mpll}. Next, in Section~\ref{sec:glue}, we apply the gluing procedure for both $(v_{q},\rho_{q})$. The gluing step for the Euler equation is similar to the approach in \cite{BDLSV19}, but with a distinct gluing parameter $\tau_q$. Additionally, we introduce the gluing step for the transport equation to ensure that the associated error term is supported on disjoint temporal intervals.
In Section~\ref{sec:defq+1}, we define the new iteration $(v_{q+1},\rho_{q+1})$ and provide inductive estimates. Specifically, we first construct the perturbation $(w_{q+1},\theta_{q+1})$ for the transport equation to cancel the glued stress term $\overline M_q$. Subsequently, we construct $\overline w_{q+1}$ for the Euler equation to cancel the glued stress term $\mathring{ \overline R}_q$ as well as the low frequency part term of $w_{q+1}\otimes w_{q+1}$.
In Section~\ref{sec:defrmq+1}, we define the new stress terms $(\mathring{R}_{q+1},M_{q+1})$ and establish the required estimates. Finally, in Section~\ref{estimate_energy}, we derive the energy estimates, completing the proof.

\subsection{Choice of parameters and mollification}\label{sec:mpll}
In the sequel, additional parameters will be indispensable and their value has to be carefully chosen in order to respect all the compatibility conditions appearing in the estimates below. First, for any fixed $\beta\in(0,\frac13), \tilde{\beta}\in[\beta,\frac{1-\beta}{2})$, $\tilde{\beta}(b+1)+\beta b<1$, and  $\alpha>0$ small enough we have
\begin{align}
    \frac{\lambda_q\tilde{\delta}_q^{1/2}\delta_{q+1}}{\tilde{\delta}_{q+1}^{1/2}\lambda_{q+1}}\leq \frac{\delta_{q+2}}{\lambda_{q+1}^{8\alpha}},\ \ \frac{\lambda_q\tilde{\delta}_q^{1/2}\delta_{q+1}^{1/2}}{\lambda_{q+1}}\leq \frac{\delta_{q+2}^{1/2}\tilde{\delta}_{q+2}^{1/2}}{\lambda_{q+1}^{8\alpha}}.\label{bd:para1/3}
\end{align}
To see this, we take logarithms and obtain 
\begin{align}
    1-b-\tilde{\beta}-2\beta b+\tilde{\beta} b<-2\beta b^2-8\alpha b,\ \ 1-b-\tilde{\beta}-\beta b<-\beta b^2-\tilde{\beta}b^2-8\alpha b,\notag
\end{align}
i.e.
\begin{align}
  8\alpha b+(b-1)(2\beta b+\tilde{\beta}-1)<0,\ \   8\alpha b+(b-1)(\tilde{\beta}(b+1)+\beta b-1)<0,\notag
\end{align}
which are satisfied by the choice of $b$ and by choosing $\alpha>0$ small enough.

 Finally, we increase $a$ such that  \eqref{ieq:ab2} holds. In the sequel, we increase $a$ to  absorb the various implicit and universal constants in the following estimates.

We take a sufficiently small $\alpha\in(0,1)$  such that \eqref{bd:para1/3} holds, and define $ l >0$ as
\begin{align}
\label{e:ell_def}
 l = \frac{\tilde{\delta}_{q+1}^{1/2}}{\tilde{\delta}_q^{1/2 }\lambda_q^{1+\frac{3\alpha}{2}}}\leq \frac{\delta_{q+1}^{1/2}}{{\delta}_q^{1/2 }\lambda_q^{1+\frac{3\alpha}{2}}}
 ,\ \ 
\end{align}
where we used  $\beta\leq \tilde{\beta}$ to deduce $ \tilde{\delta}_{q+1}^{1/2}\tilde{\delta}_q^{-1/2 }\leq \delta_{q+1}^{1/2}{\delta}_q^{-1/2 }$.

Then we replace the pair $(v_q,\rho_q,\RR_q,M_q)$ with a mollified pair $( v_{ l},\rho_l,\mathring{R}_{ l},M_l)$.
We define
\begin{align*}
 v_{ l}:=& v_q \ast_x \phi_{ l},\ \
\mathring{R}_{ l}:= \mathring R_q \ast_x  \phi_{ l}  -(v_q\mathring\otimes v_q)\ast_x \phi_{ l}  +  v_{  l}\mathring\otimes  v_{ l},\\
 \rho_{ l}:=& \rho_q \ast_x \phi_{ l},\ \
{M}_{ l}:= M_q \ast_x  \phi_{ l}  +(v_q\rho_q)\ast_x \phi_{ l}  -  v_{  l}  \rho_{ l},
\end{align*}
which obey \eqref{e:euler_trans} and \eqref{e:euler_reynolds}  for a suitable $\pi_l$. Here $\phi_l:=\frac{1}{l^3}\phi(\frac{\cdot}{l})$ is a family of standard radial mollifiers on $\mathbb{R}^3$. Since the mollification does not depend on time, we still have $\rho_l=1$  on $[0,T_q]$.  Since  $\rho_l=\rho_q=1$ on the interval $[0,T_q]$, we have $M_l=v_q\ast_x \phi_{ l}  -  v_{  l} =0$ on the interval $[0,T_q]$.

In order to bound $\mathring{R}_{ l}$ and ${M}_{ l}$, we use the basic mollification estimate, Lemma \ref{p:CET} and the fact that $\lambda_{q}^{-3/2}\leq l\leq \lambda_{q}^{-1}$ (by choosing $\tilde{\beta}(b-1)+\frac32\alpha<\frac12$) to  obtain for $N\in\mN_0$
\begin{align}
\norm{ v_{ l}-v_q}_{C^0}&\lesssim  l\|v_q\|_{C^1}\les l\lambda_q\delta_q^{1/2}\lesssim\delta_{q+1}^{1/2}\lambda_q^{-\alpha}\,,\label{e:v:ell:0}\ \
\norm{ v_{ l}}_{C^{N+1}} \lesssim   l^{-N}\|v_q\|_{C^1}\lesssim\delta_q^{1/2}\lambda_q l^{-N}, \\
\norm{ \rho_{ l}-\rho_q}_{C^0}&\lesssim  l\|\rho_q\|_{C^1}\les l\lambda_q\tilde{\delta}_q^{1/2}\lesssim\tilde{\delta}_{q+1}^{1/2}\lambda_q^{-\alpha}\,,\label{e:rho:ell:0}\ \
\norm{ \rho_{ l}}_{C^{N+1}} \lesssim   l^{-N}\|\rho_q\|_{C^1}\lesssim\tilde{\delta}_q^{1/2}\lambda_q l^{-N}, \\
\|\mathring{R}_{ l}\|_{C^{N+\alpha}}&\lesssim  l^{-N-\alpha}\| \mathring R_q\|_{C^0}+ l^{2-N-\alpha}\|v_q\|_{C^1}^2\lesssim  l^{-N-\alpha}\lambda_q^{-3\alpha}\delta_{q+1}\lesssim\delta_{q+1} l^{-N+\alpha}, \label{e:R:ell}\\
\|M_{ l}\|_{C^{N+\alpha}}&\lesssim  l^{-N-\alpha}\| M_q\|_{C^0}+ l^{2-N-\alpha}\|v_q\|_{C^1}\|\rho_q\|_{C^1}\lesssim  l^{-N-\alpha}\lambda_q^{-3\alpha}{\delta}_{q+1} ^{1/2}\tilde{\delta}_{q+1}^{1/2} \lesssim{\delta}_{q+1} ^{1/2}\tilde{\delta}_{q+1}^{1/2} l^{-N+\alpha}. \label{e:M:ell}
\end{align}

Similarly, we use the same computation as for \eqref{e:v:ell:0} to obtain
		\begin{align}\label{energy:vq-vl}
			\left| \int_{\mathbb{T}^3} |v_q|^2- |v_l|^2 \dif x \right| &=\left| \int_{\mathbb{T}^3} |v_q|^2 *\phi_{l}- |v_{l}|^2 \dif x \right| 
			\lesssim \left\| |v_q|^2 *\phi_{l}- |v_{l}|^2\right \|_{C^0} 
			\lesssim  l^2 \|v_q\|_{C^1}^2\lesssim  \delta_{q+1}l^{\alpha}.
		\end{align}

\subsection{The gluing procedure}\label{sec:glue}
To achieve the regime $\beta+2\tilde{\beta}<1$ and  the Onsager regime up to $1/3$, we employ the gluing technique introduced by Isett in \cite{Ise18}. In Section \ref{sec:gluing:v}, we glue together the exact solutions to the Euler equations similarly  to \cite{BDLSV19}. This ensures that the associated error term $\mathring{\overline R}_q$ is supported on disjoint temporal intervals. Subsequently, in Section \ref{sec:gluing:rho}, we glue together the exact solutions of the transport equations, ensuring that the corresponding error term $\overline M_q$ is also localized on pairwise disjoint temporal intervals. To this end, we define the length of the time intervals  as
$$
 \tau_q:=  \frac{l^{2\alpha}\tilde{\delta}_{q+1}^{1/2}}{\tilde{\delta}_q^{1/2}\lambda_q{\delta}_{q+1}^{1/2}}.
 $$

Then we define 
\begin{align}
 t_i = i  \tau_q, \qquad  I_i = [ t_i + \frac{ \tau_q}{3},  t_i + \frac{2 \tau_q}{3}] \cap [0,1], \qquad  J_i = ( t_i - \frac{ \tau_q}{3},  t_i+ \frac{ \tau_q}{3}) \cap [0,1].
\label{eq:ti:Ii}
\end{align}

\subsubsection{Gluing step for $(v_l,\mathring{R}_l)$}\label{sec:gluing:v} 
The gluing step for the Euler equations in \eqref{eq:euler} is analogous to that in \cite{BDLSV19}, with the exception of the choice of the parameters  $\tau_q$ and $l$. In this section, we will only utilize the facts that 
\begin{align}
     \tau_q\leq  \frac{l^{2\alpha}}{\delta_q^{1/2}\lambda_q},\ \ \tau_q\delta_{q+1}^{1/2} l^{-1-\alpha/2}= (l\lambda_q)^{\frac{3\alpha}{2}}\leq 1,\label{papa:gluing}
\end{align}
 which hold by using again that $ \tilde{\delta}_{q+1}^{1/2}\tilde{\delta}_q^{-1/2 }\leq \delta_{q+1}^{1/2}{\delta}_q^{-1/2 }$ because $\beta\leq \tilde{\beta}$.
 Here we note that the right-hand side of the first  term corresponds to the definition of $\tau_q$ in \cite{BDLSV19}. We will repeat the procedures as in \cite{BDLSV19} and provide  some details of the  proof in Appendix \ref{app:sec:glu}. 
 
First, we recall that by \eqref{e:v:ell:0} and \eqref{papa:gluing}, $ \tau_q$ obeys the CFL-like condition:
\begin{equation}
\tau_q \norm{ v_{ l}}_{C^{1+\alpha}} \lesssim   \tau_q\delta_q^{{1}/{2}} \lambda_q  l^{-\alpha}  \lesssim  l^{\alpha} \ll   1.\label{e:CFL}
\end{equation}

Therefore, for each $i$, we can uniquely  solve the Euler equations locally  on  time interval $[ t_i - \tau_q, t_i+\tau_q]$  with initial datum $ v_{ l}(\cdot, t_i)\in C^{1+\alpha}$:
\begin{align}
\label{eq:Euler:i}
\partial_t  v_i+  v_i \cdot \nabla  v_i +\nabla  \pi_i &=0,\notag\\
\div  v_i &= 0,\\
 v_i(\cdot, t_i)&= v_{ l}(\cdot, t_i).\notag
\end{align}

Then we aim to establish the estimates on  $ v_i$, the differences between $ v_i$ and $ v_l$ in H\"older space and in the some  negative order spaces. For this purpose we apply the Biot-Savart operator to define 
\begin{align*}
  z_i 
  = (-\Delta)^{-1} \curl  v_i ,\ \ z_l =  (-\Delta)^{-1} \curl  v_l.
\end{align*}
It follows that $\div  z_i =\div  z_l= 0$ and $\curl  z_i =  v_i,\curl  z_l =  v_l$.

Then by the same argument as in \cite[Corollary 3.2. Proposition 3.3, Proposition 3.4]{BDLSV19}  we have
\bp\label{prop:gluedv}
For all $t\in [ t_i- \tau_q, t_i+ \tau_q]$ and all $N\in \mN_0$
\begin{align}
\norm{ v_i(t)}_{C^{N+1+\alpha}} 
&\les  \tau_q^{-1}  l^{-N+\alpha},
\label{eq:v_i:C^N}\\
 \norm{( v_i- v_{ l})(t)}_{C^{N+\alpha}}+\tau_q\norm{(\partial_t +  v_{ l} \cdot \nabla) ( v_i- v_{ l})(t)}_{C^{N+\alpha}}  &\lesssim   \tau_q\delta_{q+1} l^{-N-1+\alpha}\,, \label{e:z_diff_k}\\
\norm{ (z_i- z_l)(t)}_{C^{N+\alpha}} +\tau_q\norm{(\partial_t +  v_{ l} \cdot \nabla) ( z_i- z_l)(t)}_{C^{N+\alpha}} &\lesssim   \tau_q\delta_{q+1} l^{-N+\alpha}\,.  \label{e:z_diff}
\end{align}
\ep
Here we remark that although the choice of parameters are different from \cite{BDLSV19}, all the estimates   can be established  analogously by using \eqref{papa:gluing}.  We  put the proof in  Appendix \ref{app:a:glu} for completeness.

Then like the usual gluing step,  we introduce the time-dependent cut-off functions  $\{ \chi_i\}_i$, which is a partition of unity in time for $[0,1]$ with the property that $\supp  \chi_i\cap \supp  \chi_{i+2}=\emptyset$ and moreover for $N\in\mN_0$
\begin{align}
\supp  \chi_i&\subset [ t_i-\frac{2 \tau_q}{3}, t_i+\frac{2 \tau_q}{3}] \, ,\ \  \ 
 \chi_i=1\textrm{ on }( t_i - \frac{ \tau_q}{3},  t_i+ \frac{ \tau_q}{3})  ,  \ \
\norm{\partial_t^N  \chi_i}_{C^0} \lesssim  \tau_q^{-N} \notag.
\end{align}

The glued stress term can be defined following the same argument as in \cite[Section 4.2]{BDLSV19}. 
We glue the constructed  exact solutions $ v_{i}$ together and define 
\begin{align}
\label{eq:bar:v_q:def}
\overline v_q(x,t):= \sum_i  \chi_i(t)  v_i(x,t) \, .
\end{align}

We note that the glued vector field $\overline  v_q$  is also divergence-free, as the cutoffs $ \chi_i$ only depend on time. Furthermore, for $t\in J_i$, we have $\mathring{\overline{R}}_q=0$, and
for  $t \in  I_i$ we  define
 \begin{align}
\mathring{\overline{R}}_q&=\partial_t \chi_i\mathcal{R}( v_i- v_{i+1})- \chi_i(1- \chi_i)( v_i- v_{i+1})\mathring{\otimes} ( v_i- v_{i+1}) \, ,\label{eq:bar:R_q:def}
\end{align}
where we used the inverse divergence operator $\RSZ$ introduced  in Section \ref{tamr}.
 By construction, $\mathring{\overline{R}}_q$ is traceless and symmetric,  and we know that $(\overline  v_q, \mathring{\overline{R}}_q)$ solves \eqref{e:euler_reynolds} for a suitable $\overline \pi_q$.

To finish this section, it  remains to estimate the glued velocity field and Reynolds stress defined in \eqref{eq:bar:v_q:def} and \eqref{eq:bar:R_q:def}.  We define $ \overline z_q:=(-\Delta)^{-1} \curl  \overline v_q$. Then the desired estimates are obtained by  the same argument as in \cite[Proposition 4.2, Proposition 4.3, Proposition 4.4]{BDLSV19}. We put the details  in Appendix \ref{app:a:glu}.
\bp\label{prop:gluedr}
For $N\in\mN_0$, 
\begin{align}
\norm{\overline{v}_q- v_{ l}}_{C^{N+\alpha}}+ l^{-1}\norm{\overline{z}_q- z_{ l}}_{C^{N+\alpha}} \les \tau_q\delta_{q+1} l^{-1-N+\alpha}&\les \delta_{q+1} ^{1/2}l^{-N+\alpha}, \label{e:vq:vell:additional} \\
 \left| \int_{\mathbb{T}^3} |\overline v_q|^2- |v_l|^2 \dif x \right|  &\lesssim \delta_{q+1} l^{\alpha} \, , \label{bd: energy vq-vl}\\
 \norm{\mathring{\overline{ R}}_{q}}_{C^{N+\alpha}}+\tau_q\norm{(\partial_t + \overline v_q\cdot \nabla) \mathring{\overline{ R}}_{q}}_{C^{N+\alpha}} &\lesssim \delta_{q+1}l^{-N+\alpha} \label{e:Rq:1}.
\end{align}
\ep

In particular, from the bounds \eqref{e:v:ell:0}, \eqref{e:vq:vell:additional}   and  the choice of parameters in \eqref{papa:gluing}  we  obtain  for all $N\in\mN_0$

\begin{align}
\norm{\overline v_q}_{C^{N+1}} &\lesssim \norm{\overline v_q-v_l}_{C^{N+1}} +\norm{v_l}_{C^{N+1}}\les  \tau_q\delta_{q+1} l^{-2-N+\alpha}+ \delta_{q}^{1/2} \lambda_q  l^{-N}\lesssim \tau_q^{-1}l^{2\alpha-N}.\label{e:vq:1}
\end{align}
Here  the estimate  is slightly different from \cite[(4.7)]{BDLSV19} due to the the different choice of parameters. However, the new bound is  still suitable for our proof. 
\subsubsection{Gluing step for $(\rho_l,M_l)$}\label{sec:gluing:rho}

Now we glue the exact solutions $\rho_i$ of the transport equations at the same times $t_i$  defined above, to make sure the glued stress $\overline M_q$ is located in disjoint time intervals $I_i$.

For each $i$, we solve the following transport equations on $[ t_i - \tau_q, t_i+\tau_q]$:
\begin{align}
\partial_t  \rho_i+  \overline  v_q \cdot\nabla  \rho_i&=0,\notag\\
 \rho_i(\cdot, t_i)&= \rho_l(\cdot, t_i).\notag
\end{align}
By \eqref{e:vq:1} we have
\begin{equation}
\tau_q \norm{ \overline  v_q}_{C^{1+\alpha}} \lesssim  l^{\alpha} \ll   1.\notag
\end{equation}
 Then using  the  estimates for the transport  equations in 
 Proposition \ref{esti:transport}, the bounds   in \eqref{e:rho:ell:0},  \eqref{e:vq:1} and interpolation we have that for all $t\in [t_i-\tau_q,t_i+\tau_q]$ and $N\geq 1$
\begin{align}
\norm{\rho_i(t)}_{C^{N+\alpha}}& \les \norm{\rho_l(t_i)}_{C^{N+\alpha}} +\tau_q\|\overline  v_q\|_{C^{N+\alpha}}\norm{\rho_l}_{C^{1}}\notag\\ &\les \tilde{\delta}_{q}^{1/2} \lambda_{q}  l^{1-N-\alpha}+\tau_q \tau_q^{-1}\tilde{\delta}_{q}^{1/2} \lambda_{q}  l^{1-N}\les \tilde{\delta}_{q}^{1/2} \lambda_{q}  l^{1-N-\alpha}.
\label{eq:rho_i:C^N}
\end{align}

Moreover, we have the following result.
\bp
For $t\in [t_i-\tau_q,t_i+\tau_q]$ and all $N\geq0$,
\begin{align}
    \norm{(\rho_i-\rho_{l})(t)}_{C^{N+\alpha}}+\tau_q\norm{(\partial_t + v_l \cdot \nabla) (\rho_i-\rho_{l})(t)}_{C^{N+\alpha}} &\lesssim  \tau_q{\delta}_{q+1}^{1/2}\tilde{\delta}_{q+1}^{1/2}l^{-N-1+\alpha}\,.\label{e:rho_diff_k}
\end{align}
\ep
\begin{proof}
 It is easy to see that  $\rho_l - \rho_i$ obeys
\begin{align}
\partial_t (\rho_{l} - \rho_i)   + v_{l} \cdot \nabla (\rho_{l} - \rho_i)  
= -(v_l - \overline v_{q}) \cdot \nabla \rho_i-\div {M}_{l}
\label{e:rho_diff_evo}
\end{align}
with initial condition  $0$ at $t = t_i$. Then by the estimates \eqref{e:M:ell}, \eqref{e:vq:vell:additional}, \eqref{eq:rho_i:C^N} and the definition of $\tau_q$  we obtain for all $N\geq 0$ and $t\in [t_i-\tau_q,t_i+\tau_q]$
\begin{align*}
\norm{(\partial_t + v_l \cdot \nabla) (\rho_i-\rho_{l})(t)}_{C^{N+\alpha}} &\lesssim\|v_l - \overline v_{q}\|_{C^{N+\alpha}}\| \nabla \rho_i\|_{C^\alpha}+\|v_l - \overline v_{q}\|_{C^\alpha}\| \nabla \rho_i\|_{C^{N+\alpha}}+\| {M}_{l}\|_{C^{N+1+\alpha}}\notag\\ 
&\lesssim \tilde{\delta}_{q}^{1/{2}} \lambda_{q}\tau_q{\delta}_{q+1}l^{-N-1}+{\delta}_{q+1}^{1/2}\tilde{\delta}_{q+1}^{1/2}l^{-N-1+\alpha}\les{\delta}_{q+1}^{1/2}\tilde{\delta}_{q+1}^{1/2}l^{-N-1+\alpha}\,.
\end{align*}
 Then by   standard estimates for the transport  equations in Proposition \ref{esti:transport}, it follows immediately that  for $t\in [t_i-\tau_q,t_i+\tau_q]$
 \begin{align*}
 \norm{(\rho_i-\rho_{l})(t)}_{C^{\alpha}}&\lesssim\left| \int_{t_0}^t \norm{(\partial_t + v_l \cdot \nabla) (\rho_i-\rho_{l})(s)}_{C^{\alpha}}\dif s\right|\lesssim  \tau_q{\delta}_{q+1}^{1/2}\tilde{\delta}_{q+1}^{1/2}l^{-1+\alpha},
 \end{align*}
 and for $N\geq1$
  \begin{align*}
    \norm{(\rho_i-\rho_{l})(t)}_{C^{N+\alpha}}&\lesssim \left|\int_{t_0}^t \norm{(\partial_t + v_l \cdot \nabla) (\rho_i-\rho_{l})}_{C^{N+\alpha}}+\tau_q\|v_l\|_{C^{N+\alpha}}\norm{(\partial_t + v_l \cdot \nabla) (\rho_i-\rho_{l})}_{C^{1}}\dif s\right|\notag\\ 
    &\lesssim \tau_q{\delta}_{q+1}^{1/2}\tilde{\delta}_{q+1}^{1/2}l^{-N-1+\alpha}+\tau_q^2\lambda_q\delta_q^{1/2}l^{1-N-\alpha}{\delta}_{q+1}^{1/2}\tilde{\delta}_{q+1}^{1/2}l^{-2+\alpha}\lesssim  \tau_q{\delta}_{q+1}^{1/2}\tilde{\delta}_{q+1}^{1/2}l^{-N-1+\alpha}\,.
\end{align*}
\end{proof}

 We also introduce vector potentials  using the inverse divergence operator $\mathcal{R}_1=\nabla\Delta^{-1}$. More precisely, as $\int_{\mathbb{T}^3}\rho_i\dif x=\int_{\mathbb{T}^3}\rho_l\dif x=1$,  defining
\begin{align*}
  y_i:=\nabla\Delta^{-1}(\rho_i-1),\ \ y_l:=\nabla\Delta^{-1}(\rho_l-1),
\end{align*}
we have
\begin{align*}
 \rho_i=\div y_i+1,\ \ \rho_l=\div y_l+1,\ \ \curl y_i=\curl y_l=0.
\end{align*}
 With this notation, now we need to estimate  $y_{l}-y_i$ in some H\"older spaces. To achieve this, first we establish the following analytic identities, which can be used to derive the equation for $y_{l}-y_i$ by applying $\mathcal{R}_1$ on  both sides of \eqref{e:rho_diff_evo}. This also allows us to utilize the basic estimates for transport equations  in Appendix \ref{esti:transport} in this context.
\bl \label{lem;identity} For any smooth vector fields $z,v:\mR^3\to\mR^3$ satisfying $\div v=0$, and  smooth function $\rho:\mR^3\to\R$, we have
 \begin{align*}
    (\curl z)\cdot\nabla \rho =\div (z\times\nabla\rho),\\     v\cdot\nabla(\div z)=\div (v\cdot\nabla z-z\cdot\nabla v),\\
    \curl(v\cdot\nabla z)=-\div((z\times \nabla)v)+v\cdot\nabla (\curl z).
\end{align*}
 \el
The proof of this Lemma is provided in Appendix \ref{app:sec:glu}. Then we have the following result.
\bp
For $t\in [t_i-\tau_q,t_i+\tau_q]$, and all $N\geq0$
\begin{align}
   \|(y_i-y_l)(t) \|_{C^{N+\alpha}} +\tau_q\norm{(\partial_t + v_l \cdot \nabla) (y_i-y_{l})(t)}_{C^{N+\alpha}}  \lesssim\tau_q{\delta}_{q+1}^{1/2} \tilde{\delta}_{q+1}^{1/2} l^{-N+\alpha}.\label{bd:yi-yl}
\end{align}
\ep
\begin{proof}
    
 By lemma \ref{lem;identity}, we use the fact that $\div v_l=0$ to deduce 
\begin{align*}
    (v_l - \overline v_{q}) \cdot \nabla \rho_i=\div \((z_l - \overline z_{q})\times\nabla \rho_i\),\\
     v_{l} \cdot \nabla (\rho_{l} - \rho_i) =\div \( v_{l}\cdot\nabla (y_{l} - y_i)- (y_{l} - y_i)\cdot\nabla v_{l}\).
\end{align*}
 Consequently, from \eqref{e:rho_diff_evo} one deduces that
\begin{align}
\div\(\partial_t (y_{l} - y_i)   + v_{l} \cdot \nabla (y_{l} - y_i)\)
=\div\((y_{l} - y_i) \cdot \nabla v_{l} -(z_l - \overline z_{q})\times \nabla\rho_i-{M}_{l}\).
\notag
\end{align}
Taking gradient on both side of the above equation, and using the identity $\nabla\div=\Delta+\curl\curl$ we obtain
\begin{align}
\Delta\(\partial_t (y_{l} - y_i)   + v_{l} \cdot \nabla (y_{l} - y_i)\)
&=\nabla\div\((y_{l} - y_i) \cdot \nabla v_{l} -(z_l - \overline z_{q} )\times \nabla\rho_i-{M}_{l}\)\notag\\ 
&-\curl\curl\(\partial_t (y_{l} - y_i)   + v_{l} \cdot \nabla (y_{l} - y_i)\).
\notag\end{align}
Since $\curl y_l=\curl y_i=0,\div v_l=0$, by Lemma \ref{lem;identity} we have 
\begin{align*}
    \curl\( v_{l} \cdot \nabla (y_{l} - y_i)\)=-\div\(((y_{l} - y_i)\times \nabla)v_l\),
\end{align*}
which implies that 
\begin{align*}
   \partial_t (y_{l} - y_i)   + v_{l} \cdot \nabla (y_{l} - y_i)
&=\Delta^{-1}\nabla\div\((y_{l} - y_i) \cdot \nabla v_{l} -(z_l - \overline z_{q}) \times \nabla\rho_i-{M}_{l}\)\\ 
&+\Delta^{-1}\curl\div\(((y_{l} - y_i)\times \nabla)v_l\).
\end{align*}
As a result, for $N\in\mN_0$
\begin{align*}
\norm{(\partial_t + v_l \cdot \nabla) (y_i-y_{l})}_{C^{N+\alpha}}
&\les\|y_{l} - y_i\|_{C^{\alpha}} \|\nabla v_{l} \|_{C^{N+\alpha}}+\|y_{l} - y_i\|_{C^{N+\alpha}} \|\nabla v_{l} \|_{C^{\alpha}}\\
&+\|z_l - \overline z_{q}\|_{C^{N+\alpha}}\|\nabla\rho_i\|_{C^{\alpha}}+\|z_l - \overline z_{q}\|_{C^{\alpha}}\|\nabla\rho_i\|_{C^{N+\alpha}}+\|{M}_{l}\|_{C^{N+\alpha}}.
\end{align*}

Then together with the estimates of  $v_l,M_l,z_l - \overline z_{q}$ and $\rho_i$ in  \eqref{e:v:ell:0}, \eqref{e:M:ell},  \eqref{e:vq:vell:additional}, and \eqref{eq:rho_i:C^N} respectively we deduce that
\begin{align*}
\norm{(\partial_t + v_l \cdot \nabla) (y_i-y_{l})}_{C^{\alpha}}&\lesssim  {\delta}_{q}^{1/2} \lambda_{q} l^{-\alpha}\|y_i-y_l \|_{C^{\alpha}}+\tilde{\delta}_{q}^{1/2} \lambda_{q}\tau_q\delta_{q+1} + {\delta}_{q+1}^{1/2} \tilde{\delta}_{q+1}^{1/2} l^{\alpha}\\
&\lesssim  {\delta}_{q}^{1/2} \lambda_{q}l^{-\alpha} \|y_i-y_l \|_{C^{\alpha}}+ {\delta}_{q+1}^{1/2} \tilde{\delta}_{q+1}^{1/2} l^{\alpha}.
\end{align*}
Using \eqref{e:CFL} and Proposition \ref{esti:transport} we obtain for $t\in[t_i-\tau_q,t_i+\tau_q]$
\begin{align*}
    \|(y_i-y_l) (t)\|_{C^{\alpha}}\lesssim |t-t_i|{\delta}_{q+1}^{1/2} \tilde{\delta}_{q+1}^{1/2} l^{\alpha}+\left|\int_{t_i}^t {\delta}_{q}^{1/2} \lambda_{q}l^{-\alpha}\|(y_i-y_l)(s) \|_{C^{\alpha}}d s\right|,
\end{align*}
which by Gronwall's inequality implies that for $t\in[t_i-\tau_q,t_i+\tau_q]$
\begin{align*}
   \|(y_i-y_l)(t) \|_{C^{\alpha}} +\tau_q\norm{(\partial_t + v_l \cdot \nabla) (y_i-y_{l})(t)}_{C^{\alpha}}  \lesssim\tau_q{\delta}_{q+1}^{1/2} \tilde{\delta}_{q+1}^{1/2} l^{\alpha}.
\end{align*}
Finally, commuting the derivatives in $N+\alpha,N\geq0$ with $\partial_t + v_l \cdot \nabla$ as in the proof of \cite[Proposition 3.4]{BDLSV19} we obtain  for $t\in[t_i-\tau_q,t_i+\tau_q]$
\begin{align*}
   \|(y_i-y_l)(t) \|_{C^{N+\alpha}} +\tau_q\norm{(\partial_t + v_l \cdot \nabla) (y_i-y_{l})(t)}_{C^{N+\alpha}}  \lesssim\tau_q{\delta}_{q+1}^{1/2} \tilde{\delta}_{q+1}^{1/2} l^{-N+\alpha}.
\end{align*}
\end{proof}

As before, we  glue the constructed $ \rho_{i}$ together using the same time-dependent cut-off functions $ \chi_i$:
\begin{align}
\label{eq:bar:rho_q:def}
\overline \rho_q(x,t):= \sum_i   \chi_i(t)  \rho_i(x,t),
\end{align}
which inherits the identity $\int_{\mathbb{T}^3}\overline \rho_q\dif x=1$ from $\rho_i$.  Moreover, we denote 
\begin{align}
    i_q:=\max\{i:t_i<T_q\},\label{def:iq}
\end{align} where we recall that $t_i=i\tau_q$ and $T_q$ is defined in \eqref{def:tq}. For any $i\leq i_q, \rho_l(t_i)=1$, which implies that $\rho_i(t)=1$ on the interval $[t_i-\tau_q,t_i+\tau_q]$ by solving the transport equation and using the fact that $\div\overline v_q=0$. By the definition of the gluing solution, we know that $\overline \rho_q(t)=1$ on $[0,T_q-\tau_q]\subset[0,t_{i_q}]$.

By the definition of the cutoff functions, on every $ J_i$ interval we have $\overline  \rho_q =  \rho_i$, so  $(\overline v_q,\overline  \rho_q)$ is an exact solution of the transport equations \eqref{e:euler_trans}. On the other hand, on every interval $ I_i$  we have
\begin{align*}
\overline  \rho_q =  \chi_i  \rho_i + (1- \chi_i)  \rho_{i+1},
\end{align*} 
which leads to 
\begin{align*}
\partial_t\overline  \rho_q+\div( \overline v_q  \overline \rho_q) = 
\partial_t \chi_i( \rho_i- \rho_{i+1}).
\end{align*}
So for all $t \in I_i$ we   define
 \begin{align}
\overline{M}_q&=-\partial_t \chi_i(y_i- y_{i+1}), \label{eq:bar:M_q:def}
\end{align}
which together with \eqref{eq:bar:R_q:def} implies  that $(\overline  v_q, \overline\rho_q, \mathring{\overline{R}}_q,\overline{M}_q)$ solves \eqref{e:euler_trans}-\eqref{e:euler_reynolds} on ${\T} \times [0,1]$. 
Since $\rho_i=1$ on $[t_i-\tau_q,t_i+\tau_q]$ for any $i\leq i_q$, we obtain  $y_i=0$  on these intervals.  Then we have  $\supp (\overline M_q) \subset {\T} \times \cup_{i\geq i_q}   I_i$.

Now we have the following estimates on the glued solution and stress term. 

\begin{proposition}
    For $N\in\mN_0$, we have 
    \begin{align}
\norm{\overline{\rho}_q- \rho_{ l}}_{C^{N+\alpha}}& \lesssim  \tau_q{\delta}_{q+1}^{1/2} \tilde{\delta}_{q+1}^{1/2} l^{-1-N+\alpha}\lesssim\tilde{\delta}_{q+1}^{1/2}  l^{-N+\alpha},\label{e:rhoq:vell}  \\
\norm{\overline \rho_q}_{C^{N+1}} &\les \tilde{\delta}_{q}^{1/2} \lambda_q  l^{-N},\label{e:rhoq:1} \\
 \norm{\overline{ M}_{q}}_{C^{N+\alpha}}+\tau_q\norm{(\partial_t + \overline v_q\cdot \nabla)\overline{ M}_{q}}_{C^{N+\alpha}} &\lesssim {\delta}_{q+1}^{1/2} \tilde{\delta}_{q+1}^{1/2}l^{-N+\alpha}.\label{e:Mq:1}
\end{align}

\end{proposition} 

\begin{proof}
By \eqref{e:rho:ell:0}, \eqref{e:rho_diff_k}  and the fact that 
$ \tau_q\delta_{q+1} ^{1/2}l^{-1}\leq 1$ in \eqref{papa:gluing}, we obtain for all $N\geq 0$
\begin{align*}
\norm{\overline{\rho}_q- \rho_{ l}}_{C^{N+\alpha}} &\lesssim \sum_i\chi_i\norm{{\rho}_i- \rho_{ l}}_{C^{N+\alpha}} \lesssim  \tau_q{\delta}_{q+1}^{1/2} \tilde{\delta}_{q+1}^{1/2} l^{-1-N+\alpha}\lesssim\tilde{\delta}_{q+1}^{1/2}  l^{-N+\alpha},\\
\norm{\overline \rho_q}_{C^{N+1}} &\lesssim \norm{\overline{\rho}_q- \rho_{ l}}_{C^{N+1}}+\norm{ \rho_{ l}}_{C^{N+1}}\les \tau_q{\delta}_{q+1}^{1/2} \tilde{\delta}_{q+1}^{1/2} l^{-2-N+\alpha}+ \tilde{\delta}_{q}^{1/2} \lambda_q  l^{-N}\notag\\ 
&\les ((l\lambda_q)^{3\alpha}+1)\tilde{\delta}_{q}^{1/2} \lambda_q  l^{-N}\les \tilde{\delta}_{q}^{1/2} \lambda_q  l^{-N}.
\end{align*}

On the other hand, for the glued stress $\overline{ M}_{q}$,  by \eqref{eq:bar:M_q:def}, the property of the cut-off functions and the estimate in \eqref{bd:yi-yl} we have for $t\in I_i$,
\begin{align*}
    \norm{\overline{ M}_{q}(t)}_{C^{N+\alpha}}\les\|  \partial_t \chi_i\|_{C_t^0}\|(y_i- y_{i+1})(t)\|_{C^{N+\alpha}}\les\tau_q^{-1}\tau_q{\delta}_{q+1}^{1/2} \tilde{\delta}_{q+1}^{1/2} l^{-N+\alpha}\les{\delta}_{q+1}^{1/2} \tilde{\delta}_{q+1}^{1/2} l^{-N+\alpha}.
\end{align*}
 Next we have  for $t\in I_i$,
  \begin{align*}
      (\partial_t +  v_l\cdot \nabla)\overline{ M}_{q}=-\partial^2_t \chi_i(y_i- y_{i+1})-\partial_t \chi_i (\partial_t +  v_l\cdot \nabla)(y_i- y_{i+1}),
  \end{align*}
which  by \eqref{bd:yi-yl} and the property of the cut-off functions implies that
\begin{align*}
    \norm{(\partial_t +  v_l\cdot \nabla)\overline{ M}_{q}(t)}_{C^{N+\alpha}}\les \tau_q^{-2}\tau_q{\delta}_{q+1}^{1/2} \tilde{\delta}_{q+1}^{1/2} l^{-N+\alpha}+\tau_q^{-1}{\delta}_{q+1}^{1/2} \tilde{\delta}_{q+1}^{1/2} l^{-N+\alpha}\les\tau_q^{-1}{\delta}_{q+1}^{1/2} \tilde{\delta}_{q+1}^{1/2} l^{-N+\alpha}.
\end{align*}
Then using the fact that $\delta_{q+1}^{1/2}\tau_q l^{-1}\leq1$ in \eqref{papa:gluing} and \eqref{e:vq:vell:additional}  we deduce that
\begin{align*}
    &\norm{(\partial_t + \overline v_q\cdot \nabla)\overline{ M}_{q}}_{C^{N+\alpha}}\les   \norm{(\partial_t +  v_l\cdot \nabla)\overline{ M}_{q}}_{C^{N+\alpha}}+\norm{(\overline v_q-v_l)\cdot\nabla \overline{ M}_{q}}_{C^{N+\alpha}}\\
    &\qquad\les   \norm{(\partial_t +  v_l\cdot \nabla)\overline{ M}_{q}}_{C^{N+\alpha}}+\norm{\overline v_q-v_l}_{C^{N+\alpha}}\norm{\nabla \overline{ M}_{q}}_{C^{\alpha}}+\norm{\overline v_q-v_l}_{C^{\alpha}}\norm{\nabla \overline{ M}_{q}}_{C^{N+\alpha}}\\
    &\qquad\lesssim\tau_q^{-1}{\delta}_{q+1}^{1/2} \tilde{\delta}_{q+1}^{1/2} l^{-N+\alpha}+ \tau_q\delta_{q+1}^{3/2}\tilde{\delta}_{q+1}^{1/2}l^{-N-2+2\alpha}\les\tau_q^{-1}{\delta}_{q+1}^{1/2} \tilde{\delta}_{q+1}^{1/2}l^{-N+\alpha}.
\end{align*}
\end{proof}

\subsection{The perturbation procedure}\label{sec:defq+1}
 To proceed the procedure, inspired by Isett's work \cite{Ise18}, we proceed with the construction of the perturbation   $(w_{q+1}+\overline w_{q+1},\theta_{q+1})$, where the support of $ (w_{q+1},\theta_{q+1})$ and that of $\overline w_{q+1}$ are disjoint. This can be achieved by choosing building blocks being disjoint, as detailed in Appendix \ref{sec:Mikado}. 
 Here  the perturbation $(w_{q+1},\theta_{q+1})$ is used to cancel the low frequency part of the stress term $\overline M_q$. However, it brings a new term $w_{q+1}\otimes w_{q+1}$ into the Euler equation.  On the way, the perturbation $\overline w_{q+1}$ is used to cancel the  stress term $\mathring{\overline R}_q$ and the low frequency part  from $w_{q+1}\otimes w_{q+1}$. 
To achieve this,  we need the the Mikado flows recalled in  Appendix \ref{sec:Mikado}. Then we define  the perturbation $(w_{q+1},\theta_{q+1})$ in Section \ref{sec:Onsager:principal:corrector:1} and  define  the perturbation $\overline w_{q+1}$ in Section \ref{sec:Onsager:principal:corrector:2}. The  corresponding estimates are shown in Section \ref{sec:Onsager:velocity:inductive}.

\subsubsection{The construction of the perturbation $(w_{q+1},\theta_{q+1})$}\label{sec:Onsager:principal:corrector:1}
In this section, we aim to construct the perturbation $(w_{q+1},\theta_{q+1})$, which is used to cancel the stress term $\overline {M}_{q}$. 
Recall that $\overline{M}_q$ has support in ${\T} \times   \cup_{i\geq i_q}   I_i$, where $ I_i$ is  defined in \eqref{eq:ti:Ii} and  $i_q$ is defined in \eqref{def:iq}. Then we  define a family of smooth temporal cut-off functions  $\{  \eta_i\}_{i\geq i_q}$  with the following properties: 
\begin{enumerate}
\item  $0 \leq  \eta_i\leq 1$, $ \eta_i \equiv 1$  on $  I_i$ ,   $\supp ( \eta_i)  \subset  ( t_i+\frac{1}{6}\tau_q, t_i+\frac{5}{6}\tau_q)$, 
\item $\norm{\partial_t^n  \eta_i}_{C_t} \les\tau_q^{-n}$,  for all $n\geq 0$.
\end{enumerate}

Since $  \eta_i \equiv 1$ on $I_i$, $  \eta_i   \eta_j \equiv 0$ for $i\neq j$, and $\supp (\overline M_q) \subset {\T} \times \cup_{i\geq i_q}   I_i$, we have that 
\begin{align}
\sum_i   \eta_i^2 \overline M_q = \overline M_q \, .
\label{eq:etaM:partition}
\end{align}

Then we define the flow maps $ \Phi_i$ for the velocity field $ \overline v_q $ as the solution of the transport equation
\begin{align}
(\partial_t +   \overline v_q  \cdot \nabla)  \Phi_i &=0 \, ,\label{eq:Phi:i:def}\\  
 \Phi_i\left(x, t_i\right) &= x ,\notag
\end{align}
 for all $t \in (t_i-\frac{\tau_q}{3},t_i+\frac{4\tau_q}{3})$. In the following we use the notation $D_{t,q}=\partial_t +   \overline v_q  \cdot \nabla_x$.

We have the following estimates for the flow maps.
\bp\label{prop:flowmap}
For all $t\in  (t_i-\frac{\tau_q}{3},t_i+\frac{4\tau_q}{3})$ and $N\geq 0$, we have
\begin{align}
 \norm{\nabla  \Phi_i(t) - {\rm Id}}_{C^0} &\les l^{\alpha},
 \label{eq:Phi:i:bnd:a} \\
 \norm{(\nabla  \Phi_i)^{-1}(t)}_{C^{N}}+ \norm{\nabla  \Phi_i(t)}_{C^{N}}  
 &\les  l^{-N},
 \label{eq:Phi:i:bnd:b} \\
 \norm{D_{t,q} (\nabla  \Phi_i)^{-1}(t)}_{C^{N}} + \norm{D_{t,q} \nabla  \Phi_i(t)}_{C^{N}} & \les \tau_q^{-1}l^{-N}.
 \label{eq:Phi:i:bnd:c}
\end{align}
\ep
The proof of this Proposition is given in Appendix \ref{sec:est:ampl}.

Next for  $t\in \supp (  \eta_i)$ we  define the vector field 
\begin{align}
 M_{q,i}(t) =\nabla\Phi_i(t)\, \((\frac34,0,0)^T+\frac{\overline M_q(t)}{\delta_{q+1} ^{1/2}\tilde{\delta}_{q+1} ^{1/2}l^{\alpha/2}}\).
\label{eq:tilde:M:q:i:def}
\end{align}
 \br
Here we introduce  additional translation in the definition of $M_{q,i}$, which on an one hand ensures that $M_{q,i}$ remains in the annulus $\overline B_1(0)\backslash B_{\frac12}(0)$ such that Lemma~\ref{l:linear_algebra2} is applicable. On the other hand,  this translation also  guarantees that the supplementary term in \eqref{eq:wthe:q+1:is:good} below maintains divergence-free as the cut-offs $\eta_i $ are only functions of $t$.
\er
 Furthermore, using the estimates on $\overline M_q, \nabla\Phi_i$ in \eqref{e:Mq:1} and \eqref{eq:Phi:i:bnd:a} respectively, we  deduce that for $t\in\supp( \eta_i)$
\begin{align}
\norm{ M_{q,i}(t)-(\frac34,0,0)^T }_{C^0}&\les \norm{\nabla\Phi_i-\Id}_{C^0}+ \frac{\norm{\nabla\Phi_i\overline M_q}_{C^0}}{\delta_{q+1} ^{1/2}\tilde{\delta}_{q+1} ^{1/2}l^{\alpha/2}}\les l^{\alpha/2} \leq \frac14,\label{eq:tilde:M:q:i:bnd:1}
\end{align}
where we choose $a$ large enough to absorb the universal constant. 

Thus, due to \eqref{eq:tilde:M:q:i:bnd:1}, for $\xi\in\Lambda^1\cup\Lambda^2$, we have $\frac12\leq|M_{q,i}|\leq1$, then we apply  Lemma~\ref{l:linear_algebra2} and define the amplitude functions 
\begin{align}
  A_{(\xi,i)}=l^{\alpha/4} \delta_{q+1}^{1/2} \eta_i \,  \Gamma_{\xi} ^{1/2}( M_{q,i}),
 \label{eq:A:xi:i:def}\ \
\tilde{A}_{(\xi,i)} =l^{\alpha/4} \tilde{\delta}_{q+1}^{1/2} \eta_i \,  \Gamma_{\xi} ^{1/2}( M_{q,i}),
\end{align}
where  $\Gamma_{ \xi}$ are the functions from Lemma~\ref{l:linear_algebra2}. 
In particular, by the definition of the cut-off functions and $i_q$ in \eqref{def:iq}, we have $\eta_i=0$ on $[0,T_q-\tau_q]$. Hence it follows that $\tilde{A}_{(\xi,i)}=0$ on $[0,T_q-\tau_q]$.

Then we have the following estimates on the amplitude functions. 
\bp\label{prop:est:A,tilA} For $\xi\in\Lambda,i=1,2$ and $N\in\mN_0$, we have
\begin{align}
\norm{  A_{(\xi,i)}}_{C^N} +  \tau_q \norm{D_{t,q}  A_{(\xi,i)}}_{C^N} \les \delta_{q+1}^{1/2}  l^{\alpha/4-N},
\label{eq:Onsager:A:xi:CN}\\
\norm{  \tilde{A}_{(\xi,i)}}_{C^N} +  \tau_q \norm{D_{t,q}  \tilde{A}_{(\xi,i)}}_{C^N} \les \tilde{\delta}_{q+1}^{1/2}  l^{\alpha/4-N}.
\label{eq:Onsager:tA:xi:CN}
\end{align}
\ep
The proof of this Proposition is given in  Appendix \ref{sec:est:ampl}.

  Now we proceed with the construction of the principle part of $(w_{q+1},\theta_{q+1})$. To this end, we  employ the   Mikado flows which are defined  as in Appendix \ref{sec:Mikado} with $ \lambda=\lambda_{q+1}$, i.e.
\begin{align*}
     W_{( \xi)}(x)=  W_{ \xi, \lambda_{q+1}}(x),\ \ \Theta_{( \xi)}(x)=  \Theta_{ \xi, \lambda_{q+1}}(x)  \, . 
\end{align*}
Here for two index sets $ \Lambda^1,\Lambda^2$ as in Lemma~\ref{l:linear_algebra2}, we use the notation   $ \Lambda^i =  \Lambda^1$ for $i$ odd, and $ \Lambda^i =  \Lambda^2$ for $i$ even. With this notation, we now use the amplitudes defined in \eqref{eq:A:xi:i:def} to  define the  principal part as 
\begin{align}
w_{q+1}^{(p)}(x,t)& = \sum_{i} \sum_{ \xi \in  \Lambda^i}  A_{(  \xi,i)}(x,t) (\nabla  \Phi_i(x,t))^{-1}  W_{(\xi)}( \Phi_i(x,t)) \,\label{eq:Onsager:w:q+1:p} ,\\
\theta_{q+1}^{(p)}(x,t) &= \sum_{i} \sum_{ \xi \in  \Lambda^i} \tilde{ A}_{(  \xi,i)}(x,t)  \Theta_{(\xi)}( \Phi_i(x,t)) \, .
\label{eq:Onsager:the:q+1:p}
\end{align}
 We then remark that the vector field  $ U_{i,\xi}=(\nabla  \Phi_i)^{-1}  W_{(\xi)}( \Phi_i)$ satisfies the following Lie-advection identity:
\begin{align}
D_{t,q}  U_{i,\xi} = ( U_{i,\xi} \cdot \nabla)  \overline v_q = (\nabla   \overline v_q)^T  U_{i,\xi} \, .
\label{eq:Onsager:Lie:advect}
\end{align}

Using  \eqref{eq:Mikado:2},   \eqref{eq:Mikado:4}, the identity \eqref{eq:etaM:partition}, and the fact that the $\eta_i$ have mutually disjoint supports, we obtain
\begin{align}
&w_{q+1}^{(p)} \theta_{q+1}^{(p)} 
= \sum_i \sum_{ \xi \in  \Lambda^i}   A_{( \xi,i)}\tilde{ A}_{( \xi,i)} (\nabla  \Phi_i)^{-1}  ( W_{( \xi)}( \Phi_i) )  ( \Theta_{(  \xi)}( \Phi_i) )  \notag\\
&= \sum_i l^{\alpha/2}\delta_{q+1} ^{1/2}\tilde{\delta}_{q+1}^{1/2}\eta_i ^2  (\nabla  \Phi_i)^{-1}    \sum_{\xi \in \Lambda^i}  \Gamma_\xi( M_{q,i})  \left( (  W_{( \xi)} \Theta_{( \xi)})( \Phi_i)  \right) \notag\\
&= \sum_i l^{\alpha/2}  \delta_{q+1} ^{1/2}\tilde{\delta}_{q+1}^{1/2}\eta_i ^2 (\nabla  \Phi_i)^{-1}   M_{q,i} + \sum_i \sum_{ \xi \in  \Lambda^i}   A_{( \xi,i)}\tilde{ A}_{( \xi,i)} (\nabla   \Phi_i)^{-1} \left( \left( \Proj_{\neq 0}( W_{( \xi)} \Theta_{( \xi)}) \right)( \Phi_i) \right)     \notag\\
&=\sum_i l^{\alpha/2}  \delta_{q+1} ^{1/2}\tilde{\delta}_{q+1}^{1/2}\eta_i ^2 (\frac34,0,0)^T+\overline M_q + \sum_i \sum_{ \xi \in  \Lambda^i}   A_{( \xi,i)}\tilde{ A}_{( \xi,i)}(\nabla  \Phi_i)^{-1} \left( (\Proj_{\geq \frac{ \lambda_{q+1}}{2}}( W_{( \xi)} \Theta_{( \xi)}) )( \Phi_i) \right) ,
\label{eq:wthe:q+1:is:good}
\end{align}
where we recall the notation $\Proj_{\neq 0} f(x) = f(x) - \fint_{\T} f(y) \dif  y$.
For the last term we   used the fact that $ W_{( \xi)}\Theta_{( \xi)}$ is $(\mathbb{T}/  \lambda_{q+1})^3$-periodic. Here we remark that the divergence of  the first term in the last identity is zero since the cut-off functions $\eta_i$ only depend on $t$.

To define the incompressibility corrector, we note that
\begin{align*}
(\nabla  \Phi_i)^{-1} ( W_{( \xi)} ( \Phi_i) ) = \curl\left( (\nabla  \Phi_i)^T  V_{( \xi)}( \Phi_i) \right).
\end{align*}
Then we define
\begin{align}
w_{q+1}^{(c)}(x,t) = \sum_i \sum_{ \xi \in  \Lambda^i} \nabla  A_{( \xi,i)}(x,t) \times \left( (\nabla  \Phi_i(x,t))^T  V_{( \xi)} ( \Phi_i(x,t) )  \right),
\label{eq:Onsager:w:q+1:c}
\end{align}
and the total velocity increment $w_{q+1}$   as
\begin{align}
w_{q+1}: = w_{q+1}^{(p)} + w_{q+1}^{(c)} = \curl \( \sum_i \sum_{ \xi \in  \Lambda^i}  A_{( \xi,i)} \, (\nabla  \Phi_i)^T ( V_{( \xi)}( \Phi_i) )\),
\label{eq:Onsager:w:q+1}
\end{align}
which is automatically incompressible. In particular, by  definition we have 
 $\supp w_{q+1}\subset \cup_{i} \supp ( \eta_i)  \subset  \cup_{i} ( t_i+\frac{1}{6}\tau_q, t_i+\frac{5}{6}\tau_q)$.

We recall that  the perturbation $\theta_{q+1}$ needs to be mean-zero, so we define the corresponding corrected perturbation
\begin{align*}
\theta_{q+1}^{(c)}:&=-\mathbb{P}_{0}\theta_{q+1}^{(p)},\ \ \ {\rm and}\ \
    \theta_{q+1}:=\theta_{q+1}^{(p)}+\theta_{q+1}^{(c)}.
\end{align*}
Then since 
 $\tilde{A}_{(\xi,i)}=0$ on $[0,T_q-\tau_q]$, by  definition we know  $\theta_{q+1}=0$ on $[0,T_q-\tau_q]$.

\subsubsection{The construction of the perturbation $\overline{w}_{q+1}$}
\label{sec:Onsager:principal:corrector:2}
In this section, we construct the new perturbation $\overline{w}_{q+1}$ which is used to cancel  the stress term $\mathring{\overline R}_q$ and the low frequency part  from $w_{q+1}\otimes w_{q+1}$.  We first note that the support of $\mathring{\overline R}_q$ is located in ${\T} \times   \cup_{i}   I_i$, where $ I_i$ is  defined in \eqref{eq:ti:Ii} and the support of $ w_{q+1}$ is located in $\cup_{i}   ( t_i+\frac{1}{6}\tau_q, t_i+\frac{5}{6}\tau_q)$.
Then we define a family of smooth  cutoff functions   $\{  \overline\eta_i \}_{i\geq0}$ with the following properties:

\begin{enumerate}
\item  $0 \leq  \overline \eta_i\leq 1$, $\overline \eta_i \equiv 1$,  on ${\T} \times   ( t_i+\frac{1}{6}\tau_q, t_i+\frac{5}{6}\tau_q)$ ,   $\supp (\overline \eta_i) \subset  {\T} \times  (t_i-\frac{\tau_q}{3},t_i+\frac{4\tau_q}{3})$, 
 \item  $ \overline\eta_i \, \overline \eta_j \equiv 0$ for every  $i \neq j$,
\item   for all $t\in [0,1]$, 
\begin{align}
    c_{\overline \eta} \leq \sum_i \int_{\mathbb{T}^3}  \overline\eta_i^2(x,t) \dif x,\label{bd:sumetageq}
\end{align}  where $ c_{ \overline\eta}>0$ is a universal constant,  
\item $\norm{\partial_t^n \overline \eta_i}_{C^m} \les\tau_q^{-n}$,  for all $n,m\geq 0$.
\end{enumerate}
 In contrast to \cite[Lemma 5.3]{BDLSV19}, our construction requires the cutoff function to satisfy $\overline \eta_i \equiv 1$ over an extended domain  to ensure that the support  covers not only the support of  $\mathring{\overline R}_q$, but also that of $ w_{q+1}$.  The existence of such cut-off functions follows by using the following domains:
\begin{align*}
    \overline O_i&=\{(x,t):t_i+\frac{\tau_q}{12}(\sin(2\pi x_1)+\frac12)\leq t\leq t_{i+1}+\frac{\tau_q}{12}(\sin(2\pi x_1)-\frac12)\}.
\end{align*}
Then the desired cut-off functions are defined as
	\begin{align*}
		 \overline\eta_i := \mathbf{1}_{\overline O_i}*_t\varphi_{\epsilon_0\tau_q} *_x\phi_{\epsilon_0},
	\end{align*}
	where $\phi_{\epsilon_0}:=\frac{1}{\epsilon_0^d}\phi(\frac{\cdot}{{\epsilon_0}})$ is a family of standard mollifiers on $\mathbb{R}^3$, and $\varphi_{\epsilon_0\tau_q}:=\frac{1}{{\epsilon_0\tau_q}}\varphi(\frac{\cdot}{{\epsilon_0\tau_q}})$ is a family of standard mollifiers with support in $(0,1)$. 
	Then by choosing $\epsilon_0>0$ small enough, the desired conclusions hold  by a similar argument as in \cite[Lemma 5.3]{BDLSV19}.

Before defining the second  amplitude function, we define the   low frequency of the product $w_{q+1}^{(p)}\otimes w_{q+1}^{(p)}$  by
\begin{align*}
    R^{(1)}_{q}:= \sum_{i} \sum_{ \xi \in  \Lambda^i} A^2_{(\xi,i)}(\nabla \Phi_i)^{-1}\, \xi\otimes \xi\, (\nabla   \Phi_i)^{-T}.
\end{align*}
In fact, by  the fact that the building blocks have mutually disjoint support~\eqref{eq:Mikado:2}, the fact that the $\eta_i$ have mutually disjoint supports, and \eqref{eq:Mikado:3} we have 
\begin{align}
   w_{q+1}^{(p)}\otimes w_{q+1}^{(p)}-R^{(1)}_{q}=\sum_i \sum_{ \xi \in \Lambda^i}  A_{(  \xi,i)}^2 (\nabla   \Phi_i)^{-1} \( \left( \Proj_{\neq 0}( W_{( \xi)} \otimes  W_{( \xi)}) \right)( \Phi_i) \)  (\nabla  \Phi_i)^{-T}.\label{eq:wq+1wq+1}
\end{align}
By the Leibniz rule, the fact that $\eta_i$ have disjoint supports, the estimate for $A_{(\xi,i)}$ in \eqref{eq:Onsager:A:xi:CN} and the estimates for $\nabla\Phi_i$ in \eqref{eq:Phi:i:bnd:b}, \eqref{eq:Phi:i:bnd:c}, we have for $N\geq0$
\begin{align}
    \norm{ R^{(1)}_{q}}_{C^N}+ \tau_q\norm{ D_{t,q}R^{(1)}_{q}}_{C^N}\lesssim \delta_{q+1}l^{\alpha/2-N}.\label{bd:rq(1)}
\end{align}

To prescribe the energy profile, we  define the energy gap
	\begin{align}\label{def:energy gap}
		 \Upsilon_q(t):=\frac13 \left(e(t)-\frac{\delta_{q+2}}{2}-  \| \overline{v}_q (t)\|_{L^2}^2-\int_{\mathbb{T}^3}\tr  R^{(1)}_{q}(t,x)\dif x  \right),
	\end{align}
	and then by \eqref{bd:sumetageq} we decompose $\Upsilon_q$ by 
	\begin{align}\label{def:rho q,i}
		\Upsilon_{q,i}(t,x):=\frac{\overline\eta_i^2(t,x)}{\sum_j \int_{\mathbb{T}^3}\overline \eta_j^2(t,y) \dif y} \Upsilon_q(t) .
	\end{align}
	From the construction, it follows that $\sum_i \int_{\mathbb{T}^3} \Upsilon_{q,i}\dif x =\Upsilon_q$ for all $t\in[0,1]$. 

\begin{proposition}
		For any $t\in[0,1]$ and $N\geq0$ we have 
		\begin{align}
			\frac{\delta_{q+1}}{6\lambda_{q}^{\alpha/3}}\leq \Upsilon_q(t)  \leq \delta_{q+1}, \label{estimate:rho1}
			\ \ \Upsilon_{q,i}(t)&\leq \frac{\delta_{q+1}}{c_{\overline{\eta}}},\\
            \|\Upsilon_{q,i}\|_{C^N}+\tau_q\|\partial_t \Upsilon_{q,i}\|_{C^N} &\lesssim \delta_{q+1}.\label{estimate:rho2}
		\end{align}
	\end{proposition}
	\begin{proof}
    First, by the estimate for the energy gaps in  \eqref{estimate:energy}, \eqref{energy:vq-vl},  \eqref{bd: energy vq-vl} and the bound in \eqref{bd:rq(1)}  we have  for some $C>0$
	\begin{align*}
			\frac{\delta_{q+1}}{2\lambda_{q}^{\alpha/3}}&\leq \frac{\delta_{q+1}}{ \lambda_{q}^{\alpha/3}}-\frac{\delta_{q+2}}{2}-C\delta_{q+1} l^{\alpha/2} \leq 3 \Upsilon_q(t)
			\\&=  e(t  )- \| v_q (t  )\|_{L^2}^2 -\frac{\delta_{q+2}}{2} -\int_{\mathbb{T}^3}\tr  R^{(1)}_{q}(t,x)\dif x 
			\\ &\quad +  \left[ \| v_q (t  )\|_{L^2}^2- \| v_l (t  )\|_{L^2}^2\right] + \left[ \| v_l(t  )\|_{L^2}^2- \| \overline{v}_q (t  )\|_{L^2}^2\right]
			\\ & \leq \delta_{q+1}+C\delta_{q+1} l^{\alpha/2}\leq  2\delta_{q+1},
	\end{align*}
		where we choose $0<\alpha<6\beta b(b-1)$ small enough to deduce that $\delta_{q+2}\leq \frac12\delta_{q+1}\lambda_{q}^{-\alpha/3}$. We also used  $l\lambda_{q}\ll1$ and  choose   $a$ large enough to absorb the universal constant. Then together with the definition of the cut-off functions, we obtain the second inequality in \eqref{estimate:rho1}.
        
    By the properties of the cutoff functions and \eqref{estimate:rho1}, 
        the  bound for the first term in \eqref{estimate:rho2} follows. Now we derive the time regularity of $\Upsilon_{q,i}$. Since $(\overline{v}_q, \mathring{\overline{R}}_q)$ obeys \eqref{e:euler_reynolds},  by the basic energy estimate, the bounds in \eqref{e:Rq:1} and \eqref{e:vq:1} we have  for $t\in[0,1]$
        \begin{align*}
          \left|  \partial_t \| \overline{v}_q(t)\|_{L^2}^2\right|\lesssim \left|\int_{\mT^d}\mathring{\overline{R}}_q(t,x)\cdot\nabla\overline{v}_q(t,x)\dif x\right|\lesssim \delta_{q+1}\tau_q^{-1}l^\alpha.
        \end{align*}
        
Since $\int_{\mathbb{T}^3}  \overline v_q\cdot\nabla R^{(1)}_{q}(t,x)\dif x=0 $, by  \eqref{e:Rq:1}  we have  for $t\in[0,1]$
\begin{align*}
   \left|  \int_{\mathbb{T}^3}\tr  \partial_tR^{(1)}_{q}(t,x)\dif x\right|=   \left| \int_{\mathbb{T}^3}\tr D_{t,q}R^{(1)}_{q}(t,x)\dif x\right|\les   \|  D_{t,q}R^{(1)}_{q}\|_{C^0}\lesssim \tau_q^{-1}\delta_{q+1}l^\alpha,
\end{align*}
which implies that 
        \begin{align}
		\|\partial_t\Upsilon_q\|_{C_t^0}\lesssim|e'(t)|+ \left|  \partial_t \| \overline{v}_q(t)\|_{L^2}^2\right|+ \left|  \int_{\mathbb{T}^3}\tr  \partial_tR^{(1)}_{q}(t,x)\dif x\right|\lesssim \tau_q^{-1}\delta_{q+1}l^\alpha,\label{bd:parttrhoq}
	\end{align}
    where we used that $|e'(t)|\leq1$ and choose $a$ large enough.
    Then by the Leibniz rule  and the properties of the cutoff functions it is easy to see that
		$\|\partial_t[\frac{\overline\eta_i^2(t,x)}{\sum_j \int_{\mathbb{T}^3}\overline \eta_j^2(t,y) \dif y}]\|_{C^N}\lesssim \tau_q^{-1}$.  Thus, the bound for the second term in \eqref{estimate:rho2} follows  by applying the Leibniz rule again.
	\end{proof}

Since $ \overline \eta_i \equiv 1$ on ${\T} \times    ( t_i+\frac{1}{6}\tau_q, t_i+\frac{5}{6}\tau_q)$, $  \overline  \eta_i    \overline \eta_j \equiv 0$ for $i\neq j$, and since $\supp (\Rbar_q) ,\supp R_q^{(1)}
\subset {\T} \times \cup_i   ( t_i+\frac{1}{6}\tau_q, t_i+\frac{5}{6}\tau_q)$, we have that 
\begin{align}
\sum_i  \overline   \eta_i^2 \Rbar_q = \Rbar_q \, ,\ \ \ \sum_i   \overline  \eta_i^2  R^{(1)}_{q} =  R^{(1)}_{q}\, .
\label{eq:eta:partition}
\end{align}
Now we define the symmetric tensor
\begin{align}
 R_{q,i} = \nabla   \Phi_i\(\Id- \frac{\overline  \eta_i  ^2(\Rbar_q+ \mathring{R}^{(1)}_{q})}{\Upsilon_{q,i}}\)\, \nabla   \Phi_i^T =\nabla  \Phi_i\(\Id-\frac{\sum_j \int_{\mathbb{T}^3}\overline \eta_j^2(t,y) \dif y} {\Upsilon_q(t)}(\Rbar_q+ \mathring{R}^{(1)}_{q})\)\, \nabla   \Phi_i^T
\label{eq:tilde:R:q:i:def}
\end{align}
for all $(x,t) \in \supp ( \overline \eta_i)$. Here $\mathring{R}^{(1)}_{q}$ means the trace-free part of ${R}^{(1)}_{q}$.  By the bounds  in \eqref{e:Rq:1}, \eqref{eq:Phi:i:bnd:a}, \eqref{bd:rq(1)} and \eqref{estimate:rho1}, we have that on $ \supp(\overline \eta_i)$ \,
\begin{align*}
\norm{ R_{q,i}(t) - \Id }_{C^0} &\lesssim\|\nabla   \Phi_i \nabla   \Phi_i^T -\Id\|_{C^0}+\|  \nabla  \Phi_i\, \frac{\sum_j \int_{\mathbb{T}^3}\overline \eta_j^2(t,y) \dif y} {\Upsilon_q(t)}(\Rbar_q+ \mathring{R}^{(1)}_{q}) \nabla   \Phi_i^T \|_{C^0} \\
&\les l^{\alpha/2} \lambda_{q}^{\alpha/3}\leq  1/2 \,   ,
\end{align*}
where we used the fact that $l\lambda_q\ll 1$, and choose  $a$ large enough to absorb the universal constant.

 One last important property of the stress $ R_{q,i}$ is obtained by recalling \eqref{eq:eta:partition}
\begin{align}
  \sum_i\Upsilon_{q,i}(\nabla  \Phi_i)^{-1}  R_{q,i} (\nabla  \Phi_i)^{-T} = \sum_i \Upsilon_{q,i} \Id - \Rbar_q -\mathring{ R}^{(1)}_{q}.
  \label{eq:tilde:R:q:i:identity}
\end{align}

Thus, since $  R_{q,i}$ obeys the conditions of Lemma~\ref{l:linear_algebra} on $\supp (\overline\eta_i)$, for $\xi\in\overline{\Lambda}^1\cup\overline{\Lambda}^2$, we can define the amplitude functions as
\begin{align}
  a_{(\xi,i)} = \Upsilon_{q,i}^{1/2} \,  \gamma_{\xi} ( R_{q,i}) ,
 \label{eq:a:xi:i:def}
\end{align}
where the $\gamma_{ \xi}$ are the functions from Lemma~\ref{l:linear_algebra}. 

Then we have the following estimates for the amplitude functions $  a_{(\xi,i)}$.
\bp\label{lem:eq:Onsager:a:xi:CN} For $\xi\in\overline{\Lambda}^1\cup\overline{\Lambda}^2,i=1,2$ and  $N\in\mN_0$, we have
\begin{align}
\norm{  a_{(\xi,i)}}_{C^N} +  \tau_q \norm{D_{t,q}  a_{(\xi,i)}}_{C^N} \les \delta_{q+1}^{1/2}  l^{-N}.
\label{eq:Onsager:a:xi:CN}
\end{align}
\ep
We put the proof of this proposition in Appendix \ref{sec:est:ampl}.

Now we  define the  principle part of the perturbation $\overline w_{q+1}$. We  shall use the index sets $\overline  \Lambda^1,\overline  \Lambda^2$ from Lemma~\ref{l:linear_algebra}. We use the notation   $ \overline \Lambda^i = \overline  \Lambda^1$ for $i$ odd, and $ \overline \Lambda^i = \overline \Lambda^2$ for $i$ even. We use the amplitude functions $a_{(\xi,i)}$ in \eqref{eq:a:xi:i:def} to define
\begin{align}
\overline w_{q+1}^{(p)}(x,t) = \sum_{i} \sum_{ \xi \in \overline  \Lambda^i}  a_{(  \xi,i)}(x,t) (\nabla  \Phi_i(x,t))^{-1}  W_{(\xi)}( \Phi_i(x,t)) \, .
\label{eq:Onsager:obarw:q+1:p}
\end{align}

Using  \eqref{eq:Mikado:2},  \eqref{eq:Mikado:4}, \eqref{eq:tilde:R:q:i:identity},
\eqref{eq:wq+1wq+1} above 
and the fact that the $\overline \eta_i$ have mutually disjoint supports, we obtain
\begin{align}
( w_{q+1}^{(p)}&+ \overline w_{q+1}^{(p)}) \otimes ( w_{q+1}^{(p)}+ \overline w_{q+1}^{(p)})+ \Rbar_q=\overline w_{q+1}^{(p)} \otimes \overline  w_{q+1}^{(p)} + w_{q+1}^{(p)} \otimes   w_{q+1}^{(p)} + \Rbar_q \notag\\
&= \sum_i \sum_{ \xi \in \overline  \Lambda^i}   a_{( \xi,i)}^2 (\nabla  \Phi_i)^{-1} \( ( W_{( \xi)}( \Phi_i) ) \otimes ( W_{(  \xi)}( \Phi_i) ) \) (\nabla  \Phi_i)^{-T} + w_{q+1}^{(p)} \otimes   w_{q+1}^{(p)}  + \Rbar_q\notag\\
&= \sum_i \Upsilon_{q,i}   (\nabla  \Phi_i)^{-1}    \sum_{\xi \in \overline  \Lambda^i}  \gamma_\xi^2( R_{q,i})  \( (  W_{( \xi)} \otimes  W_{( \xi)})( \Phi_i)  \) (\nabla  \Phi_i)^{-T} + w_{q+1}^{(p)} \otimes   w_{q+1}^{(p)}  + \Rbar_q\notag\\
&= \sum_i  \Upsilon_{q,i} (\nabla  \Phi_i)^{-1}   R_{q,i} (\nabla  \Phi_i)^{-T} + w_{q+1}^{(p)} \otimes   w_{q+1}^{(p)}  + \Rbar_q\notag\\
&\qquad+ \sum_i \sum_{ \xi \in  \overline  \Lambda^i}  a_{(  \xi,i)}^2 (\nabla   \Phi_i)^{-1} \( \left( \Proj_{\neq 0}( W_{( \xi)} \otimes  W_{( \xi)}) \right)( \Phi_i) \)  (\nabla  \Phi_i)^{-T}   \notag\\
&=\sum_i\Upsilon_{q,i}  \Id+\frac13\tr R_q^{(1)}\Id+ \sum_i \sum_{ \xi \in \Lambda^i}  A_{(  \xi,i)}^2 (\nabla   \Phi_i)^{-1} \( \left( \Proj_{\neq 0}( W_{( \xi)} \otimes  W_{( \xi)}) \right)( \Phi_i) \)  (\nabla  \Phi_i)^{-T}   \notag\\ 
&\qquad+ \sum_i \sum_{ \xi \in   \overline \Lambda^i}  a_{(  \xi,i)}^2 (\nabla  \Phi_i)^{-1} \(\left(\Proj_{\neq0}( W_{( \xi)} \otimes  W_{( \xi)}) \right)( \Phi_i) \)  (\nabla  \Phi_i)^{-T}.
\label{eq:w:q+1:is:good}
\end{align}

Similar to \eqref{eq:Onsager:w:q+1:c}, we define the incompressibility corrector
\begin{align}
 \overline w_{q+1}^{(c)}(x,t) = \sum_i \sum_{ \xi \in   \overline \Lambda^i} \nabla  a_{( \xi,i)}(x,t) \times \( (\nabla  \Phi_i(x,t))^T ( V_{( \xi)} ( \Phi_i(x,t) )  \).
\label{eq:Onsager:barw:q+1:c}
\end{align}
Then the total velocity increment $\overline w_{q+1}$ is defined as
\begin{align}
 \overline w_{q+1} = \overline  w_{q+1}^{(p)} + \overline  w_{q+1}^{(c)} = \curl\(\sum_i \sum_{ \xi \in  \overline  \Lambda^i}  a_{( \xi,i)} \, (\nabla  \Phi_i)^T ( V_{( \xi)}( \Phi_i) ) \) 
\label{eq:Onsager:barw:q+1},
\end{align}
which is automatically incompressible. 
 Finally we define
\begin{align}
v_{q+1} =  \overline v_q + w_{q+1}+ \overline  w_{q+1},\  \ \rho_{q+1}=\overline \rho_q+\theta_{q+1}.
\label{eq:Onsager:v:q+1:def}
\end{align}

Moreover, recalling  the fact that $\tau_q\leq \lambda_q^{-1}\delta_q^{-1/2}\leq \tilde{\delta}_{q+1}^{1/2}$ due to $\tilde{\beta} b+\beta<1$,  we have $\theta_{q+1}=\overline \rho_q-1=0$ on $[0,T_{q+1}]\subset[0,T_q-\tau_q]$,  hence $\rho_{q+1} = 1$ on $[0,T_{q+1}]\subset[0,T_q-\tau_q]$.
 Moreover, it is easy to see that $\theta_{q+1}$ is mean-zero. Consequently, we have  $\int\rho_{q+1}\dif x=\int\overline \rho_{q}\dif x=1$.

To conclude this section, we note that 
\begin{align}
    \overline w_{q+1}\theta_{q+1}^{(p)}=0\label{eq:barwthe},
\end{align}
since the building blocks have disjoint supports.

\subsubsection{Estimate of perturbations}
\label{sec:Onsager:velocity:inductive}
In this section we establish the desired estimates on the perturbations and derive the estimates  \eqref{e:v_q_0}, \eqref{e:v_q_inductive_est}  and \eqref{e:v_diff_prop_est} at the level of $q+1$.

First we estimate the principle part of the perturbations $w_{q+1}^{(p)}, \theta_{q+1}^{(p)}$ and $\overline w_{q+1}^{(p)}$ as defined in \eqref{eq:Onsager:w:q+1:p}, \eqref{eq:Onsager:the:q+1:p} and \eqref{eq:Onsager:obarw:q+1:p} respectively.
From the estimates for $\nabla\Phi_i$, $A_{(  \xi,i)}$ and $\tilde{A}_{(  \xi,i)}$ given in \eqref{eq:Phi:i:bnd:a}, \eqref{eq:Onsager:A:xi:CN},  the estimates for the building blocks in \eqref{eq:Mikado:bounds}, and \eqref{eq:Onsager:tA:xi:CN} respectively,  and the fact that the $\overline\eta_i$ have disjoint supports we have for some $M>1$
\begin{align}
\norm{w_{q+1}^{(p)}}_{C^0} &\lesssim \sum_{i} \sum_{ \xi \in  \Lambda^i}  \| A_{(  \xi,i)}(\nabla  \Phi_i)^{-1}  W_{(\xi)}( \Phi_i)  \|_{C^0}\leq \frac{ M}{8} \delta_{q+1}^{1/2},\label{eq:Liverpool:1}\\
\| \theta_{q+1}\|_{C^0} \les \| \theta_{q+1}^{(p)}\|_{C^0} &\les \sum_{i} \sum_{ \xi \in  \Lambda^i}  \|\tilde{A}_{(  \xi,i)}\Theta_{(\xi)}( \Phi_i)\|_{C^0}\les\tilde{\delta}_{q+1}^{1/2}l^{\alpha/4}\leq \frac{ 1}{2} \delta_{q+1}^{1/2},\label{bdLtheq+1c0}
\end{align}
where in the last inequality we choose $a$ large enough to absorb the universal constant. Here and in the following the sum over$i$ is finite, since by definition the amplitude functions have disjoint supports for different $i$.

When considering the  $ C^1$-norm, from the bounds for $\nabla\Phi_i$, $A_{(  \xi,i)}$, and $\tilde{A}_{(  \xi,i)}$ in  \eqref{eq:Phi:i:bnd:b}, \eqref{eq:Onsager:A:xi:CN}, \eqref{eq:Onsager:tA:xi:CN} respectively, from \eqref{eq:Mikado:bounds} we  lose a factor   $l^{-1}$ from the gradient, and lose  a factor   $\lambda_{q+1}$ from the gradient of $W_{(\xi)}$, i.e.
\begin{align}
 \norm{w_{q+1}^{(p)}}_{C^1} &\leq  \sum_{i} \sum_{ \xi \in  \Lambda^i} \| A_{(  \xi,i)}(\nabla  \Phi_i)^{-1} W_{(\xi)}( \Phi_i)  \|_{C^1}\lesssim\delta_{q+1}^{1/2}( \lambda_{q+1}+ l^{-1}) \leq \frac{ M}{8} \delta_{q+1}^{1/2}\lambda_{q+1},\notag\\
   \| \theta_{q+1}\|_{C^1} &\les \| \theta_{q+1}^{(p)}\|_{C^1} \les \sum_{i} \sum_{ \xi \in  \Lambda^i}  \|\tilde{A}_{(  \xi,i)}\Theta_{(\xi)}( \Phi_i)\|_{C^1}\les \tilde{\delta}_{q+1}^{1/2}l^{\alpha/4}( \lambda_{q+1}+ l^{-1}) \leq  \frac{ 1}{2} \delta_{q+1}^{1/2}\lambda_{q+1},\label{bdLtheq+1c1}
\end{align}
where we  used that by choosing $\alpha>0$  small enough 
\begin{align}
\frac{l^{-1}}{\lambda_{q+1}} = \frac{\tilde{\delta}_q^{1/2} \lambda_q^{1+ \frac{3\alpha}{2}}}{\tilde{\delta}_{q+1}^{1/2} \lambda_{q+1}} = \frac{\lambda_q^{1-\tilde{\beta} + \frac{3\alpha}{2}}}{\lambda_{q+1}^{1-\tilde{\beta}}} \leq  \lambda_{q}^{\frac{3\alpha}{2} - (b-1)(1-\tilde{\beta})} \leq \lambda_{q}^{- \frac{(b-1)(1-\tilde{\beta})}{2}}  \ll 1 \, .
\label{eq:Onsager:ell:gap}
\end{align}
We also  choose $a$ large enough to absorb the universal constant in the last inequality.

Then we note that the definition of the perturbation $\overline w_{q+1}^{(p)}$ is similar to that of $w_{q+1}^{(p)}$ by replacing $\Lambda^i, A_{(\xi,i)}$ by  $\overline \Lambda^i, a_{(\xi,i)}$ respectively. Then by the estimate of $a_{(\xi,i)}$ in \eqref{eq:Onsager:a:xi:CN} and a similar calculation we have for $j=0,1$
\begin{align}
    \norm{\overline w_{q+1}^{(p)}}_{C^j}
    \leq \frac{ M}{8} \delta_{q+1}^{1/2}\lambda_{q+1}^j.\label{eq:barLiverpool:bar1}
\end{align}

For the incompressibility correctors $w_{q+1}^{(c)}$ and $\overline w_{q+1}^{(c)}$ defined in  \eqref{eq:Onsager:w:q+1:c} and \eqref{eq:Onsager:barw:q+1:c} respectively,  we see that from  the estimates for $\nabla\Phi_i$, $A_{(  \xi,i)}$ and $V_{(\xi)}$   in  \eqref{eq:Phi:i:bnd:b}, \eqref{eq:Onsager:A:xi:CN} and \eqref{eq:Mikado:bounds} respectively, we obtain that for $j=0,1$ and $M>1$
\begin{align}
\norm{w_{q+1}^{(c)}}_{C^j}&\lesssim \sum_i \sum_{ \xi \in  \Lambda^i} \norm{\nabla  A_{( \xi,i)} \times \left( (\nabla  \Phi_i)^T ( V_{( \xi)} ( \Phi_i )  \right)}_{C^j}\lesssim \delta_{q+1}^{1/2}\frac{ l^{-1} }{\lambda_{q+1}}( \lambda_{q+1}+ l^{-1})^j\leq \frac{M}{8} \delta_{q+1}^{1/2} \lambda_{q+1}^{j},\notag
\end{align}
where we used \eqref{eq:Onsager:ell:gap}. Similarly,
\begin{align}
    \norm{\overline w_{q+1}^{(c)}}_{C^j} \lesssim \delta_{q+1}^{1/2}\frac{ l^{-1} }{\lambda_{q+1}}( \lambda_{q+1}+ l^{-1})^j
   \leq \frac{M}{8}  \delta_{q+1}^{1/2} \lambda_{q+1}^{j}.\label{eq:Liverpool:bar2}
\end{align}
Combining all the bounds above, together with \eqref{e:v:ell:0},   \eqref{e:vq:vell:additional}, \eqref{e:vq:1} we have that
\begin{align*}
    \norm{v_{q+1}-v_q}_{C^0}   & \leq \norm{w_{q+1}+\overline w_{q+1}}_{C^0} +\norm{\overline v_{q}-v_l}_{C^0} +\norm{v_{l}-v_q}_{C^0}   \leq  \frac{M}{2} \delta_{q+1}^{1/2}+M' \delta_{q+1}^{1/2} \lambda_q^{-\alpha}\leq  M \delta_{q+1}^{1/2},\\
     \norm{v_{q+1}}_{C^1}    & \leq \norm{w_{q+1}+\overline w_{q+1}}_{C^1} +\norm{\overline v_{q}}_{C^1}   \leq  \frac{M}{2} \delta_{q+1}^{1/2}\lambda_{q+1}+M' \tau_q^{-1}l^{2\alpha}\leq  M \delta_{q+1}^{1/2}\lambda_{q+1},
\end{align*}
where $ M'$ is a universal constant and we note that by choosing  $a$ large enough, we have
\begin{align*}
    \tau_ql^{-2\alpha}\delta_{q+1}^{1/2}\lambda_{q+1}=  \lambda_{q}^{(b-1)(1-\tilde{\beta})}\gg 1.
\end{align*}
We also choose $a$ large enough to absorb the constant. Then \eqref{e:v_q_inductive_est} and  \eqref{e:v_diff_prop_est}  are satisfied  for  $v_{q+1}$, and, moreover, the bound  \eqref{e:v_q_0}  holds for  $v_{q+1}$.

By combining the above estimates with \eqref{e:rho:ell:0},   \eqref{e:rhoq:vell} and \eqref{e:rhoq:1}, we deduce that 
\begin{align*}
    \norm{\rho_{q+1}-\rho_q}_{C^0}
    \leq  \norm{\theta_{q+1}}_{C^0}+ \norm{\overline\rho_{q}-\rho_l}_{C^0}+ \norm{\rho_{l}-\rho_q}_{C^0}
    \leq  \frac12 \tilde{\delta}_{q+1}^{1/2}+M'\tilde{\delta}_{q+1}^{1/2} \lambda_q^{-\alpha}\leq   \tilde{\delta}_{q+1}^{1/2},\\
     \norm{\rho_{q+1}}_{C^1}
    \leq  \norm{\theta_{q+1}}_{C^1}+ \norm{\overline\rho_{q}}_{C^1}
    \leq  \frac12 \tilde{\delta}_{q+1}^{1/2}\lambda_{q+1}+M'\tilde{\delta}_{q}^{1/2}\lambda_q\leq   \tilde{\delta}_{q+1}^{1/2}\lambda_{q+1},
\end{align*}
where  $ M'$ is a universal constant and  we choose $a$ large enough to absorb the constant. 
Then \eqref{e:v_q_inductive_est}  and \eqref{e:v_diff_prop_est} are satisfied for  $\rho_{q+1}$, and, moreover,  the bound \eqref{e:v_q_0}   holds for  $\rho_{q+1}$.

\subsection{The estimate of the stress terms $(M_{q+1},\mathring{R}_{q+1})$
}\label{sec:defrmq+1}
In this section we aim to  establish the desired estimates on the  stress terms. We recall that the stress term $\overline{M}_q$ in the transport equation is canceled by the perturbation $(w_{q+1},\theta_{q+1})$, while  it brings a new term $w_{q+1}\otimes w_{q+1}$ into the Euler equations. Then  the perturbation $\overline w_{q+1}$ is used to cancel the stress term $\mathring{\overline R}_q$ and the low frequency part from $w_{q+1}\otimes w_{q+1}$. This procedure together with the  definition of new stress terms $M_{q+1}$ and  $\RR_{q+1}$is stated in Section \ref{sec:def:r,m}. Then we establish the desired estimates in Section \ref{sec:bd:mq+1} and \ref{sec:bd:rq+1} for  $M_{q+1}$ and $\RR_{q+1}$ respectively.

\subsubsection{The definition of stress terms $M_{q+1}$ and $\RR_{q+1}$} \label{sec:def:r,m}

In order to define $M_{q+1}$, we notice  that $(\overline v_q,\overline \rho_q,\overline M_q)$ obeys the transport equation \eqref{e:euler_trans}. Using the definition of the perturbation  $(v_{q+1},\rho_{q+1})$, along with  the fact that $$\div((w_{q+1}+\overline w_{q+1})\theta_{q+1}^{(c)})=\theta_{q+1}^{(c)}\div (w_{q+1}+\overline w_{q+1})=0,$$  and \eqref{eq:barwthe}, we have that
\begin{align}
-\div M_{q+1} 
&=D_{t,q} \theta_{q+1}+ \div((w_{q+1} +\overline w_{q+1})\theta_{q+1} -\overline M_q)+(w_{q+1}+\overline w_{q+1})\cdot  \nabla  \overline \rho_{q}\notag\\
&=D_{t,q} \theta_{q+1}+ \div(w_{q+1}\theta_{q+1}^{(p)} -\overline M_q)+(w_{q+1}+\overline w_{q+1})\cdot  \nabla  \overline \rho_{q}.\notag
\end{align}
Using the inverse divergence operator $\RSZ_1$ introduced in Section \ref{tamr}, and  \eqref{eq:wthe:q+1:is:good}  we define the transport error, Nash error and oscillation error respectively as
\begin{align}
  M_{ tr}: &= \RSZ_1 ( D_{t,q} \theta_{q+1} ) ,\ \  M_{ Nash}: = \RSZ_1 ((w_{q+1}+\overline w_{q+1})\cdot  \nabla  \overline \rho_{q} ) ,\notag\\
  M_{ osc} :&= \sum_i \sum_{ \xi \in  \Lambda^i}\RSZ_1\div\left(  A_{(  \xi,i)}\tilde{A}_{(  \xi,i)} (\nabla  \Phi_i)^{-1} \left( (\Proj_{\geq \frac{ \lambda_{q+1}}{2}}( W_{( \xi)} \Theta_{( \xi)}) )( \Phi_i) \right)\right)+w_{q+1}^{(c)}\theta_{q+1}^{(p)} .\notag
\end{align}
  Then 
\begin{align}
-M_{q+1}: =  M_{tr} +   M_{ osc} +   M_{ Nash}   \, .
\label{eq:M:Onsager:q+1:B}
\end{align}

Since $ \overline \rho_{q}-1=\theta_{q+1}=\theta^{(p)}_{q+1}=\tilde{A}_{(  \xi,i)}=0$ on 
$[0,T_{q+1}]$,  it follows that  $M_{q+1}=0$ on $[0,T_{q+1}]$.

In order to define $\RR_{q+1}$, recalling that  $(\overline v_q,\Rbar_q)$ solves the system \eqref{e:euler_reynolds}, we write
\begin{align}
\div \RR_{q+1} - \nabla p_{q+1}+\nabla \overline  p_q
&=D_{t,q} (w_{q+1}+\overline w_{q+1})+ \div((w_{q+1}+\overline w_{q+1})\otimes (w_{q+1}+\overline w_{q+1})+ \Rbar_q)\notag\\ 
&\qquad+(w_{q+1}+\overline w_{q+1})\cdot  \nabla  \overline v_{q}.\notag
\end{align}
Then using the inverse divergence operator $\RSZ$ introduced in Section \ref{tamr}, we define the  transport error and Nash error,  respectively, as 
\begin{align*}
  R_{ tr} :&= \RSZ ( D_{t,q} (w_{q+1}+\overline w_{q+1}) ) ,\ \  R_{ Nash} := \RSZ ( (w_{q+1}+\overline w_{q+1})\cdot  \nabla  \overline v_{q} ).
\end{align*}
For the oscillation error,  we use \eqref{eq:w:q+1:is:good} to define
\begin{align*}
     R_{ osc}:&=\mathcal{R}\div \left( \sum_i \sum_{ \xi \in \Lambda^i}  A_{(  \xi,i)}^2 (\nabla   \Phi_i)^{-1} \( \left( \Proj_{\neq 0}( W_{( \xi)} \otimes  W_{( \xi)}) \right)( \Phi_i) \)  (\nabla  \Phi_i)^{-T} \right) (:= R_{ osc,1}) \notag\\ 
&+ \mathcal{R}\div \left( \sum_i \sum_{ \xi \in   \overline \Lambda^i}  a_{(  \xi,i)}^2 (\nabla  \Phi_i)^{-1} \( \left(\Proj_{\neq0}( W_{( \xi)} \otimes  W_{( \xi)}) \right)( \Phi_i) \)  (\nabla  \Phi_i)^{-T} \right) (:= R_{ osc,2})  \\
 &+(w_{q+1}^{(c)}+\overline w_{q+1}^{(c)})\mathring{\otimes}(w_{q+1}^{(c)}+\overline w_{q+1}^{(c)})+2(w_{q+1}^{(p)}+\overline w_{q+1}^{(p)})\mathring{\otimes}_s(w_{q+1}^{(c)}+\overline w_{q+1}^{(c)}) (:= R_{ osc,3}) ,
\end{align*}
where  we use the notation $a\otimes_sb=\frac12(a\otimes b+b\otimes a)$, and $\mathring{\otimes}_s$ is defined as the trace-free part of the symmetric tensor. Here the divergence of first two terms in the last inequality in  \eqref{eq:w:q+1:is:good}  can be written as pressure terms, so we put them into $\nabla p_{q+1}$.

Then the Reynolds stress at the level of $q+1$ is defined by
\begin{align}
\RR_{q+1}: =  R_{tr} +   R_{ osc} +   R_{ Nash}   \, .
\label{eq:R:Onsager:q+1:B}
\end{align}

To estimate the above stress terms, we will apply the stationary phase bounds in Proposition \ref{prop:phase bounds}  to the building blocks defined in Appendix \ref{sec:Mikado}. More precisely, by \eqref{eq:Phi:i:bnd:b}, \eqref{eq:Mikado:bounds} and a similar argument as in \cite[Section 6.6.1]{BV}, we have for $a \in C^\infty(\mathbb{T}^d;\mR),b \in C^\infty(\mathbb{T}^d;\mR^d)$

\begin{align}
\norm{\RSZ ( b  ( W_{(\xi)} \circ\Phi_i )) }_{C^\alpha}  & +\lambda_{q+1}\norm{\RSZ ( b  ( V_{(\xi)} \circ\Phi_i ) )}_{C^\alpha}   +\norm{\RSZ \left( b \; \left( ( \Proj_{\geq \frac{\lambda_{q+1}}{2}}(\phi_{(\xi)}^2 ) ) \circ\Phi_i \right) \right) }_{C^\alpha}   \notag\\ 
&\qquad\qquad\les \frac{  \norm{b}_{C^0}}{\lambda_{q+1}^{1-\alpha}} +   \frac{ \norm{b}_{C^{m+\alpha}}+\norm{b}_{C^0} l^{-m-\alpha}}{\lambda_{q+1}^{m-\alpha}},\label{eq:Mikado:stationary:phase:2}\\
\norm{\RSZ_1 ( a  ( W_{(\xi)} \circ\Phi_i ) )}_{C^\alpha}  & +\norm{\RSZ_1 \left( a \; \left( ( \Proj_{\geq \frac{\lambda_{q+1}}{2}}(\phi_{(\xi)}^2 ) ) \circ\Phi_i \right) \right) }_{C^\alpha}   
\les \frac{  \norm{a}_{C^0}}{\lambda_{q+1}^{1-\alpha}} +   \frac{ \norm{a}_{C^{m+\alpha}}+\norm{a}_{C^0} l^{-m-\alpha}}{\lambda_{q+1}^{m-\alpha}}.\label{eq:Mikado:stationary:phase:3}
\end{align}

\subsubsection{The estimate for $M_{q+1}$}\label{sec:bd:mq+1}
In this section, we aim to  show the estimate of  $M_{q+1}$, which is defined in \eqref{eq:M:Onsager:q+1:B} by estimating the three errors above separately.

First, for the transport error $M_{tr}$, recalling that $\Phi_i$ satisfies the equation  \eqref{eq:Phi:i:def},   the definition of the perturbations $\theta_{q+1}$, and the fact the $\theta^{(c)}_{q+1}$ is a function of time, we  write
\begin{align}
\RSZ_1 ( D_{t,q} \theta_{q+1} ) = \RSZ_1 ( D_{t,q} \theta^{(p)}_{q+1} ) &= \sum_{i} \sum_{\xi \in \Lambda^i}  \RSZ_1 \left( ( D_{t,q} \tilde{ A}_{( \xi,i)} )  \Theta_{( \xi)}( \Phi_i) \right).\notag
\end{align}
By the estimate of the derivatives of $\tilde{A}_{(\xi,i)}$ in \eqref{eq:Onsager:tA:xi:CN} we have
\begin{align}
   \norm{D_{t,q}  \tilde{A}_{(\xi,i)}}_{C^{m+\alpha}} \les\tilde{\delta}_{q+1}^{1/2}\tau_q^{-1}   l^{-m-\alpha},
\notag
\end{align}
which together with \eqref{eq:Mikado:stationary:phase:3} implies that
\begin{align*}
\norm{M_{ tr}}_{C^\alpha} &\les\norm{\RSZ_1 ( D_{t,q} \theta^{(p)}_{q+1} ) }_{C^\alpha}\les \frac{\tilde{\delta}_{q+1}^{1/2}  \tau_q^{-1}}{\lambda_{q+1}^{1-\alpha}} \left( 1 + \frac{ l^{-m-\alpha}}{\lambda_{q+1}^{m-1}} \right) \lesssim \frac{\delta_{q+1}^{1/2} \tilde{\delta}_{q}^{1/2} \lambda_q}{\lambda_{q+1}^{1-3\alpha}},
\end{align*}
where we used  \eqref{eq:Onsager:ell:gap} to deduce $(l \lambda_{q+1})^{-1} \leq \lambda_q^{-\frac{(b-1)(1-\tilde{\beta})}{2}}$, and then choose $m$ large enough such that $(l\lambda_{q+1})^{-m}l^{-\alpha}\lambda_{q+1}\lesssim 1$ in the last inequality.  Here and in the following the sum over $i$ is finite, since by the definition the amplitude functions have disjoint supports for different $i$.

For the oscillation error $M_{osc}$, we recall the identity  $$\nabla(\Proj_{\geq \frac{\lambda_{q+1}}{2}}( \phi_{(\xi)}^2 ) )( \Phi_i)=(\nabla \Phi_i)^{T} \((\nabla\Proj_{\geq \frac{\lambda_{q+1}}{2}}( \phi_{(\xi)}^2 ) )( \Phi_i)\).$$ Therefore, we have
\begin{align}
 \div &\left( A_{(\xi,i)}\tilde{A}_{(\xi,i)}(\nabla \Phi_i)^{-1} \left( (\Proj_{\geq \frac{\lambda_{q+1}}{2}}(W_{(\xi)}\Theta_{(\xi)}) )( \Phi_i) \right)    \right)  
\notag\\
&= \div \left( A_{(\xi,i)}\tilde{A}_{(\xi,i)} (\nabla \Phi_i)^{-1} \xi\left( (\Proj_{\geq \frac{\lambda_{q+1}}{2}}(\phi_{(\xi)}^2) )( \Phi_i) \right)     \right)  \notag\\
&= (\Proj_{\geq \frac{\lambda_{q+1}}{2}}(\phi_{(\xi)}^2))( \Phi_i)
\div \left( A_{(\xi,i)}\tilde{A}_{(\xi,i)}(\nabla \Phi_i)^{-1}    \xi \right)\notag\\ 
&\quad+ A_{(\xi,i)}\tilde{A}_{(\xi,i)}\left( (\nabla\Proj_{\geq \frac{\lambda_{q+1}}{2}}(\phi_{(\xi)}^2) )( \Phi_i) \right) ^T(\nabla \Phi_i) (\nabla \Phi_i)^{-1}   \xi,    \notag
\end{align}
 which   by the fact that $(\xi \cdot \nabla) \Proj_{\geq \frac{\lambda_{q+1}}{2}}( \phi_{(\xi)}^2 )  = 0$ implies
\begin{align*}
     M_{ osc} =\sum_i \sum_{\xi \in \Lambda^i} \RSZ_1 \left(  (\Proj_{\geq \frac{\lambda_{q+1}}{2}}(\phi_{(\xi)}^2) )(\Phi_i)  \div \left(  A_{(\xi,i)}\tilde{A}_{(\xi,i)}  (\nabla  \Phi_i)^{-1}   \xi  \right) \right)  +w_{q+1}^{(c)}\theta_{q+1}^{(p)} .
\end{align*}
By the estimate for the derivatives of $A_{(\xi,i)},\tilde{A}_{(\xi,i)}$ and $(\nabla \Phi_i)^{-1}$ in \eqref{eq:Onsager:A:xi:CN}, \eqref{eq:Onsager:tA:xi:CN} and \eqref{eq:Phi:i:bnd:b} respectively, we have
\begin{align*}
  \norm{  \div(  A_{( \xi,i)}\tilde{A}_{(\xi,i)}   (\nabla  \Phi_i)^{-1}  \xi) }_{C^{m+\alpha}}\lesssim\delta_{q+1}^{1/2}\tilde{\delta}_{q+1}^{1/2} l^{-1-m-\alpha},
\end{align*}
 which together with  \eqref{eq:Mikado:stationary:phase:3}, the estimate of  $w_{q+1}^{(c)}$ in \eqref{eq:Liverpool:bar2}, the estimate of $\theta_{q+1}^{(p)}$ in \eqref{bdLtheq+1c0}, \eqref{bdLtheq+1c1} and interpolation implies that
 \begin{align*}
\norm{M_{osc}}_{C^\alpha} \les \frac{\delta_{q+1}^{1/2}\tilde{\delta}_{q+1}^{1/2}l^{-1}}{\lambda_{q+1}^{1-\alpha}}  \left( 1 + \frac{ l^{-m-\alpha}}{\lambda_{q+1}^{m-1}} \right) +\norm{w_{q+1}^{(c)}}_{C^\alpha}\norm{\theta_{q+1}^{(p)} }_{C^\alpha}\les  \frac{\delta_{q+1}^{1/2}\tilde{\delta}_{q+1}^{1/2}l^{-1}}{\lambda_{q+1}^{1-2\alpha}}  \les \frac{\delta_{q+1}^{1/2}\tilde{ \delta}_q^{1/2}\lambda_q }{\lambda_{q+1}^{1-4\alpha}  } . 
\end{align*}
Here we used again the  bound $(l \lambda_{q+1})^{-1} \leq \lambda_q^{-\frac{(b-1)(1-\tilde{\beta})}{2}}$ by \eqref{eq:Onsager:ell:gap}, and then choose $m$ large enough.

The Nash error is written  as $M_{ Nash} = \RSZ_1 ( w_{q+1}\cdot  \nabla \overline \rho_{q} )+\RSZ_1 ( \overline w_{q+1}\cdot  \nabla \overline \rho_{q} ),$
where 
\begin{align*}
      \RSZ_1 ( w_{q+1}\cdot  \nabla \overline \rho_{q} )&=\sum_{i} \sum_{ \xi \in  \Lambda^i}\RSZ_1(  ( \nabla   \overline \rho_{q} )^T   A_{(  \xi,i)}(\nabla  \Phi_i)^{-1}  W_{(\xi)}( \Phi_i))\\
      &+\RSZ_1  \left(  ( \nabla   \overline \rho_{q} )^T\nabla A_{( \xi,i)}\times( (\nabla  \Phi_i)^T ( V_{( \xi)} ( \Phi_i )  )\right),
     \end{align*}
 and   the term $\RSZ_1 ( \overline w_{q+1}\cdot  \nabla \overline \rho_{q} )$ can be written in a similar form with $A_{(\xi,i)},\Lambda^i$ replaced by $a_{(\xi,i)},\overline\Lambda^i$ respectively. 
 
     By the estimate of the derivatives of $A_{(\xi,i)}$, $\nabla   \overline \rho_{q}$, and $(\nabla \Phi_i)^{-1}$ in \eqref{eq:Onsager:A:xi:CN}, \eqref{e:rhoq:1}, and \eqref{eq:Phi:i:bnd:b} respectively, we have
     \begin{align*}
       \|  ( \nabla   \overline \rho_{q} )^T  A_{(  \xi,i)}(\nabla  \Phi_i)^{-1}\|_{C^{m+\alpha}}+ l\|    ( \nabla \overline \rho_{q} )^T\nabla A_{( \xi,i)}\times ((\nabla  \Phi_i)^T \cdot)\|_{C^{m+\alpha}}\lesssim \delta_{q+1}^{1/2} \tilde{\delta}_q^{1/2} \lambda_q l^{-m-\alpha},
     \end{align*}
    at which point we apply \eqref{eq:Mikado:stationary:phase:3} to  deduce that
\begin{align}
\norm{ \RSZ_1 ( w_{q+1}\cdot  \nabla \overline \rho_{q} )}_{C^\alpha}\les
   \frac{\delta_{q+1}^{1/2}\tilde{\delta}_q^{1/2} \lambda_q}{\lambda_{q+1}^{1-\alpha}} \left( 1 +\frac{ l^{-1}}{\lambda_{q+1}}+ \frac{ l^{-1-m-\alpha}}{\lambda_{q+1}^{m}} \right) \lesssim \frac{\delta_{q+1}^{1/2} \tilde{\delta}_{q}^{1/2} \lambda_q}{\lambda_{q+1}^{1-\alpha}}.\notag
\end{align}
Here we used again the  bound $(l \lambda_{q+1})^{-1} \leq \lambda_q^{-\frac{(b-1)(1-\tilde{\beta})}{2}}$ by \eqref{eq:Onsager:ell:gap}, and then choose $m$ large enough.

The second term  $\RSZ_1 ( \overline w_{q+1}\cdot  \nabla \overline \rho_{q} )$ is bounded by the same argument, but we omit the details. Furthermore, we have
\begin{align*}
     \| M_{ Nash}\|_{C^\alpha} \lesssim \norm{ \RSZ_1 ( w_{q+1}\cdot  \nabla \overline \rho_{q} )}_{C^\alpha}+\norm{ \RSZ_1 ( \overline w_{q+1}\cdot  \nabla \overline \rho_{q} )}_{C^\alpha}\les\frac{\delta_{q+1}^{1/2} \tilde{\delta}_{q}^{1/2} \lambda_q}{\lambda_{q+1}^{1-\alpha}}.\notag
\end{align*}

The above bounds together with \eqref{bd:para1/3} immediately imply the desired estimate \eqref{e:R_q_inductive_est} for  $M_{q+1}$:
\begin{align}
\norm{M_{q+1}}_{C^\alpha} \les \frac{\delta_{q+1}^{1/2} \tilde{\delta}_{q}^{1/2} \lambda_q}{\lambda_{q+1}^{1-4\alpha}}\leq\frac12\frac{\delta_{q+2}^{1/2}\tilde{\delta}_{q+2}^{1/2}}{\lambda_{q+1}^{3\alpha}},\notag
\end{align}
where the extra power of  $\lambda_{q+1}^{-\alpha}$ is used to absorb the implicit constant by choosing $a$ sufficiently large.

Additionally, when considering the   $ C^1$-norm, we lose a factor   $l^{-1}$ in the estimates of the derivatives of $A_{(\xi,i)},\tilde{A}_{(\xi,i)}$, $\nabla   \overline \rho_{q}$, and $(\nabla \Phi_i)^{-1}$ as presented in \eqref{eq:Onsager:A:xi:CN}, \eqref{eq:Onsager:tA:xi:CN}, \eqref{e:rhoq:1}, and \eqref{eq:Phi:i:bnd:b} respectively. Moreover, we lose a factor   $\lambda_{q+1}$ from the estimates of $W_{(\xi)},V_{(\xi)}$ in \eqref{eq:Mikado:bounds}. Consequently, we have 
\begin{align*}
    \norm{M_{q+1}}_{C^1} \les \frac{\delta_{q+1}^{1/2} \tilde{\delta}_{q}^{1/2} \lambda_q}{\lambda_{q+1}^{1-4\alpha}}\lambda_{q+1}\leq\frac12\frac{\delta_{q+2}^{1/2}\tilde{\delta}_{q+2}^{1/2}}{\lambda_{q+1}^{3\alpha-1}}.
\end{align*}

\subsubsection{The estimate for $\RR_{q+1}$}\label{sec:bd:rq+1}
 With the  new Reynolds stress $\RR_{q+1}$ established, we   shall see that it satisfies the estimate \eqref{e:R_q_inductive_est} at level $q+1$ by individually estimating the three errors mentioned above.

First, we consider the transport error $R_{tr}$. Using the definition of $w_{q+1}^{(p)}$ in \eqref{eq:Onsager:w:q+1:p}, and the  Lie-advection identity \eqref{eq:Onsager:Lie:advect} we rewrite the transport stress as 
\begin{align}
\RSZ \( D_{t,q} w^{(p)}_{q+1} \) &= \sum_{i} \sum_{\xi \in \Lambda^i} \RSZ\(  A_{( \xi,i)} (\nabla \overline v_q)^T (\nabla  \Phi_i)^{-1}  W_{( \xi)}( \Phi_i) \)  +  \RSZ \( ( D_{t,q}  A_{( \xi,i)} )  (\nabla  \Phi_i)^{-1}  W_{( \xi)}( \Phi_i) \).\notag
\end{align}
By the estimate of the derivatives of $A_{(\xi,i)}$, $\overline   v_q$, and $(\nabla \Phi_i)^{-1}$ in \eqref{eq:Onsager:A:xi:CN}, \eqref{e:vq:1} and \eqref{eq:Phi:i:bnd:b} respectively we have
\begin{align*}
    \| A_{( \xi,i)} (\nabla \overline  v_q)^T (\nabla  \Phi_i)^{-1} \|_{C^{m+\alpha}}\lesssim \delta_{q+1}^{1/2}  \tau_q^{-1}  l^{-m-\alpha},\ \
    \|( D_{t,q}  A_{( \xi,i)})  (\nabla  \Phi_i)^{-1}\|_{C^{m+\alpha}}\lesssim  \delta_{q+1}^{1/2}  \tau_q^{-1} l^{-m-\alpha}.
\end{align*}
Then by \eqref{eq:Mikado:stationary:phase:2} we obtain
\begin{align*}
\norm{\RSZ ( D_{t,q} w^{(p)}_{q+1} )}_{C^\alpha} &\les  \frac{\delta_{q+1}^{1/2}  \tau_q^{-1}}{\lambda_{q+1}^{1-\alpha}}\left( 1 + \frac{ l^{-m-\alpha}}{\lambda_{q+1}^{m-1}} \right) \lesssim \frac{\delta_{q+1} \tilde{\delta}_{q}^{1/2} \lambda_q}{ \tilde{\delta}_{q+1}^{1/2}\lambda_{q+1}^{1-3\alpha}},
\end{align*}
where we used the fact that  $(l \lambda_{q+1})^{-1} \leq \lambda_q^{-\frac{(b-1)(1-\tilde{\beta})}{2}}$ by \eqref{eq:Onsager:ell:gap}, and the last inequality was obtained by taking the parameter $m$  sufficiently large (in terms of $\tilde{\beta}$ and $b$).  Here and in the following the sum over $i$ is finite, since by the definition the amplitude functions have disjoint supports for distinct $i$

For the incompressibility corrector defined in \eqref{eq:Onsager:w:q+1:c}, since $\Phi_i$ obeys \eqref{eq:Phi:i:def}, we have 
\begin{align*}
   D_{t,q} w_{q+1}^{(c)}= \sum_i \sum_{ \xi \in  \Lambda^i}  D_{t,q} \nabla  A_{( \xi,i)}\times \left( (\nabla  \Phi_i)^T  V_{( \xi)} ( \Phi_i)  \right)+\sum_i \sum_{ \xi \in  \Lambda^i} \nabla  A_{( \xi,i)} \times \left( ( D_{t,q} \nabla  \Phi_i)^T  V_{( \xi)} ( \Phi_i )  \right).
\end{align*}
 By the estimates for $A_{(\xi,i)}$, $\overline   v_q$ in \eqref{eq:Onsager:A:xi:CN}, \eqref{e:vq:1}  respectively, we have 
\begin{align*}
      \| D_{t,q} \nabla  A_{( \xi,i)}\|_{C^{m+\alpha}}&\lesssim   \| D_{t,q}   A_{( \xi,i)}\|_{C^{m+1+\alpha}}+  \|\overline v_q \|_{C^{1+\alpha}}\|A_{( \xi,i)}\|_{C^{m+1+\alpha}}+  \|\overline v_q \|_{C^{m+1+\alpha}}\|A_{( \xi,i)}\|_{C^{1+\alpha}}\\
      &\lesssim \delta_{q+1}^{1/2}  \tau_q^{-1} l^{-1-m-\alpha},
 \end{align*}
which together with the estimate for $(\nabla \Phi_i)^{T}$ in  \eqref{eq:Phi:i:bnd:b} implies that 
\begin{align*}
   \| D_{t,q} \nabla  A_{( \xi,i)}\times \left( (\nabla  \Phi_i)^T \cdot\right)\|_{C^{m+\alpha}}+ \| \nabla  A_{( \xi,i)}\times \left( (D_{t,q} \nabla  \Phi_i)^T \cdot\right)\|_{C^{m+\alpha}}\les  \delta_{q+1}^{1/2}  \tau_q^{-1} l^{-1-m-\alpha},
\end{align*}
while we gain a factor   $\lambda_{q+1}^{-1}$ from $V_{(\xi)}$ compared with $W_{(\xi)}$. Then by applying \eqref{eq:Mikado:stationary:phase:2} we have 
\begin{align*}
    \norm{\RSZ ( D_{t,q}w^{(c)}_{q+1} )}_{C^\alpha}\les  \frac{\delta_{q+1}^{1/2}  \tau_q^{-1}}{\lambda_{q+1}^{1-\alpha}} \frac{l^{-1}}{\lambda_{q+1}}\left( 1 + \frac{ l^{-m-\alpha}}{\lambda_{q+1}^{m-1}} \right)\lesssim\frac{\delta_{q+1} \tilde{\delta}_{q}^{1/2} \lambda_q}{ \tilde{\delta}_{q+1}^{1/2}\lambda_{q+1}^{1-3\alpha}},
\end{align*}
where we  used the fact that  $(l \lambda_{q+1})^{-1} \leq \lambda_q^{-\frac{(b-1)(1-\tilde{\beta})}{2}}$ by \eqref{eq:Onsager:ell:gap}, and the last inequality was obtained by taking the parameter $m$  sufficiently large. 

Then   compared to  $\RSZ ( D_{t,q}  w^{(p)}_{q+1} )$ and  $\RSZ ( D_{t,q}  w^{(c)}_{q+1} )$ , we obverse that the definition of $\RSZ ( D_{t,q} \overline w^{(p)}_{q+1} )$ and $\RSZ ( D_{t,q} \overline w^{(c)}_{q+1} )$  can be obtained by replacing $\Lambda^i, A_{(\xi,i)}$ by  $\overline \Lambda^i, a_{(\xi,i)}$ respectively. Then by the estimate of $a_{(\xi,i)}$ in \eqref{eq:Onsager:a:xi:CN} we have the same bound as above:
\begin{align*}
\norm{\RSZ ( D_{t,q}\overline w^{(p)}_{q+1} )}_{C^\alpha}+\norm{\RSZ ( D_{t,q}\overline w^{(c)}_{q+1} )}_{C^\alpha} \lesssim \frac{\delta_{q+1} \tilde{\delta}_{q}^{1/2} \lambda_q}{ \tilde{\delta}_{q+1}^{1/2}\lambda_{q+1}^{1-3\alpha}}.
\end{align*}

To address the oscillation error $R_{osc}$, we begin by bounding $R_{osc,1}$, and the second term $R_{osc,2}$ can be bounded by a similar argument. We recall the identity    $$\nabla(\Proj_{\geq \frac{\lambda_{q+1}}{2}}( \phi_{(\xi)}^2 ) )( \Phi_i)=(\nabla \Phi_i)^{T} \((\nabla\Proj_{\geq \frac{\lambda_{q+1}}{2}}( \phi_{(\xi)}^2 ) )( \Phi_i)\).$$ Therefore, we have
\begin{align*}
& \div \left( A_{(\xi,i)}^2 (\nabla \Phi_i)^{-1} \left( (\Proj_{\geq \frac{\lambda_{q+1}}{2}}(W_{(\xi)} \otimes W_{(\xi)}) )(\Phi_i) \right)  (\nabla \Phi_i)^{-T}   \right)  
\notag\\
&= \div \left( A_{(\xi,i)}^2 (\nabla \Phi_i)^{-1} (\xi \otimes \xi) (\nabla \Phi_i)^{-T} \left( (\Proj_{\geq \frac{\lambda_{q+1}}{2}}(\phi_{(\xi)}^2) )( \Phi_i) \right)     \right)  \notag\\
&=\left((\Proj_{\geq \frac{\lambda_{q+1}}{2}}(\phi_{(\xi)}^2) )(\Phi_i) \right)    \div \left( A_{(\xi,i)}^2 (\nabla \Phi_i)^{-1}    (\xi \otimes \xi) (\nabla \Phi_i)^{-T}\right)\\
&\qquad +  A_{(\xi,i)}^2 (\nabla \Phi_i)^{-1} (\xi \otimes \xi) (\nabla \Phi_i)^{-T} (\nabla \Phi_i)^{T} \left((\nabla \Proj_{\geq \frac{\lambda_{q+1}}{2}}(\phi_{(\xi)}^2) )( \Phi_i) \right),
\end{align*}
where the last term equals to 0 since $(\xi \cdot \nabla) \Proj_{\geq \frac{\lambda_{q+1}}{2}}( \phi_{(\xi)}^2 )  = 0$.
By the estimate of the derivatives of $A_{(\xi,i)}$ in \eqref{eq:Onsager:A:xi:CN} and $(\nabla \Phi_i)^{-1}$ in \eqref{eq:Phi:i:bnd:b} we have
\begin{align*}
  \|  \div(  A_{( \xi,i)}^2 (\nabla  \Phi_i)^{-1} ( \xi\otimes  \xi) (\nabla  \Phi_i)^{-T} )\|_{C^{m+\alpha}}\lesssim\delta_{q+1} l^{-1-m-\alpha}.
\end{align*}
Then we  apply \eqref{eq:Mikado:stationary:phase:2} to obtain
\begin{align*}
\norm{R_{osc,1}}_{C^\alpha} \les \frac{\delta_{q+1} l^{-1}}{\lambda_{q+1}^{1-\alpha}}  \left( 1 + \frac{ l^{-m-\alpha}}{\lambda_{q+1}^{m-1}} \right) \les \frac{\delta_{q+1} \tilde{\delta}_{q}^{1/2} \lambda_q}{ \tilde{\delta}_{q+1}^{1/2}\lambda_{q+1}^{1-3\alpha}} \, ,
\end{align*}
where we  used the fact that  $(l \lambda_{q+1})^{-1} \leq \lambda_q^{-\frac{(b-1)(1-\tilde{\beta})}{2}}$ by \eqref{eq:Onsager:ell:gap}, and then choose $m$  large enough. 

By replacing $\Lambda^i, A_{(\xi,i)}$ with  $\overline \Lambda^i, a_{(\xi,i)}$ and applying a similar argument as before, we obtain
\begin{align*}
\norm{R_{osc,2}}_{C^\alpha}\les \frac{\delta_{q+1} \tilde{\delta}_{q}^{1/2} \lambda_q}{ \tilde{\delta}_{q+1}^{1/2}\lambda_{q+1}^{1-3\alpha}} \, . 
\end{align*}
The last term $R_{osc,3}$ is estimated by using the bounds for  the perturbations in \eqref{eq:Liverpool:1}-\eqref{eq:Liverpool:bar2} and  the definition of $ l$ in \eqref{e:ell_def} as
\begin{align*}
      \norm{R_{osc,3}}_{C^\alpha} \les   \|w_{q+1}^{(c)}+\overline w_{q+1}^{(c)}\|_{C^\alpha}\(\|w_{q+1}^{(c)}+\overline w_{q+1}^{(c)}\|_{C^\alpha}+\|w_{q+1}^{(p)}+\overline w_{q+1}^{(p)}\|_{C^\alpha}\)\les\frac{\delta_{q+1}l^{-1} }{\lambda_{q+1}^{1-2\alpha}}\les\frac{\delta_{q+1} \tilde{\delta}_{q}^{1/2} \lambda_q}{ \tilde{\delta}_{q+1}^{1/2}\lambda_{q+1}^{1-4\alpha}}.
\end{align*}

In the end, we only need to estimate the Nash error $R_{ Nash} = \RSZ ( w_{q+1}\cdot  \nabla  \overline v_{q} )+ \RSZ (\overline w_{q+1}\cdot  \nabla  \overline v_{q} )$. For the first term, due to the definition of the perturbation $w_{q+1}$, we have
\begin{align*}
      \RSZ ( w_{q+1}\cdot  \nabla \overline  v_{q} )&=\sum_{i} \sum_{ \xi \in  \Lambda^i}\RSZ (  ( \nabla \overline  v_{q})^T   A_{(  \xi,i)}(\nabla  \Phi_i)^{-1}  W_{(\xi)}( \Phi_i))\\
      &+\RSZ  \left(  ( \nabla \overline  v_{q})^T\nabla  A_{( \xi,i)}\times( (\nabla  \Phi_i)^T V_{( \xi)} ( \Phi_i )  \right).
     \end{align*}
     By the estimate of the derivatives of $A_{(\xi,i)}$, $\nabla \overline  v_q$, and $(\nabla \Phi_i)^{-1}$ in \eqref{eq:Onsager:A:xi:CN}, \eqref{e:vq:1} and \eqref{eq:Phi:i:bnd:b} respectively, we have
     \begin{align*}
       \|  ( \nabla  \overline  v_{q})^T   A_{(  \xi,i)}(\nabla  \Phi_i)^{-1}\|_{C^{m+\alpha}}+ l\|    ( \nabla  \overline  v_{q})^T\nabla  A_{( \xi,i)}\times( (\nabla  \Phi_i)^T \cdot)\|_{C^{m+\alpha}}\lesssim  \delta_{q+1}^{1/2}  \tau_q^{-1} l^{-m-\alpha}.
     \end{align*}
     Then we use \eqref{eq:Mikado:stationary:phase:2} to show that
\begin{align}
    \|  \RSZ ( w_{q+1}\cdot  \nabla \overline  v_{q} )\|_{C^0} \lesssim  \frac{\delta_{q+1}^{1/2}  \tau_q^{-1} }{ \lambda_{q+1}^{1-\alpha}} \left( 1 +\frac{ l^{-1}}{\lambda_{q+1}}+ \frac{ l^{-1-m-\alpha}}{\lambda_{q+1}^{m}} \right) \lesssim \frac{\delta_{q+1} \tilde{\delta}_{q}^{1/2} \lambda_q}{ \tilde{\delta}_{q+1}^{1/2}\lambda_{q+1}^{1-3\alpha}},\notag
\end{align}
where we  used the fact that  $(l \lambda_{q+1})^{-1} \leq \lambda_q^{-\frac{(b-1)(1-\tilde{\beta})}{2}}$ by \eqref{eq:Onsager:ell:gap}, and the last inequality was obtained by taking the parameter $m$  sufficiently large. 

Then the expression of  $\RSZ ( \overline w_{q+1}\cdot  \nabla \overline  v_{q} )$ is similar to $\RSZ (  w_{q+1}\cdot  \nabla \overline  v_{q} )$ with $\Lambda^i, A_{(\xi,i)}$ replaced  by  $\overline \Lambda^i, a_{(\xi,i)}$ respectively. Then by the estimate of $a_{(\xi,i)}$ in \eqref{eq:Onsager:a:xi:CN} we deduce
\begin{align*}
\norm{\RSZ ( \overline w_{q+1}\cdot  \nabla \overline  v_{q} )}_{C^\alpha} \lesssim \frac{\delta_{q+1} \tilde{\delta}_{q}^{1/2} \lambda_q}{ \tilde{\delta}_{q+1}^{1/2}\lambda_{q+1}^{1-3\alpha}}.
\end{align*}
  By all the bounds above together with \eqref{bd:para1/3} we have
\begin{align}
\norm{\RR_{q+1}}_{C^\alpha} \les  \frac{\delta_{q+1} \tilde{\delta}_{q}^{1/2} \lambda_q}{ \tilde{\delta}_{q+1}^{1/2}\lambda_{q+1}^{1-4\alpha}}\leq \frac{\delta_{q+2}}{\lambda_{q+1}^{3\alpha}},
\notag
\end{align}
where  the extra power of  $\lambda_{q+1}^{-\alpha}$ is used to absorb the implicit constant by choosing $a$ sufficiently large.

\subsection{Estimates on the energy}\label{estimate_energy}
	To conclude the proof of Proposition~\ref{prop:1/3}, we have to check that the iterative condition on the energy \eqref{estimate:energy} holds at the level $q+1$.  To this end, by the definition of $v_{q+1}$ in \eqref{eq:Onsager:v:q+1:def}, we first have
		\begin{align}
				  e(t )-\frac{\delta_{q+2}}{2}-\|v_{q+1}(t )\|_{L^2}^2&= e(t )-\frac{\delta_{q+2}}{2}-\|\overline{v}_q(t )\|_{L^2}^2-
				  \int_{\mathbb{T}^3} |w_{q+1}^{(p)}+ \overline w_{q+1}^{(p)}|^2(t )\dif x \notag
				\\& \quad -     \int_{\mathbb{T}^3}\[|w_{q+1}^{(c)}+ \overline w_{q+1}^{(c)}|^2+2(w_{q+1}^{(p)}+ \overline w_{q+1}^{(p)})\cdot (w_{q+1}^{(c)}+ \overline w_{q+1}^{(c)})\](t ) \dif x\notag\\ 
                &\quad-\int_{\mathbb{T}^3}2\overline{v}_q\cdot (w_{q+1}+\overline w_{q+1}) (t )\dif x.\label{eq:e-vq+1}
		\end{align}
		For the integrand in the first line on the right hand side, by taking the trace on  both sides of \eqref{eq:w:q+1:is:good} and using the fact that $ \mathring{\overline{R}}_q $ is traceless, we deduce that 
        \begin{align}
|w_{q+1}^{(p)}&+ \overline w_{q+1}^{(p)}|^2\notag\\
&=3\sum_i\Upsilon_{q,i}  +\tr R_q^{(1)}+ \sum_i \sum_{ \xi \in \Lambda^i} \tr\left[ A_{(  \xi,i)}^2 (\nabla   \Phi_i)^{-1} \( \left( \Proj_{\neq 0}( W_{( \xi)} \otimes  W_{( \xi)}) \right)( \Phi_i) \)  (\nabla  \Phi_i)^{-T}  \right] \notag\\ 
&\qquad+ \sum_i \sum_{ \xi \in   \overline \Lambda^i} \tr\left[  a_{(  \xi,i)}^2 (\nabla  \Phi_i)^{-1} \( \left(\Proj_{\neq0}( W_{( \xi)} \otimes  W_{( \xi)}) \)( \Phi_i) \right)  (\nabla  \Phi_i)^{-T}\right].\notag
\end{align}
		By integrating on both sides, together with the fact that $\sum_i  \int_{\mathbb{T}^3} \Upsilon_{q,i} \dif x =\Upsilon_q$, the definitions in \eqref{def:energy gap}, \eqref{def:rho q,i}, and the definition  of the building blocks we obtain for $t\in[0,1]$
		\begin{align*}
		  	\int_{  \mathbb{T}^3}	|(w_{q+1}^{(p)}+ \overline w_{q+1}^{(p)})(t )|^2 \dif x&=     e(t )-\frac{\delta_{q+2}}{2}-\|\overline{v}_{q}(t )\|_{L^2}^2
			\\+   \sum_i \sum_{\xi \in \Lambda^i} & \int_{  \mathbb{T}^3} \tr \left[ A_{(\xi,i)}^2 (\nabla \Phi_i)^{-1} \xi \otimes \xi  (\nabla \Phi_i)^{-T}  \right]\left( ( \mathbb{P}_{\neq0}(\phi_{(\xi)}^2))(\Phi_i) \right) (t  ) \dif x
            \\+   \sum_i \sum_{\xi \in \overline\Lambda^i} & \int_{  \mathbb{T}^3} \tr \left[ a_{(\xi,i)}^2 (\nabla \Phi_i)^{-1} \xi \otimes \xi  (\nabla \Phi_i)^{-T}  \right]\left( ( \mathbb{P}_{\neq0}(\phi_{(\xi)}^2))(\Phi_i) \right) (t  ) \dif x.
		\end{align*}
		Then we use the estimates for $a_{(\xi,i)}$, $A_{(\xi,i)}$,  and $(\nabla \Phi_i)^{-1}$ in \eqref{eq:Onsager:a:xi:CN},  \eqref{eq:Onsager:A:xi:CN} and \eqref{eq:Phi:i:bnd:b} respectively, to obtain for any $t\in[0,1] $
        \begin{align*}
    \| \tr[ a_{( \xi,i)}^2 (\nabla  \Phi_i)^{-1} ( \xi\otimes  \xi) (\nabla  \Phi_i)^{-T}] \|_{C^{m}}+\| \tr[ A_{( \xi,i)}^2 (\nabla  \Phi_i)^{-1} ( \xi\otimes  \xi) (\nabla  \Phi_i)^{-T}] \|_{C^{m}}\lesssim\delta_{q+1} l^{-m},
\end{align*}
which together with  \eqref{esti:integral} with $N=1$ implies that 
		\begin{align*}
         &\left| e(t )-\frac{\delta_{q+2}}{2}-\|\overline{v}_{q}(t )\|_{L^2}^2 -\int_{  \mathbb{T}^3}	|(w_{q+1}^{(p)}+ \overline w_{q+1}^{(p)})(t )|^2 \dif x\right|\lesssim \frac{\delta_{q+1}l^{-1}}{\lambda_{q+1}} \lesssim\frac{\delta_{q+1} \tilde{\delta}_{q}^{1/2} \lambda_q}{ \tilde{\delta}_{q+1}^{1/2}\lambda_{q+1}^{1-2\alpha}},
		\end{align*}
		 where we recall that the sum over $i$ is finite, since by the definition the amplitude functions have disjoint supports for different $i$.
		
		Then the second term on the right hand side of \eqref{eq:e-vq+1} is  estimated by the bounds for the perturbations in \eqref{eq:Liverpool:1}-\eqref{eq:Liverpool:bar2}. We have for $t\in[0,1]$
		\begin{align*}
		  	&\left| 	     \int_{\mathbb{T}^3}\[|w_{q+1}^{(c)}+ \overline w_{q+1}^{(c)}|^2+2(w_{q+1}^{(p)}+ \overline w_{q+1}^{(p)})\cdot (w_{q+1}^{(c)}+ \overline w_{q+1}^{(c)})\](t ) \dif x\right| \\
            &\qquad\les \|w_{q+1}^{(c)}+ \overline w_{q+1}^{(c)}\|_{C^0}^2+\|w_{q+1}^{(p)}+ \overline w_{q+1}^{(p)}\|_{C^0}\|w_{q+1}^{(c)}+ \overline w_{q+1}^{(c)}\|_{C^0}\lesssim    \frac{\delta_{q+1}l^{-1}}{\lambda_{q+1}}
		\lesssim  \frac{\delta_{q+1} \tilde{\delta}_{q}^{1/2} \lambda_q}{ \tilde{\delta}_{q+1}^{1/2}\lambda_{q+1}^{1-2\alpha}} .
		\end{align*}
		
		For the last term on the right hand side of \eqref{eq:e-vq+1},  we recall that $w_{q+1}+\overline w_{q+1}$ can be written as
		\begin{align*}
			w_{q+1}+\overline w_{q+1}  = \curl \( \sum_i \sum_{\xi \in \Lambda^i}  A_{(\xi,i)} \, (\nabla \Phi_i)^T (V_{(\xi)}( \Phi_i) ) +\sum_i \sum_{ \xi \in  \overline  \Lambda^i}  a_{( \xi,i)} \, (\nabla  \Phi_i)^T ( V_{( \xi)}( \Phi_i) ) \). 
		\end{align*}
		Then we use integration by parts,  the estimates for $\overline v_q$ in \eqref{e:vq:1},  the estimates for $a_{(\xi,i)}$, $A_{(\xi,i)}$,  and $(\nabla \Phi_i)^{-1}$ in \eqref{eq:Onsager:a:xi:CN},  \eqref{eq:Onsager:A:xi:CN} and \eqref{eq:Phi:i:bnd:b} respectively, and the bound on the building block in \eqref{eq:Mikado:bounds}  to obtain  for any $t\in[0,1]$
		\begin{align*}
			&       \left|\int_{\mathbb{T}^3}\overline{v}_q\cdot (w_{q+1}+\overline w_{q+1}) (t )\dif x  \right| 
			\\& \lesssim \( \sum_i \sum_{\xi \in \Lambda^i} \| A_{(\xi,i)} \, (\nabla \Phi_i)^T (V_{(\xi)}(\Phi_i))\| _{C^0}+\sum_i \sum_{\xi \in \overline\Lambda^i} \| a_{(\xi,i)} \, (\nabla \Phi_i)^T (V_{(\xi)}(\Phi_i))\| _{C^0}\)\|\overline v_q\|_{C^1}
			\\ &\lesssim \frac{ \delta_{q+1}^{1/2}\tau_q^{-1}}{\lambda_{q+1}} \lesssim\frac{\delta_{q+1} \tilde{\delta}_{q}^{1/2} \lambda_q}{ \tilde{\delta}_{q+1}^{1/2}\lambda_{q+1}^{1-2\alpha}}. 
		\end{align*}
		Combining the above estimates  we obtain by \eqref{bd:para1/3} that
        \begin{align*}
	\left|  	e(t)-\frac{\delta_{q+2}}{2}-\|v_{q+1}(t)\|_{L^2}^2  \right|  \lesssim\frac{\delta_{q+1} \tilde{\delta}_{q}^{1/2} \lambda_q}{ \tilde{\delta}_{q+1}^{1/2}\lambda_{q+1}^{1-2\alpha}}\leq \delta_{q+2}\lambda_{q+1}^{-\alpha},
		\end{align*}
        where we the extra power of $\lambda_{q+1}^{-\alpha}$ is used to absorb the constant by choosing $a$ large enough. Then we obtain \eqref{estimate:energy}  at the level $q+1$ by choosing $a$ large enough again:
	\begin{align*}
		\delta_{q+2}\lambda_{q+1}^{-\alpha/3} \leq\frac{\delta_{q+2}}{2}- \delta_{q+2}\lambda_{q+1}^{-\alpha} \leq    e(t )- \|v_{q+1}(t )\|_{L^2}^2 \leq\frac{ \delta_{q+2}}2+\delta_{q+2}\lambda_{q+1}^{-\alpha}\leq \delta_{q+2}.
	\end{align*}
	
This completes the proof of Proposition \ref{prop:1/3}.

\section{Construction of non-unique solutions in $L_t^1W^{1,s}$ scales}\label{cogpss2}
In this section, our primary objective is to establish the non-uniqueness of stochastic Lagrangian trajectories for solutions to the Navier-Stokes or Euler equations that possess a specific level of Sobolev regularity, as stated in Theorem \ref{thm:nonuniode+sde}. By applying the superposition principle, it suffices to demonstrate the non-uniqueness of the corresponding Fokker-Planck equations, which is exactly the claim of Theorem \ref{thm:non_pde_23}.
More precisely, for any triple  $(p,r,s)\in \mathcal{A}$, we construct a solution $v\in L^r_tL^p\cap L^2_tL^2\cap L^1_tW^{1,s}\cap C_tL^1$ to the Navier-Stokes equations \eqref{eq:ns} such that the related advection-diffusion  equations \eqref{eq:fpe} admit two positive solutions in the space $L^r_tL^p\cap L^2_tL^2\cap C_tL^1$  with initial data $\rho_0 = 1$.  To achieve this, we apply the convex integration method again. However, in this section, in contrast to Section \ref{sec:Euler:C1/3}, we  will add more intermittency in the building blocks in both the spatial and temporal directions, and handle  the advection-diffusion  equations \eqref{eq:fpe} and the Navier-Stokes equations \eqref{eq:ns} at different frequency scales. 

Without loss of generality, we set $T=1$. The convex integration iteration is indexed
by a parameter $q \in \mathbb{N}_0$. We consider an increasing sequence $\{\lambda_q\}_{q\in\mathbb{N}_0}\subset\mathbb{N}$ which diverges to $\infty$, and a sequence $\{\delta_q\}_{q\in\mathbb{N}_0}\subset (0,1]$
 which is decreasing to 0. 
 Let
$$\lambda_{q }=a^{(b^{q })},q\geq0,\ \ \ \delta_q=\frac1{48^2}\lambda_2^{2\beta}\lambda_q^{-2\beta},q\geq2,\ \ \delta_0=1,\ \ \delta_1=\frac1{48^2}.$$
where $\beta>0$ will be chosen sufficiently small and $a,b$ will be chosen sufficiently large.  In addition, we use that $\sum_{q\geq1} \delta_q^{1/2}\leq \frac1{48}(1+ \sum_{q\geq2}a^{(2-q)b\beta})\leq \frac1{48}(1+\frac{1}{1-a^{-b\beta}})< \frac1{16}$
which boils down to
\begin{align}\label{ieq:2ab2}
a^{b\beta}>2,
\end{align}
assumed from now on.

At each step $q$, a pair $(v_q,\rho_q, \mathring{R}_q,M_q)$ is constructed solving the following system  on $[0,1]$:
\begin{align}\label{eq:2qth}
\partial_t \rho_q-\kappa\Delta \rho_{q}+\div(v_q \rho_q)&=-\div M_q,\notag\\
\partial_tv_q+\div(v_q\otimes v_q)-\nu\Delta v_q+\nabla\pi_{q}
&=\div\mathring{R}_q,\ \
\div v_q=0.
\end{align}
 where $\kappa,\nu\in[0,1]$, $\mathring{R}_q$ is  a trace-free symmetric matrix and $M_q$ is some vector field.

As before, to handle the initial condition, we
 let $T_q := \frac13 - \sum_{ 1\leq r\leq q} \delta_r^{1/2}\in(0,\frac13]$. Here we define $\sum_{1\leq r\leq 0}:=0.$ 

Under the above assumptions, our main iteration reads as follows:

\bp\label{prop:case2}
Let $d\geq2$. For  any triple $(p,r,s)\in \mathcal{A}$, there exists a choice of parameters $a,b,\beta$ such
that the following holds: Let $(v_q,\rho_q, \mathring{R}_q, M_q)$ be a solution to \eqref{eq:2qth} satisfying $\int\rho_q\dif x=1$,
\begin{align}\label{bd:2vql2}
\|v_q\|_{L^2_tL^2}\leq  {C}_vC_0\sum_{m=1}^q\delta_m^{1/2},\ \ \|\rho_q\|_{L^2_tL^2}\leq C_\rho  C_0\sum_{m=0}^{q+1}\delta_m^{1/2}
\end{align}
for some universal constants $C_0,C_v, C_\rho\geq1$, and
\begin{align}
\|v_q\|_{C_{t,x}^2}\leq C_0\lambda_q^{4d+3},\label{bd:2vqc1}\ \ \|\rho_{q}\|_{C_{t,x}^1}\leq C_0\lambda_q^{3d+2},&\\
\|\mathring{R}_q\|_{L^1_tL^1}\leq C_0^2\delta_{q+1},\ \ \|M_q\|_{L^1_tL^1}\leq C_0^2\delta^2_{q+2},&\label{bd:2rql1}\\
 \rho_q-1= v_q=\mathring{R}_q=M_q=0\ {\rm on}\ [0,T_q].&\label{bd:2mq=rhoq=0}
\end{align}
Then there exists 
$(v_{q+1},\rho_{q+1}, \mathring{R}_{q+1}, M_{q+1})$ which solves \eqref{eq:2qth} and satisfies \eqref{bd:2vql2}-\eqref{bd:2mq=rhoq=0} at the level $q+1$ and
\begin{align}
\|v_{q+1}-v_{q}\|_{L^2_tL^2}\leq
 {C}_vC_0\delta_{q+1}^{1/2}\label{bd:2vq+1-vql2},\ \ \|\rho_{q+1}-\rho_{q}\|_{L^2_tL^2}\leq
C_\rho  C_0 \delta_{q+2}^{1/2}.
\end{align}
Moreover,
\begin{align}
\|v_{q+1}-v_{q}\|_{L_t^rL^p}\leq \delta_{q+1}^{1/2},\ \ \|v_{q+1}-v_{q}\|_{C_tL^1}\leq \delta_{q+1}^{1/2},\label{bd:2vq+1-vqlpr}\ \
\|v_{q+1}-v_{q}\|_{L_t^1W^{1,s}}\leq \delta_{q+1}^{1/2},\\
 \|\rho_{q+1}-\rho_{q}\|_{L_t^rL^p}\leq \delta_{q+2}^{1/2},\ \ \|\rho_{q+1}-\rho_{q}\|_{C_tL^1}\leq \delta_{q+2}^{1/2},\label{bd:2rhoq+1-rhoql1}\ \
    \inf_{t\in [0,1]}(\rho_{q+1} - \rho_q) \geq -\delta_{q+2}^{1/2}.
\end{align}
 Here $C_0$ is determined by the choice of the starting iterations, and $C_v,C_\rho$ are two constants determined by the generalized Holder inequality for $v_q$, $\rho_q$ respectively, and other implicit geometrical constants in the proof.
\ep
 Here we remark that all the parameters are independent of the choices of $\kappa$ and $\nu$,  the extra power in  the bound on $M_q$ in \eqref{bd:2rql1} is used to absorb the universal constant to avoid exponential explosion during the iterative process, see \eqref{bd:2wq+1l2} for the details.

The proof of Proposition \ref{prop:case2} is presented in Section \ref{proof:prop2} below. With  Proposition \ref{prop:case2} in hand, the proof of Theorem \ref{thm:non_pde_23} follows by a similar argument as for Theorem \ref{thm:con:1/3}.
\begin{proof}[Proof of Theorem \ref{thm:non_pde_23}]
 For any triple $(p,r,s)\in \mathcal{A}$, without loss of generality we  assume that $p>1$ and  $T=1$.
As before, we intend to start the iteration from 
$$\rho_0(t,x)=1+\frac{\sin \pi x_1}4\chi_0(t),\  v_0=0,\ \mathring{R}_0=0,\ M_0(t,x)=(\partial_t\chi_0(t)+\kappa\chi_0(t)\pi^2)\frac{\cos \pi x_1}{4\pi}(1,0,...,0).$$
where $x=(x_1,...,x_d)$, $\chi_0$ is a smooth function with $\chi_0(t)=0$ on $[0,\frac13]$, $\chi_0(t)=1$ on $[\frac23,1]$.
By choosing $C_0$ large enough, we have for $\kappa\in[0,1]$
\begin{align*}
   \| \rho_0\|_{L^2_tL^2}+\|\rho_0\|_{C_{t,x}^1}\lesssim 1\leq C_0,\ \ \|M_0\|_{L^1_tL^1}\lesssim 1\leq\frac1{48^4} C_0^2.
\end{align*}
 Then \eqref{bd:2vql2}-\eqref{bd:2mq=rhoq=0} are satisfied.

Next, we use Proposition \ref{prop:case2} to build inductively $(v_q,\rho_q, \mathring{R}_q, M_q)$ for every $q \geq 1$.   By \eqref{bd:2vq+1-vql2}-\eqref{bd:2rhoq+1-rhoql1}, the sequence $\{v_q\}_{q\in \N}$ is
Cauchy in $$ C([0,1];L^1)\cap L^r([0,1];L^p)\cap L^2([0,1]\times\mT^d)\cap L^1([0,1];W^{1,s})$$ and the sequence $\{\rho_q\}_{q\in \N}$ is
Cauchy in $$  L^r([0,1];L^p)\cap  L^2([0,1]\times\mT^d)\cap C([0,1];L^1).$$
We then denote by $(v,\rho)$ the limit, where $v$ is also divergence-free.   For the case $r=\infty$, we have $v,\rho \in C([0,1];L^p)$. Clearly by \eqref{bd:2rql1}, $(\rho,v)$  solves \eqref{eq:fpe} and  \eqref{eq:ns}. Then from the fact $\int \rho_q\dif x=1$ we deduce that  $\int \rho\dif x=1$.
 \eqref{bd:2mq=rhoq=0} ensures that $\rho(t) \equiv 1$ for every $t$ sufficiently close to 0.

 Moreover, $\rho$ is non-negative on $\mathbb{T}^d$ by \eqref{ieq:2ab2} and \eqref{bd:2rhoq+1-rhoql1}:
\begin{align*}
\inf_{t\in [0,1]}\rho \geq \inf_{t\in [0,1]} \rho_0 +
\sum_{q=0}^\infty \inf_{t\in [0,1]}(\rho_{q+1}- \rho_q) \geq\frac34-\sum_{q=0}^\infty \delta_{q+1}^{1/2}\geq \frac12,
\end{align*}
and $\rho$  does not coincide with the solution which is constantly equal to 1, since  by \eqref{ieq:2ab2} and \eqref{bd:2rhoq+1-rhoql1}
\begin{align*}
  \|\rho - 1\|_{C_tL^1}\geq \|1 - \rho_0\|_{ C_tL^1} -
\sum_{q=0}^\infty \|\rho_{q+1}-\rho_q\|_{C_t L^1}\geq \frac1{16}- 
\sum_{q=0}^\infty\delta_{q+1}^{1/2}> 0.
\end{align*}
\end{proof}

Then the non-uniqueness of stochastic Lagrangian trajectories follows from Theorem \ref{thm:non_pde_23} and the  superposition principle.

 \begin{proof}[Proof of Theorem \ref{thm:nonuniode+sde}]
  We only prove the case where $\kappa>0$, while the case $\kappa=0$ follows similarly and more easily.
 For any triple $(p,r,s)\in \mathcal{A}$, without loss of generality we  assume that $p>1$ and  $T=1$. By Theorem \ref{thm:non_pde_23}, there exists  $v\in L^r([0,1];L^p)\cap L^2([0,1]\times \mT^d)\cap L^1([0,1]; W^{1,s})\cap C([0,1];L^1)$  and a non-constant  positive density $\rho  \in  L^r([0,1];L^p)\cap L^2([0,1]\times\mT^d)\cap C([0,1];L^1)$ 
satisfying \eqref{eq:fpe} and \eqref{eq:ns}. If $r=\infty$, we have additionally $v\in C([0,1];L^p)$.

 Moreover,  by \eqref{bd:2vqc1}, \eqref{bd:2vq+1-vql2} and interpolation we conclude that $v\in L^{2(1+\epsilon)}([0,1]\times\mT^d)$ for some $\epsilon>0$ small enough. Then,   it follows that
\begin{align*}
    \int_0^1\int_{\mT^d}|v(s,x)|^{1+\epsilon}\rho(s,x)\dif x\dif s\leq \|v\|_{L^{2+2\epsilon}_{t,x}}^{1+\epsilon}\|\rho\|_{L^2_{t,x}}<\infty,
\end{align*}
and that $t \to \rho(t,x)\dif  x$ is weakly continuous on $[0,1]$ since $\rho\in C([0,1];L^1)$.  Using the superposition principle (see \cite[Section 7.2]{Tre14}) for  \eqref{eq:fpe},  there exists  a probability measure $\mathbf{Q}$ on $C([0,1];\mT^d) $
 equipped  with its Borel $\sigma$-algebra and its natural filtration generated by the canonical process $\Pi_t, t \in [0,1]$, defined by $$\Pi_t(\omega) := \omega(t),\ \  \omega\in C([0,1];\mT^d),$$
which is a martingale solution associated to diffusion operator $$L:=\kappa\Delta+v\cdot \nabla.$$
 More precisely, for every smooth function $f$ on $\mT^d$, the process
\begin{align*}
     f( \Pi_t) -f(\Pi_0)-\int^t_0 Lf( \Pi_s)\dif s
\end{align*}
is a $\bQ$-martingale with respect to the natural filtration with
 the initial law $\mathbf{Q}\circ \Pi^{-1}_0=\cL^d$. Here  $\cL^d$  denotes the Lebesgue measure on the torus.  Moreover,  for every $t \in [0,1]$, it holds that $\rho(t) \cL^d=\mathbf{Q}\circ \Pi^{-1}_t$.
Since  $\overline\rho\equiv1$ is also a solution satisfying all the above conditions, we have another martingale solution  $ \overline{\mathbf{Q}}$ 
satisfying $ \cL^d=\overline{\mathbf{Q}}\circ \Pi^{-1}_t$.

Now we define $\{\mathbf{Q}^x\}_{x\in\mT^d},  \{ \overline{\mathbf{Q}}^x\}_{x\in\mT^d}$ as the  regular conditional probabilities with respect to
$\Pi_0$. Then for a.e. $x\in\mT^d$, $\mathbf{Q}^x,  \overline{\mathbf{Q}}^x$ are both  martingale solutions  to \eqref{eq:sde} with  initial condition $x$. 
We define 
\begin{align*}
    A(v):=\{x\in\mT^d:\mathbf{Q}^x,  \overline{\mathbf{Q}}^x\ {\rm are\ two\ distinct\  martingale\ solutions\ associated\ to}\ L\}.
\end{align*}
Then we prove that the Lebesgue measure of $A(v)$ is positive.
We assume by contradiction that  the Lebesgue measure is zero. In this case, we have $\mathbf{Q}^x= \overline{\mathbf{Q}}^x$ for a.e. $x\in\mT^d$. Consequently, we have for any smooth function $f$ and $t\in[0,1]$
\begin{align*}
    \int_{\mT^d}f\rho(t)\dif x=\int f(\Pi_t)\dif \mathbf{Q}= \int_{\mT^d}\int f(\Pi_t)\dif \mathbf{Q}^x\dif x= \int_{\mT^d}\int f(\Pi_t)\dif \mathbf{\overline Q}^x\dif x=\int_{\mT^d}f\dif x,
\end{align*}
which leads to a  contradiction as $\rho$ is non-constant.

 Since $|v|,\rho  \in L^2([0,1]\times\mathbb{T}^d)$, we have
\begin{align*}
  \int_{\mT^d}\bE^{x}\int_0^1|v(s,\Pi_s)|\dif s\dif x= \int\int_0^1|v(s,\Pi_s)|\dif s\dif \mathbf{Q} = \int_0^1 \int_{\mT^d}|v(s,x)|\rho(s,x)\dif x\dif s<\infty,
\end{align*}
which implies that  $\bE^{x}[\int_0^1|v(s,\Pi_s)|\dif s]<\infty$ for a.e. $x\in\mT^d$.  Similarly, we have $\mathbf{ \overline E}^{x}[\int_0^1|v(s,\Pi_s)|\dif s] <\infty$ for a.e. $x\in\mT^d$. 
 

\end{proof}

\section{Proof of Proposition \ref{prop:case2}}\label{proof:prop2}
The proof is also based on   convex integration schemes. Compared with Section \ref{sec:proof13}, here we use two distinct scales for two equations during the iteration. We begin by fixing the parameters and then proceed with a mollification step in Section~\ref{sec:choicepara2}. 
Section~\ref{sec:per} presents the perturbation construction of $(w_{q+1}{+\overline{w}_{q+1}},\theta_{q+1})$ and the new iteration $(v_{q+1},\rho_{q+1})$. Here, the perturbation $(w_{q+1},\theta_{q+1})$ is designed to cancel the  stress term  $ M_q$ in the advection-diffusion equations, while  the perturbation $ w_{q+1}$ is  to cancel the stress term $\mathring{R}_q$ in the fluid equations.  Here, unlike  the previous case, the term $w_{q+1}\otimes w_{q+1}$ is automatically small in $L^1$ space due to the choice of smaller scales for the advection-diffusion equation.  We emphasize that the method employed in Section \ref{sec:proof13} can not be  applied directly in this context, since the additional intermittency  in the building blocks introduces extra oscillation errors (see Section \ref{sec:idea2} for more explanation). Then we establish the inductive estimates. Finally, in Section~\ref{sec:error}, we define the new stress components $(\mathring{R}_{q+1},M_{q+1})$ and establish the  inductive estimates respectively.

\subsection{Choice of parameters and mollification}\label{sec:choicepara2}
In the sequel, additional parameters will be indispensable and their value have to be carefully chosen to respect all the compatibility conditions appearing in the estimates below. First, for a sufficiently small $\alpha\in(0,1)$ to be chosen, we take $l:=\lambda_{q+1}^{-\frac{3\alpha}{2}}\lambda_{q}^{-2d-\frac32}$ and have 
\begin{align}
    l^{-1}\leq\lambda_{q+1}^{2\alpha},\  \ l\lambda_{q}^{4d+3}\ll\lambda_{q+1}^{-\alpha}\leq \delta_{q+3}^2,\ \ \lambda_{q}^{4d+3}\leq\lambda_{q+1}^\alpha\label{para22}
\end{align}
 provided $\alpha{b}\geq  4d+3,\alpha>4\beta b^2$. 

For fixed $d\geq2$  and $(p,r,s)\in \mathcal{A},$ without loss of generality, we assume $p>1$. Then we  introduce a large constant  $N:=N(p,r,s,d)$ to be chosen in Lemma \ref{lem:para2} below. In the sequel, we also need
$$\alpha b\geq 4d+3,\ \ \alpha>4\beta b^2,\ \
(12d+43)\alpha<\frac1{2N}.$$
 The above can be obtained by choosing $\alpha>0$ small such that $(12d+43)\alpha<\frac1{2N}$, and choosing $b\in \mathbb{N}$ large enough such that $b>\frac{4d+3}{\alpha}$, and finally choosing $\beta>0$ small such that $\alpha>4\beta b^2$.

Then we increase $a$ such that \eqref{ieq:2ab2} holds. In the sequel, we also increase $a$  to absorb  various implicit and universal constants in the subsequent estimates.

Now we replace $(v_{q},\rho_q)$ by a mollified field $(v_l,\rho_l)$, and 
define
\begin{align}
v_l=(v_{q }*_x\phi_l)*_t\varphi_l,\ \ 
 \rho_l=(\rho_{q}*_x\phi_l)*_t\varphi_l,\notag
\end{align}
where $\phi_l:=\frac{1}{l^d}\phi(\frac{\cdot}{l})$ is a family of standard  radial  mollifiers on $\mathbb{R}^d$, and $\varphi_l:=\frac{1}{l}\varphi(\frac{\cdot}{l})$ is a family of standard  radial  mollifiers with support in $(0,1)$.  For the mollification around $t = 0$,  since $v_q,\rho_q,\mathring{R}_q$ and $M_q$ are constants around $t = 0$, see \eqref{bd:2mq=rhoq=0}, we  can directly extend these definitions  to $t\leq 0$ by their values at $t=0$.

By straightforward calculations and \eqref{eq:2qth} we obtain 
\begin{align}
\partial_t \rho_l-\kappa \Delta\rho_l+\div(v_l \rho_l)&=-\div M_l,\notag\\
\partial_tv_l-\nu\Delta v_l+\div (v_l\otimes v_l)+\nabla\pi_l
&=\div\mathring{R}_l,\ \ \div v_l=0\label{eq:2v_l}
\end{align}
 for some  suitable $\pi_l$, where
 \begin{align}
M_l&:=(M_{q }*_x\phi_l)*_t\varphi_l-v_l \rho_l+(v_{q } \rho_{q})*_x\phi_l*_t\varphi_l,\notag\\
\mathring{R}_l&:=(\mathring{R}_{q }*_x\phi_l)*_t\varphi_l+v_l\mathring{\otimes}v_l-(v_{q }\mathring{\otimes}v_{q })*_x\phi_l*_t\varphi_l.
\end{align}

Moreover, since $l\leq \frac12\delta_{q+1}^{1/2}=\frac{T_q-T_{q+1}}2$,  by \eqref{bd:2mq=rhoq=0} we know that $ \rho_l-1={v}_{ {l}}=\mathring{ {R}}_{ {l}}= M_l=0$ on $[0,\frac{T_q+T_{q+1}}2]$.

 Then by the mollification estimates in Lemma \ref{p:CET}, the space-time embedding $W^{d+1+\epsilon,1}\subset L^\infty$ and the bounds \eqref{bd:2vqc1}, \eqref{bd:2rql1}, \eqref{para22} we obtain
 \begin{align}
 \|\mathring{ {R}}_{ {l}}\|_{L_t^1L^1}\leq \|\mathring{R}_q\|_{L^1_tL^1}+Cl^{2}\|v_q\|_{C_{t,x}^1}^2\leq  C_0^2\delta_{q+1}+CC_0^2l^{2}\lambda_q^{8d+6}\leq 2C_0^2\delta_{q+1},\label{bd:2rll1l1}\\
  \|M_{ {l}}\|_{L_t^1L^1}\leq \|M_q\|_{L^1_tL^1}+Cl^{2}\|v_q\|_{C_{t,x}^1}\|\rho_q\|_{C_{t,x}^1}\leq  C_0^2\delta_{q+2}^2+CC_0^2l^{2}\lambda_q^{7d+5}\leq 2C_0^2\delta_{q+2}^2,\label{bd:2mll1l1}
   \end{align} 
   and for $N\geq 0$,
    \begin{align}
 \|\mathring{ {R}}_{ {l}}\|_{C_{t,x}^N}&\lesssim  {l}^{-d-1-\epsilon-N}\|\mathring{R}_q\|_{L^1_tL^1}+l^{2-N}\|v_q\|_{C_{t,x}^1}^2\lesssim {l}^{-d-1-\epsilon-N}+l^{2-N}\lambda_q^{8d+6}\lesssim  {l}^{-d-2-N},\label{bd:2rlcn}\\
    \|M_l\|_{C_{t,x}^N}&\lesssim l^{-d-1-\epsilon-N} \|M_{q }\|_{L_t^1L^1}+l^{2-N}\|v_q\|_{C_{t,x}^1}\|\rho_q\|_{C_{t,x}^1}\lesssim l^{-d-1-\epsilon-N}+l^{2-N}\lambda_q^{7d+5}\lesssim l^{-d-2-N},\label{bd:2mlcn}
 \end{align} 
 where we use $l^{-1}$  to absorb the implicit  constants by choosing $a$ large enough.
 
\subsection{The construction of  perturbations and inductive estimates}\label{sec:per}

 As outlined in Section \ref{sec:idea}, we proceed with the construction of the perturbation $(w_{q+1}+\overline w_{q+1},\theta_{q+1})$, ensuring that the supports of $ (w_{q+1},\theta_{q+1})$  and $\overline w_{q+1}$ are disjoint. The perturbation $(w_{q+1},\theta_{q+1})$ is designed to cancel the stress terms $ M_q$.
 At the same time, the perturbations $ w_{q+1}$ are exclusively utilized to cancel the  stress terms $\mathring{R}_q$.  In contrast to the previous case, this extra product $ w_{q+1}\otimes w_{q+1}$ is already sufficiently small  since we use  two distinct scales for two equations during the iteration (see the inductive condition \eqref{bd:2rql1}). 
 
\subsubsection{Construction of $w_{q+1}$ and $\theta_{q+1}$}\label{sec:2defq+1}
Let us now proceed with the construction of the perturbation $(w_{q+1}, \theta_{q+1})$  by employing  the building blocks and temporal jets introduced in Section \ref{gij}. 

For given $p,r,s$ and $d$, we choose the parameters $ {\lambda}, {r}_{\perp}, {r}_{\parallel}, {\eta}$ as follows:

\bl\label{lem:para2}
Let $d\geq2$,   $(p,r,s)\in \mathcal{A}$ and $p>1$. 
There exists a choice of parameters $ {\lambda}, {r}_{\perp}, {r}_{\parallel}, {\eta}$ and $N=N(p,r,s,d)\in\mathbb{N}$, such that  
\begin{align}
    {\lambda}^{-1}\ll  {r}_{\perp}\ll  {r}_{\parallel}\ll 1,\ \ {\eta} ^{-1}\leq  {\lambda}^d,\label{bd:lempara2}
\end{align}
 and
 \begin{align}
  \max\{   \frac{ {r}_\perp}{ {r}_\parallel},\ \  {r}_\perp^{-1} {\lambda}^{-1},\ \  {\lambda}  {r}_\perp^{\frac{d-1}s-\frac{d-1}{2}}  {r}_\parallel^{\frac1s-\frac12} {\eta}^{\frac12},\ \  {r}_\perp^{\frac{d-1}{2}}  {r}_\parallel^{\frac12} {\eta}^{-\frac12},\ \  {r}_\perp^{\frac{d-1}{p}-\frac{d-1}{2}}  {r}_\parallel^{\frac{1}{p}-\frac12} {\eta}^{\frac1r-\frac12}\}\leq {\lambda}^{-\frac1N}.\label{paralam}
 \end{align}

\el
\begin{proof}
    For any $(p,r,s)\in \mathcal{A}$ and $p>1$, it is easy to see that there exists $k\in\mQ,k>1$ such that 
    $$\frac1r-\frac12+k(\frac{1}{p}-\frac12)>0, \ \ \frac1d<\frac1{2k}-\frac12+\frac1s.$$
    Then there exists  $M\in\N$ such that $\frac1r-\frac12+k(\frac{1}{p}-\frac12)>\frac1M$,  $\frac1d<(1-\frac1M)(\frac1{2k}-\frac12+\frac1s)$. 
     We choose $N\in\mathbb{N}$ large enough such that $N\geq \max\{M,\frac{2kM^2}{(M-1)d},\frac{4Mk}{d(M-1)(k-1)}\}$ and  $\frac1d-(1-\frac1M)(\frac1{2k}-\frac12+\frac1s)\leq-\frac{2}{dN}$.  Then we define $ {\eta}:= {\lambda}^{-\frac{d}{k}(1-\frac1M)}, {r}_\perp:= {\lambda}^{-1+\frac1M}\gg  {\lambda}^{-1}, {r}_\parallel:= {\lambda}^{-1+\frac1{M}+\frac1N}$ and have
        $$\frac{ {r}_\perp}{ {r}_\parallel}=  {\lambda}^{-\frac1N},\  {r}_\perp^{-1} {\lambda}^{-1}= {\lambda}^{-\frac1M}\leq {\lambda}^{-\frac1N},\  {r}_\perp^{\frac{d-1}{2}}  {r}_\parallel^{\frac12} {\eta}^{-\frac12}\leq  {\lambda}^{-(1-\frac1M)\frac d2(1-\frac1k)+\frac1{N}}\leq {\lambda}^{-\frac1N},$$
    $$ {\lambda}  {r}_\perp^{\frac{d-1}s-\frac{d-1}{2}}  {r}_\parallel^{\frac1s-\frac12} {\eta}^{\frac12}\leq {\lambda}^{1-(1-\frac1M)d(\frac1{2k}-\frac12+\frac1s)+\frac1N}\leq {\lambda}^{-\frac1N}, $$ 
    $$ {r}_\perp^{\frac{d-1}{p}-\frac{d-1}{2}}  {r}_\parallel^{\frac{1}{p}-\frac12} {\eta}^{\frac1r-\frac12}\leq  {\lambda}^{-(1-\frac1M)\frac dk(\frac1r-\frac12+k(\frac{1}{p}-\frac12))+\frac{1}{N}}\leq  {\lambda}^{-(1-\frac1M)\frac dk\frac1M+\frac{1}{N}}\leq {\lambda}^{-\frac1N}.$$
\end{proof}
Then we choose ${\lambda}:=\lambda_{q+1 }$, and $ {r}_{\perp}, {r}_{\parallel}, {\eta}$ in terms of $\lambda_{q+1}$ according to Lemma \ref{lem:para2}. Moreover, we define
\begin{align}
    {\lambda}:=\lambda_{q+1 },\ \  \overline{\mu}:=  {r}_\parallel^{-\frac12} {r}_\perp^{-\frac{d-1}{2}} {\lambda}_{q+1}^{\frac1{2N}}\leq  \lambda_{q+1}^{\frac d2},\ \  {\mu}:=  {r}_\parallel^{-\frac12} {r}_\perp^{-\frac{d-1}{2}} \leq  \lambda_{q+1}^{\frac d2}, \ \  {\sigma}:= {\lambda}_{q+1}^{\frac1{2N}}.\label{def:2musigma}
\end{align}  
 It is required that $b$ is a multiple of $M$ to ensure that $ {\lambda}  {r}_{\perp}=a^{(b^{q+1})/M}\in\mathbb{N}$, where $M$ is given in Lemma \ref{lem:para2}.

Next, using the building blocks introduced in Section \ref{gij}, we define the perturbation $(w_{q+1},\theta_{q+1})$ similarly to that in \cite[Section 5.3]{BCDL21}. 
Let $\chi\in C_c^\infty(-\frac34,\frac34)$  be a non-negative function such that
$\sum_{n\in\mathbb{Z}} \chi(t - n) = 1$ for every $t \in \R$. Let $\tilde{\chi} \in C_c^\infty(-\frac45,\frac45)$ be a non-negative function satisfying $\tilde{\chi} = 1$ in
$[-\frac34,\frac34]$ and $\sum_{n\in\mathbb{Z}} \tilde{\chi}(t - n) \leq2$.

We fix a parameter $\zeta = 20/\delta^2_{q+3}$ and consider two disjoint sets $\Lambda^1, \Lambda^2$ introduced in Lemma \ref{l:linear_algebra2} with $d\geq2$. Next,  we use the notation $ \Lambda^i =  \Lambda^1$ for $i$ odd, and $ \Lambda^i =  \Lambda^2$ for $i$ even.
In the following we abuse the notation and define for $n\in\N$
\begin{align*}
    W_{(\xi,g)}(x,t):=W_{(\xi)}(x,(\frac{n}{\zeta})^{1/2}H_{(\xi)}(t)).
\end{align*} 
 Similarly, we could define  $\Theta_{(\xi,g)},V_{(\xi,g)}
 $ and all other terms appearing in Section \ref{sec:bbte}. Now by the identity \eqref{eq:ptthe+}, the definition of $H_{(\xi)}(t)$ in \eqref{eq:parth} and the choice of $\mu$  in \eqref{def:2musigma} we have
 \begin{align}
     \partial_t\Theta_{(\xi,g)}+(\frac{n}{\zeta})^{1/2}g_{(\xi)}\div (W_{(\xi,g)}\Theta_{(\xi,g)})=0.\label{eq:2ptthe+n}
 \end{align}

As the next step, we define the principle part of the perturbation $w_{q+1}$ by
\begin{align}
w_{q+1}^{(p)}:=\sum_{ n\geq3}\tilde{\chi}(\zeta|M_l|-n)\(\frac{n}{\zeta}\)^{1/2}\sum_{\xi\in\Lambda^n}W_{(\xi,g)}g_{(\xi)}.\notag
\end{align}
We remark that here and in the following the first sum runs for $n$ in the range 
\begin{align}
    3\leq n\leq 1+\zeta|M_l|\leq 1+l^{-2/3-d-1-\epsilon}\|M_l\|_{L^1_tL^1}\leq 1+Cl^{-d-2}
\label{bd:2|ml|}\end{align}
 because of  the bounds  \eqref{para22}, \eqref{bd:2rll1l1} and the space-time Sobolev embedding  $W^{d+1+\epsilon,1}\subset L^\infty$.  Here we choose $a$ large enough to absorb the embedding constant.

 Moreover, we define the incompressibility corrector by
\begin{align}
w_{q+1}^{(c)}:=\sum_{n\geq3}\sum_{\xi\in\Lambda^n}\(-\tilde{\chi}(\zeta|M_l|-n)&(\frac{n}{\zeta})^{1/2}\frac{1}{(n_*\lambda_{q+1})^2}\nabla\Phi_{(\xi,g)}\xi\cdot\nabla\psi_{(\xi,g)}\notag\\ 
&+\nabla( \tilde{\chi}(\zeta|M_l|-n))(\frac{n}{\zeta})^{1/2}:V_{(\xi,g)}\)g_{(\xi)}.\notag
\end{align}
Here we denote $(\nabla( \tilde{\chi}(\zeta|M_l|-n)):{V}_{({\xi,g})})^i:=\sum_{j=1}^d\partial_j( \tilde{\chi}(\zeta|M_l|-n)){V}_{({\xi,g})}^{ij},\ i=1,2,...,d$.

By the identity \eqref{divOmega} we have
\begin{align}
w_{q+1}^{(p)}+w_{q+1}^{(c)}&=\sum_{n\geq3}\sum_{\xi\in\Lambda^n}\div\(\tilde{\chi}(\zeta|M_l|-n)(\frac{n}{\zeta})^{1/2} V_{(\xi,g)}\)g_{(\xi)}.\label{eq:2wq+1p+wq+1c}
\end{align}
Since $V_{(\xi,g)}$ is skew-symmetric, we obtain
$$\div ( w_{q+1}^{(p)}+w_{q+1}^{(c)})=0.$$

Then we  define the principle part of the perturbation $\theta_{q+1}$ by
\begin{align}
\theta_{q+1}^{(p)}:&=\sum_{n\geq3}\chi(\zeta|M_l|-n)(\frac{n}{\zeta})^{1/2}\sum_{\xi\in\Lambda^n}\Gamma_{\xi}\(\frac{M_l}{|M_l|}\)\Theta_{(\xi,g)}g_{(\xi)},
\end{align}
where $\Gamma_{\xi}$ is introduced in Lemma \ref{lem:geo:tran}.
  $\theta_{q+1}^{(p)}$ is non-negative since all the components are non-negative.
 
The mean corrector is defined  by 
\begin{align}
\theta_{q+1}^{(c)}:&=-\mathbb{P}_{0}\theta_{q+1}^{(p)},\notag
\end{align}
where we recall that $\mathbb{P}_{0}f=\fint_{\mT^d} f\dif x.$  

Using the fact that $g_{(\xi)}$ have disjoint support, the identity \eqref{eq:intwthe} above, the geometry Lemma \ref{l:linear_algebra2} and the fact that $\chi\tilde{\chi}=\chi$ we obtain 
\begin{align}
    w_{q+1}^{(p)}\theta_{q+1}^{(p)}&=\sum_{n\geq3}\chi(\zeta|M_l|-n)\frac{n}{\zeta}\sum_{\xi\in\Lambda^n}\Gamma_{\xi}\(\frac{M_l}{|M_l|}\)W_{(\xi,g)}\Theta_{(\xi,g)}g_{(\xi)}^2\notag\\ 
    &=\sum_{n\geq3}\sum_{\xi\in\Lambda^n}\chi(\zeta|M_l|-n)\frac{n}{\zeta}\Gamma_{\xi}\(\frac{M_l}{|M_l|}\)\mP_{\neq0}(W_{(\xi,g)}\Theta_{(\xi,g)})g_{(\xi)}^2\notag\\ 
    &\ \ \ \   +\sum_{n\geq3}\sum_{\xi\in\Lambda^n}\chi(\zeta|M_l|-n)\frac{n}{\zeta}\Gamma_{\xi}\(\frac{M_l}{|M_l|}\)\xi (g_{(\xi)}^2-1)+\sum_{n\geq3}\chi(\zeta|M_l|-n)\frac{n}{\zeta}\frac{M_l}{|M_l|}.\label{2wq+1ptheq+1p}
\end{align}
We observe that in  \eqref{2wq+1ptheq+1p}, there is an undesirable term of the form $(\cdot) (g_{(\xi)}^2-1)$ arising from the temporal intermittency. To deal with this term, we define the temporal corrector as follows:
\begin{align*}
 \theta_{q+1}^{(o)}:&=-\sigma^{-1}\sum_{n\geq3}\sum_{\xi\in\Lambda^n}h_{(\xi)}\div  \(\chi(\zeta|M_l|-n)\frac{n}{\zeta}\Gamma_{\xi}\(\frac{M_l}{|M_l|}\)\xi\).\notag
\end{align*}

Recalling the definition of $h_{(\xi)}(t)$ in \eqref{eq:parth}  we have
\begin{align}
  \partial_t\theta^{(o)}_{q+1}&+ \sum_{n\geq3}\sum_{\xi\in\Lambda^n}(g_{(\xi)}^2-1)\div\(\chi(\zeta|M_l|-n)\frac{n}{\zeta}\Gamma_{\xi}\(\frac{M_l}{|M_l|}\)\xi\)\notag\\ 
  &=-\sigma^{-1}\sum_{n\geq3}\sum_{\xi\in\Lambda^n}h_{(\xi)}\partial_t\div \(\chi(\zeta|M_l|-n)\frac{n}{\zeta}\Gamma_{\xi}\(\frac{M_l}{|M_l|}\)\xi\).\label{eq:2parttheo}
\end{align}

Finally, the  perturbations  are defined as
\begin{align}
w_{q+1}:=w_{q+1}^{(p)}+w_{q+1}^{(c)},\ \ \theta_{q+1}:=\theta_{q+1}^{(p)}+\theta_{q+1}^{(c)}+    \theta_{q+1}^{(o)},\notag
\end{align}
where  $w_{q+1}$ is  mean-zero and divergence-free, and $\theta_{q+1}$ is mean-zero.
Since $M_l(t)=0$ for $t\in [0,\frac{T_q+T_{q+1}}2]$, by definition we know that  $w_{q+1}(t)=\theta_{q+1}(t)=0$ for $t\in [0,\frac{T_q+T_{q+1}}2]$. 

 The new scalar $\rho_{q+1}$ is defined by
\begin{align*}
    \rho_{q+1}:=\theta_{q+1}+\rho_l,
\end{align*}
which satisfies $\int_{\mT^d}\rho_{q+1}\dif x=1$. 
Consequently,  since $\rho_l=1$ for $t\in [0,\frac{T_q+T_{q+1}}2]$, we have $\rho_{q+1}(t)=1$ for $t\in [0,\frac{T_q+T_{q+1}}2]$.

\subsubsection{Estimates of $ w_{q+1}$}
First we establish the estimate of the amplitude functions  defined in Section \ref{sec:2defq+1}.
\bp\label{lem:2chi}
 For $ N\in\mathbb{N}_0$ we have
\begin{align*}
    \sum_{n\geq3}\|\chi(\zeta|M_l|-n)\|_{C_{t,x}^N}+ \sum_{n\geq3}\|\tilde{\chi}(\zeta|M_l|-n)\|_{C_{t,x}^N}&\lesssim l^{- (d+4)N-(d+2)},\\
     \sum_{n\geq3}\sum_{\xi\in\Lambda^n}\norm{\chi(\zeta|M_l|-n)\Gamma_\xi\(\frac{M_l}{|M_l|}\)}_{C_{t,x}^N}& \lesssim l^{-(2d+8)N-(d+2)},\\
   \(\frac{n}{\zeta}\)^N1_{ \{\chi(\zeta|M_l|-n)>0\}}+\(\frac{n}{\zeta}\)^N1_{ \{ \tilde{\chi}(\zeta|M_l|-n)>0\}}&\lesssim l^{-N(d+2)}.
\end{align*}
\ep
 We give the  proof of this lemma in Appendix \ref{sec:est:ampl}.
 
Recalling that $w_{q+1}=w_{q+1}^{(p)}+w_{q+1}^{(c)}$ is defined in Section \ref{sec:2defq+1}, we first estimate ${w}^{(p)}_{q+1}$ in the $L^2_tL^2$-norm. Similar as in \eqref{2wq+1ptheq+1p} we have
\begin{align}
|w_{q+1}^{(p)}|^2&\lesssim \sum_{n\geq3}\tilde{\chi}(\zeta|M_l|-n)\frac{n}{\zeta}\sum_{\xi\in\Lambda^n}|W_{(\xi,g)}|^2g_{(\xi)}^2.\notag
\end{align}
Then by the improved H\"older inequality in Lemma \ref{ihiot}, the estimate of $\tilde{\chi}$ in Proposition \ref{lem:2chi}, the bounds \eqref{int4} and \eqref{paralam} we have
\begin{align}
\|w_{q+1}^{(p)}(t)\|^2_{L^2}
&\lesssim\sum_{n\geq3}\left\|\tilde{\chi}(\zeta|M_l(t)|-n)\frac{n}{\zeta}\right\|_{L^1}\sum_{\xi\in\Lambda^n}\|W_{(\xi,g)}\|_{C_tL^2}^2g_{(\xi)}^2(t)\notag\\ 
&\quad+(r_\perp\lambda_{q+1})^{-1}\left\|\tilde{\chi}(\zeta|M_l|-n)\frac{n}{\zeta}\right\|_{C_{t,x}^1}\sum_{\xi\in\Lambda^n}\|W_{(\xi,g)}\|_{C_tL^2}^2g_{(\xi)}^2(t)\notag\\
&\lesssim\(\big{\|}\sum_{n\geq3}\tilde{\chi}(\zeta|M_l(t)|-n)(|M_l(t)|+\zeta^{-1}
)\big{\|}_{L^1}+l^{-3d-8}\lambda_{q+1}^{-\frac{1}{N}}\)\sum_{\xi\in\Lambda^1\cup\Lambda^2}g_{(\xi)}^2(t)\notag\\ 
&\lesssim (\|M_l(t)\|_{L^1}+\delta^2_{q+2})\sum_{\xi\in\Lambda^1\cup\Lambda^2}g_{(\xi)}^2(t),\notag
\end{align}
where we used the facts that $\tilde{\chi}$ is non-negative, $\sum_{n\in\mathbb{Z}} \tilde{\chi}(t - n) \leq2$, and used conditions on the parameters to have $(6d+16)\alpha-\frac{1}{N}<-\alpha<-4\beta b$. Then we  apply the improved H\"older inequality  of Lemma \ref{ihiot} again in time. Together with the bounds on $g_{(\xi)}$ in \eqref{bd:gwnp}, $M_l$ in \eqref{bd:2mll1l1}, \eqref{bd:2mlcn} and the choice of parameters in \eqref{para22} we obtain
\begin{align}
    \|w_{q+1}^{(p)}\|^2_{L_t^2L^2}&\lesssim (\|M_l\|_{L_t^1L^1}+\delta^2_{q+2}+\sigma^{-1}\|M_l\|_{C_{t,x}^1})\|\sum_{\xi\in\Lambda^1\cup\Lambda^2}g_{(\xi)}^2\|_{L_t^1}\notag\\
    &\lesssim C_0^2(\delta^2_{q+2}+\lambda_{q+1}^{(2d+6)\alpha-\frac1{2N}})\lesssim C_0^2(\delta^2_{q+2}+ \lambda_{q+1}^{-\alpha})\lesssim \frac{C_0^2}{16}\delta^2_{q+2},\label{bd:2wq+1pl2}
\end{align}
where we used  conditions on the parameters to have $(2d+8)\alpha<\frac{1}{2N}$.

For the general $L^u_tL^m$-norm with $u,m\in[1,\infty]$, by the estimates for the building blocks in \eqref{int2}-\eqref{int4}  and the estimates for $\tilde{\chi}$ in Proposition \ref{lem:2chi} we obtain
\begin{align}
\|w_{q+1}^{(p)}\|_{L^u_tL^m}&\lesssim\sum_{n\geq3}\sum_{\xi\in\Lambda^n}\left\|\tilde{\chi}(\zeta|M_l|-n)\(\frac{n}{\zeta}\)^{1/2}\right\|_{C_{t,x}^0}\|W_{(\xi,g)}\|_{C_tL^m}\|g_{(\xi)}\|_{L_t^u}\notag\\ 
&\lesssim l^{-2d-4}r_\perp^{\frac{d-1}{m}-\frac{d-1}{2}} r_\parallel^{\frac{1}{m}-\frac12}\eta^{\frac1u-\frac12}
,\label{bd:2wq+1plp}
\end{align}
\begin{align}
\|w_{q+1}^{(c)}&\|_{L^u_tL^m}\lesssim\sum_{n\geq3}\sum_{\xi\in\Lambda^n}\left\|\tilde{\chi}(\zeta|M_l|-n)\(\frac{n}{\zeta}\)^{1/2}\right\|_{C^1_{t,x}}\notag\\ 
&\qquad\qquad\qquad\times\(\frac{1}{\lambda_{q+1}^2}\|\nabla\Phi_{(\xi,g)} \xi\cdot\nabla\psi_{(\xi,g)} \|_{L^m}+\|V_{(\xi,g)}\|_{L^m}\)\|g_{(\xi)}\|_{L_t^u}\notag\\
&\lesssim l^{-3d-8}r_\perp^{\frac{d-1}{m}-\frac{d-1}{2}} r_\parallel^{\frac{1}{m}-\frac12}\(\frac{r_\perp}{ r_\parallel}+\lambda_{q+1}^{-1}\)\eta^{\frac1u-\frac12}\lesssim l^{-3d-8}r_\perp^{\frac{d-1}{m}-\frac{d-1}{2}} r_\parallel^{\frac{1}{m}-\frac12}\eta^{\frac1u-\frac12}\lambda_{q+1}^{-\frac1N}.\label{bd:2wq+1clp}
\end{align}

Combining these with the choice of parameters in \eqref{para22},  \eqref{paralam}, and the bound \eqref{bd:2wq+1pl2} we obtain 
\begin{align}
\|w_{q+1}\|_{L^2_tL^2}&\lesssim \frac{C_0}4\delta_{q+2}+l^{-3d-8}\lambda_{q+1}^{-\frac1N}
\leq \frac{3C_0}8\delta_{q+2}^{1/2},\label{bd:2wq+1l2}
\end{align} 
where we used   conditions on the parameters to have
${(6d+16)\alpha-\frac{1}{N}}<-\alpha<-\beta b $. 
We also selected $\delta_{q+2}^{1/2}$ to be small enough by choosing $a$  large enough to absorb the universal constant.

By the above bounds \eqref{bd:2wq+1plp}, \eqref{bd:2wq+1clp} and  the choice of parameters in  \eqref{para22}, \eqref{paralam} we have 
\begin{align}
     \|w_{q+1}\|_{L^r_tL^p}&\lesssim  l^{-3d-8}r_\perp^{\frac{d-1}{p}-\frac{d-1}{2}} r_\parallel^{\frac{1}{p}-\frac12}\eta^{\frac1r-\frac12}\lesssim \lambda_{q+1}^{(6d+16)\alpha-\frac{1}{N}}\lesssim\lambda_{q+1}^{-\alpha},\label{bd:wq+1lrp} \\
      \|w_{q+1}\|_{C_tL^1}& \lesssim  l^{-3d-8}r_\perp^{\frac{d-1}{2}} r_\parallel^{\frac12}\eta^{-\frac12}\lesssim \lambda_{q+1}^{(6d+16)\alpha-\frac{1}{N}}\lesssim\lambda_{q+1}^{-\alpha}, \label{bd:wq+1ctl1}
\end{align}
where we used  conditions on the parameters to have
${(6d+17)\alpha<\frac{1}{N}}$. 

Next, we estimate the $C_{t,x}^2$-norm. By the fact that 
\begin{align*}
    \partial_t(V_{(\xi,g)}(t))=\(\frac{n}{\zeta}\)^{1/2}g_{(\xi)}\(\partial_t V_{(\xi)}\)\(\(\frac n{\zeta}\)^{1/2}H_{(\xi)}(t)\),
\end{align*}
and the estimates for  the building blocks in \eqref{int4}, \eqref{bd:gwnp}, the identity  \eqref{eq:2wq+1p+wq+1c}  and the estimates for the amplitude functions in Proposition \ref{lem:2chi} we have 
\begin{align}
    \|w_{q+1}\|_{C_{t,x}^2}&\lesssim\sum_{n\geq3}\sum_{\xi\in\Lambda^n}\|\tilde{\chi}(\zeta|M_l|-n)(\frac{n}{\zeta})^{1/2}\|_{C_{t,x}^{3}}\(\|g_{(\xi)}\nabla V_{(\xi,g)}\|_{C_{t,x}^2} +\|g_{(\xi)}V_{(\xi,g)}\|_{C_{t,x}^2}\)\notag\\
&\lesssim  \sum_{n\geq3}\|\tilde{\chi}(\zeta|M_l|-n)(\frac{n}{\zeta})^{3/2}\|_{C_{t,x}^{3}} \lambda_{q+1}^2\mu^2 r_\parallel^{-\frac12}r_\perp^{-\frac{d-1}{2}}\sigma^2\eta^{-\frac52}\lesssim\lambda_{q+1}^{(12d+36)\alpha+4d+\frac{5}{2}}\leq \lambda_{q+1 }^{4d+3},\label{bd:2wq+1ctx2}
\end{align}
where we used \eqref{bd:lempara2}  \eqref{def:2musigma}, and  conditions on the parameters to have $(12d+36)\alpha<\frac12$.

We conclude this part with estimates in $W^{1,s}$-norms. By the estimates for the building blocks in \eqref{int4},  \eqref{bd:gwnp},  and the estimate for the amplitude functions in Proposition \ref{lem:2chi} we obtain
\begin{align}
\|w_{q+1}\|_{L^1_tW^{1,s}}&\lesssim\sum_{n\geq3}\sum_{\xi\in\Lambda^n}\|\tilde{\chi}(\zeta|M_l|-n)(\frac{n}{\zeta})^{1/2}\|_{C_{t,x}^{2}}\|V_{(\xi,g)}\|_{C_tW^{2,s}}\|g_{(\xi)}\|_{L_t^1}\notag\\
&\lesssim l^{-4d-12} \lambda_{q+1} r_\perp^{\frac{d-1}{s}-\frac{d-1}{2}} r_\parallel^{\frac{1}{s}-\frac12}\eta^{\frac12}\lesssim \lambda_{q+1}^{(8d+24)\alpha-\frac1N} \lesssim \lambda_{q+1}^{-\alpha},\label{bd:2wq+1w1p}
\end{align}
where  we used  \eqref{para22}, \eqref{paralam} and  conditions on the parameters to have  $(8d+25)\alpha<\frac1N.$ 

\subsubsection{Estimates of $\theta_{q+1}$}
 Recall that $\theta_{q+1}$ is defined in Section \ref{sec:2defq+1}.
We first estimate $\theta^{(p)}_{q+1}$ in $L^2_tL^2$-norm by a similar argument as in \eqref{bd:2wq+1pl2}. Noting  the fact that $\Gamma_{\xi}$ are uniformly bounded,    we have
\begin{align}
|\theta_{q+1}^{(p)}|^2
&\lesssim \sum_{n\geq3}\chi(\zeta|M_l|-n)\frac{n}{\zeta}\sum_{\xi\in\Lambda^n}\left|\Gamma_{\xi}\(\frac{M_l}{|M_l|}\)\Theta_{(\xi,g)}\right|^2g_{(\xi)}^2
\lesssim \sum_{n\geq3}\chi(\zeta|M_l|-n)\frac{n}{\zeta}\sum_{\xi\in\Lambda^n}\left|\Theta_{(\xi,g)}\right|^2g_{(\xi)}^2.\notag
\end{align}

 Then by the same argument as in \eqref{bd:2wq+1pl2}, we have for some ${C}_\rho\geq1$
\begin{align}
    \|\theta_{q+1}^{(p)}\|^2_{L_t^2L^2}&\lesssim (\|M_l\|_{L_t^1L^1}+\delta_{q+2}+\sigma^{-1}\|M_l\|_{C_{t,x}^1})\left\|\sum_{\xi\in\Lambda^1\cup\Lambda^2}g_{(\xi)}^2\right\|_{L_t^1}\notag\\
    &\lesssim C_0^2(\delta_{q+2}+\lambda_{q+1}^{(2d+6)\alpha-\frac1{2N}})\lesssim C_0^2(\delta_{q+2}+ \lambda_{q+1}^{-\alpha})
   \leq \frac{C_\rho^2C_0^2}{4}\delta_{q+2},\label{bd:2theq+1pl2}
\end{align}
where we used \eqref{para22} and  conditions on the parameters to have $(2d+7)\alpha<\frac{1}{2N}$.

 For the general $L^u_tL^m$-norm with $m,u\in[1,\infty]$, by the same argument as \eqref{bd:2wq+1plp}, we have
\begin{align}
\|\theta_{q+1}^{(p)} \|_{L^u_tL^m}
&\lesssim l^{-2d-4}r_\perp^{\frac{d-1}{m}-\frac{d-1}{2}} r_\parallel^{\frac{1}{m}-\frac12}\eta^{\frac1u-\frac12}
\label{bd:2theq+1plp}.
\end{align}
Moreover, by the estimate for $h_{(\xi)}$ in \eqref{bd:gwnp} and Proposition \ref{lem:2chi}  we have
\begin{align}
\|\theta_{q+1}^{(c)} \|_{C_t}&\lesssim \|\theta_{q+1}^{(p)} \|_{C_tL^1}
\lesssim l^{-2d-4}r_\perp^{\frac{d-1}{2}} r_\parallel^{\frac12}\eta^{-\frac12}\lesssim 
\lambda_{q+1}^{(4d+8)\alpha-\frac1N}\lesssim\lambda_{q+1}^{-\alpha},\label{bd:2theq+1clp}\\
 \|   \theta_{q+1}^{(o)}\|_{C_tW^{1,\infty}}&\lesssim \sigma^{-1}\sum_{n\geq3}\sum_{\xi\in\Lambda^n}\|h_{(\xi)}\|_{L_t^\infty}\norm{\chi(\zeta|M_l|-n)\frac{n}{\zeta}\Gamma_{\xi}\(\frac{M_l}{|M_l|}\)\xi}_{W^{2,\infty}}\notag\\ 
 &\lesssim\sigma^{-1} l^{-6d-20}\lesssim \lambda_{q+1}^{(12d+40)\alpha-\frac1{2N}}\lesssim \lambda_{q+1}^{-\alpha}.\label{bd:2theq+1olp}
\end{align}
where we used \eqref{paralam}, \eqref{def:2musigma} and the  conditions on the parameters to deduce
$(12d+41)\alpha<\frac1{2N}$ and choose $a$ large enough to absorb the universal constant.
Then together with the bound \eqref{bd:2vqc1} and \eqref{para22} we show that
\begin{align}
\|\rho_{q+1}-\rho_q\|_{L^2_tL^2}&\leq\|\theta_{q+1}\|_{L^2_tL^2}+l\|\rho_q\|_{C_{t,x}^1}\leq \frac{C_\rho  C_0}2\delta_{q+2}^{1/2}+C\lambda_{q+1}^{-\alpha}+C_0l\lambda_{q}^{3d+2}
\leq C_\rho  C_0\delta_{q+2}^{1/2},\notag
\end{align}
once we choose $a$ large enough in order to absorb the universal constant. 
Then we have \eqref{bd:2vq+1-vql2} and then \eqref{bd:2vql2} for $\rho_{q+1}$.

 Furthermore, together with the  bounds \eqref{bd:2vqc1} and \eqref{bd:2theq+1plp}-\eqref{bd:2theq+1olp} we have 
\begin{align*}
    \|\rho_{q+1}-\rho_q\|_{L^r_tL^p}&\leq\|\theta_{q+1}\|_{L^r_tL^p}+l\|\rho_q\|_{C_{t,x}^1}\lesssim l^{-2d-4}r_\perp^{\frac{d-1}{p}-\frac{d-1}{2}} r_\parallel^{\frac{1}{p}-\frac12}\eta^{\frac1r-\frac12}+\lambda_{q+1}^{-\alpha}+C_0l\lambda_{q}^{3d+2}\notag\\ 
&\lesssim \lambda_{q+1}^{(4d+8)\alpha-\frac1N}+\lambda_{q+1}^{-\alpha}\lesssim 
\lambda_{q+1}^{-\alpha}\leq \delta_{q+2}^{1/2},\notag
\end{align*}
where we used \eqref{para22}, \eqref{paralam} and  conditions on the parameters to have
$(4d+9)\alpha<\frac1{N}$ and   choose $a$ large enough to absorb the universal constant. Consequently, we obtain the first bound in \eqref{bd:2rhoq+1-rhoql1}.

The second bound in \eqref{bd:2rhoq+1-rhoql1} can be derived by using  the above bounds again:
\begin{align*}
    \|\rho_{q+1}-\rho_q\|_{C_tL^1}&\leq \|\theta_{q+1}\|_{C_tL^1}+l\|\rho_q\|_{C_{t,x}^1}\lesssim l^{-2d-4}r_\perp^{\frac{d-1}{2}} r_\parallel^{\frac12}\eta^{-\frac12}+\lambda_{q+1}^{-\alpha}+C_0l\lambda_{q}^{3d+2}\notag\\ 
&\lesssim \lambda_{q+1}^{(4d+8)\alpha-\frac1N}+\lambda_{q+1}^{-\alpha}\lesssim 
\lambda_{q+1}^{-\alpha}\leq \delta_{q+2}^{1/2},
\end{align*}
where we use  conditions on the parameters to have
$(4d+9)\alpha<\frac1{N}$  choose $a$ large enough to absorb the universal constant.

 Since $\theta_{q+1}^{(p)}$ is non-negative,  by the estimates for  $\rho_q,\theta_{q+1}^{(c)}$ and $\theta_{q+1}^{(o)}$ in \eqref{bd:2vqc1}, \eqref{bd:2theq+1clp} and \eqref{bd:2theq+1olp} respectively we obtain 
\begin{align}
    \inf_{t\in [0,1]}(\rho_{q+1} - \rho_q)&\geq-\|\theta_{q+1}^{(c)}\|_{C_{t}}-\|\theta_{q+1}^{(o)}\|_{C_{t,x}^0}-\|\rho_l-\rho_q\|_{C_{t,x}^0}\geq -C\lambda_{q+1}^{-\alpha}-lC_0\lambda_q^{3d+2}\geq -\delta_{q+2}^{1/2},\notag
\end{align}
which yields the last bound in \eqref{bd:2rhoq+1-rhoql1}. Here we used \eqref{para22}  and choose $a$ large enough to absorb the universal constant.

Now we  estimate  $ \theta_{q+1}$ in $C_{t,x}^1$-norm. By the estimates for the building blocks and the amplitude functions in \eqref{int4theta}, \eqref{bd:gwnp} and Proposition \ref{lem:2chi} respectively, we have
\begin{align*}
\|\theta_{q+1}^{(c)}\|_{C_{t,x}^1}&\lesssim\|\theta_{q+1}^{(p)} \|_{C_{t,x}^1}\lesssim\sum_{n\geq3}\sum_{\xi\in\Lambda^n}\norm{\chi(\zeta|M_l|-n)(\frac{n}{\zeta})^{1/2}\Gamma_{\xi}\(\frac{M_l}{|M_l|}\)}_{C_{t,x}^1}\|\Theta_{(\xi,g)} \|_{C_{t,x}^1}\|g_{(\xi)}\|_{C_t^1}\notag\\ 
&\lesssim l^{-4d-12}\lambda_{q+1}\mu r_\parallel^{-\frac12}r_\perp^{-\frac{d-1}{2}}\sigma\eta^{-2}\lesssim
\lambda_{q+1}^{(8d+24)\alpha+3d+\frac{3}{2}},\notag\\
     \|   \theta_{q+1}^{(o)}\|_{C_{t,x}^1}&\lesssim \sigma^{-1}\sum_{n\geq3}\sum_{\xi\in\Lambda^n}\|h_{(\xi)}\|_{C_t^1}\norm{\chi(\zeta|M_l|-n)\frac{n}{\zeta}\Gamma_{\xi}\(\frac{M_l}{|M_l|}\)\xi}_{W^{2,\infty}}\lesssim \lambda_{q+1}^{(12d+40)\alpha+d}.
\end{align*}
 where we used  \eqref{bd:lempara2} and \eqref{def:2musigma}. By choosing $(8d+24)\alpha<1/2$ we deduce 
\begin{align}
\|\rho_{q+1}\|_{C_{t,x}^1}\leq\|\rho_{l}\|_{C_{t,x}^1}+\|\theta_{q+1}\|_{C_{t,x}^1}\leq C_0\lambda_{q}^{3d+2}+ \lambda_{q+1}^{3d+2}\leq C_0\lambda_{q+1}^{3d+2}.\notag
\end{align}
which implies  \eqref{bd:2vqc1} for $\rho_{q+1}$.

 Finally, we consider the bound of $\theta_{q+1}^{(p)}$ in $W^{1,s}$-norm which will be used below. By the estimates on the building blocks in \eqref{int4theta}, \eqref{bd:gwnp} and on the amplitude functions in Proposition \ref{lem:2chi} we obtain
\begin{align}
\|\theta_{q+1}^{(p)} \|_{L^1_tW^{1,s}}&\lesssim\sum_{n\geq3}\sum_{\xi\in\Lambda^n}\norm{\chi(\zeta|M_l|-n)(\frac{n}{\zeta})^{1/2}\Gamma_{\xi}\(\frac{M_l}{|M_l|}\)}_{C_{t,x}^1}\|\Theta_{(\xi,g)} \|_{C_tW^{1,s}}\|g_{(\xi)}\|_{L_t^1}\notag\\ 
&\lesssim l^{-4d-12}\lambda_{q+1}r_\parallel^{\frac{1}{s}-\frac12}r_\perp^{\frac{d-1}{s}-\frac{d-1}{2}}\eta^{\frac12}
\lesssim\lambda_{q+1}^{(8d+24)\alpha-\frac1N}\lesssim\lambda_{q+1}^{-\alpha} ,\label{bd:2theq+1w1s}
\end{align}
where we  used  \eqref{paralam} and  conditions on the parameters to have $(8d+25)\alpha<\frac1N$.

\subsubsection{Construction of the perturbation $\overline w_{q +1}$}\label{sec:2defvq+12}
Let us now proceed with the construction of the perturbation $\overline w_{q+1 }$  by employing  the generalized intermittent jets and temporal jets introduced in Section \ref{gij}. 

As the next step, we shall define certain amplitude functions used in the definition of the perturbations  $\overline w_{q +1}$. For any $\xi\in \overline{\Lambda}$, we define
\begin{align}
A:=2\sqrt{ {l}^2+|\mathring{ {R}}_{ {l}}|^2}
,\ \
a_{({\xi})}:=A^{1/2}\gamma_{{\xi}}(\Id-\frac{\mathring{ {R}}_{ {l}}}{A}),\notag
\end{align}
where $\gamma_{ {\xi}}$ is introduced in Lemma \ref{lem:3}.
Since we have
$$|\Id-\frac{\mathring{ {R}}_{ {l}}}{A}-\Id|\leq1/2,$$
by Lemma \ref{lem:3}
 it  follows that
\begin{align}
A\ \Id-\mathring{ {R}}_{ {l}}=\sum_{{\xi}\in \overline {\Lambda}}A\gamma_{{ \xi}}^2(\Id-\frac{\mathring{ {R}}_{ {l}}}{A}){ \xi}\otimes{\xi}=\sum_{{ \xi}\in \overline {\Lambda}}a_{({\xi})}^2{\xi}\otimes{ \xi}.\label{**}
\end{align}
We  have the following estimate for the amplitude function.
\bp\label{prop:bd:2axi}
For $\xi\in \overline{\Lambda}$ and $ N\in\mathbb{N}_0$ we have
\begin{align}
\|a_{({ \xi})}\|_{C_{t,x}^N}\lesssim  {l}^{-2d-3-(d+4)N}.\label{bd:2axi}
\end{align}
\ep
 The proof of this proposition is given in Appendix \ref{sec:est:ampl}.
 
Now  recalling the temporal jets  $H_{(\xi)}(t)$  in Section \ref{sec:tj} we define 
$$\overline {W}_{( {\xi,g} )}(x,t)= \overline {W}_{({ {\xi}})}(x, {H}_{(\xi)}(t)),$$
and similarly define $ \overline {V}_{( {\xi,g})},\overline {\varphi}_{( {\xi,g})}$ and other terms appearing in Section \ref{sec:bbns}. 

With these preparations in hand, we define the principal part
\begin{align}
\overline w_{q+1 }^{(p)}:&=\sum_{{ {\xi}}\in \overline {\Lambda}}a_{({ {\xi}})} \overline {W}_{( {\xi,g})} {g}_{(\xi)},\notag
\end{align}
and by the identity \eqref{**} and the fact that $g_{(\xi)}$ have disjoint support we have
\begin{align}
(&\overline w_{q+1}^{(p)}+w_{q+1}^{(p)})\otimes( \overline w_{q+1}^{(p)}+w_{q+1}^{(p)})+\mathring{{R}}_{{l}}=
\overline w_{q+1}^{(p)}\otimes \overline w_{q+1}^{(p)}+w_{q+1}^{(p)}\otimes  w_{q+1}^{(p)}+\mathring{{R}}_{{l}}\notag\\ &=\sum_{\xi\in\overline{\Lambda}}a_{({\xi})}^2\overline{W}_{({\xi,g})}\otimes \overline{W}_{({\xi,g})}g_{(\xi)}^2-\sum_{\xi\in\overline{\Lambda}}a_{({\xi})}^2{\xi}\otimes{\xi}+w_{q+1}^{(p)}\otimes  w_{q+1}^{(p)}+A \Id\notag\\ 
&=\sum_{{\xi}\in\overline{\Lambda}}a_{({\xi})}^2\mathbb{P}_{\neq0}(\overline{W}_{({\xi,g})}\otimes \overline{W}_{({\xi,g})})g_{(\xi)}^2+w_{q+1}^{(p)}\otimes  w_{q+1}^{(p)},\label{2.12}
\end{align}
where we use the notation $\mathbb{P}_{\neq0}f:=f-\fint_{\mathbb{T}^d}f\dif x$.

We define the incompressibility corrector by
\begin{align}
\overline w_{q+1}^{(c)}&:=\sum_{{\xi}\in\overline{\Lambda}}-\frac{1}{(\overline n_*\lambda_{q+1})^2}a_{({\xi})}\nabla\overline{\Phi}_{({\xi,g})}{\xi}\cdot\nabla\overline{\psi}_{({\xi,g})}g_{({\xi})}+\nabla a_{({\xi})}:\overline{V}_{({\xi,g})}g_{({\xi})}.\label{def:wq+1c}
\end{align}
Here $(\nabla a_{({\xi})}:\overline{V}_{({\xi,g})})^i=\sum_{j=1}^d\partial_ja_{({\xi})}\overline{V}_{({\xi,g})}^{ij},\ i=1,2,...,d$.

By identity \eqref{divOmegans} we have
\begin{align}
\overline w_{q+1}^{(p)}+\overline w_{q+1}^{(c)}
&=\sum_{{\xi}\in\overline{\Lambda}}(a_{({\xi})}\div \overline{V}_{({\xi,g})}+\nabla a_{({\xi})}:\overline{V}_{({\xi,g})})g_{({\xi})}=\sum_{{\xi}\in\overline{\Lambda}}\div (a_{({\xi})}\overline{V}_{({\xi,g})})g_{({\xi})}.\label{wq+12p+wq+12c}
\end{align}
Since $a_{({\xi})}\overline{V}_{({\xi,g})}$ is skew-symmetric, we obtain
$$\div ( \overline w_{q+1}^{(p)}+\overline w_{q+1}^{(c)})=0.$$

Next we looking back at \eqref{2.12}, we see that there are still two bad terms that we need to deal with (in fact, we need to deal with the divergence of the these terms).  To this end, we   introduce the temporal corrector
\begin{align}\label{2.14}
 \overline w_{q+1}^{(t)}&:=\frac1{\overline{\mu}}\sum_{{\xi}\in\overline{\Lambda}}\mathbb{P}_{\neq0}\mathbb{P}_{\rm H}(a_{({\xi})}^2\overline{\psi}^2_{({\xi,g})}\overline{\phi}^2_{({\xi,g})}{\xi})g_{({\xi})},
\end{align}
where $\mathbb{P}_{\rm H}$ is the Helmholtz projection. By a direct computation and \eqref{*5} we obtain
\begin{align}
&\partial_t  \overline w_{q+1}^{(t)}+\sum_{{\xi}\in\overline{\Lambda}}\mathbb{P}_{\neq0}(a_{({\xi})}^2g_{({\xi})}^2 \div(\overline{W}_{({\xi,g})}\otimes \overline{W}_{({\xi,g})}))\notag
\\
&\qquad= \frac{1}{\overline{\mu}}\sum_{{\xi}\in\overline{\Lambda}}\mathbb{P}_H\mathbb{P}_{\neq0}\partial_t(a_{({\xi})}^2g_{({\xi})}\overline{\psi}^2_{({\xi,g})}\overline{\phi}^2_{({\xi,g})}{\xi})
-\frac{1}{\overline{\mu}}\sum_{{\xi}\in\overline{\Lambda}}\mathbb{P}_{\neq0}( a^2_{({\xi})}g_{({\xi})}\partial_t(\overline{\psi}^2_{({\xi,g})}\overline{\phi}^2_{({\xi,g})}{\xi}))\notag
\\&\qquad=(\mathbb{P}_H-\Id)\frac{1}{\overline{\mu}}\sum_{{\xi}\in\overline{\Lambda}}\mathbb{P}_{\neq0}\partial_t(a_{({\xi})}^2g_{({\xi})}\overline{\psi}^2_{({\xi,g})}\overline{\phi}^2_{({\xi,g})}{\xi})
+\frac{1}{\overline
{\mu}}\sum_{{\xi}\in\overline{\Lambda}}\mathbb{P}_{\neq0}(\partial_t (a^2_{({\xi})}g_{({\xi})})\overline{\psi}^2_{({\xi,g})}\overline{\phi}^2_{({\xi,g})}{\xi}).\label{int1}
\end{align}
Note that the first term on the right hand side can be viewed as a pressure term. 

 Similarly as before, to handle the undesirable term of the form $(\cdot) (g_{(\xi)}^2-1)$ in \eqref{2.12},  we define a new perturbation term
\begin{align}
\overline  w_{q +1}^{(o)}:&=-\sigma^{-1}\mathbb{P}_{\rm{H}}\mathbb{P}_{\neq0}\sum_{\xi\in\overline \Lambda}h_{(\xi)}\div(a_{(\xi)}^2\xi\otimes\xi),\notag
\end{align}
which implies that 
\begin{align}
 \partial_t \overline  w_{q +1}^{(o)}  +\sum_{\xi\in\overline{\Lambda}}(g_{(\xi)}^2-1)\div(a_{({\xi})}^2{\xi}\otimes{\xi})=\nabla p-\sigma^{-1} \mathbb{P}_{\rm{H}}\mathbb{P}_{\neq0}\sum_{\xi\in\overline \Lambda}h_{(\xi)}\partial_t\div(a_{(\xi)}^2\xi\otimes\xi).\label{eq:partbarwq+1o}
\end{align}
Here $p$ denotes some pressure term. 
 
Finally, the total perturbation $\overline w_{q+1 }$  and the  new velocity $v_{q+1 }$ are defined by
$$\overline w_{q+1 }:=\chi_{q}\overline w_{q +1}^{(p)}+\chi_{q}\overline w_{q+1 }^{(c)}+\chi_{q}^2\overline w_{q +1}^{(t)}+\chi_{q}^2\overline w_{q +1}^{(o)},\ \ v_{q +1}:= {v}_{ {l}}+w_{q +1}+\overline w_{q +1},$$
which are mean-zero and divergence-free.
Here $\chi_q$ is a smooth time cut off satisfying $\chi_{q}=0$ for $t\leq T_{q+1}$, $\chi_{q}=1
$ for $t\geq \frac{T_q+T_{q+1}}{2}$ and  for $n\geq0$
\begin{align}
   \| \chi_q\|_{C^n_t}\lesssim (T_q-T_{q+1})^{-n}\lesssim {l}^{-n}.\label{bd:2chiqcn}
\end{align}
Then we obtain $v_{q +1}=0$  on $[0,T_{q+1}]$.

\subsubsection{Estimates of $\overline w_{q+1}$}

 Recall that $\overline w_{q+1}$ is  defined in Section \ref{sec:2defvq+12}.  
First we  estimate $\overline w^{(p)}_{q +1}$ in $L^2_tL^2$-norm by applying the improved H\"older inequality from Lemma \ref{ihiot} in the spatial  variable. By the definition of $a_{(\xi)}$, the bounds on $\overline W_{(\xi)}$ and $a_{(\xi)}$ in \eqref{int4ns}  and \eqref{bd:2axi}
we obtain 
\begin{align}
\|\overline w_{q+1 }^{(p)} (t)\|^2_{L^2}&\lesssim\sum_{{\xi}\in \overline{\Lambda}}\(\|a_{({ \xi})} (t)\|_{L^2}^2+( {r}_\perp\lambda_{q +1})^{-1}\|a_{({ \xi})} \|_{C_{t,x}^1}^2\)\| \overline{W}_{({\xi,g})} \|^2_{C_tL^2} {g}_{(\xi)}^2(t)\notag\\
&\lesssim\sum_{{ {\xi}}\in \overline {\Lambda}}(\|\mathring{ {R}}_{ 
{l}}(t)\|_{L^1}+ {l}+{\lambda}_{q +1}^{(12d+28)\alpha-\frac{1}{N}}) {g}_{(\xi)}^2(t)\lesssim (\|\mathring{ {R}}_{ {l}}(t)\|_{L^1}+\delta_{q+1}) \sum_{{ {\xi}}\in \overline {\Lambda}}{g}_{(\xi)}^2(t),\notag
\end{align}
where we used the choice of parameters in \eqref{paralam} and  conditions on the parameters to have $(12d+28)\alpha-\frac{1}{N}<-\alpha<-2\beta $. Then we again apply the improved H\"older inequality  in Lemma \ref{ihiot} in time. By the bounds on  $g_{(\xi)}$ and $\mathring{R}_l$ in \eqref{bd:gwnp}, \eqref{bd:2rll1l1} and \eqref{bd:2rlcn} we obtain
for some universal constant $ {C}_v\geq1$
\begin{align}
    \|\overline w_{q +1}^{(p)}\|_{L^2_tL^2}^2&\lesssim  (\|\mathring{R}_q\|_{L^1_tL^1}+\delta_{q+1}+ {\sigma}^{-1}\|\mathring{R}_l\|_{C_{t,x}^1})\| \sum_{{ {\xi}}\in \overline {\Lambda}}{g}_{(\xi)}^2\|_{L_t^1}\lesssim C_0^2(\delta_{q+1}+ {\sigma}^{-1} {l}^{-d-3})\notag \\
    &\lesssim C_0^2(\delta_{q+1}+\lambda_{q +1}^{(2d+6)\alpha-\frac1{2N}})
    \leq\frac{ {C}_v^2C_0^2}{16}\delta_{q+1},\label{bd:2wq+12pl2}
\end{align}
where we used  conditions on the parameters to have $(2d+6)\alpha-\frac{1}{2N}<-\alpha<-2\beta$.

Then we turn to bound the perturbations $\overline w_{q+1}$ in general $L^m_tL^n$-norm with $m\in[1,\infty],n\in(1,\infty)$. Recalling  the choice of parameters in \eqref{paralam} and \eqref{def:2musigma}, the estimates on the building blocks and the amplitude functions $a_{(\xi)}$ in  \eqref{int2ns}-\eqref{int4ns}, \eqref{bd:gwnp} and \eqref{bd:2axi} respectively, we have that
\begin{align}
\|\overline w_{q+1 }^{(p)} \|_{L^m_tL^n}&\lesssim \sum_{{ \xi}\in \overline{\Lambda}}\|a_{({ \xi})} \|_{C_{t,x}^0}\| \overline{W}_{({ \xi,g})} \|_{C_tL^n}\| {g}_{(\xi)}\|_{L_t^m}\lesssim  {l}^{-2d-3} {r}_\perp^{\frac{d-1}{n}-\frac{d-1}{2}}  {r}_\parallel^{\frac{1}{n}-\frac12} {\eta}^{\frac1m-\frac12},\label{bd:2wq+12plp}
\end{align}
\begin{align}
\|\overline w_{q+1 }^{(c)}& \|_{L^m_tL^n}\lesssim\sum_{{ \xi}\in \overline{\Lambda}}\frac{1}{\lambda_{q +1}^2}\|a_{({ \xi})} \nabla \overline{\Phi}_{({ \xi,g})} { \xi}\cdot\nabla \overline{\psi}_{({ \xi,g})} \|_{C_tL^n}\| {g}_{(\xi)}\|_{L_t^m}+\sum_{{ \xi}\in \overline{\Lambda}}\|\nabla a_{({\xi})} : \overline{V}_{({ \xi,g})} \|_{C_tL^n}\| {g}_{(\xi)}\|_{L_t^m}\notag\\
&\lesssim  {l}^{-3d-7} {r}_\perp^{\frac{d-1}{n}-\frac{d-1}{2}}  {r}_\parallel^{\frac{1}{n}-\frac12}(\frac{ {r}_\perp}{  {r}_\parallel}+ \frac1{\lambda_{q+1}}) {\eta}^{\frac1m-\frac12}\lesssim  {l}^{-3d-7} {r}_\perp^{\frac{d-1}{n}-\frac{d-1}{2}}  {r}_\parallel^{\frac{1}{n}-\frac12}
{\eta}^{\frac1m-\frac12}{\lambda}_{q+1}^{-\frac1{N}},\label{bd:2wq+12clp}
\end{align}
\begin{align}
\|\overline w_{q+1 }^{(t)} \|_{L^m_tL^n}&\lesssim\overline  {\mu}^{-1}\sum_{{ \xi}\in \overline{\Lambda}}\|a_{({ \xi})} \|_{C_{t,x}^0}^2\| \overline{\psi}_{({\xi,g})} \|_{C_tL^{2n}}^2\|\overline {\phi}_{({ \xi,g})} \|^2_{C_tL^{2n}}\| {g}_{(\xi)}\|_{L_t^m}\notag\\
&\lesssim  {l}^{-4d-6} \overline{\mu}^{-1} {r}_\perp^{\frac{d-1}{n}-d+1} {r}_\parallel^{\frac{1}{n}-1} {\eta}^{\frac1m-\frac12}\lesssim  {l}^{-4d-6} {r}_\perp^{\frac{d-1}{n}-\frac{d-1}2} {r}_\parallel^{\frac{1}{n}-\frac12} {\eta}^{\frac1m-\frac12} {\lambda}_{q+1}^{-\frac1{2N}},\label{bd:2wq+12tlp}
\end{align}
and together with the boundedness of ${h}_{(\xi)}$ in \eqref{bd:gwnp}, we obtain
\begin{align}
\|\overline w_{q +1}^{(o)} \|_{C_tW^{1,\infty}}&\lesssim {\sigma}^{-1}\| {h}_{(\xi)}\|_{L_t^\infty}\|a_{(\xi)}^2\|_{C_tW^{2+\alpha,\infty}}\lesssim  {\sigma}^{-1} {l}^{-6d-15}\lesssim \lambda_{q+1 }^{(12d+30)\alpha-\frac1{2N}}\lesssim \lambda_{q +1}^{-\alpha},\label{bd:2wq+12olp}
\end{align}
where we used  conditions on the parameters to have $(12d+31)\alpha<\frac1{2N}$. 
Combining the choice of parameters in \eqref{paralam}, the bound \eqref{bd:2wq+12pl2} and the fact that $\|\chi_q\|_{L_t^\infty}\leq1$  we obtain
\begin{align}
\|\overline w_{q+1 }\|_{L^2_tL^2}&\leq\frac{ {C}_vC_0}4\delta_{q+1}^{1/2}+Cl^{-4d-6}{\lambda}_{q+1}^{-\frac1{2N}}+ C\lambda_{q +1}^{-\alpha}\leq\frac{3 {C}_vC_0}8\delta_{q+1}^{1/2},\label{bd:barwq+1l2}
\end{align}
where we used the conditions on the parameters to have $(8d+12)\alpha-\frac1{2N}<-\alpha<-\beta $. In the last inequality we choose $a$ large enough to absorb the universal constant.

 The above inequality together with the bounds in \eqref{bd:2vqc1}, \eqref{para22} and \eqref{bd:2wq+1l2} yields:
\begin{align*}
\|v_{q+1}-v_{q }\|_{L^2_tL^2}&\leq \|w_{q+1}\|_{L^2_tL^2}+ \|\overline w_{q+1}\|_{L^2_tL^2}+l\|v_{q }\|_{C_{t,x}^1}
\leq \frac34{C}_vC_0\delta_{q+1}^{1/2}+lC_0\lambda_{q }^{4d+3}\leq 
{C}_vC_0\delta_{q+1}^{1/2}.
\end{align*}
Then \eqref{bd:2vq+1-vql2} and \eqref{bd:2vql2} holds for $v_{q+1}$.

Similarly, recalling the choice of  parameters  in  \eqref{paralam}  we have the following bounds:
\begin{align}
    \|\overline w_{q+1 }\|_{L^r_tL^p}&\lesssim l^{-4d-6} {r}_\perp^{\frac{d-1}{p}-\frac{d-1}{2}}  {r}_\parallel^{\frac{1}{p}-\frac12} {\eta}^{\frac1r-\frac12}+\lambda_{q+1 }^{-\alpha}\lesssim \lambda_{q +1}^{(8d+12)\alpha-\frac1N}+\lambda_{q+1 }^{-\alpha}\lesssim \lambda_{q+1 }^{-\alpha},\label{bd:barwq+1lrp}\\
    \|\overline w_{q+1 }\|_{C_tL^1}& \lesssim  \|\overline w_{q+1 }\|_{C_tL^{1+\epsilon}}\lesssim  l^{-4d-6} {r}_\perp^{\frac{d-1}{2}}  {r}_\parallel^{\frac12} {\eta}^{-\frac12}\lambda_{q+1}^{d\epsilon}+{\sigma}^{-1} {l}^{-6d-15}
   \notag\\\ 
   &\lesssim \lambda_{q +1}^{(8d+13)\alpha-\frac1N}+ \lambda_{q+1 }^{(12d+30)\alpha-\frac1{2N}}\lesssim \lambda_{q+1 }^{-\alpha},\label{bd:barwq+1ctl1} 
    \end{align}
where we  choose $\epsilon>0$ small enough such that $d\epsilon<\alpha$ and  choose the parameters such that $(12d+31)\alpha<\frac1{2N}$. 
Together with the bounds on $v_q, w_{q+1}$ in \eqref{bd:2vqc1}, \eqref{bd:wq+1lrp} and \eqref{bd:wq+1ctl1} respectively we obtain 
\begin{align*}
     \|v_{q+1}-v_{q }\|_{L^r_tL^p}&\lesssim    \|w_{q+1}\|_{L^r_tL^p}+  \|\overline w_{q+1}\|_{L^r_tL^p}+   l\|v_{q }\|_{C_{t,x}^1}\lesssim \lambda_{q+1}^{-\alpha}+lC_0\lambda_{q }^{4d+3}\lesssim \lambda_{q+1}^{-\alpha}\leq 
\delta_{q+1}^{1/2},\\
  \|v_{q+1}  -v_{q }\|_{C_tL^1}& \lesssim    \|w_{q+1}\|_{C_tL^1}+  \|\overline w_{q+1}\|_{C_tL^1}+   l\|v_{q }\|_{C_{t,x}^1}\lesssim \lambda_{q+1}^{-\alpha}+lC_0\lambda_{q }^{4d+3}\lesssim \lambda_{q+1}^{-\alpha}\leq 
\delta_{q+1}^{1/2},
\end{align*}
where in the last inequality we choose $a$ large enough to absorb the universal constant.  Then we obtain the first two bounds in \eqref{bd:2vq+1-vqlpr}.

Next we estimate the $C_{t,x}^2$-norm. Taking into account the fact that  $$\partial_t\( \overline {V}_{( {\xi,g})}(t)\)= {g}_{(\xi)}(t)(\partial_t \overline  {V}_{( {\xi})})( {H}_{(\xi)}(t)),$$ 
 and using the estimates on the building blocks $ \overline {\psi}_{({ {\xi}})}, \overline {\phi}_{({ {\xi}})}, \overline {V}_{({ {\xi}})},{g}_{(\xi)}$ in  \eqref{int2ns}-\eqref{int4ns}, \eqref{bd:gwnp} respectively, the estimates on $a_{({ {\xi}})}$ in \eqref{bd:2axi} and  the identity \eqref{wq+12p+wq+12c} we have
\begin{align}
\|\overline w_{q +1}^{(p)}+\overline w_{q +1}^{(c)}\|_{C_{t,x}^2}&\lesssim\sum_{{ {\xi}}\in \overline {\Lambda}} \|a_{({ {\xi}})}\|_{C_{t,x}^{3}}(\| {g}_{(\xi)}\nabla \overline {V}_{({ {\xi,g}})}\|_{C_{t,x}^{2}} +\| {g}_{(\xi)} \overline {V}_{({ {\xi,g}})}\|_{C_{t,x}^{2}})\notag\\
&\lesssim l^{-5d-15}\lambda_{q +1}^2 \overline{\mu}^2  {r}_\parallel^{-\frac12} {r}_\perp^{-\frac{d-1}{2}} {\sigma}^2 {\eta}^{-\frac52}\lesssim{\lambda}_{q +1}^{(10d+30)\alpha+4d+\frac{5}{2}},\notag
\end{align}

\begin{align}
\|\overline w_{q +1}^{(t)}\|_{C_{t,x}^2}&\lesssim\frac{1}{\overline {\mu}}\sum_{{ {\xi}}\in\overline  {\Lambda}}\sum_{i=0}^2\|a^2_{({ {\xi}})}\overline  {\psi}^2_{({ {\xi,g}})} \overline {\phi}^2_{({ {\xi,g}})} {g}_{(\xi)}\|_{C^i_tW^{2-i+\alpha,\infty}}\notag\\
&\lesssim  {l}^{-6d-14}\lambda_{q+1 }^{2+\alpha}\overline {\mu}  {r}_\parallel^{-3} {r}_\perp^{-d+3} {\sigma}^2 {\eta}^{-\frac52}\lesssim \lambda_{q +1}^{(12d+29)\alpha+4d+\frac{5}{2}},\notag
\end{align}
and by the estimates on ${h}_{(\xi)},a_{(\xi)}$ in \eqref{bd:gwnp}, \eqref{bd:2axi}, we have
\begin{align}
\|\overline w_{q +1}^{(o)} \|_{C_{t,x}^2}&\lesssim {\sigma}^{-1}\sum_{\xi\in\overline \Lambda}\| {h}_{(\xi)}\|_{C^2}\|\div(a_{(\xi)}^2\xi\otimes\xi)\|_{C_{t,x}^{ 2+\alpha}}\lesssim  {\eta}^{-2} {l}^{-7d-19}\lesssim \lambda_{q +1}^{(14d+38)\alpha+2d}.\notag
\end{align}
 Thus by  using \eqref{bd:2chiqcn} and $(12d+33)\alpha<\frac12$ we  obtain 
$$\|\overline w_{q +1}\|_{C_{t,x}^2}\leq (\|\chi_q\|_{C_t^1}^2 +\|\chi_q\|_{C_t^2})\lambda_{q+1 }^{(12d+29)\alpha+4d+\frac{5}{2}}\leq \lambda_{q+1 }^{4d+3}.$$

  Combining the bounds on $ w_{q +1}$ in \eqref{bd:2wq+1ctx2} we get  \eqref{bd:2vqc1} for $v_{q+1}$.
$$\|v_{q+1 }\|_{C_{t,x}^2}\leq \|v_q\|_{C_{t,x}^2}+\| w_{q +1}\|_{C_{t,x}^2}+\|\overline w_{q +1}\|_{C_{t,x}^2}\leq  C_0\lambda_q^{4d+3} +2\lambda_{q+1 }^{4d+3}\leq C_0\lambda_{q+1 }^{4d+3}.$$

We conclude this part with the bound in $W^{1,s}$-norm. By the  estimates for the building blocks $ \overline {\psi}_{({ {\xi}})}, \overline {\phi}_{({ {\xi}})}$, $\overline {V}_{({ {\xi}})},{g}_{(\xi)}$ in  \eqref{int2ns}-\eqref{int4ns}, \eqref{bd:gwnp} respectively and the estimates for $a_{({ \xi})}$ in \eqref{bd:2axi}   it follows  that 
\begin{align}
\|\overline w_{q +1}^{(p)} +\overline w_{q+1 }^{(c)}\|_{L^1_tW^{1,s}}&\lesssim\sum_{{ {\xi}}\in \overline {\Lambda}}\|a_{({ {\xi}})}\|_{C_{t,x}^{2}}\| \overline {V}_{( {\xi,g})}\|_{C_tW^{2,s}}\|g_{(\xi)}\|_{L_t^1}\notag\\
&\lesssim   {l}^{-4d-11}  {r}_\perp^{\frac{d-1}{s}-\frac{d-1}{2}}  {r}_\parallel^{\frac{1}{s}-\frac12}\lambda_{q+1 } {\eta}^{\frac12}\lesssim \lambda_{q+1 }^{(8d+22)\alpha-\frac1N}\lesssim\lambda_{q+1 }^{-\alpha},\notag
\end{align}
\begin{align}
\|\overline w_{q+1 }^{(t)} \|_{L^1_tW^{1,s}}&\lesssim\frac{1}{ \overline{\mu}}\sum_{{ {\xi}}\in \overline {\Lambda}}\|a_{({ {\xi}})}^2\|_{C_{t,x}^{1}}\| \overline {\psi}^2_{({ {\xi,g}})}\overline  {\phi}^2_{({ {\xi,g}})}\|_{C_tW^{1,s}}\|g_{(\xi)}\|_{L_t^1}\lesssim  {l}^{-5d-10}\overline {\mu}^{-1} {r}_\perp^{\frac{d-1}{s}-d+1} {r}_\parallel^{\frac{1}{s}-1}\lambda_{q +1} {\eta}^{\frac12}\notag\\
&\lesssim  {l}^{-5d-10} {r}_\perp^{\frac{d-1}{s}-\frac{d-1}2} {r}_\parallel^{\frac{1}{s}-\frac12}\lambda_{q+1 } {\eta}^{\frac12}\lesssim \lambda_{q +1}^{(10d+20)\alpha-\frac1N}\lesssim\lambda_{q+1 }^{-\alpha},\notag
\end{align}
where we  used  \eqref{paralam} and  conditions on the parameters to have $(10d+21)\alpha<\frac1N$. Then since $\|\chi_q\|_{L_t^\infty}\leq1$, together with the the estimate on $\overline w_{q +1}^{(o)}$ in \eqref{bd:2wq+12olp}  we obtain
\begin{align}
     \|\overline {w}_{q+1}\|_{L^1_tW^{1,s}}\lesssim \lambda_{q+1}^{-\alpha}.\label{bd:barwq+1w1s}
\end{align}

 Hence together with  \eqref{para22}  and  the bound on ${w}_{q+1}$ in \eqref{bd:2wq+1w1p}  we deduce
\begin{align*}
\|v_{q+1}-v_{q }\|_{L^1_tW^{1,s}}&\leq \|{w}_{q+1}\|_{L^1_tW^{1,s}}+ \|\overline {w}_{q+1}\|_{L^1_tW^{1,s}}+l\|v_{q }\|_{C_{t,x}^2}\leq \lambda_{q+1}^{-\alpha} +lC_0\lambda_{q }^{4d+3}\leq \delta_{q+1}^{1/2},
\end{align*}
which implies the last estimate in  \eqref{bd:2vq+1-vqlpr}. Here in the last inequality we choose $a$ large enough to absorb the universal constant.

\subsection{The estimates of the stress terms}\label{sec:error}
 In this section we complete the proof of Proposition \ref{prop:case2}  by proving the remaining estimates on the stress terms in \eqref{bd:2rql1}. The stress term $M_l$ will be canceled by the  perturbation $(w_{q+1},\theta_{q+1})$, as will be showed in Section \ref{sec:defmq+1}. The estimate of the new stress term $M_{q+1}$ will be estimated in Section \ref{sec:estmq+1}. The  stress term $\mathring{R}_l$ will be canceled by the  perturbation $\overline w_{q+1}$ will be shown in Section \ref{sec:defrq+12}. The estimate of the new stress term $\mathring{R}_{q+1}$ is contained in Section \ref{sec:estrq+12}.
\subsubsection{Construction of the  stress term $M_{q+1}$}\label{sec:defmq+1}
First we recall that the supports of $g_{(\xi)}$ are disjoint for different $\xi\in \Lambda^1\cup\Lambda^2\cup\overline\Lambda$. As a result,  based on the definitions of the perturbations, we have that 
$$( \chi_q\overline w_{q+1}^{(p)}+ \chi_q\overline w_{q+1}^{(c)}+ \chi_q^2\overline w_{q+1}^{(t)})\theta^{(p)}_{q+1}=0.$$

Then together with the identity   $\div( v_{q+1}\theta_{q+1}^{(c)})=\theta_{q+1}^{(c)}\div v_{q+1}=0$ we obtain that

\begin{align*}
 -\div &M_{q+1}\notag\\ 
 &=\partial_t\theta_{q+1}+\div (w_{q+1}^{(p)}\theta_{q+1}^{(p)}-M_l)(:=\div M_{osc})\\
 &-\kappa\Delta \theta_{q+1}+\div\(v_l\theta_{q+1}+(w_{q+1}+\overline w_{q+1})(\rho_l+\theta_{q+1}^{(o)})+(w_{q+1}^{(c)}+ \chi_q^2\overline w_{q+1}^{(o)})\theta_{q+1}^{(p)}\),(:=\div M_{lin})
\end{align*}
where we use the inverse divergence operator $\mathcal{R}_1$ defined in Section \ref{tamr} to define the linear error by
\begin{align*} 
M_{lin}:&=-\kappa\mathcal{R}_1\Delta \theta_{q+1}+v_l\theta_{q+1}+(w_{q+1}+\overline w_{q+1})(\rho_l+\theta_{q+1}^{(o)})+(w_{q+1}^{(c)}+ \chi_q^2\overline w_{q+1}^{(o)})\theta_{q+1}^{(p)}.
\end{align*}

To define the oscillation error $M_{osc}$, by the identities \eqref{2wq+1ptheq+1p} and \eqref{eq:2parttheo} we have
\begin{align*}
 \div M_{osc}&=\partial_t\theta_{q+1}+\div (\theta_{q+1}^{(p)}w_{q+1}^{(p)}-M_l)\notag\\ 
 &=\mathbb{P}_{\neq0} \(\partial_t\theta^{(p)}_{q+1}+\div (\theta_{q+1}^{(p)}w_{q+1}^{(p)}-M_l)+\partial_t\theta^{(o)}_{q+1}\)\\
 &=\sum_{n\geq3}\sum_{\xi\in\Lambda^n}\mathbb{P}_{\neq0}\(\partial_t\[\chi(\zeta|M_l|-n)(\frac{n}{\zeta})^{\frac12}\Gamma_{\xi}\(\frac{M_l}{|M_l|}\)g_{(\xi)}\]\Theta_{(\xi,g)}\)(:=\div M_{osc,t})\notag\\ 
 &+\mathbb{P}_{\neq0}\[\nabla[\chi(\zeta|M_l|-n)\frac{n}{\zeta}\Gamma_{\xi}\(\frac{M_l}{|M_l|}\)]g_{(\xi)}^2\mP_{\neq0}(W_{(\xi,g)}\Theta_{(\xi,g)})\](:=\div M_{osc,x})\\
 &+\mathbb{P}_{\neq0}\[\chi(\zeta|M_l|-n)(\frac{n}{\zeta})^{\frac12}g_{(\xi)}\Gamma_{\xi}\(\frac{M_l}{|M_l|}\)(\partial_t\Theta_{(\xi,g)}+(\frac{n}{\zeta})^{\frac12}g_{(\xi)}\div(W_{(\xi,g)}\Theta_{(\xi,g)}))\]\\
 &+\div \(\sum_{n\geq3}\chi(\zeta|M_l|-n)\frac{n}{\zeta}\frac{M_l}{|M_l|}-M_l\)(:=\div M_{osc,c})\notag\\
 &-\sigma^{-1}\sum_{n\geq3}\sum_{\xi\in\Lambda^n}h_{(\xi)}\partial_t\div \[\chi(\zeta|M_l|-n)\frac{n}{\zeta}\Gamma_{\xi}\(\frac{M_l}{|M_l|}\)\xi\](:=\div M_{osc,o}),\notag
\end{align*}
where the last third term equals to 0 by \eqref{eq:2ptthe+n}. Now by  using the inverse divergence operators $\mathcal{R}_1,\mathcal{B}_1$ introduced in Section \ref{tamr}  we define $M_{osc}:=M_{osc,t}+M_{osc,x}+M_{osc,c}+M_{osc.o}$, where 
\begin{align*}
   M_{osc,t}&:=\sum_{n\geq3}\sum_{\xi\in\Lambda^n}\mathcal{R}_1\(\partial_t[\chi(\zeta|M_l|-n)(\frac{n}{\zeta})^{1/2}\Gamma_{\xi}\(\frac{M_l}{|M_l|}\)g_{(\xi)}]\Theta_{(\xi,g)}\),\\
   M_{osc,x}&:=\sum_{n\geq3}\sum_{\xi\in\Lambda^n}\mathcal{B}_1\(\nabla[\chi(\zeta|M_l|-n)\frac{n}{\zeta}\Gamma_{\xi}\(\frac{M_l}{|M_l|}\)],\mathbb{P}_{\neq0}(W_{(\xi,g)}\Theta_{(\xi,g)})\)g_{(\xi)}^2,\\
   M_{osc,c}&:=\sum_{n\geq3}\chi(\zeta|M_l|-n)\frac{n}{\zeta}\frac{M_l}{|M_l|}- M_l,\\
   M_{osc,o}&:=-\sigma^{-1}\sum_{n\geq3}\sum_{\xi\in\Lambda^n}h_{(\xi)}\partial_t\(\chi(\zeta|M_l|-n)\frac{n}{\zeta}\Gamma_{\xi}\(\frac{M_l}{|M_l|}\)\xi\).
\end{align*}
 
Then we have
$-{M}_{q+1}:=M_{osc}+M_{lin}.$ 
It is easy to see that ${M}_{q+1}(t)=0$ on $[0,T_{q+1}]$, since 
$M_l(t)= w_{q+1}(t)=\overline w_{q+1}(t)=0$ on $[0,T_{q+1}]$.

\subsubsection{Estimates of the $M_{q+1}$}\label{sec:estmq+1}
To conclude the proof of Proposition \ref{prop:case2} we verify the last bound in \eqref{bd:2rql1} for $M_{q+1}$ by estimating each term in the definition of $M_{q+1}$ separately, as done previously.

For the oscillation error, we consider $M_{osc,t}$ first. By the improved H\"older inequality in  Lemma \ref{ihiot}, the estimates for  amplitude functions in Proposition \ref{lem:2chi}, and the estimates for the building blocks in \eqref{int4theta}, \eqref{bd:gwnp} we obtain
\begin{align}
    \|M_{osc,t}\|_{L^1_tL^1}
    &\lesssim \sum_{n\geq3}\sum_{\xi\in\Lambda^n}\norm{\chi(\zeta|M_l|-n)(\frac{n}{\zeta})^{1/2}\Gamma_{\xi}\(\frac{M_l}{|M_l|}\)}_{C_{t,x}^1}\|\Theta_{(\xi,g)}\|_{C_tL^1}\|g_{(\xi)}\|_{W_t^{1,1}}\notag\\
    &\lesssim l^{-4d-12}r_\perp^{\frac{d-1}2}r_\parallel^{\frac12}\eta^{-\frac12}\sigma\lesssim \lambda_{q+1}^{(8d+24)\alpha-\frac 1{2N}}\lesssim \lambda_{q+1}^{-\alpha},\notag
\end{align}
where we used the choice of parameters in 
 \eqref{para22}, \eqref{paralam}, \eqref{def:2musigma}, and   conditions on the parameters to have $(8d+25)\alpha<\frac 1{2N}$.

For the second term  $M_{osc,x}$, we obverse that $W_{(\xi)}\Theta_{(\xi)}$ is $(\mathbb{T}/r_\perp\lambda_{q+1})^d$-periodic. So, together with the bounds for the amplitude functions in Proposition \ref{lem:2chi}, and the estimates for the building blocks in \eqref{int4}, \eqref{int4theta}, \eqref{bd:gwnp}, we apply  Theorem \ref{bb1_1} to obtain
\begin{align}
\|M_{osc,x} \|_{L^1_tL^1}&\lesssim\sum_{n\geq3}\sum_{\xi\in\Lambda^n}\norm{\mathcal{B}_1\(\nabla[\chi(\zeta|M_l|-n)\frac{n}{\zeta}\Gamma_{\xi}\(\frac{M_l}{|M_l|}\)],\mathbb{P}_{\neq0}(W_{(\xi,g)}\Theta_{(\xi,g)})\)}_{C_tL^1}\norm{g_{(\xi)}^2}_{L_t^1}\notag\\
&\lesssim\sum_{n\geq3}\sum_{\xi\in\Lambda^n}\norm{\chi(\zeta|M_l|-n)\frac{n}{\zeta}\Gamma_{\xi}\(\frac{M_l}{|M_l|}\)}_{C_{t,x}^2}\|W_{(\xi,g)}\Theta_{(\xi,g)}\|_{C_tL^{1}}(r_\perp\lambda_{q+1})^{-1}\notag\\
&\lesssim l^{-6d-20}(r_\perp\lambda_{q+1})^{-1}\|\Theta_{(\xi,g)}\|_{C_tL^2}\|W_{(\xi,g)}\|_{C_tL^{2}}\lesssim \lambda_{q+1}^{(12d+40)\alpha-\frac{1}{N}}
\lesssim \lambda_{q+1}^{-\alpha},\notag
\end{align}
where we used \eqref{para22}, \eqref{paralam} and  conditions on the parameters to have $(12d+41)\alpha<\frac{1}{N}$.

For  the third term  $M_{osc,c}$, by a similar argument as in \cite[(31)]{BCDL21}  we have
 \begin{align}
    \left| M_{osc,c}\right|&\leq \left|  \sum_{n=-1}^{2}\chi(\zeta|M_l|-n)M_l\right|+\left|\sum_{n\geq3}\chi(\zeta|M_l|-n)(\frac{n}{\zeta}\frac{M_l}{|M_l|}- M_l)\right|\notag\\
     &\leq \frac{3}{\zeta}+\sum_{n\geq3}\chi(\zeta|M_l|-n)\left|\frac{n}{\zeta}- |M_l|\right|\leq \frac{3}{20}\delta^2_{q+3}+\frac{1}{20}\delta^2_{q+3}
     \leq \frac15C_0^2\delta^2_{q+3}.\notag
\end{align}
 By the bounds on $h_{(\xi)}$ in \eqref{bd:gwnp}  and the bounds on the amplitude functions in Proposition \ref{lem:2chi}  we have
\begin{align}
   \|  M_{osc,o}\|_{L^1_tL^1}&\lesssim \sigma^{-1}\sum_{n\geq3}\sum_{\xi\in\Lambda^n}\|h_{(\xi)}\|_{L_t^\infty}\|\chi(\zeta|M_l|-n)\frac{n}{\zeta}\Gamma_{\xi}\xi\|_{C_{t,x}^1}\notag\\ 
   &\lesssim \sigma^{-1}l^{-4d-12}\lesssim \lambda_{q+1}^{(8d+24)\alpha-\frac1{2N}}\lesssim \lambda_{q+1}^{-\alpha},\notag
\end{align}
where we used  conditions on the parameters to have $(8d+25)\alpha<\frac{1}{2N}$.

 Now we turn to the linear error $M_{lin}$, where we note that all terms have already been estimated in the previous sections. More precisely, by the estimates on $\theta_{q+1}^{(p)}$ and $\theta_{q+1}^{(o)}$ in \eqref{bd:2theq+1w1s}  and \eqref{bd:2theq+1olp}   respectively, we have for $0\leq \kappa\leq1$
\begin{align}
    \|\kappa\mathcal{R}_1\Delta \theta_{q+1}\|_{L^1_tL^1}\lesssim \|\theta_{q+1}^{(p)}\|_{L^1_tW^{1,s}}+\|\theta_{q+1}^{(o)}\|_{L_t^1W^{1,\infty}}\lesssim \lambda_{q+1}^{-\alpha}.\notag
\end{align}
 Moreover, according to  the estimates of 
$v_l,\rho_l$ in \eqref{bd:2vqc1},  the estimates of  $\theta_{q+1}$ in \eqref{bd:2theq+1plp}-\eqref{bd:2theq+1olp},  the estimates of  $w_{q+1}$ in \eqref{bd:wq+1ctl1}, and   the estimates of  $\overline w_{q+1}$ in  \eqref{bd:barwq+1ctl1} we obtain
\begin{align}
\|v_l\theta_{q+1}&+(w_{q+1}+\overline w_{q+1})(\rho_l+\theta_{q+1}^{(o)})+(w_{q+1}^{(c)}+ \chi_q^2\overline w_{q+1}^{(o)})\theta_{q+1}^{(p)}\|_{L^1_tL^1}\notag\\ 
&\leq \|v_l\|_{C_{t,x}^0}\|\theta_{q+1}\|_{C_tL^1}+(\|w_{q+1}\|_{C_tL^1}+\|\overline w_{q+1}\|_{C_tL^1})(\|\rho_l\|_{C_{t,x}^0}+\|\theta_{q+1}^{(o)}\|_{C_{t,x}^0})\notag\\ 
&\quad+ (\|w_{q+1}^{(c)}\|_{L^2_tL^2}+\|\overline w_{q+1}^{(o)}\|_{C_{t,x}^0})\| \theta_{q+1}^{(p)}\|_{L^2_tL^2}\notag\\
&\lesssim C_\rho 
C_0(\lambda_{q }^{4d+3}+1)l^{-6d-20}(r_\perp^{\frac{d-1}2}r_\parallel^{\frac12}\eta^{-\frac12}\lambda_{q+1}^{d\epsilon}+\sigma^{-1})\lesssim C_\rho C_0\lambda_{q+1}^{(12d+42)\alpha-\frac1{2N}}\lesssim C_0\lambda_{q+1}^{-\alpha},\notag
\end{align}
where we choose $\epsilon>0$ small enough such that $d\epsilon<\alpha$. We also used the bounds in \eqref{para22}, \eqref{paralam}, and   conditions on the parameters to have $(12d+43)\alpha<\frac{1}{2N}$. In the last inequality we choose $a$ large enough to absorb the universal constant.

Summarizing all the bounds above we obtain \eqref{bd:2rql1} for $M_{q+1}$:
\begin{align*}
    \|M_{q+1}\|_{L^1_tL^1}\leq \frac15C_0^2\delta^2_{q+3}+CC_0^2\lambda_{q+1}^{-\alpha}\leq C_0^2\delta^2_{q+3},\end{align*}
where we choose $a$ large to absorb the universal constant.

\subsubsection{Construction of the Reynolds stress $\mathring{R}_{q+1 }$}\label{sec:defrq+12}
From \eqref{eq:2qth},  \eqref{eq:2v_l} and the definition of  the perturbations $w_{q+1},\overline w_{q+1}$ we obtain

\begin{align}
&\ \ \ \div\mathring{R}_{q +1}-\nabla\pi_{q +1}+\nabla {\pi}_{ {l}}\notag\\
&=\partial_t(\chi_{q}\overline w_{q +1}^{(p)}+\chi_{q}\overline w_{q +1}^{(c)}+w_{q+1})-\nu\Delta (w_{q +1}+ \overline w_{q +1})+\div\( {v}_{ {l}}\otimes (w_{q+1 }+\overline w_{q+1})+ (w_{q+1 }+\overline w_{q+1})\otimes  {v}_{ {l}}\)\notag\\
&\quad+\div\((\chi_{q}\overline w_{q+1 }^{(c)}+\chi_{q}^2\overline w_{q +1}^{(t)}+\chi_{q}^2\overline w_{q+1 }^{(o)}+w_{q+1}^{(c)})\otimes (w_{q +1}+\overline w_{q +1})\)\notag\\ 
&\quad+\div\((\chi_{q}\overline w_{q+1 }^{(p)}+w_{q+1 }^{(p)})\otimes (\chi_{q}\overline w_{q+1 }^{(c)}+\chi_{q}^2\overline w_{q+1 }^{(t)}+\chi_{q}^2\overline w_{q +1}^{(o)}+w_{q+1}^{(c)})\)\notag\\
&\quad+\partial_t (\chi_{q}^2\overline w_{q+1 }^{(t)}+\chi_{q}^2\overline w_{q+1 }^{(o)})+\div\((\chi_{q}\overline w_{q+1 }^{(p)}+w_{q+1 }^{(p)})\otimes (\chi_{q}\overline w_{q+1 }^{(p)}+w_{q+1 }^{(p)})+\mathring{ {R}}_{ {l}}\),\notag
\end{align}
where by using the inverse divergence operator $\mathcal{R}$ introduced in Section \ref{tamr} we define
\begin{align*}
 {R}_{lin}:&=\mathcal{R}\partial_t(\chi_{q}\overline w_{q +1}^{(p)}+\chi_{q}\overline w_{q +1}^{(c)}+w_{q+1})-\nu\mathcal{R}\Delta (w_{q +1}+ \overline w_{q +1})\\
 &\quad+{v}_{ {l}}\mathring{\otimes}  (w_{q+1 }+\overline w_{q+1})+ (w_{q+1 }+\overline w_{q+1})\mathring{\otimes}  {v}_{ {l}}.  \\
 {R}_{cor}:&=(\chi_{q}\overline w_{q+1 }^{(c)}+\chi_{q}^2\overline w_{q +1}^{(t)}+\chi_{q}^2\overline w_{q+1 }^{(o)}+w_{q+1}^{(c)})\mathring{\otimes} (w_{q +1}+\overline w_{q +1})\notag\\ 
 &\quad+(\chi_{q}\overline w_{q+1 }^{(p)}+w_{q+1 }^{(p)})\mathring{\otimes} (\chi_{q}\overline w_{q+1 }^{(c)}+\chi_{q}^2\overline w_{q+1 }^{(t)}+\chi_{q}^2\overline w_{q +1}^{(o)}+w_{q+1}^{(c)}).
\end{align*}

In order to define the remaining oscillation error in the last line,  first by the definition of $\chi_q$ and the fact that  $\mathring{ {R}}_{ {l}}=w_{q+1 }=w^{(p)}_{q+1 }=0 $ for $t\in[0,\frac{T_q+T_{q+1}}{2}]$, we know that  $\mathring{ {R}}_{ {l}}=\chi_{q}\mathring{ {R}}_{ {l}},w^{(p)}_{q+1 }=\chi_{q}w^{(p)}_{q+1 }$ and $w_{q+1 }=\chi_{q}w_{q+1 }$. Then   we apply the identities  \eqref{2.12}, \eqref{int1}-\eqref{eq:partbarwq+1o}  to obtain
\begin{align*}
\partial_t& (\chi_{q}^2\overline w_{q+1 }^{(t)}+\chi_{q}^2\overline w_{q+1 }^{(o)})+\div\((\chi_{q}\overline w_{q+1 }^{(p)}+w_{q+1 }^{(p)})\otimes (\chi_{q}\overline w_{q+1 }^{(p)}+w_{q+1 }^{(p)})+\mathring{ {R}}_{ {l}}\)\\
&= \chi_{q}^2\partial_t\overline w_{q+1 }^{(t)}+\chi_{q}^2\sum_{{\xi}\in\overline{\Lambda}}\div\(a_{({\xi})}^2\mathbb{P}_{\neq0}(\overline{W}_{({\xi,g})}\otimes \overline{W}_{({\xi,g})})g_{(\xi)}^2\)\\
&\quad+\chi_{q}^2\partial_t\overline w_{q+1 }^{(o)}+\chi_{q}^2\sum_{\xi\in\overline{\Lambda}}\div\(a_{({\xi})}^2(g_{(\xi)}^2-1){\xi}\otimes{\xi}\)\notag\\ 
&\quad+\chi_{q}^2\nabla A+\div(w_{q+1}\otimes w_{q+1})+(\chi_{q}^2)'(\overline w_{q+1 }^{(t)}+\overline w_{q+1 }^{(o)})
\notag\\
&=\sum_{{ {\xi}}\in\overline  {\Lambda}}\chi_{q}^2 \mathbb{P}_{\neq0}\left(\nabla a_{({ {\xi}})}^2 {g}_{(\xi)}^2\mathbb{P}_{\neq0} ( \overline {W}_{({ {\xi,g}})}\otimes  \overline {W}_{({ {\xi,g}})})\right)+\frac{1}{\overline {\mu}}\chi_{q}^2\sum_{{ {\xi}}\in\overline  {\Lambda}}\mathbb{P}_{\neq0}\left(\partial_t(a_{({ {\xi}})}^2 {g}_{(\xi)})\overline  {\phi}_{({ {\xi,g}})}^2 \overline {\psi}_{({ {\xi,g}})}^2{ {\xi}}\right)\\
&\quad-\sigma^{-1}\chi_{q}^2 \mathbb{P}_{\rm H}\mathbb{P}_{\neq0}\sum_{\xi\in\overline \Lambda}h_{(\xi)}\partial_t\div(a_{(\xi)}^2\xi\otimes\xi) +\partial_t(\chi_{q}^2)(\overline w_{q+1 }^{(t)}+\overline w_{q+1 }^{(o)})+ \nabla {p}_1+\div(w_{q+1}\otimes w_{q+1}).\notag 
\end{align*}
 Here $p_1$ denotes the pressure term.

Therefore using the inverse divergence operators  
  $\mathcal{R}$, $\mathcal{B}$  introduced in Section \ref{tamr} we have
  
\begin{align*}
 {R}_{osc}&:=\sum_{{ {\xi}}\in\overline {\Lambda}} {g}_{(\xi)}^2\mathcal{B}\left(\nabla a_{({ {\xi}})}^2\chi_{q}^2,\mathbb{P}_{\neq0}(\overline {W}_{({ {\xi,g}})}\otimes \overline {W}_{({ {\xi,g}})})\right)(:= {R}_{osc,x})\\
&\ \ \ \ +\frac{1}{ \overline {\mu}}\sum_{{ {\xi}}\in \overline{\Lambda}}\mathcal{R}\left(\partial_t(a_{({ {\xi}})}^2 {g}_{(\xi)})\chi_{q}^2 \overline{\phi}_{({ {\xi,g}})}^2 \overline{\psi}_{({ {\xi,g}})}^2{ {\xi}}\right)(:= {R}_{osc,t})\\
&\ \ \ \  +(\chi_{q}^2)'\mathcal{R}(\overline w_{q+1 }^{(t)}+\overline w_{q+1 }^{(o)}) -\sigma^{-1}\chi_{q}^2\mathbb{P}_{\neq0}\sum_{\xi\in\overline \Lambda}h_{(\xi)}\partial_t(a_{(\xi)}^2\xi\mathring{\otimes}\xi)(:= {R}_{osc,o})\notag\\
&\ \ \ \ +w_{q+1}\mathring{\otimes} w_{q+1}.
\end{align*}

 Then  the Reynolds stress at the level $q+1 $ is given by  $\mathring{R}_{q +1}= {R}_{lin}+ {R}_{cor}+ {R}_{osc}.$
It is easy to see that $\mathring{R}_{q +1}$ is a trace-free and symmetric matrix satisfying $\mathring{R}_{q +1}=0$ on $[0,T_{q+1}]$.

\subsubsection{Estimate of $\mathring{R}_{q +1}$}\label{sec:estrq+12}
We estimate each term in the definition of $\mathring{R}_{q+1 }$ separately.

For the linear error $ {R}_{lin}$, by the estimates for the building blocks in  \eqref{int4ns}, \eqref{bd:gwnp},  the estimate for $a_{(\xi)}$ in \eqref{bd:2axi}  and the identity \eqref{wq+12p+wq+12c} we obtain  for $\epsilon>0$ small enough
\begin{align}
\|\mathcal{R}&\partial_t(\chi_{q}\overline w_{q +1}^{(p)}+\chi_{q}\overline w_{q +1}^{(c)}) \|_{L^1_tL^1}\lesssim\sum_{{ {\xi}}\in \overline{\Lambda}} \|\chi_q\|_{C_t^1}\| a_{({ {\xi}})}\|_{C_{t,x}^1}(\|\partial_t( {g}_{(\xi)} {V}_{({ {\xi,g}})})\|_{L^1_tL^{1+\epsilon}}+\| {g}_{(\xi)} {V}_{({ {\xi,g}})}\|_{L^1_tL^{1+\epsilon}})\notag\\ 
&\lesssim  {l}^{-3d-8} {r}_\perp^{\frac{d-1}{2}}  {r}_\parallel^{
\frac12}(\frac{ {r}_\perp \overline{\mu}}{ {r}_\parallel}+\frac{\sigma\eta^{-\frac12}}{\lambda_{q+1}}) r_{\perp}^{-d\epsilon}\lesssim  {l}^{-3d-8}(\frac{ {r}_\perp}{ {r}_\parallel} {\lambda}_{q+1}^{\frac1{2N}}+ {r}_\perp^{\frac{d-1}{2}}  {r}_\parallel^{\frac12}\eta^{-\frac12})\lambda_{q+1}^{d\epsilon}\notag\\
&\lesssim
{\lambda}_{q +1}^{(6d+17)\alpha-\frac{1}{2N}}\lesssim\lambda_{q+1 }^{-\alpha},
\end{align}
where we  note that $\mathcal{R}\div $ is not $L^p$ bounded for $p>1$, and choose $\epsilon>0$ small enough such that $d \epsilon<\alpha$. We also used  the choice of parameters in \eqref{para22}, \eqref{paralam}, \eqref{def:2musigma} and used the  conditions on the parameters to have $(6d+18)\alpha<\frac{1}{2N}$.

The definition of $w_{q+1}$  is analogous to that of $\overline w_{q+1}$, allowing us to bound it in a similar manner.  By the estimates on the building blocks in   \eqref{int4}, \eqref{bd:gwnp},  and the estimate on $\tilde{\chi}$ in Proposition \ref{lem:2chi} we obtain  for $\epsilon>0$ small enough

\begin{align}
\|\mathcal{R}\partial_tw_{q+1}\|_{L^1_tL^1}
&\lesssim\sum_{n\geq3}\sum_{\xi\in\Lambda^n}\|\tilde{\chi}(\zeta|M_l|-n)(\frac{n}{\zeta})^{\frac12}\|_{C_{t,x}^1}\(\|\partial_t(g_{(\xi)}V_{(\xi,g)})\|_{L_t^1L^{1+\epsilon}}+\|g_{(\xi)}V_{(\xi,g)}\|_{L_t^1L^{1+\epsilon}}\)\notag\\
&\lesssim l^{-3d-8}(\frac{r_\perp}{r_\parallel}+r_\perp^{\frac{d-1}{2}} r_\parallel^{\frac12}\eta^{-\frac12})\lambda_{q+1 }^{d\epsilon} 
\lesssim\lambda_{q+1}^{(6d+17)\alpha-\frac{1}{N}}
\lesssim \lambda_{q+1}^{-\alpha},\notag
\end{align}
 where we choose $\epsilon>0$ small enough such  that $d\epsilon<\alpha$ and  used  conditions on the parameters to have $(6d+18)\alpha<\frac{1}{N}$.

Using \eqref{para22}, the bounds in  \eqref{bd:barwq+1w1s},  \eqref{bd:2wq+1w1p} and the boundedness property of the inverse divergence operator in Theorem \ref{bb} we have for $0\leq\nu\leq1$
 \begin{align*}
   \nu\|\mathcal{R}\Delta (w_{q +1}+ \overline w_{q +1})\|_{L^1_tL^1}\lesssim\|w_{q+1 }\|_{L^1_tW^{1,s}}+\|\overline w_{q+1 }\|_{L^1_tW^{1,s}}\lesssim \lambda_{q +1}^{-\alpha}.
 \end{align*}

By the estimates of $v_q$ in \eqref{bd:2vqc1}, and the estimates of the components of $w_{q+1 },\overline w_{q+1}$ in \eqref{bd:2wq+12plp}-\eqref{bd:2wq+12olp}, \eqref{bd:2wq+1plp} and \eqref{bd:2wq+1clp} we deduce
\begin{align} 
\|{v}_{ {l}}\otimes &(w_{q+1 }+\overline w_{q+1})+ (w_{q+1 }+\overline w_{q+1})\otimes  {v}_{ {l}}\|_{L^1_tL^1}\lesssim\| {v}_{ {l}}\|_{C_{t,x}^0}(\|\overline w_{q+1}\|_{C_tL^{1+\epsilon}}+\| w_{q+1}\|_{C_tL^{1+\epsilon}})\notag\\
&\lesssim C_0\lambda_{q}^{4d+3} {l}^{-6d-15}( {r}_\perp^{\frac{d-1}{1+\epsilon}-\frac{d-1}2} {r}_\parallel^{\frac{1}{1+\epsilon}-\frac12} {\eta}^{-\frac12}+ {\sigma}^{-1})\lesssim C_0\lambda_{q+1 }^{(12d+30)\alpha-\frac{1}{2N}+\epsilon d}\lesssim C_0\lambda_{q+1 }^{-\alpha},\notag
\end{align}
where we used  \eqref{para22}, \eqref{paralam} and choose $\epsilon>0$ small enough such that $ d \epsilon<\alpha$. We also used the  conditions on the parameters to have $(12d+32)\alpha<\frac{1}{2N}$.

The corrector error is estimated using the choice of parameters  \eqref{para22}, \eqref{paralam} and the bounds on perturbations  in \eqref{bd:2wq+12pl2}-\eqref{bd:2wq+12olp}, \eqref{bd:2wq+1pl2}-\eqref{bd:2wq+1clp} as
\begin{align}
\| {R}_{cor} \|_{L^1_tL^1}&\leq\|\chi_{q}\overline  w_{q+1 }^{(c)} +\chi_{q}^2\overline w_{q +1}^{(t)}+\chi_{q}^2\overline w_{q+1 }^{(o)}+w_{q+1}^{(c)} \|_{L^2_tL^{2}}\notag\\ 
&\quad\quad\times (\|w_{q +1} \|_{L^2_tL^{2}}+\|\overline w_{q +1} \|_{L^2_tL^{2}}
+\|w_{q +1}^{(p)} \|_{L^2_tL^{2}}+\|\overline w_{q +1}^{(p)} \|_{L^2_tL^{2}})\notag\\
&\lesssim (l^{-4d-6}\lambda_{q+1 }^{-\frac1{2N}}+\lambda_{q+1 }^{-\alpha})
C_v C_0\lesssim C_0 \lambda_{q+1 }^{-\alpha}
,\notag
\end{align}
where we used  conditions on the parameters to have $(8d+13)\alpha<\frac{1}{2N}$. We choose $a$ large enough to absorb the constant  $C_v$ and the implicit constant.

Now we consider the oscillation term.
In order to bound the first term  $ {R}_{osc,x}$, by the boundedness property of inverse divergence operator in Theorems \ref{bb} and \ref{bb1}, the estimates for the building blocks in  \eqref{int4ns}, \eqref{bd:gwnp}  and   the estimates for $a_{(\xi)}$ in \eqref{bd:2axi}  we have 
\begin{align}
\| {R}_{osc,x} \|_{L^1_tL^1}
&\lesssim \sum_{{ {\xi}}\in\overline  {\Lambda}}\|a_{({ {\xi}})}^2\|_{C_{t,x}^2}\|\mathcal{R}( \overline {W}_{({ {\xi,g}})}\otimes  \overline {W}_{({ {\xi,g}})})\|_{C_tL^{1+\epsilon}}\| {g}_{(\xi)}^2\|_{L_t^1}
\notag\\
&\lesssim  {l}^{-6d-14}( {r}_\perp\lambda_{q +1})^{-1} {r}_\perp^{(d-1)(\frac{1}{1+\epsilon}-1)} {r}_\parallel^{\frac{1}{1+\epsilon}-1}\lesssim \lambda_{q+1 }^{(12d+28)\alpha-\frac{1}{N}+d\epsilon}\lesssim\lambda_{q +1}^{-\alpha},\notag
\end{align}
where we  choose $\epsilon>0$ small enough such that $d\epsilon<\alpha$. We also used \eqref{para22}, \eqref{paralam}, \eqref{def:2musigma} and  conditions on the parameters to have $(12d+30)\alpha<\frac{1}{N}$. 

For the second term ${R}_{osc,t}$ we use the estimates for the building blocks $g_{(\xi)},\overline {\phi}_{({ {\xi}})}, \overline {\psi}_{({ {\xi}})}$ and  for $a_{(\xi)}$ in \eqref{bd:gwnp},  \eqref{int2ns}, \eqref{int3ns}, 
and \eqref{bd:2axi} respectively to deduce
\begin{align}
\|  {R}_{osc,t}\|_{L^1_tL^1}&\lesssim \frac{1}{ \overline {\mu}}\sum_{{ {\xi}}\in\overline  {\Lambda}}\|a_{({ {\xi}})}^2\|_{C_{t,x}^1}\| {g}_{(\xi)}\|_{W^{1,1}}\| \overline {\phi}_{({ {\xi,g}})}^2  \overline {\psi}_{({ {\xi,g}})}^2 {\xi}\|_{C_tL^{1+\epsilon}}\notag\\ 
&\lesssim   {l}^{-5d-10} \overline {\mu}^{-1} {r}_\perp^{(d-1)(\frac{1}{1+\epsilon}-1)} {r}_\parallel^{\frac{1}{1+\epsilon}-1}\sigma\eta^{-\frac12}\notag\\ 
&\lesssim  {l}^{-5d-10}\lambda_{q +1}^{d\epsilon} {r}_\perp^{\frac{d-1}{2}} {r}_\parallel^{\frac12}\eta^{-\frac12}\sigma\lesssim\lambda_{q+1 }^{(10d+21)\alpha-\frac1{2N}}\lesssim\lambda_{q+1 }^{-\alpha},\notag
\end{align}
where we  choose $\epsilon>0$ small enough such that $d\epsilon<\alpha$ and used \eqref{para22}, \eqref{paralam} and  condition on the parameters to have $(10d+22)\alpha<\frac1{2N}$.

We continue with $ {R}_{osc,o}$. By  the bound on the cut-off function in \eqref{bd:2chiqcn}, the bounds on temporal jets $g_{(\xi)},h_{(\xi)}$  in \eqref{bd:gwnp}, the bounds on perturbations  $\overline w_{q+1 }^{(t)},\overline w_{q+1 }^{(o)}$  in  \eqref{bd:2wq+12tlp}, \eqref{bd:2wq+12olp} and the bounds on $a_{(\xi)}$ in \eqref{bd:2axi} we have
\begin{align}
   \|  {R}_{osc,o}\|_{L^1_tL^1}&\lesssim\|\chi_{q}^2\|_{C_t^1}(\|\overline w_{q+1 }^{(t)}\|_{C_tL^{1+\epsilon}}+\|\overline w_{q +1}^{(o)}\|_{C_tL^{1+\epsilon}})+\sigma^{-1}\sum_{\xi\in\overline \Lambda}\|h_{(\xi)}\|_{L_t^\infty}\|a_{(\xi)}^2\|_{C_{t,x}^1}\notag\\ 
   &\lesssim  {l}^{-4d-7} {r}_\perp^{\frac{d-1}{1+\epsilon}-\frac{d-1}2} {r}_\parallel^{\frac{1}{1+\epsilon}-\frac12} {\eta}^{-\frac12}+{l}^{-6d-16} \sigma^{-1}\lesssim C_0^2\lambda_{q +1}^{(12d+32)\alpha-\frac1{2N}}\lesssim C_0^2\lambda_{q +1}^{-\alpha},\notag
\end{align}
where we   choose $\epsilon>0$ small enough such that $d\epsilon<\alpha$. We used \eqref{para22}, \eqref{paralam}, \eqref{def:2musigma} and   conditions on the parameters to have $(12d+33)\alpha<\frac{1}{2N}$. We choose $a$ large enough to absorb the universal constant.

By \eqref{bd:2wq+1l2} we have that 
\begin{align*}
 \|w_{q+1}\mathring{\otimes} w_{q+1}\|_{L^1_tL^1}\leq    \frac{C_0^2}2
 \delta_{q+2}.
\end{align*}

Summarizing all the estimates above we obtain the first term in \eqref{bd:2rql1}:
$$\|\mathring{R}_{q+1} \|_{L^1_tL^1}\leq  \frac{C_0^2}2
 \delta_{q+2}+C C_0^2 \lambda_{q +1}^{-\alpha}\leq C_0^2\delta_{q+2}.$$
 Here we used the condition  \eqref{para22} to deduce the last inequality.

Thus we finish the proof of Proposition \ref{prop:case2}.

\noindent{\bf Acknowledgement.}  We are very grateful to Jan Burczak, L\'aszl\'o Sz\'ekelyhidi Jr. and Bian Wu for directing our attention to the  seminal works \cite{BSW23} and \cite{Row24}. These papers establish, through distinct approaches, the non-uniqueness of solutions to  transport equations and their associated Lagrangian trajectories.

\appendix
 \renewcommand{\appendixname}{Appendix~\Alph{section}}
  \renewcommand{\theequation}{A.\arabic{equation}}
\section{Some technical tools}\label{app:b}
We collect some technical tools used in the construction of convex integration schemes.
  \subsection{Inverse divergence operators}\label{tamr}
We first recall the following inverse divergence operator $\mathcal{R}$ as in \cite[Appendix B.2]{CL22}, which acts on vector fields $v$ with $\int_{\mathbb{T}^d} v\dif x = 0$ as
$$(\mathcal{R}v)_{ij}=\mathcal{R}_{ijk}v_k,$$
where

$$\mathcal{R}_{i j k}=\frac{2-d}{d-1} \Delta^{-2} \partial_{i} \partial_{j} \partial_{k}-\frac{1}{d-1} \Delta^{-1} \partial_{k} \delta_{i j}+\Delta^{-1} \partial_{i} \delta_{j k}+\Delta^{-1} \partial_{j} \delta_{i k}.$$
Then $\mathcal{R}v(x)$ is a symmetric trace-free matrix for each $x \in \mathbb{T}^d$, and $\mathcal{R}$ is a right
inverse of the $\div$ operator, i.e. $\div(\mathcal{R}v) = v$.
 Here we use the notation  $\mathcal{R}v:=\mathcal{R}(v-\int v\dif x)$ for a general vector field. In the following we use $C^\infty(\mathbb{T}^d;\mathbb{R}^d)$ to denote the space of smooth functions from $\mathbb{T}^d$ to $\mathbb{R}^d$, and we use $C_0^\infty(\mathbb{T}^d;\mathbb{R}^d)$ to denote the subspace of functions with zero spatial mean. Similarly, we define $C^\infty(\mathbb{T}^d;\mathbb{R}^{d\times d})$ and $C_0^\infty(\mathbb{T}^d;\mathbb{R}^{d\times d})$.   By $\mathcal{S}_0^{d\times d}$ we denote the space of symmetric trace-free matrices.
\bt$($\cite[Theorem B.3]{CL22}$)$\label{bb}
 Let $1\leq p \leq\infty$. For any vector field 
$f \in C_0^\infty(\mathbb{T}^d; \mathbb{R}^d)$, $\sigma\in\mathbb{N}$,
$$\|\mathcal{R}f(\sigma\cdot)\|_{ L^p}\lesssim\sigma^{-1} \| f\|_{ L^p}.$$
\et

We also introduce the bilinear version $\mathcal{B}: C^\infty(\mathbb{T}^d; \mathbb{R}^d) \times C_0^\infty(\mathbb{T}^d; \mathbb{R}^{d\times d})\to C^\infty(\mathbb{T}^d; \mathcal{S}_0^{d\times d})$  by $$(\mathcal{B}(v,A))_{ij}=v_m\mathcal{R}_{ijk}A_{mk}-\mathcal{R}(\partial_iv_m\mathcal{R}_{ijk}A_{mk}).$$

\bt$($\cite[Theorem B.4]{CL22}$)$\label{bb1}
Let $1 \leq p \leq\infty$. For any $v \in C^\infty(\mathbb{T}^d;\mathbb{R}^d)$ and $A\in C_0^\infty(\mathbb{T}^d; \mathbb{R}^{d\times d})$, we have $\div(\mathcal{B}(v,A))=vA-\int_{\mathbb{T}^d}vA\dif x,$ and
$$\|\mathcal{B}(v,A)\|_{L^p}\lesssim\|v\|_{C^1}\|\mathcal{R}A\|_{L^p}.$$
\et

We also need to define the inverse divergence operator acting on  scalars. We define $\mathcal{R}_1:=\nabla\Delta^{-1}$ as a right
inverse of the $\div$ operator, i.e. $\div(\mathcal{R}_1v) = v$ for scalars $v$ with $\int_{\mathbb{T}^d} v\dif x = 0$.   Here we use the notation  $\mathcal{R}_1v:=\mathcal{R}_1(v-\int v\dif x)$ for a general scalar function $v$. Then  since $\mathcal{R}_1$ is a Calder\'on-Zygmund operator, we have 
\bt\label{bb_1}
 Let $1\leq p \leq\infty$. For any vector field 
$f \in C_0^\infty(\mathbb{T}^d; \mathbb{R})$, $\sigma\in\mathbb{N}$,
$$\|\mathcal{R}_1f(\sigma\cdot)\|_{ L^p}\lesssim\sigma^{-1} \| f\|_{ L^p}.$$
\et
We introduce the the bilinear version $\mathcal{B}_1: C^\infty(\mathbb{T}^d; \mathbb{R}) \times C_0^\infty(\mathbb{T}^d; \mathbb{R})\to C^\infty(\mathbb{T}^d; \mathbb{R}^{d})$  by $$\mathcal{B}_1(v,f)=v\mathcal{R}_1f-\mathcal{R}_1\(\nabla v\cdot\mathcal{R}_{1}f+\int vf\dif x\).$$

\bt$($\cite[Lemma 3.3]{BCDL21}$)$\label{bb1_1}
Let $1 \leq p \leq\infty$. For any $v \in C^\infty(\mathbb{T}^d; \mathbb{R})$ and $f\in C_0^\infty(\mathbb{T}^d; \mathbb{R})$, we have $\div(\mathcal{B}_1(v,f))=vf-\int_{\mathbb{T}^d}vf\dif x,$ and  for $\sigma\in\mathbb{N}$,
$$\|\mathcal{B}_1(v,f(\sigma\cdot))\|_{W^{k,p}}\lesssim\sigma^{k-1}\|v\|_{C^{k+1}}\|f\|_{W^{k,p}}.$$
\et

 We recall the following estimates  stationary phase bounds from \cite[Lemma 2.2]{DSJ17}, which are useful for estimating the oscillation errors.
 \bp\label{prop:phase bounds}
		Let $ \lambda\xi\in\mZ^d, N \geq 1$, $a \in C^\infty(\mathbb{T}^d;\mR),b \in C^\infty(\mathbb{T}^d;\mR^d)$, $\Phi \in C^\infty(\mathbb{T}^d;\mathbb{R}^d)$ be smooth functions and assume that there exists a constant $\hat C\geq1$ such that
		$ \hat C^{-1} \leq |\nabla \Phi| \leq \hat C$
		holds on $\mathbb{T}^d$. Then
		\begin{align}\label{esti:integral}
			\left| \int_{\mathbb{T}^d} a(x)e^{i\lambda \xi \cdot \Phi(x)} \dif x\right| \lesssim \frac{\|a\|_{C^N}+\|a\|_{C^0} \|\nabla \Phi\|_{C^N}}{\lambda^N},
		\end{align}
		and for the operators $\mathcal{R}$ and $\mathcal{R}_1$ defined above, we have for $\alpha \in (0,1)$
		\begin{align*}
			\left\|\mathcal{R}\left(b(x)e^{i\lambda \xi\cdot \Phi(x)}\right)\right\|_{C^\alpha}\lesssim \frac{\|b\|_{C^0}}{\lambda^{1-\alpha}}+\frac{\|b\|_{C^{N+\alpha}}+\|b\|_{C^0}\|\nabla \Phi\|_{C^{N+\alpha}}}{{\lambda^{N-\alpha}}},\\
			\left\|\mathcal{R}_1\left(a(x)e^{i\lambda \xi\cdot \Phi(x)}\right)\right\|_{C^\alpha}\lesssim \frac{\|a\|_{C^0}}{\lambda^{1-\alpha}}+\frac{\|a\|_{C^{N+\alpha}}+\|a\|_{C^0}\|\nabla \Phi\|_{C^{N+\alpha}}}{{\lambda^{N-\alpha}}},
		\end{align*}
		where the implicit constants depend on $\hat C$, $\alpha$ and $N$, but not on the frequency $\lambda$.
	\ep

\subsection{Commutator estimate}
We recall   the following  commutator estimate which can be seen as a generalization of \cite[Proposition A.2]{BDLSV19}:
\begin{lemma}\label{p:CET}
Let $f,g\in C^{\infty}({\T}\times[0,1])$ and $\psi$ a standard radial smooth and compactly supported kernel. For any $r\geq 0$ and $\theta_1,\theta_2 \in (0,1]$ we have the estimate
\begin{align*}
   \|(f*\psi_ l)( g*\psi_l)-(fg)*\psi_ l\|_{C^r} \les l^{\theta_1+\theta_2 -r}  \|f\|_{C^{\theta_1}} \|g\|_{C^{\theta_2}} \, ,\notag
\end{align*}
where the implicit constant  depends only on $r$ and $\psi$.
\end{lemma}

\subsection{Estimates for transport equations}
	Now we recall some standard estimates for solutions to the transport equation:
	\begin{align}\label{eq:phi}
		\partial_tf+v\cdot\nabla f&=g,
		\\ f(0)&=f_0,\notag
	\end{align}
    where $v$ is  a given smooth vector field. 
	 We have the following proposition.
	\begin{proposition}\label{esti:transport}\cite[Proposition B.1]{BDLSV19}
		Assume $t\|v\|_{C^1}\leq 1$. Then, any solution f of \eqref{eq:phi} satisfies
        \begin{align*}
            \|f(t)\|_{C^0}&\leq\|f_0\|_{C^0}+\int_{0}^{t}\|{g(\tau)}\|_{C^0} \dif \tau ,\\
            \|f(t)\|_{C^\alpha}&\leq e^\alpha\left( \|f_0\|_{C^\alpha}+\int_{0}^{t}\|{g(\tau)}\|_{C^\alpha} \dif \tau \right),
        \end{align*}
		for all $\alpha\in[0,1)$. More generally, for any $N\geq1$ and $\alpha\in[0,1)$
        \begin{align*}
            \|f(t)\|_{C^{N+\alpha}}\lesssim &\|f_0\|_{C^{N+\alpha}}+|t|\|{v}\|_{C^{N+\alpha}} \|f_0\|_{C^1} +\int_{0}^{t} \left( \|g(\tau)\|_{C^{N+\alpha}}+(t-\tau)\|{v}\|_{C^{N+\alpha}} \|g(\tau)\|_{C^1} \right) \dif \tau,
        \end{align*}
		where the implicit constant depends on $N$ and $\alpha$.
        
		 Consequently, the  flow $\Phi$ of $v$ starting at time $0$ (i.e. $\frac{\dif }{\dif t}\Phi=v(\Phi(t),t)$ and $\Phi(0)={\rm Id}$) satisfies
		\begin{align*}
			\| \nabla \Phi(t)-\mathrm{Id}\|_{C^0}&\lesssim 
            |t|\|v\|_{C^1} \, ,	\\ \|  \Phi(t)\|_{C^N}\lesssim &\, |t|\|v\|_{C^N}, \, N\geq 2 .
		\end{align*}
	\end{proposition}

\subsection{Improved H\"older inequality on $\mathbb{T}^d$}
This lemma improves the usual H\"older inequality by using the decorrelation between frequencies, which  establish the desired estimates.
\bl$($\cite[Theorem B.1]{CL22}$)$\label{ihiot}
Let $d\geq2,p \in [1, \infty]$ and $a, f : \mathbb{T}^d \to \mathbb{R}$ be smooth functions. Then for any $\sigma\in\mathbb{N}$,
$$| \|a f(\sigma\cdot)\|_{L^p}-\|a\|_{L^p}\|f\|_{L^p} |\lesssim\sigma^{-1/p}\|a\|_{C^1}\|f\|_{L^p}.$$
\el

\renewcommand{\appendixname}{Appendix~\Alph{section}}
  \renewcommand{\theequation}{B.\arabic{equation}}
 
  \section{Building blocks and auxiliary estimates in Section \ref{sec:proof13}}
  \label{s:appA.3}
In this section, we first  introduce the building blocks used in the convex integration method. Then we provide some auxiliary estimates in  the gluing   and perturbation steps in Section \ref{sec:proof13}.

 \subsection{Mikado flows} 
\label{sec:Mikado}
In this section we revisit the Mikado flows introduced in \cite{DSJ17} and extend their applicability to transport equations.

First we introduce the following geometric lemmas in  $\mathbb{R}^d$, where $d\geq2$. The first lemma is to show that every vector in the annulus $ \overline B_1(0)\backslash B_{\frac12}(0)$ can
 be written as a positive  combination of vectors in some  subset 
of $\mathbb{S}^{d-1}\cap \mathbb{Q}^d $, which generalizes  \cite[Lemma 3.1]{BCDL21} about the sphere $\partial B_1(0)$ to the annulus.  The proof follows from a similar argument as in \cite[Lemma 3.1]{BCDL21}.
\begin{lemma}\label{lem:geo:tran}
   \label{l:linear_algebra2}
   Let ${B}_{r}(0)$ denote the  ball of radius $r$ around 0 in $\mathbb{R}^d$.
There exists  a finite set $\Lambda\in \mathbb{S}^{d-1}\cap \mathbb{Q}^d $ and a non-negative  $C^\infty$-function  $\Gamma_{\xi}: \overline B_1(0)\backslash B_{\frac12}(0)\to\R $ such
that for every $\frac12\leq |R|\leq1$
$$R=\sum_{\xi\in\Lambda}\Gamma_\xi(R)\xi.$$
\end{lemma}
\begin{proof}
  For each vector $v\in\overline B_1(0) \backslash B_{\frac12}(0)$, we consider a collection $\Lambda(v) =\{\xi_1(v),...,\xi_d(v)\}\subset \partial B_1(0)$
 of linearly independent unit vectors in $\mQ^d$ with the property that the $d$-dimensional  open symplex $\Sigma(v)$ with vertices $0,2\xi_1(v),...,2\xi_d(v)$ contains $v$. Since $\{ \Sigma(v):
 v \in\overline B_1(0) \backslash B_{\frac12}(0) \} $ is an open cover of  $\overline B_1(0) \backslash B_{\frac12}(0)$, we consider a finite sub-cover and the corresponding collections $\Lambda_i = \Lambda(v_i),\ i=1,...,N$. We set
$$\Lambda:=\cup^N_{ j=1}\Lambda_j.$$
 For each fixed $j$, each vector $R \in \overline B_1(0) \backslash B_{\frac12}(0)$ can be written in a unique way as linear
 combination of the vectors in $\Lambda_j=\{\xi_{j,i}\}_{i=1,...,d}$. If denotes  $b_{j,i}(R)$ the corresponding
 coefficients which depend linearly on $R$, then the latter are all strictly
 positive if $R$ belongs to $\Sigma(v_j)$. Then we consider a partition of unity $\chi_j$ on $\overline B_1(0) \backslash B_{\frac12}(0)$ associated to this cover and for every $\xi_{j,i} = \xi \in \Lambda$ we set
$$ \Gamma_\xi(R) := \chi_j(R)b_{j,i}(R).$$
 The coefficients $ \Gamma_\xi$ are then smooth nonnegative functions of $R$. 
\end{proof}
Then the next  lemma shows that  symmetric matrices in $\overline{B}_{\frac12}(\mathrm{Id})$  can
 be written as a  combination of some  symmetric tensors.
\begin{lemma}\label{l:linear_algebra}
    \label{lem:3} \cite[Lemma 4.2]{CL22}
Let $\overline{B}_{\frac12}(\mathrm{Id})$ denote the closed ball of radius $\frac12$ around the identity matrix $\mathrm{Id}$, in the
space of $d\times d$ symmetric matrices. There exists a finite set $\overline{\Lambda}\in \mathbb{S}^{d-1}\cap \mathbb{Q}^{d}$ such that for each $\xi\in{\overline{\Lambda}}$ there exists
a $C^\infty$-function $\gamma_{{\xi}}$: $\overline{B}_{\frac12}(\mathrm{Id})\to \mathbb{R}$ such that for every symmetric matrix satisfying $|R-\mathrm{Id}|\leq 1/2$
$$R=\sum_{{{\xi}}\in\overline{\Lambda}}\gamma_{{\xi}}^2(R)({{\xi}}\otimes{{\xi}}).$$
\end{lemma}

Then by taking suitable rational rotations, there are four disjoint sets  $\Lambda^1, \Lambda^2,\overline\Lambda^1,\overline \Lambda^2$ such that both $\Lambda^1$ and $\Lambda^2$  enjoy the property of Lemma \ref{l:linear_algebra2}, and  similarly  both  $\overline\Lambda^1$ and $\overline\Lambda^2$  enjoy the property of Lemma \ref{l:linear_algebra}.
 For  convenience, in the following we set $\Lambda:=\Lambda^1\cup \Lambda^2\cup\overline\Lambda^1\cup\overline \Lambda^2$.

Now we apply the two lemmas in dimension 3 and obtain the set $\Lambda$. 
For each $\xi\in\Lambda$ we define $A^i_\xi\in \mathbb{S}^{2}\cap \mathbb{Q}^3,\ i=1,2$ such that $\{\xi, A^i_\xi,i=1,2\} $
form an orthonormal basis in $\mathbb{R}^3$. We label by $n_*$ the smallest
natural number such that
$$\{n_*\xi, n_*A^i_\xi,i=1,2\}\subset\mathbb{Z}^3.$$

Let ${\Phi} : \mathbb{R}^{2} \to \mathbb{R}$ be a smooth function with support in a ball of radius $\epsilon_\Lambda$, where $\epsilon_\Lambda>0$ will be chosen later in terms of $\Lambda$. We normalize ${\Phi}$ such that ${\phi}  = -\Delta{\Phi} $ obeys
$$\int_{\mathbb{R}^{2}}{\phi}^2(x_1,x_2)\dif x_1\dif x_2=1.$$
By definition we know $\int_{\mathbb{R}^{2}}{\phi}  \dif x=0$.

We periodize them so that they are viewed as periodic functions on $\mathbb{T}^{2}$.
Consider a large  parameter  $\lambda \in\mN$. For every ${\xi}\in{\Lambda}$ we introduce
$${\Phi} _{({\xi})}(x):=
{\Phi} ({n}_*{\lambda}(x-\alpha_{{\xi}})\cdot A^1_{{\xi}},{n}_*{\lambda}(x-\alpha_{{\xi}})\cdot A^{2}_{{{\xi}}}),$$
$${\phi} _{({\xi})}(x):=
{\phi} ({n}_*{\lambda}(x-\alpha_{{\xi}})\cdot A^1_{{{\xi}}},{n}_*{\lambda}(x-\alpha_{{\xi}})\cdot A^{2}_{{{\xi}}}),$$
where $\alpha_{{\xi}}\in\mathbb{Q}^2$ are shifts  to ensure that $\{{\phi} _{({\xi)}}\}_{{\xi}\in{\Lambda}}$  have mutually disjoint supports by choosing $\epsilon_\Lambda>0$ correspondingly. For the existence of such shifts $\alpha_{\xi}$, we refer to \cite[Section 6.4]{BV}.

By construction, we know that $\xi \cdot \nabla \Phi_{(\xi)} = \xi\cdot\nabla\phi_{(\xi)}=0$ and $(n_*\lambda)^2\phi_{(\xi)}=-\Delta\Phi_{(\xi)}$. 

With this notation, the  Mikado flows $W_{(\xi)} \colon {\T} \to \R^3,\Theta_{(\xi)} \colon {\T} \to \R$ are defined as
\begin{align}
W_{(\xi)}(x) := W_{\xi,\lambda}(x) := \xi \, \phi_{(\xi)}(x),\ \ \Theta_{(\xi)}(x) := \Theta_{\xi,\lambda}(x) :=  \phi_{(\xi)}(x).
\label{eq:Mikado:def}
\end{align}
Since $\xi \cdot \nabla \phi_{(\xi)} = 0$, we immediately deduce that 
\begin{align}
 \div W_{(\xi)} = 0, \ \  \div \left(W_{(\xi)} \otimes W_{(\xi)} \right) = 0 \ \ \mbox{and} \ \ \div \left(W_{(\xi)}\Theta_{(\xi)} \right) = 0. \notag
\end{align}
By construction, the functions $W_{(\xi)}$ have zero mean on ${\T}$ and are in fact $(\mathbb{T}/{\lambda})^3$-periodic. Moreover, by our choice of $\alpha_\xi$ we have that 
\begin{align}
W_{(\xi)} \otimes W_{(\xi')} \equiv 0,\ \ W_{(\xi)}\Theta_{(\xi')} \equiv 0 \qquad \mbox{whenever} \qquad \xi \neq \xi' \in   \Lambda\,,
\label{eq:Mikado:2}
\end{align}
 and the normalization of $\phi_{(\xi)}$ ensures that 
\begin{align}
\fint_{\T} W_{(\xi)}(x) \otimes W_{(\xi)}(x)\,\dif x=\xi\otimes \xi,\ \ \fint_{\T} W_{(\xi)}(x) \Theta_{(\xi)}(x)\,\dif x=\xi\,. 
\label{eq:Mikado:3}
\end{align}
Lastly, using~\eqref{eq:Mikado:3}, the definition of the functions $\gamma_\xi,\Gamma_\xi$ in Lemma~\ref{l:linear_algebra2}, Lemma~\ref{l:linear_algebra} and the $L^2$ normalization of the functions $\phi_{(\xi)}$ we have that 
\begin{align}
 \sum_{\xi \in\overline \Lambda^i} \gamma_{\xi}^2(R) \fint_{\T} W_{(\xi)}(x) \otimes W_{(\xi)}(x) \dif x = R \, ,\ \   \sum_{\xi \in \Lambda^i} \Gamma_{\xi}(M) \fint_{\T} W_{(\xi)}(x)\Theta_{(\xi)}(x) \dif x = M,
 \label{eq:Mikado:4}
\end{align}
for every $i=1,2$, any symmetric matrix $R \in  \overline B_{\frac12}(\Id)$, and $M\in\overline  B_{1}(0)\backslash B_{\frac12}(0)$.

To define the incompressibility corrector in Section~\ref{sec:Onsager:principal:corrector:1}-Section~\ref{sec:Onsager:principal:corrector:2},   we note that $W_{(\xi)}$ can be written as $\curl V_{(\xi)}$ similar as \cite[(6.31)]{BV}, where we define 
\begin{align}
V_{(\xi)}: =  \frac{1}{(n_* \lambda)^2} \nabla \Phi_{(\xi)} \times \xi.
\notag
\end{align}

 With this notation we have the bounds for $N\geq 0$
\begin{align}
\norm{W_{(\xi)}}_{C^N} + \norm{\Theta_{(\xi)}}_{C^N} + \lambda\norm{V_{(\xi)}}_{C^N} \les\lambda^{N}.
\label{eq:Mikado:bounds}
\end{align}

\subsection{The estimates in  gluing steps}\label{app:sec:glu}
   In this section, we provide some estimates on the glued solutions for the Euler equations, which follow by a similar argument as in \cite[Section3, Section4]{BDLSV19}.  Here we give the detailed calculations taking into account  the different definitions of the parameters.
\begin{proof}[Proof of Proposition \ref{prop:gluedv}]\label{app:a:glu}
First by \cite[Proposition 3.1]{BDLSV19}, \eqref{e:v:ell:0} and \eqref{papa:gluing} we have for $t\in [ t_i- \tau_q, t_i+ \tau_q],N\geq0$
\begin{align}
    \norm{ v_i(t)}_{C^{N+1+\alpha}} \les \| v_{ l}\|_{C^{N+1+\alpha}}\les\delta_q^{1/2} \lambda_q l^{-N-\alpha}
\les  \tau_q^{-1}  l^{-N+\alpha}.\label{bd:vivlc1}
\end{align}

We note that $v_l-v_i$ obeys 
\begin{align*}
    (\partial_t+v_l\cdot\nabla)(v_l-v_i)=-(v_l-v_i)\cdot\nabla v_i-\nabla (\pi_l-\pi_i)+\div \mathring{R}_l,
\end{align*}
and 
\begin{align*}
 \nabla(   \pi_l-\pi_i)=\nabla\Delta^{-1}\div \((v_l-v_i)\cdot \nabla v_l+(v_l-v_i)\cdot \nabla v_i\)+\nabla\Delta^{-1}\div \div \mathring{R}_l.
\end{align*}

 Then by the transport estimate in 
 Proposition \ref{esti:transport}, \eqref{e:R:ell} and \eqref{bd:vivlc1} we have  for $t\in [ t_i- \tau_q, t_i+ \tau_q]$
\begin{align*}
\| (v_l-v_i)(t)\| _{C^{\alpha}} & \lesssim\left|\int_{t_i}^{t}\;\| (v_l-v_i)\cdot\nabla v_i(s)\| _{C^{\alpha}}+\| \nabla(\pi_l-\pi_i)(s)\| _{C^{\alpha}}+\|\mathring{R}_l\| _{C^{1+\alpha}}\mathrm{d}s\right|\\
 & \lesssim\tau_{q}\| \mathring{R}_l\| _{C^{1+\alpha}}+\left|\int_{t_i}^{t}\;\| \left(v_l-v_i\right)(s)\| _{C^{\alpha}}\left(\| v_i(s)\| _{C^{1+\alpha}}+\| v_l\| _{C^{1+\alpha}}\right)\mathrm{d}s\right|\\
 & \lesssim\tau_{q}\delta_{q+1}l^{-1+\alpha}+\left|\int_{t_i}^{t}\| \left(v_l-v_i\right)(s)\| _{C^{\alpha}}\tau_q^{-1}l^\alpha\mathrm{d}s\right|.
\end{align*}

 By Gronwall's inequality  we conclude that for $t\in [ t_i- \tau_q, t_i+ \tau_q]$
 \begin{align*}
\| (v_l-v_i)(t)\| _{C^{\alpha}} 
 & \lesssim\tau_{q}\delta_{q+1}l^{-1+\alpha}.
\end{align*}
 Moreover, by \eqref{e:R:ell} and \eqref{bd:vivlc1}  we have  for $t\in [ t_i- \tau_q, t_i+ \tau_q]$
\begin{align*}
\| \nabla(\pi_l-\pi_i)(t)\| _{C^{\alpha}} &  \lesssim\|(v_l-v_i)\cdot \nabla v_l(t)
\| _{C^{\alpha}}+\|(v_l-v_i)\cdot \nabla v_i(t)
\| _{C^{\alpha}}+\|\mathring{R}_l\| _{C^{1+\alpha}}\\
&\lesssim\tau_{q}\delta_{q+1}l^{-1+\alpha}\tau_{q}^{-1}l^\alpha+\delta_{q+1}l^{-1+\alpha}\lesssim \delta_{q+1}l^{-1+\alpha},
\end{align*}
and then we obtain for $t\in [ t_i- \tau_q, t_i+ \tau_q]$
\begin{align*}
\norm{(\partial_t +  v_{ l} \cdot \nabla) ( v_i- v_{ l})(t)}_{C^{\alpha}}  &  \lesssim\|\nabla v_i(v_l-v_i)
(t)\| _{C^{\alpha}}+\|\mathring{R}_l\| _{C^{1+\alpha}}+\delta_{q+1}l^{-1+\alpha}\\ & \lesssim\tau_{q}\delta_{q+1}l^{-1+\alpha}\tau_{q}^{-1}l^\alpha+\delta_{q+1}l^{-1+\alpha}\lesssim \delta_{q+1}l^{-1+\alpha}.
\end{align*}

Let $\theta$ be a multi-index with $\left|\theta\right|=N$, then by  \eqref{e:R:ell} and \eqref{bd:vivlc1} we have for $t\in [ t_i- \tau_q, t_i+ \tau_q]$
\begin{align*}
\| \partial^{\theta}\nabla(\pi_l-\pi_i)(t)\| _{C^{\alpha}} & \lesssim\| \mathring{R}_l\| _{C^{1+N+\alpha}}+\| (v_l-v_i)(t)\| _{C^{\alpha}}(\| v_i(s)\| _{C^{1+N+\alpha}}+\| v_l\| _{C^{1+N+\alpha}})\\
 & \qquad+\| (v_l-v_i)(t)\| _{C^{N+\alpha}}(\| v_i(t)\| _{C^{1+\alpha}}+\| v_l\| _{C^{1+\alpha}})\\
 & \lesssim\delta_{q+1}l^{-N-1+\alpha}+\tau_{q}\delta_{q+1}l^{-1-N+\alpha}\tau_q^{-1}l^\alpha+\| v_l-v_i\| _{C^{N+\alpha}}\tau_q^{-1}l^\alpha\\
 & \lesssim\delta_{q+1}l^{-N-1+\alpha}+\| v_l-v_i\| _{C^{N+\alpha}}\tau_q^{-1}l^\alpha.
\end{align*}
Moreover, by a similar calculation we have for $t\in [ t_i- \tau_q, t_i+ \tau_q]$
\begin{align*}   
\| \partial^{\theta}(\partial_t +  v_{ l} \cdot \nabla) (v_l-v_i)(t)\| _{C^{\alpha}}\lesssim\delta_{q+1}l^{-N-1+\alpha}+\| (v_l-v_i)(t)\| _{C^{N+\alpha}}\tau_q^{-1}l^\alpha.
\end{align*}
Using the transport estimate in 
 Proposition \ref{esti:transport} once again, we have  for $t\in [ t_i- \tau_q, t_i+ \tau_q]$
\begin{align*}
\| \partial^{\theta}&(v_l-v_i)(t)\| _{C^{\alpha}} \lesssim\left|\int_{t_{i}}^{t}\;\|(\partial_t +  v_{ l} \cdot \nabla) \partial^{\theta}\left(v_l-v_i\right)(s)\| _{C^{\alpha}}\mathrm{d}s\right|\\
 & \lesssim\left|\int_{t_{i}}^{t}\| [(\partial_t +  v_{ l} \cdot \nabla) ,\partial^{\theta}](v_l-v_i)(s)\| _{C^{\alpha}}+\delta_{q+1}l^{-N-1+\alpha}+\| (v_l-v_i)(s)\| _{C^{N+\alpha}}\tau_q^{-1}l^\alpha\mathrm{d}s\right|.
\end{align*}
By interpolation and \eqref{bd:vivlc1} we have  for $t\in [ t_i- \tau_q, t_i+ \tau_q]$
\begin{align*}
\|[(\partial_t +  v_{ l} \cdot \nabla) ,\partial^{\theta}](v_l-v_i)(t)\| _{C^{\alpha}} & \lesssim\| v_l\| _{C^{1+\alpha}}\|(v_l-v_i)(t)\| _{C^{N+\alpha}}+\| v_l\| _{C^{1+N+\alpha}}\| (v_l-v_i)(t)\| _{C^{\alpha}}\\
 & \lesssim\tau_q^{-1}l^{\alpha}\|(v_l-v_i)(t)\| _{C^{N+\alpha}}+\delta_{q+1}l^{-N-1+\alpha},
\end{align*} 
which implies that for $t\in [ t_i- \tau_q, t_i+ \tau_q]$

\begin{align*}
\| \partial^{\theta}(v_l-v_i)(t)\| _{C^{\alpha}} 
 & \lesssim\left|\int_{t_{i}}^{t}l^{-N-1+\alpha}\delta_{q+1}+\|( v_l-v_i)(s)\| _{C^{N+\alpha}}\tau_q^{-1}l^\alpha\mathrm{d}s\right|.
\end{align*}
Using  Gronwall's inequality, we obtain  \eqref{e:z_diff_k}: for $t\in [ t_i- \tau_q, t_i+ \tau_q]$
\begin{align*}
\|( v_l-v_i)(t)\| _{C^{N+\alpha}}+\tau_q\norm{(\partial_t +  v_{ l} \cdot \nabla) ( v_i- v_{ l})(t)}_{C^{N+\alpha}} & \lesssim\tau_q\delta_{q+1}l^{-N-1+\alpha}.
\end{align*}

    By a similar calculation as in \cite[(3.20)]{BDLSV19},  $z_i-z_l$ obeys
\begin{align*}
 (\partial_t + v_l \cdot\nabla) (z_i - z_l) 
 &= \Delta^{-1}\curl \div \mathring{R}_l + \Delta^{-1} \nabla \div \( ((z_i - z_l) \cdot \nabla) v_l \) \notag\\
 &\qquad + \Delta^{-1} \curl \div \( ((z_i - z_l) \times \nabla)v_l + ((z_i - z_l) \times \nabla) v_i^T \) \, .
\end{align*}
Consequently, by interpolation,  \eqref{e:R:ell} and \eqref{bd:vivlc1} we have for $t\in [ t_i- \tau_q, t_i+ \tau_q]$ 
\begin{align*}
   &\| (\partial_t + v_l \cdot\nabla) (z_i - z_l)(t) \|_{C^{N+\alpha}}\\
   &\les (\|v_i(t)\|_{C^{N+1+\alpha}}+\|v_l\|_{C^{N+1+\alpha}})\|(z_i - z_l)(t)\|_{C^\alpha}\notag\\ 
   &\qquad+(\|v_i(t)\|_{C^{1+\alpha}}+\|v_l\|_{C^{1+\alpha}})\|(z_i - z_l)(t)\|_{C^{N+\alpha}}+\|\mathring{R}_l\|_{C^{N+\alpha}}\\
   &\les \tau_q^{-1}l^{-N+\alpha}\|(z_i - z_l)(t)\|_{C^\alpha}+\tau_q^{-1}l^\alpha\|(z_i - z_l)(t)\|_{C^{N+\alpha}}+\delta_{q+1}l^{-N+\alpha}.
\end{align*}
When $N=0$, we have by Gronwall's inequality that for $t\in [ t_i- \tau_q, t_i+ \tau_q]$
\begin{align*}
     \|( z_i - z_l)(t) \|_{C^{\alpha}}+\tau_q \| (\partial_t + v_l \cdot\nabla) (z_i - z_l)(t) \|_{C^{\alpha}}\les \tau_q\delta_{q+1}l^{\alpha}.
\end{align*}
By commuting the derivatives in $N+\alpha,N\geq0$ with $\partial_t + v_l \cdot \nabla$ as before, we finish the proof of \eqref{e:z_diff}.
\end{proof}

\begin{proof}[Proof of Proposition \ref{prop:gluedr}]
By \eqref{papa:gluing}, \eqref{e:z_diff_k}, \eqref{e:z_diff} and the definition of the cut-off function  we have
\begin{align*}
    \norm{\overline{v}_q- v_{ l}}_{C^{N+\alpha}} &\les \sum_i  \chi_i \norm{v_i- v_{ l}}_{C^{N+\alpha}}\lesssim  \tau_q\delta_{q+1} l^{-1-N+\alpha}\les \delta_{q+1} ^{1/2}l^{-N+\alpha},\\
      \norm{\overline{z}_q- z_{ l}}_{C^{N+\alpha}} &\les \sum_i  \chi_i \norm{z_i- z_{ l}}_{C^{N+\alpha}}\lesssim  \tau_q\delta_{q+1} l^{-N+\alpha},
\end{align*}
which implies \eqref{e:vq:vell:additional}.

To estimate the energy gap, by the same calculation as in \cite[Proposition 4.4]{BDLSV19}, we  obtain on the time interval $I_i$
\begin{align*}
    |\overline v_q|^2 - |v_l|^2 = \chi_i(| v_i|^2 - |v_l|^2) + (1 -\chi_i)(|v_{i+1}|^2 - |v_l|^2) - \chi_i(1 -\chi_i)|v_i- v_{i+1}|^2.
\end{align*}
Then since $v_i$ solves the Euler equation, $(v_l,\mathring{R}_l )$ solves \eqref{e:euler_reynolds},  by the basic energy estimate, \eqref{e:R:ell} and \eqref{bd:vivlc1}, we have 
\begin{align*}
\left|\frac{\dif}{\dif t}\(
\int_{\mT^3} |v_i|^2 - |v_l|^2  \dif x\)\right|=\left|2\int_{\mT^3}
\nabla v_l:\mathring{R}_l \dif x\right|\les \|\nabla v_l\|_{C^0}\|\mathring{R}_l \|_{C^0}\les  \tau_q^{-1}\delta_{q+1}l^\alpha,
\end{align*}
Moreover,  after integrating in
time we deduce for $t\in [t_i-\tau_q,t_i+\tau_q]$
\begin{align*}
\left|
\int_{\mT^3} (|v_i|^2 - |v_l|^2  )(t)\dif x\right|\lesssim \delta_{q+1}l^\alpha.
\end{align*}
Furthermore, together with \eqref{e:vq:vell:additional} we obtain
\begin{align*}
    \left| \int_{\mathbb{T}^3}(|\overline v_q|^2- |v_l|^2 )(t)\dif x \right|  &\lesssim \delta_{q+1} l^{\alpha} +\|(v_i- v_{i+1})(t)\|_{C^\alpha}^2\lesssim \delta_{q+1} l^{\alpha},
\end{align*}
where we rewrite $ v_i -  v_{i+1} = ( v_i -  v_{ l}) - ( v_{i+1} -  v_{ l})$.  This concludes the proof of \eqref{bd: energy vq-vl}.

 For the Reynolds stress, by the definition of the cut-off functions,  the definition in \eqref{eq:bar:R_q:def}, the choice of parameters in \eqref{papa:gluing}, the bounds \eqref{e:z_diff_k} and \eqref{e:z_diff} we have for $t\in I_i$
  \begin{align*}
      \norm{\mathring{\overline{ R}}_{q}(t)}_{C^{N+\alpha}} &\lesssim \|\partial_t \chi_i\|_{C_t^0}\|\mathcal{R}( v_i- v_{i+1})(t)\|_{C^{N+\alpha}} +\| (v_i- v_{i+1})(t)\|_{C^{N+\alpha}} \| (v_i- v_{i+1})(t)\|_{C^\alpha}\\
      &\les \tau_q^{-1}\tau_q\delta_{q+1}l^{-N+\alpha}+\tau_q^2\delta^2_{q+1}l^{-2-N+\alpha}\les 
    \delta_{q+1}l^{-N+\alpha},
  \end{align*}
  where we rewrite $ v_i -  v_{i+1} = ( v_i -  v_{ l}) - ( v_{i+1} -  v_{ l})$. By  direct calculation we have for $t\in I_i$
  \begin{align*}
(\partial_t + v_l\cdot \nabla)\mathring{\overline{R}}_{q} & =\partial_{t}^{2}\chi_{i}\mathcal{R}\curl\left(z_{i}-z_{i+1}\right)+\partial_{t}\chi_{i}\mathcal{R}\curl(\partial_t + v_l\cdot \nabla)(z_{i}-z_{i+1})\\
 & +\partial_{t}\chi_{i}[v_l\cdot\nabla,\mathcal{R}\curl](z_{i}-z_{i+1}) +\partial_{t}(\chi_{i}^{2}-\chi_{i})(v_{i}-v_{i+1})\mathring{\otimes}(v_{i}-v_{i+1})\\
 & +(\chi_{i}^{2}-\chi_{i})\left((\partial_t + v_l\cdot \nabla)(v_{i}-v_{i+1})\right)\mathring{\otimes}
\left(v_{i}-v_{i+1}\right),\\
 & +(\chi_{i}^{2}-\chi_{i})\left(v_{i}-v_{i+1}\right)\mathring{\otimes}\left((\partial_t + v_l\cdot \nabla)(v_{i}-v_{i+1})\right),
\end{align*}
which together with the definition of the cut-off functions, the choice of parameters in  \eqref{papa:gluing} and  the bounds in \eqref{e:v:ell:0},  \eqref{eq:v_i:C^N}-\eqref{e:z_diff} implies that  for $t\in I_i$
\begin{align*}
    \|(\partial_t + v_l\cdot \nabla)\mathring{\overline{R}}_{q} (t)\|_{C^{N+\alpha}}&\les \tau_q^{-2}\|(z_{i}-z_{i+1})(t)\|_{C^{N+\alpha}}+\tau_q^{-1}\|(\partial_t + v_l\cdot \nabla)(z_{i}-z_{i+1})(t)\|_{C^{N+\alpha}}\\
 & +\tau_q^{-1}\|v_l\|_{C^{1+\alpha}}\|(z_{i}-z_{i+1})(t)\|_{C^{N+\alpha}}+\tau_q^{-1}\|v_l\|_{C^{N+1+\alpha}}\|(z_{i}-z_{i+1})(t)\|_{C^{\alpha}}\\
 & +\tau_q^{-1}\|(v_{i}-v_{i+1})(t)\|_{{C^{N+\alpha}}}\|(v_i-v_{i+1})(t)\|_{C^\alpha}\\
 & +\|(\partial_t + v_l\cdot \nabla)(v_{i}-v_{i+1})(t)\|_{C^{N+\alpha}}\|
(v_{i}-v_{i+1})(t)\|_{C^\alpha} \\
&+\|(\partial_t + v_l\cdot \nabla)(v_{i}-v_{i+1})(t)\|_{C^{\alpha}}\|
(v_{i}-v_{i+1})(t)\|_{C^{N+\alpha}}\\
&\les \tau_q^{-1}\delta_{q+1}l^{-N+\alpha}+\tau_q\delta_{q+1}^2l^{-N-2+\alpha}\les \tau_q^{-1}\delta_{q+1}l^{-N+\alpha},
\end{align*}
where the commutator is bounded using \cite[Proposition D.1]{BDLSV19} as $\mathcal{R}\curl$ is a Calder\'on-Zygmund operator.

Then we use interpolation, \eqref{papa:gluing} and the estimate \eqref{e:vq:vell:additional} to deduce that for $t\in I_i$
\begin{align*}
      \|(\partial_t + \overline v_q\cdot \nabla)\mathring{\overline{R}}_{q}(t) \|_{C^{N+\alpha}}&\les   \|(\partial_t + v_l\cdot \nabla)\mathring{\overline{R}}_{q} (t)\|_{C^{N+\alpha}}\notag\\ 
      &\qquad+\|(v_l-\overline v_q)(t)\|_{C^{\alpha}}\|\mathring{\overline{R}}_{q}(t) \|_{C^{1+N+\alpha}}+\|(v_l-\overline v_q)(t)\|_{C^{N+\alpha}}\|\mathring{\overline{R}}_{q}(t) \|_{C^{1+\alpha}}\\
      &\les \tau_q^{-1}\delta_{q+1}l^{-N+\alpha}+\tau_q\delta_{q+1}^2l^{-N-2+\alpha}\les \tau_q^{-1}\delta_{q+1}l^{-N+\alpha}.
\end{align*}

\end{proof}

To conclude this section, we provide the derivation of the analytic identities in  Lemma \ref{lem;identity}.

 \begin{proof}[Proof of Lemma \ref{lem;identity}]
  The proof follows by direct calculation:
\begin{align*}
    (\curl z)\cdot\nabla \rho&=\sum_{i,j,k}\epsilon_{ijk}\partial_jz^k\partial_i \rho=\sum_{i,j,k}\(\partial_j(\epsilon_{ijk}z^k\partial_i \rho)-\epsilon_{ijk}z^k\partial_i\partial_j\rho\)\\
    &=\sum_{i,j,k}\partial_j(\epsilon_{ijk}z^k\partial_i \rho)=\div (z\times\nabla\rho),
\end{align*}
     where the second term equals to 0 since $\epsilon_{ijk}\partial_i\partial_j\rho=-\epsilon_{jik}\partial_j\partial_i\rho$. Here $\epsilon_{ijk}$ is the Levi-Civita symbol in 3D which is defined as follows:
$\epsilon_{ijk} = 1/-1$
  if $\{i, j, k\}$ are an even/odd permutation of the indices $\{1, 2, 3\}$. Otherwise, $\epsilon_{ijk} = 0$.  For the second term, since $\div v=0$, we have
     \begin{align*}
         v\cdot\nabla(\div z)&=\sum_{i,j}v^i\partial_i\partial_j z^j=\sum_{i,j}\(\partial_j (v^i\partial_i z^j)-\partial_jv^i\partial_i  z^j\)=\div (v\cdot\nabla z)-\sum_{i,j}\(\partial_i(\partial_jv^i z^j)-\partial_i\partial_jv^i z^j\)\\
         &=\div (v\cdot\nabla z-z\cdot\nabla v).
         \end{align*}
         For the last term,  since $\div v=0$, we have for $i=1,2,3$
            \begin{align*}
   [ \curl(v\cdot\nabla z)]^i&=\sum_{j,k,l}\epsilon_{ijk}\partial_j(v^l\partial_lz^k)=\sum_{j,k,l}\(\epsilon_{ijk}\partial_jv^l\partial_lz^k+\epsilon_{ijk}v^l\partial_j\partial_lz^k\)\\
    &=\sum_{j,k,l}\(\partial_l(\epsilon_{ijk}\partial_jv^lz^k)-\epsilon_{ijk}\partial_l\partial_jv^lz^k+v^l\partial_l(\epsilon_{ijk}\partial_jz^k)\)\\
    &=-[\div((z\times \nabla)v)]^i+[v\cdot\nabla (\curl z)]^i,
\end{align*}
where we use the notation $[(z\times \nabla)v]^{il}=\sum_{j,k}\epsilon_{ikj}z^k\partial_jv_l$.
 \end{proof}

\subsection{The estimates in  perturbation steps}\label{sec:est:ampl}
In this section, we show the estimates of the amplitude functions appearing in Section \ref{sec:proof13}. Before that, we give the estimates of the flow maps.
 \begin{proof}[Proof of Proposition \ref{prop:flowmap}]
Since for every $t \in  (t_i-\frac{\tau_q}{3},t_i+\frac{4\tau_q}{3})$ we have $|t-t_i| \leq 2 \tau_q$,  using the estimates \eqref{e:vq:1} and Proposition \ref{esti:transport}
we obtain
  \begin{align*}
 \norm{\nabla  \Phi_i(t) - {\rm Id}}_{C^0} &\les\tau_q\|\nabla \overline v_q\|_{C^\alpha} \les l^{\alpha},\end{align*}
 which implies the estimate in \eqref{eq:Phi:i:bnd:a}. Then it is easy to see that  $(\nabla  \Phi_i)^{-1}(t)$ is well-defined  for  $t \in  (t_i-\frac{\tau_q}{3},t_i+\frac{4\tau_q}{3})$ and
$
      \norm{(\nabla  \Phi_i)^{-1}(t)}_{C^{0}}\les 1.$

By applying Proposition \ref{esti:transport} again we have for $N\geq1$
  \begin{align}
 \norm{\nabla  \Phi_i(t)}_{C^{N}}  \les\tau_q\|\nabla \overline v_q\|_{C^{N}}  \les  l^{-N},\notag
 \end{align}
 which together with the Leibniz rule  implies that for $N\geq1$
 \begin{align*}
     \norm{(\nabla  \Phi_i)^{-1}(t)}_{C^{N}}&\lesssim \norm{(\nabla  \Phi_i)^{-1}(t)}_{C^{0}}\sum_{m=0}^{N-1}\norm{\nabla  \Phi_i(t)}_{C^{N-m}}  \norm{(\nabla  \Phi_i)^{-1}(t)}_{C^{m}}\\
     &\les \sum_{m=0}^{N-1}l^{-N+m}  \norm{(\nabla  \Phi_i)^{-1}(t)}_{C^{m}}.
 \end{align*}
 By induction we  obtain  \eqref{eq:Phi:i:bnd:b}. To derive \eqref{eq:Phi:i:bnd:c}, since $D_{t,q}\nabla  \Phi_i=-\nabla  \Phi_iD\overline v_q,$ we have for $N\geq0$
 \begin{align}
 \norm{D_{t,q} \nabla  \Phi_i(t)}_{C^{N}} & \les \|\nabla\Phi_i(t)\|_{C^0}\|\overline v_q\|_{C^{N+1}}+\|\nabla\Phi_i(t)\|_{C^{N}}\|\overline v_q\|_{C^1}\les \tau_q^{-1}l^{-N},
\end{align}
and 
\begin{align*}
    &\norm{D_{t,q} (\nabla  \Phi_i)^{-1}(t)}_{C^{N}}\\
    &
    \les \norm{D_{t,q} \nabla  \Phi_i(t)}_{C^{N}}\norm{(\nabla  \Phi_i)^{-1}(t)}_{C^{0}}^2+\norm{D_{t,q} \nabla  \Phi_i(t)}_{C^{0}}\norm{(\nabla  \Phi_i)^{-1}(t)}_{C^{N}}\norm{(\nabla  \Phi_i)^{-1}(t)}_{C^{0}}\\
    &\les \tau_q^{-1}l^{-N}.
\end{align*}
\end{proof}

\begin{proof}[Proof of Proposition \ref{prop:est:A,tilA}]
We recall the definition of $M_{q,i}$ in \eqref{eq:tilde:M:q:i:def}.
Using the Leibniz rule, the estimates on $\overline M_q, \nabla\Phi_i$ in \eqref{e:Mq:1}, and \eqref{eq:Phi:i:bnd:a}-\eqref{eq:Phi:i:bnd:c} we have for all $t \in (t_i-\frac{\tau_q}{3},t_i+\frac{4\tau_q}{3})$ and $N\geq 0$
    \begin{align*}
        \norm{ M_{q,i}(t)}_{C^N} &\les \|\nabla\Phi_i(t)\|_{C^0}\(1+\frac{\|\overline M_q\|_{C^N}}{\delta_{q+1} ^{1/2}\tilde{\delta}_{q+1} ^{1/2}l^{\alpha/2}}\)+ \|\nabla\Phi_i(t)\|_{C^N}\(1+\frac{\|\overline M_q\|_{C^0}}{\delta_{q+1} ^{1/2}\tilde{\delta}_{q+1} ^{1/2}l^{\alpha/2}}\) \les   l^{-N},\\
        \norm{D_{t,q}  M_{q,i}(t)}_{C^N}&\les \|D_{t,q} \nabla\Phi_i(t)\|_{C^0}\(1+\frac{\|\overline M_q\|_{C^N}}{\delta_{q+1} ^{1/2}\tilde{\delta}_{q+1} ^{1/2}l^{\alpha/2}}\)+ \|D_{t,q} \nabla\Phi_i(t)\|_{C^N}\(1+\frac{\|\overline M_q\|_{C^0}}{\delta_{q+1} ^{1/2}\tilde{\delta}_{q+1} ^{1/2}l^{\alpha/2}}\)\notag\\ 
        &+ \|\nabla\Phi_i(t)\|_{C^0}\(1+\frac{\|D_{t,q} \overline M_q\|_{C^N}}{\delta_{q+1} ^{1/2}\tilde{\delta}_{q+1} ^{1/2}l^{\alpha/2}}\)+ \|\nabla\Phi_i(t)\|_{C^N}\(1+\frac{\|D_{t,q} \overline M_q\|_{C^0}}{\delta_{q+1} ^{1/2}\tilde{\delta}_{q+1} ^{1/2}l^{\alpha/2}}\)\les  \tau_q^{-1} l^{-N}.
    \end{align*}
    
    By applying \cite[Proposition C.1]{TCP} we have for $N\geq1$
    \begin{align*}
        \| \Gamma_{\xi} ^{1/2}( M_{q,i})(t)\|_{C^N}\les \| \Gamma_{\xi} ^{1/2}( M_{q,i})(t)\|_{C^0}+\|\Gamma_{\xi} ^{1/2}\|_{C^1}\|M_{q,i}(t)\|_{C^N}+\|\Gamma_{\xi} ^{1/2}\|_{C^{N}}\|M_{q,i}(t)\|_{C^1}^N\les   l^{-N}.
    \end{align*}
    It is easy to see that the above is still valid for $N=0$.
    
    Then we calculate $ D_{t,q}\Gamma_{\xi} ^{1/2}( M_{q,i})=\nabla\Gamma_{\xi} ^{1/2}( M_{q,i})\cdot  D_{t,q}M_{q,i}$, by the same argument we have for $N\geq0$
    \begin{align*}
        \|\nabla \Gamma_{\xi} ^{1/2}( M_{q,i})(t)\|_{C^N}\les l^{-N},
    \end{align*}
    which together with  the Leibniz rule implies that
    \begin{align*}
         \|D_{t,q}&\Gamma_{\xi} ^{1/2}( M_{q,i})(t)\|_{C^N}\\
         &\lesssim  \|\nabla\Gamma_{\xi} ^{1/2}( M_{q,i})(t)\|_{C^0}\|  D_{t,q}M_{q,i}(t)\|_{C^N}+\|\nabla \Gamma_{\xi} ^{1/2}( M_{q,i})(t)\|_{C^{N}}\|  D_{t,q}M_{q,i}(t)\|_{C^0}\lesssim  \tau_q^{-1} l^{-N}.
    \end{align*}

    By applying the Leibniz rule again, and by the properties of the cut-off functions we have
    \begin{align*}
         \|A_{(\xi,i)}\|_{C^N}&\les l^{\alpha/4} \delta_{q+1}^{1/2}\| \eta_i\|_{C^0}\|\Gamma_{\xi} ^{1/2}( M_{q,i})\|_{C^N}+ l^{\alpha/4} \delta_{q+1}^{1/2}\| \eta_i\|_{C^N}\|\Gamma_{\xi} ^{1/2}( M_{q,i})\|_{C^0}\les \delta_{q+1}^{1/2}  l^{\alpha/4-N},\\
           \|D_{t,q}A_{(\xi,i)}\|_{C^N}&\les l^{\alpha/4} \delta_{q+1}^{1/2}\(\| D_{t,q}\eta_i\Gamma_{\xi} ^{1/2}( M_{q,i})\|_{C^N}+ \|\eta_iD_{t,q} \Gamma_{\xi} ^{1/2}( M_{q,i})\|_{C^N}\)\les  \tau_q^{-1} \delta_{q+1}^{1/2}  l^{\alpha/4-N},
    \end{align*}
     where we used $\tau_q^{-1}\leq l^{-1}$, since $\delta_{q+1}^{1/2}\leq l^{2\alpha}\lambda_q^{3\alpha/2}$ by choosing $\alpha>0$ small enough.
 The other term  $\tilde{A}_{(\xi,i)}$ can be estimated by the same calculation. So we omit the proof.
\end{proof}

  \begin{proof}[Proof of Proposition \ref{lem:eq:Onsager:a:xi:CN}]

First we estimate the derivatives of $ R_{q,i}$ defined in \eqref{eq:tilde:R:q:i:def}. More precisely, we show that for $N\geq0$ and $t\in \supp(\overline \eta_i)$
   \begin{align} 
 \norm{ R_{q,i}(t)}_{C^N} +  \tau_q \norm{D_{t,q}  R_{q,i}(t)}_{C^N} \les   l^{-N}.
 \label{eq:tilde:R:q:i:bnd}
\end{align}
  In fact, by the properties of the cutoff functions, the $C^N$-bounds in \eqref{e:Rq:1}, \eqref{bd:rq(1)}  and  \eqref{estimate:rho1},  we have for $N\geq0$
\begin{align*}
     \norm{ \frac{\sum_j \int
     \overline \eta_j^2(t,y) \dif y} {\Upsilon_q(t)}(\Rbar_q+ \mathring{R}^{(1)}_{q}) }_{C^N} \les \frac{\lambda_{q}^{\alpha/3}} {\delta_{q+1}}\norm{\Rbar_q+ \mathring{R}^{(1)}_{q} }_{C^N} \les \lambda_{q}^{\alpha/3}l^{\alpha/2-N}
 \les   l^{-N},
\end{align*}
 where we used the fact that $(l\lambda_q)^\alpha\ll1$.
Then by the Leibniz rule and  \eqref{eq:Phi:i:bnd:b}, the  bound for the first term follows. 

Then, to estimate the  material derivative,  we  have
\begin{align*}
    D_{t,q}&\( \frac{\sum_j \int
    \overline \eta_j^2 \dif x} {\Upsilon_q}(\Rbar_q+ \mathring{R}^{(1)}_{q}) \)= \partial_t \left( \frac{ \sum _j\int
    \overline \eta_j^2\dif x}{\Upsilon_q}\right) (\Rbar_q+ \mathring{R}^{(1)}_{q})+  \frac{ \sum_j \int
    \overline \eta_j^2\dif x}{\Upsilon_q} D_{t,q} (\Rbar_q+ \mathring{R}^{(1)}_{q}).
\end{align*}
		By the Leibniz rule for the derivative of the product, \eqref{estimate:rho1} and \eqref{bd:parttrhoq}  we derive
		\begin{align}
			\left| \partial_t \left( \frac{ \sum _j\int
            \overline \eta_j^2\dif x}{\Upsilon_q}\right)\right| \lesssim \frac{\sup_j\|\partial_t\overline \eta_j\|_{C^0}}{\Upsilon_q} +\left| \frac{\partial_t\Upsilon_q}{\Upsilon_q^2}\right|\lesssim \delta_{q+1}^{-1}\lambda_{q}^{\alpha/3}\tau_q^{-1}+\delta_{q+1}^{-2}\lambda_{q}^{2\alpha/3}\tau_q^{-1}\delta_{q+1}l^\alpha\leq \delta_{q+1}^{-1}\lambda_{q}^{\alpha/3}\tau_q^{-1},\notag
		\end{align}
        where we used the fact that  $(l\lambda_q)^{\alpha}\ll 1$. Then together with \eqref{estimate:rho1},  \eqref{e:Rq:1} and \eqref{bd:rq(1)} we obtain that
		\begin{equation}
				\norm{D_{t,q}\( \frac{\overline  \eta_i  ^2(\Rbar_q+ \mathring{R}^{(1)}_{q})}{\Upsilon_{q,i}}\)}_{C^N}
				\lesssim \delta_{q+1}^{-1} \lambda_{q}^{\alpha/3} \tau_q^{-1}\delta_{q+1}l^{-N+\alpha/2} + \delta_{q+1}^{-1}  \lambda_{q}^{\alpha/3} \tau_q^{-1}\delta_{q+1}l^{-N+\alpha/2}
			\lesssim
				\tau_q^{-1} l^{-N}\, ,\notag
		\end{equation}
        where we used the fact that  $(l\lambda_q)^{\alpha}\ll 1$. By applying the  Leibniz rule again, together with the estimates on $\Phi_i$ in  \eqref{eq:Phi:i:bnd:b}, \eqref{eq:Phi:i:bnd:c}  we  obtain the bound for the second term in \eqref{eq:tilde:R:q:i:bnd}.

	By applying \cite[Proposition C.1]{TCP}  and \eqref{eq:tilde:R:q:i:bnd} we have for $N\geq1$
		\begin{align}
			\|\gamma_{\xi}( R_{q,i})\|_{C^{N}}\lesssim \|\gamma_\xi\|_{C^1}\| R_{q,i}\|_{C^{N}}+  \|\gamma_\xi\|_{C^{N}}\| R_{q,i}\|^{N}_{C^{1}} \lesssim l^{-N}.\label{chain:a xi:2}
		\end{align}
By the definition of $\Upsilon_{q,i}$ in \eqref{def:rho q,i}, \eqref{estimate:rho1} and the properties of the cut-off functions  $\overline\eta_i$, we deduce that $$ \|\Upsilon_{q,i}^{1/2}\|_{C^N}\lesssim  \|\Upsilon_{q}^{1/2}\|_{C_t^0}\|\overline\eta_i\|_{C^N}\lesssim\delta_{q+1}^{1/2}.$$ Then together with \eqref{chain:a xi:2}, we obtain for $N\geq1$
		\begin{align*}
			\| a_{(\xi,i)}\|_{C^N} \lesssim  \| \Upsilon_{q,i}^{1/2}\|_{C^0} \|\gamma_{\xi}(R_{q,i})\|_{C^{N}}+\|\Upsilon_{q,i}^{1/2}\|_{C^N} \|\gamma_{\xi}(R_{q,i})\|_{C^{0}} \lesssim  \delta_{q+1}^{1/2}l^{-N}.
		\end{align*}
		It is easy to see that the bound is also valid for $N=0$.
        
		Then we calculate that
		\begin{align*}
			D_{t,q} (\Upsilon_{q,i}^{1/2})=\left[  \partial_t \left( \frac{\overline \eta_i}{(\sum_j \int
            \overline\eta_j^2\dif x)^{1/2}} \right) +
			\frac{ \overline{v}_q\cdot \nabla \overline\eta_i }{(\sum _j\int
            \overline\eta_j^2\dif x)^{1/2}} \right]  \Upsilon_q^{1/2}+ \frac{\overline\eta_i \partial_t(\Upsilon_q^{1/2})}{(\sum _j\int
            \overline\eta_j^2\dif x)^{1/2}} .
		\end{align*}

        By the properties of the cut-off functions and \eqref{e:vq:1} we have for $N\geq1$
        \begin{align*}
             \norm{\partial_t \left( \frac{\overline \eta_i}{(\sum_j \int
             \overline\eta_j^2\dif x)^{1/2}} \right)}_{C^N}&\lesssim \tau_q^{-1},\\ 
            \|\overline{v}_q\cdot \nabla \overline\eta_i \|_{C^N}\lesssim \|\overline{v}_q\|_{C^1}\| \nabla \overline\eta_i \|_{C^N}+\|\overline{v}_q\|_{C^N}\| \nabla \overline\eta_i \|_{C^0}\lesssim\tau_q^{-1}l^{-N+1}&\lesssim\tau_q^{-1}l^{-N}.
        \end{align*}
   It is easy to see that the bounds are also valid for $N=0$, since 
        \begin{align*}
            \|\overline{v}_q\cdot \nabla \overline\eta_i \|_{C^0}\lesssim \|\overline{v}_q\|_{C^1}\| \nabla \overline\eta_i \|_{C^0}\lesssim\tau_q^{-1}.
        \end{align*}
        
        By \eqref{estimate:rho1} and \eqref{bd:parttrhoq} we have $$|  \partial_t(\Upsilon_q^{1/2})|\lesssim \left| \frac{\partial_t \Upsilon_q}{\Upsilon_q^{1/2}} \right| \lesssim \tau_q^{-1} \delta_{q+1}^{1/2}l^\alpha \lambda_q^{\alpha/6}\lesssim \tau_q^{-1} \delta_{q+1}^{1/2},$$
        where we used the fact that $l\lambda_q\leq1$ in the last inequality. Then it follows that for $N\geq0$
		\begin{align*}
			\|  D_{t,q} (\Upsilon_{q,i}^{1/2})\|_{C^N} \lesssim & \tau_q^{-1}l^{-N}\|\Upsilon_q^{1/2}\|_{C^0_t}+ \tau_q^{-1} \delta_{q+1}^{1/2}\lesssim  \tau_q^{-1}\delta_{q+1}^{1/2}l^{-N}.
		\end{align*}

          Moreover, 
        \begin{align*}
            D_{t,q}(\gamma_{\xi}( R_{q,i}))=\nabla\gamma_{\xi}( R_{q,i})\cdot D_{t,q}R_{q,i}.
        \end{align*}
	By applying \cite[Proposition C.1]{TCP}  and \eqref{eq:tilde:R:q:i:bnd} we have for $N\geq1$
		\begin{align}
            \|\nabla\gamma_{\xi}( R_{q,i})\|_{C^{N}}\lesssim \|\nabla\gamma_\xi\|_{C^1}\| R_{q,i}\|_{C^{N}}+  \|\nabla\gamma_\xi\|_{C^{N}}\| R_{q,i}\|^{N}_{C^{1}} \lesssim l^{-N},\notag
		\end{align}
        which together with  \eqref{eq:tilde:R:q:i:bnd}  implies that for $N\geq1$
        \begin{align*}
            \| D_{t,q}(\gamma_{\xi}( R_{q,i}))\|_{C^N}\lesssim \|\nabla\gamma_{\xi}( R_{q,i})\|_{C^N}\| D_{t,q}R_{q,i}\|_{C^0}+\|\nabla\gamma_{\xi}( R_{q,i})\|_{C^0}\| D_{t,q}R_{q,i}\|_{C^N}\lesssim \tau_q^{-1}l^{-N}.
        \end{align*}
        It is easy to see that the bound is also valid for $N=0$.
        
In the end, by applying the chain rule again, we summarize all  the bounds above and obtain $N\geq 0$
		\begin{align*}
			\| D_{t,q} a_{(\xi,i)}\|_{C^N}&\lesssim  \|  D_{t,q} (\Upsilon_{q,i}^{1/2})\gamma_{\xi}(\tilde R_{q,i})\|_{C^{N}}+  \| \Upsilon_{q,i}^{1/2}D_{t,q} \gamma_{\xi}(\tilde R_{q,i})\|_{C^{N}} \lesssim\tau_q^{-1}   \delta_{q+1}^{1/2}l^{-N}.
		\end{align*}
	\end{proof}

\renewcommand{\appendixname}{Appendix~\Alph{section}}
  \renewcommand{\theequation}{C.\arabic{equation}}

  \section{ Building blocks and  auxiliary estimates in Section \ref{proof:prop2}}
  \label{s:appA.5}
  In this section, we first  introduce the building blocks used in the convex integration method. Then we provide some estimates on the amplitude functions  appearing in Section \ref{proof:prop2}.
\subsection{Generalized intermittent spatial-time jets}\label{gij}

\subsubsection{Building blocks for the transport equations} \label{sec:bbte}
 In this section, we present the building blocks for advection-diffusion equations, which can be viewed as a generalization of those in \cite[Section 4]{BCDL21}, incorporating more intermittency in the spatial domain. 

For parameters $\lambda,r_\perp, r_\parallel > 0 $, we assume
$$\lambda^{-1}\ll r_{\perp}\ll r_{\parallel}\ll 1,\ \ \lambda r_{\perp}\in\mathbb{N}.$$

We recall the geometric Lemma \ref{l:linear_algebra2} and the 
 two disjoint families $\Lambda^1,\Lambda^2$ discussed in that lemma for $d\geq2$.
For each $\xi\in\Lambda^1\cup \Lambda^2$ let us define $A^i_\xi\in \mathbb{S}^{d-1}\cap \mathbb{Q}^d,\ i=1,2,...,d-1$  such that $\{\xi, A^i_\xi,i=1,...,d-1\} $
form an orthonormal basis in $\mathbb{R}^d$. We label by $n_*$ the smallest
natural number such that
$$\{n_*\xi, n_*A^i_\xi,i=1,...,d-1\}\subset\mathbb{Z}^d,$$
for every $\xi\in\Lambda^1\cup \Lambda^2$.

 Now we introduce the cut-off functions used in the construction. 
 First, there is a smooth mean-zero function  $\phi : \mathbb{R}^{d-1} \to \mathbb{R}$  with support in $B(0,1)$   satisfying $\phi\equiv 1$ on $B(0,\frac13)$. Moreover, the function  $\Phi$,  defined by $\phi = -\Delta\Phi$, is a smooth function with support in  $B(0,1)$.
 In fact, let $\phi_0: \mathbb{R}^{d-1} \to \mathbb{R}$ be  a smooth function with support in $B(0,1)$ satisfying $\phi_0\equiv 1$ on $B(0,\frac13)$. By \cite[Lemma 1.12]{MB03} we can define $\Phi _0: \mathbb{R}^{d-1} \to \mathbb{R}$  by solving $\phi_0 = -\Delta\Phi_0$. We let $\Phi:=\Phi _0\phi_0$ which is a smooth function with support in a ball of radius 1. We define $\phi=-\Delta\Phi$, which is a smooth function with support in a ball of radius 1 and mean-zero.
 It is easy to check that on $B(0,\frac13)$, $\phi=-\Delta(\Phi _0\phi_0)=-\Delta\Phi _0=\phi_0=1$.

Let  $\psi : \mathbb{R}\to\mathbb{R}$  be a smooth, mean-zero function with
support in $B(0,1)$ satisfying $\psi\equiv 1$ on $B(0,\frac13)$.
Define $\phi' : \mathbb{R}^{d-1} \to \mathbb{R}$ to be a smooth non-negative function with support in $B(0,\frac13)$ satisfying $$\int_{\mathbb{R}^{d-1}}\phi'(x_1,x_2,...,x_{d-1})\dif x_1\dif x_2..\dif x_{d-1}=1,$$
and let $\psi' : \mathbb{R} \to \mathbb{R}$ be a smooth non-negative function with support in $B(0,\frac13)$ such that
$$\int_\mathbb{R}\psi'(x_d)\dif x_d=1.$$

Then it is easy to see that 
\begin{align}
    \phi\phi'=\phi',\ \ \psi\psi' =\psi'.\label{phi2}
\end{align}

We define the rescaled cut-off functions by
$$\phi _{r_{\perp}}(x_1,x_2,...,x_{d-1})=\frac1{r_{\perp}^{(d-1)/2}}\phi (\frac{x_1}{r_{\perp}},\frac{x_2}{r_{\perp}},...,\frac{x_{d-1}}{r_{\perp}}),$$
$$\Phi _{r_{\perp}}(x_1,x_2,...,x_{d-1})=\frac1{r_{\perp}^{(d-1)/2}}\Phi (\frac{x_1}{r_{\perp}},\frac{x_2}{r_{\perp}},...,\frac{x_{d-1}}{r_{\perp}}),$$
$$\psi _{r_{\parallel}}(x_d)=\frac1{r_{\parallel}^{1/2}}\psi (\frac{x_d}{r_{\parallel}}).$$
Similarly, we define the rescaled cut-off  $\phi '_{r_{\perp}},\psi'_{r_{\parallel}}$ with respect to $\phi ',\psi'$.
We periodize $\phi _{r_{\perp}},\Phi _{r_{\perp}},\psi _{r_{\parallel}},\phi '_{r_{\perp}},\psi'_{r_{\parallel}}$ so that they can be viewed as   functions on $\mathbb{T}^{d-1}$ and $\mathbb{T}$ respectively. Consider a large time oscillation parameter $\mu > 0$. For every $\xi\in\Lambda^1\cup\Lambda^2$ we introduce
\begin{align*}
    \psi _{(\xi)}(t,x):&=\psi _{r_{\parallel}}(n_*r_{\perp}\lambda(x\cdot \xi-\mu t)),\\
 \Phi _{(\xi)}(x):&=\Phi _{r_{\perp}}(n_*r_{\perp}\lambda x\cdot A^1_{\xi},...,n_*r_{\perp}\lambda x\cdot A^{d-1}_{\xi}),\\
\phi _{(\xi)}(x):&=\phi _{r_{\perp}}(n_*r_{\perp}\lambda x\cdot A^1_{\xi},...,n_*r_{\perp}\lambda x\cdot A^{d-1}_{\xi}).
\end{align*}
Here we remark that we do not need to translate the building blocks such that the supports are disjoint, since the  disjoint support property will be achieved by choosing suitable time jets in the following section.
 Similarly, we define the building blocks  $\phi '_{(\xi)},\psi'_{(\xi)}$. 

The building blocks $W_{(\xi)} :\mathbb{R} \times\mathbb{T}^d \to \mathbb{R}^d,\Theta_{(\xi)} :\mathbb{R} \times\mathbb{T}^d \to \mathbb{R}$ are defined as
$$W_{(\xi)}(t,x):=\xi\psi_{(\xi)}(t,x)\phi_{(\xi)}(x),\ \ \Theta_{(\xi)}(t,x):=\psi'_{(\xi)}(t,x)\phi'_{(\xi)}(x),$$
at which point, together with identities in \eqref{phi2} we have that
\begin{align}
\int _{\mathbb{T}^d}W_{(\xi)}\Theta_{(\xi)}\dif x=\xi,\label{eq:intwthe}\\
\partial_t\Theta_{(\xi)}+\mu r_\perp^{\frac{d-1}{2}}r_\parallel^{\frac12}\div (W_{(\xi)}\Theta_{(\xi)})=0.\label{eq:ptthe+}
\end{align}

Since $W_{(\xi)}$ is not divergence-free, inspired by \cite[Section 4.1]{CL22} we
introduce the skew-symmetric  corrector term
\begin{align}
V_{(\xi)}:=\frac{1}{(n_*\lambda)^2}(\xi\otimes\nabla\Phi_{(\xi)}-\nabla\Phi_{(\xi)}\otimes\xi)\psi_{(\xi)}.\notag
\end{align}

Then by a direct computation
\begin{align}
\div V_{(\xi)}&=\psi_{(\xi)}\phi_{(\xi)}\xi -\frac{1}{(n_*\lambda)^2}\nabla\Phi_{(\xi)}\xi\cdot\nabla\psi_{(\xi)}=W_{(\xi)} -\frac{1}{(n_*\lambda)^2}\nabla\Phi_{(\xi)}\xi\cdot\nabla\psi_{(\xi)}.\label{divOmega}
\end{align}
Finally, we obtain that for $N, M \geq 0$ and $p\in [1, \infty]$ the
following holds
\begin{align}
\|\nabla^N\partial_t^M\psi_{(\xi)}\|_{C_tL^p}\lesssim r_\parallel^{\frac{1}{p}-\frac12}(\frac{r_\perp\lambda}{r_\parallel})^N(\frac{r_\perp\lambda\mu}{r_\parallel})^M,&\label{int2}\\
\|\nabla^N\phi_{(\xi)}\|_{L^p}+\|\nabla^N\Phi_{(\xi)}\|_{L^p}\lesssim r_\perp^{\frac{d-1}{p}-\frac{d-1}2}\lambda^N,&\label{int3}\\
\|\nabla^N\partial_t^MW_{(\xi)}\|_{C_tL^p}+\lambda\|\nabla^N\partial_t^MV_{(\xi)}\|_{C_tL^p}
\lesssim r_\perp^{\frac{d-1}{p}-\frac{d-1}2}r_\parallel^{\frac{1}{p}-\frac12}\lambda^N(\frac{r_\perp\lambda\mu}{r_\parallel})^M,&\label{int4}\\
\|\nabla^N\partial_t^M\Theta_{(\xi)}\|_{C_tL^p}\lesssim r_\perp^{\frac{d-1}{p}-\frac{d-1}{2}}r_\parallel^{\frac{1}{p}-\frac12}\lambda^N(\frac{r_\perp\lambda\mu}{r_\parallel})^M,& \label{int4theta}
\end{align}
where the implicit constants may depend on $p,N$ and $M$, but are independent of $\lambda,r_\perp,r_\parallel,\mu$. These estimates can be easily deduced from the definitions.

\subsubsection{Building blocks for the Navier-Stokes or Euler equations}\label{sec:bbns}
In this section we recall the generalized intermittent jets introduced in \cite[Section3]{LZ23}.

We recall the geometric Lemma \ref{l:linear_algebra} and the set  $\overline \Lambda:=\overline \Lambda^1$  discussed in that lemma for $d\geq2$. Additionally, we use the parameters $\lambda,r_\perp, r_\parallel > 0 $ defined in the previous section.
For each $ {\xi}\in \overline{\Lambda}$ let us define $\overline A^i_{ {\xi}}\in \mathbb{S}^{d-1}\cap \mathbb{Q}^d,\ i=1,2,...,d-1$,
 such that $\{ {\xi}, A^i_{ {\xi}},i=1,...,d-1\} $ form an orthonormal basis in $\mathbb{R}^d$. 
We label by $ \overline{n}_*$ the smallest
natural number such that for every $ {\xi}\in \overline{\Lambda}$
$$\{ \overline{n}_* {\xi}, \overline {n}_*\overline A^i_{ {\xi}},i=1,...,d-1\}\subset\mathbb{Z}^d.$$

Let $\overline {\Phi} : \mathbb{R}^{d-1} \to \mathbb{R}$ be a smooth function with support in a ball of radius 1. We normalize $\overline {\Phi}$ such that $\overline {\phi}  = -\Delta\overline {\Phi} $ obeys
$$\int_{\mathbb{R}^{d-1}} \overline{\phi}^2(x_1,x_2,...,x_{d-1})\dif x_1\dif x_2..\dif x_{d-1}=1.$$
By definition we know that $\int_{\mathbb{R}^{d-1}}\overline {\phi} \dif x_1\dif x_2..\dif x_{d-1}=0$.

 Define $ \overline{\psi}  : \mathbb{R}\to\mathbb{R}$ to be a smooth, mean-zero function with
support in a ball of radius 1 satisfying 
$$\int_\mathbb{R} \overline{\psi}^2(x_d)\dif x_d=1.$$

We define the rescaled cut-off functions 
$$ \overline{\phi} _{ {r}_{\perp}}(x_1,x_2,...,x_{d-1})=\frac1{ {r}_{\perp}^{(d-1)/2}} \overline{\phi }(\frac{x_1}{ {r}_{\perp}},\frac{x_2}{ {r}_{\perp}},...,\frac{x_{d-1}}{ {r}_{\perp}}),$$
$$ \overline{\Phi} _{ {r}_{\perp}}(x_1,x_2,...,x_{d-1})=\frac1{ {r}_{\perp}^{(d-1)/2}}\overline {\Phi} (\frac{x_1}{ {r}_{\perp}},\frac{x_2}{ {r}_{\perp}},...,\frac{x_{d-1}}{ {r}_{\perp}}),$$
$$ \overline{\psi} _{ {r}_{\parallel}}(x_d)=\frac1{ {r}_{\parallel}^{1/2}} \overline{\psi}(\frac{x_d}{ {r}_{\parallel}}).$$
We periodize them so that they can be viewed as functions on $\mathbb{T}^{d-1}$ and $\mathbb{T}$ respectively.
Consider a new large time oscillation parameter $\overline {\mu} > 0$. For every $ {\xi}\in \overline{\Lambda}$ we introduce
$$ \overline{\psi} _{( {\xi})}(t,x):=
\overline {\psi} _{ {r}_{\parallel}}( \overline{n}_* {r}_{\perp} {\lambda}(x\cdot { {\xi}}-\overline  {\mu} t)),$$
 $$ \overline{\Phi} _{( {\xi})}(x):=
\overline {\Phi} _{ {r}_{\perp}}(\overline {n}_* {r}_{\perp} {\lambda} x\cdot\overline A^1_{ {\xi}},..., \overline{n}_* {r}_{\perp} {\lambda}x\cdot\overline A^{d-1}_{{ {\xi}}}),$$
$$ \overline{\phi} _{( {\xi})}(x):=
\overline {\phi} _{ {r}_{\perp}}( \overline{n}_* {r}_{\perp} {\lambda}x\cdot\overline A^1_{{ {\xi}}},..., \overline{n}_* {r}_{\perp} {\lambda}x\cdot\overline A^{d-1}_{{ {\xi}}}).$$
Here, similarly as before, we do not need to translate the building blocks since the the disjoint support property will be achieved by selecting suitable time jets in the following section.

The intermittent jets $ \overline{W}_{( {\xi})} :\mathbb{R} \times\mathbb{T}^d \to \mathbb{R}^d$ are defined as
$$\overline {W}_{( {\xi})}(t,x)= {\xi}\overline {\psi} _{( {\xi})}(t,x) \overline{\phi} _{( {\xi})}(x).$$

By definition and a basic calculation we have that
\begin{align}
\partial_t(\overline {\psi}^2_{( {\xi})} \overline{\phi}^2_{( {\xi})} {\xi})+ \overline {\mu}\div( \overline{W}_{( {\xi})}\otimes \overline {W}_{( {\xi})})=0,\label{*5}\\
\int_{\mathbb{T}^d}\overline W_{( {\xi})}\otimes\overline W_{( {\xi})}\dif x= {\xi}\otimes {\xi}.\label{***}
\end{align}

Since $\overline {W}_{( {\xi})}$ is not divergence-free, we
introduce the skew-symmetric  corrector term
\begin{align}
 \overline{V}_{({ {\xi}})}:=\frac{1}{( \overline{n}_* {\lambda})^2}({ {\xi}}\otimes\nabla\overline {\Phi}_{({ {\xi}})}-\nabla \overline{\Phi}_{({ {\xi}})}\otimes{ {\xi}}) \overline{\psi}_{({ {\xi}})}.\notag
\end{align}

Then by a direct computation
\begin{align}
\div \overline {V}_{({ {\xi}})}=\overline {W}_{({ {\xi}})} -\frac{1}{(\overline {n}_* {\lambda})^2}\nabla\overline {\Phi}_{({ {\xi}})}{ {\xi}}\cdot\nabla \overline{\psi}_{({ {\xi}})}.\label{divOmegans}
\end{align}
Finally, we obtain that for $N, M \geq 0$ and $p\in [1, \infty]$ the
following holds
\begin{align}
\|\nabla^N\partial_t^M\overline {\psi}_{({ {\xi}})}\|_{C_tL^p}&\lesssim { {r}}_\parallel^{\frac{1}{p}-\frac12}(\frac{{ {r}}_\perp{ {\lambda}}}{{ {r}}_\parallel})^N(\frac{{ {r}}_\perp{ {\lambda}}{\overline  {\mu}}}{{ {r}}_\parallel})^M,\label{int2ns}\\
\|\nabla^N\overline {\Phi}_{({ {\xi}})}\|_{L^p}+\|\nabla^N\overline {\phi}_{({ {\xi}})}\|_{L^p}&\lesssim { {r}}_\perp^{\frac{d-1}{p}-\frac{d-1}2}{ {\lambda}}^N,\label{int3ns}\\
\|\nabla^N\partial_t^M \overline{W}_{({ {\xi}})}\|_{C_tL^p}+{ {\lambda}}\|\nabla^N\partial_t^M \overline{V}_{({ {\xi}})}\|_{C_tL^p}&\lesssim { {r}}_\perp^{\frac{d-1}{p}-\frac{d-1}2}{ {r}}_\parallel^{\frac{1}{p}-\frac12}{ {\lambda}}^N(\frac{{ {r}}_\perp{ {\lambda}}{\overline  {\mu}}}{{ {r}}_\parallel})^M,\label{int4ns}
\end{align}
where the implicit constants may depend on $p,N$ and $M$, but are independent of ${ {\lambda}},{ {r}}_\perp,{ {r}}_\parallel,{\overline {\mu}}$.

\subsubsection{Temporal jets}\label{sec:tj}
 In this section, we introduce additional intermittency in the  time direction similarly as in \cite[Section 4.2]{CL21}.
For $\xi\in\Lambda^1\cup\Lambda^2\cup\overline\Lambda$, let us choose temporal functions $g_{(\xi)}(t)$ and $h_{(\xi)}
(t)$ to oscillate the
building blocks intermittently in time. Let $G \in C_c^\infty(0, 1)$ be non-negative and 
$$\ \ \int_0^1G^2(t)\dif t=1.$$

Letting $\eta >0$  be a small constant satisfying $\eta \cdot{\rm c ard}(\Lambda^1\cup\Lambda^2\cup\overline\Lambda)\ll1$, we define $\tilde{g}_{(\xi)}: \mathbb{T}\to \mathbb{R}$ as the 1-periodic extension of $\eta^{-\frac12}G(\frac {t -t_\xi}{\eta})$, where $t_\xi$ are chosen so that $\tilde{g}_{(\xi)}$ have disjoint supports for different $\xi\in \Lambda^1\cup\Lambda^2\cup\overline \Lambda$. 
We will also oscillate the perturbations at a large
frequency $\sigma\in\mathbb{N}$. So, we define
$$g_{(\xi)}(t)=\tilde{g}_{(\xi)}(\sigma t).$$
For the corrector term we define $H_{(\xi)},h_{(\xi)}:\mathbb{T}\to\mathbb{R}$ by
\begin{align}
H_{(\xi)}(t)=\int_0^{t}g_{(\xi)}(s)\dif s,\ \ h_{(\xi)}(t)=\int_0^{\sigma t}(\tilde{g}_{(\xi)}^2(s)-1)\dif s. \label{eq:parth}
\end{align}

In view of the zero-mean condition for $\tilde{g}_{(\xi)}^2(t)-1$,  we obtain  that $h_{(\xi)}$ is $\mathbb{T}/\sigma$-periodic, and for any $n\geq0,p\geq1$
\begin{align}
    \|g_{(\xi)}\|_{W_t^{n,p}}\lesssim(\frac{\sigma}{\eta})^n\eta^{\frac1p-\frac12},\  \  \|h_{(\xi)}\|_{L_t^\infty}\leq1.\label{bd:gwnp}
\end{align}

\subsection{The estimates on the amplitude functions}
  In this section, we show the estimates of the amplitude functions appearing in Section \ref{proof:prop2}.
\begin{proof}[Proof of Proposition \ref{lem:2chi}]
      First, we  estimate  $\chi(\zeta|M_l|-n)$ in $C_{t,x}^N$-norm for $N\in \mN$.  We recall that from \eqref{bd:2mlcn}
      \begin{align*}
    \|M_l\|_{C_{t,x}^N}\lesssim l^{-d-2-N}.
 \end{align*} 
 On $\supp\chi(\zeta|M_l|-n),n\geq3$, we have $|M_l|\geq \zeta^{-1}$.  Then we apply \cite[Proposition C.1]{TCP} to a smooth function $f(z)$ satisfying  $f(z)=|z|$ on $|z|\geq\zeta^{-1}$.  Since $|D^Nf(z)|\lesssim\zeta^{N-1}$ on $|z|\geq\zeta^{-1}$, we have for  $N\geq1$
  \begin{align*}
      \left\||M_l|\right\|_{C_{t,x}^N}&\lesssim   \left\||M_l|\right\|_{C_{t,x}^0}+\|Df\|_{C^0}\|M_l\|_{C_{t,x}^N}+\|Df\|_{C^{N-1}}\|M_l\|_{C_{t,x}^1}^N\\
      &\lesssim l^{-d-2-N}+  \zeta^{N-1}l^{-(d+3)N}\lesssim \zeta^{N-1}l^{-(d+3)N}. 
  \end{align*}
 Then we apply the chain rule  from \cite[Proposition 4.1]{TCP} to $f(z)=\chi(\zeta z-n),|D^mf|\lesssim \zeta^m$ to obtain
  \begin{align*}
     \|\chi(\zeta|M_l|-n)\|_{C_{t,x}^N}&\lesssim  \|\chi(\zeta|M_l|-n)\|_{C_{t,x}^0}+\|Df\|_{C^0}\norm{|M_l|}_{C_{t,x}^N}+\|Df\|_{C^{N-1}}\norm{|M_l|}_{C_{t,x}^1}^N\\
     &\lesssim \zeta^N l^{-N(d+3)}\lesssim l^{-N(d+4)},
  \end{align*}
  where we used the condition $\zeta\lesssim l^{-1}$. Then  by \eqref{bd:2|ml|} we have for $N\in\mN$
  \begin{align*}
        \sum_{n\geq3}\|\chi(\zeta|M_l|-n)\|_{C_{t,x}^N}\lesssim \sum_{n=3}^{1+Cl^{-d-2}} l^{-N(d+4)}\lesssim l^{-N(d+4)-(d+2)}.
  \end{align*}
This bound is also valid for $N=0$. The bound for $\tilde{\chi}(\zeta|M_l|-n)$ is similar to the one described above.
 
 Next we  estimate $\chi(\zeta|M_l|-n)\Gamma_\xi(\frac{M_l}{|M_l|})$  in $C_{t,x}^N$-norm. By the Leibniz rule we get 
\begin{align*}
    \left\|\frac{M_l}{|M_l|}\right\|_{C_{t,x}^N}\lesssim \sum_{m=0}^N\|M_l\|_{C_{t,x}^{N-m}}\left\|\frac{1}{|M_l|}\right\|_{C_{t,x}^{m}}.
\end{align*}

 On $\supp\chi(\zeta|M_l|-n),n\geq3$, we have $|M_l|\geq \zeta^{-1}$.   Then we apply \cite[Proposition C.1]{TCP} to a smooth function $f(z)$ satisfying  $f(z)=\frac1{|z|}$ on $|z|\geq\zeta^{-1}$. Since $|D^Nf(z)|\lesssim\zeta^{N+1}$ on $|z|\geq\zeta^{-1}$, we have for  $N\geq1$
 \begin{align*}
     \left\|\frac{1}{|M_l|}\right\|_{C_{t,x}^{m}}\lesssim  \zeta+\zeta^{2} l^{-d-2-m}+\zeta^{m+1} l^{-m(d+3)}\lesssim l^{-m(d+4)-1},
 \end{align*}
 which implies that 
 \begin{align*}
      \left\|\frac{M_l}{|M_l|}\right\|_{C_{t,x}^N}\lesssim \sum_{m=1}^N  l^{-d-2-N+m}l^{-m(d+4)-1}+  l^{-d-2-N}l^{-1}\lesssim l^{-(d+5)N-(d+3)}.
 \end{align*}
 We apply the chain rule from \cite[Proposition 4.1]{TCP} to $f(z)=\Gamma_\xi(z),|D^mf|\lesssim 1$ to obtain  for $N\geq1$
 \begin{align*}
     \left\|\Gamma_\xi\(\frac{M_l}{|M_l|}\)\right\|_{C_{t,x}^N}\lesssim \left\|\frac{M_l}{|M_l|}\right\|_{C_{t,x}^N}+\left\|\frac{M_l}{|M_l|}\right\|_{C_{t,x}^1}^N\lesssim l^{-(2d+8)N}.
 \end{align*}
 By the chain rule we get for $N\in\mN$
 \begin{align*}
     \sum_{n\geq3}\sum_{\xi\in\Lambda^n}\left\|\chi(\zeta|M_l|-n)\Gamma_\xi\(\frac{M_l}{|M_l|}\)\right\|_{C_{t,x}^N}\lesssim l^{-(2d+8)N-(d+2)}.
 \end{align*}
 This bound is also valid for $N=0$.
 
 For the last term, we have for $N\in\mN_0$
 \begin{align*}
     \(\frac{n}{\zeta}\)^N1_{ \{\chi(\zeta|M_l|-n)>0\}}+ \(\frac{n}{\zeta}\)^N1_{ \{\tilde{\chi}(\zeta|M_l|-n)>0\}}\lesssim (\|M_l\|_{C_{t,x}^0}+\zeta^{-1})^N \lesssim l^{-N(d+2)}.
 \end{align*}
\end{proof}

\begin{proof}
    [Proof of Proposition  \ref{prop:bd:2axi}]
 First we recall the definition of $A$:
\begin{align}
A:=2\sqrt{ {l}^2+|\mathring{ {R}}_{ {l}}|^2}
,\notag
\end{align}
which together with \eqref{bd:2rlcn} implies that
   $\|A\|_{C_{t,x}^0}\lesssim  {l}^{-d-2}$.

Next we estimate the $C_{t,x}^N$-norm for $N\in\mathbb{N}$. We apply the chain rule in \cite[Proposition C.1]{TCP}  to $f(z)=\sqrt{ {l}^2+z^2},|D^mf(z)|\lesssim  {l}^{-m+1}$ to obtain
\begin{align*}
\left\|\sqrt{ {l}^2+|\mathring{ {R}}_{ {l}}|^2}\right\|_{C_{t,x}^N}&\lesssim\left\|\sqrt{ {l}^2+|\mathring{ {R}}_{ {l}}|^2}\right\|_{C_{t,x}^0}+\|Df\|_{C^0}\|\mathring{ {R}}_{ {l}}\|_{C_{t,x}^N}+\|Df\|_{C^{N-1}}\|\mathring{ {R}}_{ {l}}\|_{C_{t,x}^1}^N\\
&\lesssim  {l}^{-d-2-N}+ {l}^{-N+1} {l}^{-(d+3)N},
\end{align*}
 which implies that for $N\geq1$
\begin{align}\label{A.1}
\|A\|_{C_{t,x}^N}&\lesssim\left\|\sqrt{ {l}^2+|\mathring{ {R}}_{ {l}}|^2}\right\|_{C_{t,x}^N}
\lesssim  {l}^{1-(d+4)N}.
\end{align}

Let us now estimate the $C_{t,x}^N$-norm. By the Leibniz rule we get
$$\|a_{( {\xi})}\|_{C_{t,x}^N}\lesssim\sum_{m=0}^N\left\|A^{1/2}\right\|_{C_{t,x}^m}\left\|\gamma_{ {\xi}}(\Id-\frac{\mathring{ {R}}_{ {l}}}{A})\right\|_{C_{t,x}^{N-m}}.$$
Applying \cite[Proposition C.1]{TCP} to $f(z)=z^{1/2}$, $|D^mf(z)|\lesssim|z|^{1/2-m}$, for $m=1,...,N$, and using \eqref{A.1} we obtain  for $m\geq1$
\begin{align*}
\|A^{1/2}\|_{C_{t,x}^m}&\lesssim \|A^{1/2}\|_{C_{t,x}^0}+ {l}^{-1/2}\|A\|_{C_{t,x}^m}+ {l}^{1/2-m}\|A\|_{C_{t,x}^1}^m\lesssim  {l}^{1/2-(d+4)m}.
\end{align*}

Next we estimate $\gamma_{ {\xi}}(\Id-\frac{\mathring{ {R}}_{ {l}}}{A})$. By \cite[Proposition C.1]{TCP} we need to estimate
$$\left\|\frac{\mathring{ {R}}_{ {l}}}{A}\right\|_{C_{t,x}^{N-m}}+\left\|\frac{\nabla_{t,x}\mathring{ {R}}_{ {l}}}{A}\right\|_{C_{t,x}^0}^{N-m}+\left\|\frac{\mathring{ {R}}_{ {l}}}{A^2}\right\|_{C_{t,x}^0}^{N-m}\left\|A\right\|_{C_{t,x}^1}^{N-m}.$$
We use $A\geq l$  to have that
$$\left\|\frac{\nabla_{t,x}\mathring{ {R}}_{ {l}}}{A}\right\|_{C_{t,x}^0}^{N-m}\lesssim  {l}^{-N+m} {l}^{-(d+3)(N-m)}\lesssim  {l}^{-(d+4)(N-m)},$$
and in view of $|\frac{\mathring{ {R}}_{ {l}}}{A}|\leq 1$ that
$$\left\|\frac{\mathring{ {R}}_{ {l}}}{A^2}\right\|_{C_{t,x}^0}^{N-m}\lesssim\left\|\frac{1}{A}\right\|_{C_{t,x}^0}^{N-m}\lesssim  {l}^{-N+m},$$
and by \eqref{A.1} that
$$\|A\|_{C_{t,x}^1}^{N-m}\lesssim  {l}^{-(d+3)(N-m)}.$$
Moreover, we write
\begin{align*}
\left\|\frac{\mathring{ {R}}_{ {l}}}{A}\right\|_{C_{t,x}^{N-m}}&\lesssim  \sum_{k=0}^{N-m}\|\mathring{ {R}}_{ {l}}\|_{C_{t,x}^k}\left\|\frac{1}{A}\right\|_{C_{t,x}^{N-m-k}}.
\end{align*}
Using \eqref{A.1}  and \cite[Proposition C.1]{TCP} we obtain
\begin{align*}
\left\|\frac{1}{A}\right\|_{C_{t,x}^{N-m-k}}\lesssim \left\|\frac{1}{A}\right\|_{C_{t,x}^{0}}+ {l}^{-2}\|{A}\|_{C_{t,x}^{N-m-k}}+ {l}^{-N+m+k-1}\|{A}\|_{C_{t,x}^1}^{N-m-k}\\
\lesssim  {l}^{-2} {l}^{1-(d+4)(N-m-k)}+ {l}^{-(N-m-k)-1} {l}^{-(d+3)(N-m-k)}\lesssim  {l}^{-(d+4)(N-m-k)-1}.
\end{align*}
Thus, we obtain
\begin{align*}
\|\frac{\mathring{ {R}}_{ {l}}}{A}\|_{C_{t,x}^{N-m}}&\lesssim  \sum_{k=0}^{N-m-1} {l}^{-d-2-k} {l}^{-(d+4)(N-m-k)-1}+ {l}^{-d-2-(N-m)} {l}^{-1}\lesssim  {l}^{-(d+3)-(d+4)(N-m)}.
\end{align*}
Finally, the above bounds lead to
$$\norm{\gamma_{ {\xi}}\(\Id-\frac{\mathring{ {R}}_{ {l}}}{A}\)}_{C_{t,x}^{N-m}}\lesssim  {l}^{-(d+3)-(d+4)(N-m)}.$$
Combining this with the bounds for $A^{1/2}$ above yields for $N\in\mathbb{N}$
\begin{align*}
\|a_{({ {\xi}})}\|_{C_{t,x}^N}&\lesssim {l}^{-(d+2)/2} {l}^{-(d+3)-(d+4)N}+\sum_{m=1}^{N-1} {l}^{1/2-(d+4)m} {l}^{-(d+3)-(d+4)(N-m)}+ {l}^{1/2-(d+4)N}\\
&\lesssim  {l}^{-2d-3-(d+4)N},
\end{align*}
where the final bound is also valid for $N=0$.

  \end{proof}

\end{document}